\documentclass[]{ociamthesis}
\usepackage{amsmath,amsfonts,amsthm}

\title{Metrics of Special Curvature with Symmetry}
\author{Brandon Dammerman}
\college{The Queen's College}
\degree{Doctor of Philosophy} \degreedate{May 2004}

\begin{document}
\baselineskip=18pt plus1pt
 \setcounter{secnumdepth}{2} \setcounter{tocdepth}{3}

\maketitle
\renewenvironment{acknowledgementslong}
 \thispagestyle{empty}
\begin{center}
\vspace*{1.5cm} {\Large \bfseries Acknowledgements}
\end{center}
\vspace{0.5cm}
\begin{quote}
Foremost, I would like to thank my supervisor, Dr. Andrew Dancer,
for his patient, lucid, and insightful instruction. The results of
Chapter $4$ were obtained under the guidance of Prof. K. P. Tod, and
I am very grateful to him for his help and instruction. Thanks also
to my fellow graduate students Gil Cavalcanti, Marco Gualtieri,
Steve Marshall, and Fred Witt for useful discussions.

Finally, I thank the Rhodes Trust for its generous funding and The
Queen's College for granting me a lectureship in my final year.
\end{quote}

\pagebreak
\include{abst}
\renewenvironment{abstract}  \thispagestyle{empty}
\begin{center}

  {\Large \bfseries  Metrics of Special Curvature with Symmetry}
  \end{center}
\vspace*{.45cm}
\begin{center}
Brandon Dammerman
\end{center}
\begin{center}
Oxford University D. Phil. Thesis
\end{center}
\begin{center}
May 2004
\end{center}
  \begin{center}
\vspace*{.45cm}
  {\Large \bfseries  Abstract}
  \end{center}
  \vspace{0.5cm}
\begin{quote}
~~~Various curvature conditions are studied on metrics admitting a
symmetry group.  We begin by examining a method of diagonalizing
cohomogeneity-one Einstein manifolds and determine when this method
can and cannot be used.  Examples, including the well-known Stenzel
metrics, are discussed.

~~~Next, we present a simplification of the Einstein condition on a
compact four manifold with $T^{2}$-isometry to a system of
second-order elliptic equations in two-variables with well-defined
boundary conditions.

~~~We then study the Einstein and extremal Kahler conditions on
Kahler toric manifolds.  After constructing explicitly new extremal
Kahler and constant scalar curvature metrics on
$\mathbb{P}(\mathcal{O}_{\mathbb{C}P^{1}\times\mathbb{C}P^{1}}
\oplus\mathcal{O}_{\mathbb{C}P^{1}\times\mathbb{C}P^{1}}(1,-1))$, we
demonstrate how these metrics can be obtained by continuously
deforming the Fubini-Study metric on $\mathbb{C}P^{3}$.  We also
define a generalization of Kahler toric manifolds, which we call
fiberwise Kahler toric manifolds, and construct new explicit
extremal Kahler and constant scalar curvature metrics on both
compact and non-compact manifolds in all even dimensions. We also
calculate the Futaki invariant on manifolds of this type.

~~~After describing an Hermitian non-Kahler analogue to fiberwise
Kahler toric geometry, we construct constant scalar curvature
Hermitian metrics with $J$-invariant Riemannian tensor.  In
dimension four, we write down explicitly new constant scalar
curvature Hermitian metrics with $J$-invariant Ricci tensor on
$\mathbb{C}P^{1}\times\mathbb{C}P^{1}$ and
$\mathbb{C}P^{2}\sharp\overline{\mathbb{C}P^{2}}$.

~~~Finally, we integrate the scalar curvature equation on a large
class of cohomogeneity-one metrics.

\end{quote}

\newcommand{\Ind}{\mathrm{ind}}
\newcommand{\Sym}{\mathrm{Sym}}
\newcommand{\Ker}{\operatorname{Ker}}
\newcommand{\rk}{\mathrm{rk\ }}
\newcommand{\Tr}{\mathrm{Tr}}
\newcommand{\Proof}{\noindent{\em Proof:} }
\newcommand{\IP}[1]{\langle #1\rangle}
\newcommand{\Ip}[1]{\left( #1\right)}
\newtheorem{theorem}{Theorem}[chapter]
\newtheorem{corollary}[theorem]{Corollary}
\newtheorem{lemma}[theorem]{Lemma}
\newtheorem{prop}[theorem]{Proposition}
\newtheorem{defn}[]{Definition}
\newtheorem{example}[defn]{Example}
\newtheorem{remark}[theorem]{Remark}
\newcommand{\OO}{\mathcal O}
\newcommand{\JJ}{\mathcal J}
\newcommand{\ZZ}{\mathbb Z}
\newcommand{\RR}{\mathbb R}
\newcommand{\GG}{\mathfrak{g}_e}
\newcommand{\GGG}{\mathfrak{g}}
\newcommand{\CC}{\mathbb C}
\newcommand{\SSS}{\mathbb S}
\newcommand{\PP}{\mathbb P}
\newcommand{\DD}{\mathcal{D}}
\newcommand{\HH}{\mathbb H}
\newcommand{\MM}{\mathbf{M}}
\newcommand{\im}{\mathrm{im\ }}
\newcommand{\Aa}{\underline{A}}
\newcommand{\Bb}{\underline{B}}
\newcommand{\Cc}{\underline{C}}
\newcommand{\Vv}{\underline{V}}
\newcommand{\Ww}{\underline{W}}
\renewcommand{\arraystretch}{1.2}

\newcommand{\CPT}{\CC P^3}
\newcommand{\w}{\omega}
\newcommand{\s}{\sigma}
\newcommand{\hs}{\hat\sigma}
\newcommand{\A}{\mathfrak{a}}
\newcommand{\X}{\tilde X}
\newcommand{\f}{\tilde f}
\newcommand{\lra}{\longrightarrow}
\newcommand{\ra}{\rightarrow}
\author{Brandon Dammerman}
\title{Metrics of Special Curvature with Symmetry}

\chapter{Introduction}

    Much of Riemannian geometry involves the search for distinguished
metrics within the infinite-dimensional space of smooth metrics on a
$C^{\infty}$-differentiable manifold, particularly when the manifold
is of dimension greater than or equal to four.  Here, the term
'distinguished' refers to those metrics that exhibit additional
structure or satisfy some curvature condition.  For example, one
could look for metrics of constant Riemannian or constant scalar
curvature or whose Ricci tensor satisfies the Einstein condition. In
particular, Einstein metrics are of interest to both geometers and
physicists.  Originally, the study of Einstein metrics was motivated
by theoretical physics; space-time in the vacuum state is an
Einstein manifold with zero Einstein constant.  Furthermore, in the
compact case Einstein metrics are critical points of the scalar
curvature functional. Alternatively, one could study metrics with a
compatible complex structure or metrics of special holonomy. Once a
condition has been placed upon a metric, the mathematician tries to
determine which manifolds admit metrics of that type, how abundant
those metrics are on a given manifold, and if at all possible write
down the metrics explicitly.  In this thesis, we study the Einstein
and constant scalar curvature conditions on manifolds of general
type as well as on those admitting Hermitian and Kahler structures,
in many cases deriving explicit solutions.  We place particular
emphasis on those cases in which the equations of the curvature
conditions take a simplified form. Simplification is accomplished in
two ways: we either demand that the manifold admit a large symmetry
group or we use Kaluza-Klein bundle constructions.  Both require an
understanding of Riemannian submersions

    The thesis begins with a brief review of the theory of Riemannian
submersions using the formalism developed by O'Neill in
\cite{O'Neill}.  Indeed, though the thesis considers various
curvature conditions on metrics, some of which have compatible
complex structures and some of which do not, all of the manifolds
considered can be viewed as Riemannian submersions.  Therefore, an
understanding of the curvature conditions on a Riemannian
submersion will be essential to every part of the thesis.

    In Chapter $2$, we also discuss a special kind of Riemannian submersion:
the quotient by a group action.  More precisely, we present the
equations of the Einstein and constant scalar curvature conditions
on a metric admitting \textit{orthogonally transitive} isometries,
i.e. those isometric group actions in which the distribution
orthogonal to the tangent space of the orbits of the action is
integrable.  Though not all of the metrics considered below admit
orthogonally transitive isometries, many do.  Isometry group
actions on cohomogeneity-one manifolds are automatically
orthogonally transitive, and cohomogeneity-one manifolds will be
considered in Chapters $3,$ $5-8,$ and $11$.  Manifolds of higher
cohomogeneity under orthogonally transitive actions will be
considered in Chapters $4-7$ and $9$.  In Chapter $2$, we also
present a rough topological classification of Einstein manifolds
with an orthogonally transitive group action.

    Most of the new material begins in Chapter $3$ with our investigation
of when it is possible to globally diagonalize an Einstein metric
of cohomogeneity-one under the action of a compact, connected
group $G$.  The Einstein condition on a cohomogeneity-one manifold
reduces to a system of ordinary differential equations, cf.
\cite{Bergery}.  While this is a great simplification of the
general Einstein condition, as yet there exists no classification
of cohomogeneity-one Einstein metrics even in the compact case.
Indeed, such a classification appears, at the moment, to be only a
distant goal.  There are, however, many known examples of
cohomogeneity-one Einstein metrics in all dimensions.  To the
author's knowledge, \textit{all} of the known metrics can be
written in diagonal form relative to a $G$-invariant basis.  Many
of these examples have been found either by looking at those group
actions under which every cohomogeneity-one metric is
automatically diagonalizable (called the monotypic case) or by
demanding \textit{ad hoc} that the metric be diagonal. In some
special cases, it has been observed that the Einstein condition
allows one to globally diagonalize certain cohomogeneity-one
metrics in the non-monotypic case.  For example, it is well-known
to relativists that the Einstein condition on a four manifold with
$SU(2)$-symmetry allows one to globally diagonalize the metric,
though this result does not seem to appear in the pure mathematics
literature. From that result, we devise a general procedure for
determining when the Einstein condition can be used to globally
diagonalize a cohomogeneity-one metric in the non-monotypic case.
As an illustration, we show that this method works in the case of
$SO(n+1)/SO(n-1)$, with $n\geq 3$, principal orbits (the Stenzel
metrics found in \cite{Stenzel} are of this type) but fails when
the principal orbits are isomorphic to $SU(n+1)/U(n-1)$ or
$SO(n_{1}+n_{2}+2)/(SO(n_{1})\times SO(n_{2}))$ with both $n_{1}$
and $n_{2}$ greater than or equal to two.

    After reviewing what is known about Einstein four manifolds
admitting a large isometry group, we leave the cohomogeneity-one
category to investigate compact Einstein four-manifolds with
$T^{2}$-isometry in Chapter $4$.  Though it was known previously
\cite{Kundt} that the torus action must be orthogonally transitive
in the Einstein case, little else is known about metrics of this
type in general.  To date, the only known Einstein metrics of this
type are the homogeneous metrics on $S^{4}$, $S^{2}\times S^{2}$,
and $\mathbb{C}P^{2}$; the Hermitian-Einstein metric (the Page
metric) on $\mathbb{C}P^{2}\sharp\overline{\mathbb{C}P^{2}}$; and
the Kahler-Einstein metric (found by Siu in \cite{Siu}) on
$\mathbb{C}P^{2}\sharp3\overline{\mathbb{C}P^{2}}$, though this
metric has not been written down explicitly.  After reviewing the
topological classification of smooth torus actions on compact
simply-connected four manifolds discovered by Orlik and Raymond in
\cite{OandR} and the \textit{nuts} and \textit{bolts} classification
developed by Gibbons and Hawking, we present our main result which
provides a simplification of the Einstein equations.  In general,
the Einstein condition is an overdetermined system of partial
differential equations.  In the compact $T^{2}$-invariant case, we
show that the system reduces to a set of four second-order elliptic
equations in four independent functions in two variables together
with a set of well-defined boundary conditions.  This simplification
is effected by making use of a natural complex structure on the
two-dimensional quotient space $M/T^{2}$.  While we are unable at
this time to prove the existence of any new Einstein metrics, the
reduction makes the system much more amenable to analytical methods
for determining the existence and uniqueness of $T^{2}$-invariant
Einstein metrics.  It is our hope that our work will provide a
useful first step in the eventual classification of such metrics.
Finally, we end this section with a discussion of a simplified
system in which the metric is diagonal.  In this case, more can be
said about the geometric structure of Einstein metrics.

    Chapter $5$ begins a study of curvature conditions on manifolds
with large symmetry group in the presence of an integrable complex
structure.  Since the Einstein condition appears, at the moment,
to be difficult to solve even in the relatively simple case of a
four manifold with $T^{2}$-symmetry, there is little hope of
arriving at a solution to the Einstein condition on a
$2n$-dimensional manifold with $T^{n}$-symmetry.  It seems
appropriate therefore to require that the manifold admit
additional structures which simplify the curvature equations.  One
such condition which is commonly studied is the reduction of the
holonomy group to $U(n)$ from $SO(2n)$.  In this case the manifold
will be Kahler and the various curvature conditions can be
expressed in much simplified form.  When the manifold is Kahler
toric, i.e. when there exists a $T^{n}$ action preserving the
Kahler structure, these conditions take yet simpler forms as we
will see below.

    Kahler toric manifolds are of interest because of, among other
reasons, the presence of a Legendre transformation between complex
coordinates and symplectic Darboux coordinates.  A recent theorem
of Abreu \cite{Abreu} allows one to fix a symplectic structure on
a toric manifold and parameterize all of the compatible Kahler
metrics via a symplectic potential.  It is this property, which is
absent in general Kahler geometry, that makes Kahler toric
geometry so tractable.  After reviewing the symplectic formalism,
we discuss various curvature conditions on Kahler toric manifolds.
In addition to the Kahler-Einstein condition, we also discuss the
extremal Kahler condition (introduced by Calabi in \cite{Calabi})
and the more restrictive constant scalar curvature condition on
Kahler toric manifolds,  both of which have been widely studied by
differential geometers.  While the majority of Chapter $5$ is a
review of material necessary for later chapters, we do present a
simplification of the scalar curvature equation on a Kahler toric
manifold found by Abreu in \cite{Abreu2}.

    Specific examples of Kahler toric manifolds with distinguished
metrics, particularly in real dimensions four and six, are
discussed in Chapter $6$.  A great deal is known about the
existence of Kahler-Einstein toric manifolds (see \cite{Mabuchi}),
though relatively few metrics are known explicitly.  Before
looking at explicit examples, we show how to integrate the
Einstein equations to obtain the Legendre transform of the
Monge-Ampere equations of complex coordinates and write down the
Futaki invariant, an obstruction to the existence of positive
Kahler-Einstein metrics, in symplectic coordinates.  We discuss
the Kahler-Einstein metric on
$\mathbb{C}P^{2}\sharp3\overline{\mathbb{C}P^{2}}$ which has not
been written down explicitly.  We analyze this case in detail and
propose an ansatz for the possible form of the Siu metric. Next,
we proceed to write down the Kahler-Einstein metric found by
Sakane on
$\mathbb{P}(\mathcal{O}_{\mathbb{C}P^{1}\times\mathbb{C}P^{1}}\oplus\mathcal{O}_{
\mathbb{C}P^{1}\times\mathbb{C}P^{1}}(1,-1))$ \textit{explicitly}
in symplectic coordinates.  It is our hope that analyzing the form
of known Kahler-Einstein toric metrics may lead to a method of
writing down all of Kahler-Einstein toric metrics.

    Recently, there has been a great deal of interest in finding
and classifying constant scalar curvature metrics on compact
Kahler manifolds.  For example, a recent paper by Donaldson
\cite{Donaldson} discusses the question of finding constant scalar
curvature metrics on Kahler toric manifolds and the relation to
notions of stability.  In the second part of Chapter $6$, we look
at the more general question of finding extremal Kahler metrics.
We begin by reviewing the construction of the extremal Kahler
metric conformal to the Hermitian-Einstein metric on
$\mathbb{C}P^{2}\sharp\overline{\mathbb{C}P^{2}}$ and giving a
novel interpretation of the symplectic potential of that extremal
Kahler metric.  We then find new \textit{explicit} extremal Kahler
metrics in dimension six. Specifically, we find a continuous
two-parameter family of extremal Kahler metrics on
$\mathbb{P}(\mathcal{O}_{\mathbb{C}P^{1}\times\mathbb{C}P^{1}}\oplus\mathcal{O}_{
\mathbb{C}P^{1}\times\mathbb{C}P^{1}}(1,-1))$.
Within this two-parameter family lie two one-parameter families of
constant scalar curvature metrics.  These two families intersect
at one point; this point corresponds to the Kahler-Einstein metric
found by Sakane.  Furthermore, this two-parameter family, along
with a one-parameter family of extremal Kahler metrics of
nonconstant scalar curvature on a $\mathbb{C}P^{2}$-bundle over
$\mathbb{C}P^{1}$, arise from a continuous deformation of the
Fubini-Study metric on $\mathbb{C}P^{3}$.

    Finally, we conclude Chapter $6$ with a discussion of the
connection, found by Derdzinski in \cite{Derd}, between
Hermitian-Einstein metrics and extremal Kahler metrics in
dimension four.  Using a result of LeBrun in \cite{LeBrun}, we
apply this connection to the Kahler toric case and exhibit the
equations which must be satisfied for a non-Kahler Hermitian toric
metric to be Einstein in dimension four.

    In addition to the $T^{3}$-symmetry present on any six-dimensional
Kahler toric manifold, all of the new extremal Kahler metrics we
have constructed admit additional symmetries; they are all in fact
of cohomogeneity-one under their full symmetry group.  It seems
natural, therefore, to look for explicit distinguished metrics on
Kahler manifolds with large symmetry group.  To begin such a
study, we look, in Chapter $7$, at Kahler manifolds of
cohomogeneity-$d$ under a compact semisimple group action.  We see
that the presence of a moment map and a complex structure requires
that each homogeneous principal orbit, $G/K$, becomes a Riemannian
submersion over a coadjoint orbit $G/L$ with fibers $L/K$ of
dimension $d$.  If we demand that each fiber be abelian (i.e.
isomorphic to $T^{d}$) and restrict to the monotypic case, we
recover a generalization of Kahler toric geometry which we refer
to as fiberwise Kahler toric geometry.  Metrics of this type can
be thought of as Kahler toric metrics with additional structure,
and describing a fiberwise Kahler toric metric is equivalent to
choosing a Kahler toric metric and specifying the values of a
collection of constants.

    Remarkably, metrics of the form described can be obtained by
starting with very different geometric assumptions.  In the second
construction of fiberwise Kahler toric metrics we begin with a
principal $T^{d}$-bundle over the product of Kahler-Einstein
manifolds endowed with metrics as described in \cite{WangZiller}.
Letting $\hat{g}$ denote a family of such metrics depending on $d$
variables, we construct a metric $g=\check{g}+\hat{g}$ where
$\check{g}$ is a metric on the $d$-dimensional space of variables.
Imposing a Kahler structure in a natural way, we retrieve metrics
of the same form with the same curvature equations as those
obtained by the previous construction; we also refer to these
metrics as fiberwise Kahler toric.  Indeed, the term fiberwise
Kahler toric will be used to describe the form of the metric and
not the particular construction.   We end Chapter $7$ by deriving
the equations of the Einstein and extremal Kahler conditions.

    Having defined fiberwise Kahler toric geometry, we specialize to
the $d=1$ case and begin our search for distinguished metrics. The
Kahler-Einstein condition on fiberwise Kahler toric manifolds in
the $d=1$ case has been studied in great detail (see
\cite{DandW}). However, to the author's knowledge, the extremal
Kahler and constant scalar curvature conditions have not been
studied on manifolds of this type.  In Chapter $8$, we demonstrate
how to integrate the scalar curvature equation on a fiberwise
Kahler toric manifold in the $d=1$ case.  This integration will
allow us to find extremal Kahler and constant scalar curvature
Kahler metrics. After discussing smoothness conditions, we
demonstrate that every fiberwise Kahler toric manifold admits at
least one extremal Kahler metric and most admit a multi-parameter
family of such metrics. This contrasts with the constant scalar
curvature case in which there are examples of fiberwise Kahler
toric metrics with $d=1$ which are known to admit no constant
scalar curvature metric even in real dimension four.  We are able
to prove the existence of a one-parameter family of constant
scalar curvature metrics on most (we will make this statement more
precise below) fiberwise Kahler toric metrics.  This method is
useful not only because it proves the existence of constant scalar
curvature Kahler metrics, but also because it allows us to write
the metrics down explicitly. Furthermore, it allows us to study
how, on a given fiberwise Kahler toric manifold, the constant
scalar curvature Kahler metrics sit within the larger family of
extremal Kahler metrics.  Conceivably, one could use this method
to determine all of the Kahler classes that admit constant scalar
curvature or extremal Kahler metrics and to study how those
distinguished classes sit in the space of all Kahler classes.
Unfortunately, we cannot answer all of these questions here.  We
do, however, go through many examples in detail.  After a review
of the four-dimensional case, we study the Kahler manifolds
$\mathbb{P}(\mathcal{O}_{\mathbb{C}P^{1}\times\mathbb{C}P^{1}}\oplus\mathcal{O}_{\mathbb{C}
P^{1}\times\mathbb{C}P^{1}}(p,q))$ and find extremal Kahler
metrics for all values of $p$ and $q$ in $\mathbb{Z}$.  Moreover,
we see that constant scalar curvature metrics exists when $p>0$
and $q<0$ (or vice versa) or when $p=q=0$ but not when both $p$
and $q$ are positive (or negative).

    The case of fiberwise Kahler toric manifolds with $d>1$ is more
difficult to analyze because the curvature conditions are partial,
instead of ordinary, differential equations. We cannot provide,
therefore, new explicit examples of distinguished metrics. We do,
however, integrate the equations of the Einstein condition to give
an analogue of the Monge-Ampere equations as we did in the Kahler
toric case.  Also, we demonstrate how to express the Futaki
invariant in the coordinates of fiberwise Kahler toric geometry.

    In Chapter $10$, we look for distinguished metrics which are
Hermitian but non-Kahler.  As in the second construction of $d=1$
fiberwise Kahler toric metrics, we study metrics of the form
$g=dt^{2}+\hat{g}$ where $\hat{g}$ is a family of metrics on a
principal $S^{1}$-bundle over the product of Kahler-Einstein
manifolds.  Again, we impose a natural integrable complex
structure on the manifold.  However, we no longer demand that the
metric be Kahler; instead, we demand that the Riemannian curvature
tensor be $J$-invariant, i.e.
$R(\cdot,\cdot,\cdot,\cdot)=R(J\cdot,J\cdot,J\cdot,J\cdot)$.  Note
that the $J$-invariance of the Riemannian curvature tensor implies
the $J$-invariance of the Ricci-tensor, i.e.
$r(\cdot,\cdot)=r(J\cdot,J\cdot)$, and that the Riemannian and
Ricci tensors of a Kahler manifold are automatically
$J$-invariant.  The Hermitian-Einstein condition on manifolds of
this type was studied in \cite{WangWang}.  In this chapter, we
study the constant scalar curvature condition and demonstrate,
after integrating the scalar curvature equation, a Hermitian
constant scalar curvature metric always exists when each principal
$S^{1}$-bundle is over a product of more than one Kahler-Einstein
manifold.  Therefore, most of the manifolds which do not admit
constant scalar curvature Kahler metrics will admit non-Kahler
constant scalar curvature Hermitian metrics with $J$-invariant
Riemannian curvature tensor.

    We place particular emphasis on the case in which the total
manifold is of dimension four.  In this case, the metric $\hat{g}$
is a metric on an $S^{1}$-bundle over $\mathbb{C}P^{1}$.  In the
compact case, the manifold is either $\mathbb{C}P^{2}$ or a
Hirzebruch surface
$\mathbb{P}(\mathcal{O}_{\mathbb{C}P^{1}}\oplus\mathcal{O}_{\mathbb{C}P^{1}}(q))$
where $q\in\mathbb{Z}$.  It is well-known that, for $q\neq0$, no
Hirzebruch surface admits a constant scalar curvature Kahler
metric.  Therefore, it seems natural to widen the search for
distinguished metrics. We propose that a natural generalization of
both the constant scalar curvature Kahler and Hermitian-Einstein
conditions, particularly in dimension four, is to look for
constant scalar curvature Hermitian metrics with $J$-invariant
Ricci tensor.  In dimension four, non-Kahler Hermitian metrics
with $J$-invariant Ricci tensor (but without necessarily being of
constant scalar curvature) have been studied extensively by
Apostolov and Gauduchon \cite{Apostolov} and others.  As such
metrics are plentiful (some manifolds admit infinite dimensional
families), one would like to impose addition restraints on the
curvature tensor.  An interesting additional constraint is the
constant scalar curvature condition.  To the author's knowledge,
the only known non-Kahler Hermitian metric of constant scalar
curvature with $J$-invariant Ricci tensor on a four manifold is
the Page metric on
$\mathbb{C}P^{2}\sharp\overline{\mathbb{C}P^{2}}$.  Here, we are
able to construct new one-parameter families of such metrics on
both $\mathbb{C}P^{1}\times\mathbb{C}P^{1}$ and
$\mathbb{C}P^{2}\sharp\overline{\mathbb{C}P^{2}}$.

    Both in the $d=1$ fiberwise Kahler toric and the analogous Hermitian
cases, we demonstrate how to integrate the scalar curvature
equation.  Chapter $11$ investigates when it is possible to
integrate the scalar curvature equation on a general
cohomogeneity-one manifold which may not have a compatible complex
structure.  More precisely, we study manifolds of
cohomogeneity-one under the action of a compact, connected group
$G$ with connected isotropy subgroup $K$.  Assuming that the
isotropy representation is composed of distinct $Ad(K)$-invariant
summands (the monotypic assumption) and that there exists a
summand of dimension one, we show how to integrate the scalar
curvature when the Ricci tensor is invariant under the action of a
natural rank-$2$ skew-symmetric map $\tilde{J}:TM\rightarrow TM$.
We note that in the Einstein case, the Ricci tensor \textit{must}
be invariant under this map.  The integration of the scalar
curvature equation in the Einstein case amounts to finding a
\textit{conserved quantity} in the Hamiltonian formulation of the
cohomogeneity-one Einstein equations.  The discovery of such an
integral may provide a first step in determining when the full
integration of the Einstein equations is possible.

\chapter{Submersions, Isometries, and Curvature}
    This chapter is devoted to reviewing the theory of Riemannian
submersions and deriving the curvature equations of a manifold
admitting an orthogonally transitive group action.  These
curvature equations will be used extensively in all of the
following chapters.

\section{Riemannanian Submersions}

    In our exposition of the basic principles of Riemannian submersions,
we follow closely the presentation found in \cite{Besse}.  Let
$(M,g)$ be a smooth Riemannian manifold without boundary and let
$(B, \check{g})$ be a Riemannian manifold with piecewise smooth
boundary.  A submersion, $\pi :M\rightarrow B$, defines a tangent
map $\pi_{*}:T_{x}M\rightarrow T_{b}B$ for each point $x$ of $M$
such that $\pi (x)=b$.  Denote by $\mathcal{V}$ the kernel of the
tangent map $\pi_{*}$, i.e. $\mathcal{V}$ is the distribution in
$TM$ such that, for any point $x$ in $M$, $\mathcal{V}_{x}$ is the
kernel of the map $\pi_{*}:T_{x}M\rightarrow T_{b}B$.  In this
thesis, the dimension of the distribution will be constant on all
points of $M$ which lie in the inverse image of the interior of
$B$.  However, at points in the inverse image of the boundary of
$B$, the dimension of $\mathcal{V}$ may decrease.  As $M$ is
endowed with a Riemannian metric, the orthogonal complement of
$\mathcal{V}$ in $TM$, which will be denoted $\mathcal{H}$, is
well-defined.  Just as in the case of the \textit{vertical}
distribution $\mathcal{V}$, the \textit{horizontal} distribution
$\mathcal{H}$ is of constant dimension in the inverse image of the
interior of $B$.  As a final piece of notation, we denote by
$\pi^{-1}({b})$ the inverse image of a point $b$ in $B$.  The
tangent space to each point $x$ in $\pi^{-1}({b})$ is simply the
vertical space $\mathcal{V}_{x}$.

    We say that $\pi$ is a \textbf{Riemannian submersion} if it
induces an isometry from $\mathcal{H}_{x}$ to $T_{b}M$ for each
point $x$ such that $\pi (x)=b$.  In this thesis, we frequently
make use of a special type of Riemannian submersion known as the
\textit{quotient by an isometric action}.  For ($M$, $g$) a
Riemannian manifold, let $G$ be a closed subgroup of the isometry
group of M.  If we set $B=M/G$ and let $\check{g}$ be the quotient
metric induced by $g$, then the projection from $M$ to $M/G$ is a
Riemannian submersion.  In this case, each fiber, $\pi^{-1}(b)$,
will be a homogeneous manifold isomorphic to $G/K$ where $K$ is
the isotropy subgroup of $G$.

    Given a Riemannian submersion, the various curvatures of $M$
can be expressed in terms of the curvature of the base space $B$,
the curvature of the spaces $\pi^{-1}(b)$, and two invariants
\cite{O'Neill}. The invariants, which following O'Neill we refer
to as $A$ and $T$, are defined in terms of the Levi-Civita
connections of $M$ and of the fibers $\pi^{-1}(b)$. Let $\nabla$
be the Levi-Civita connection of $g$ and let $\hat{\nabla}$ be the
collection of Levi-Civita connections associated to the metrics
$\hat{g}_{b}$ induced by $g$ on each fiber $\pi^{-1}(b)$.
Equivalently, one can define $\hat{\nabla}$ to be the vertical
projection of $\nabla$: for $X$ and $Y$ vector fields in $TM$, we
define the connection of the fiber to be $\hat{\nabla}$ such that
$\hat{\nabla}_{X}Y=\mathcal{V}\nabla_{X}Y$, where $\mathcal{V}X$
denotes the projection of the vector field $X$ onto the vertical
distribution (the horizontal projection is defined in an analogous
way).  We are now prepared to define the two invariants of a
Riemannian submersion.

\begin{defn}
For vector fields $E$ and $F$ on $M$, let $T$ be the (2,1)-tensor
field on $M$ such that
\begin{equation*}
T_{E}F=\mathcal{H}\nabla_{\mathcal{V}E}\mathcal{V}F+\mathcal{V}\nabla_{\mathcal{V}E}\mathcal{H}
F.
\end{equation*}

\end{defn}

The various identities of the invariant $T$ are collected in
\cite{Besse} and will not be reproduced here.  We note, however,
that $T$ can be thought of as the collection of the second
fundamental forms of the fibers $\pi^{-1}(b)$.  The vertical
distribution $\mathcal{V}$, since it is defined to be the kernel
of a linear map of the tangent space $TM$, is integrable.  This
implies that each of the fibers $\pi^{-1}(b)$ is a submanifold of
$M$. Therefore, one can define on each fiber the second
fundamental form and it is easily checked that $T$ restricted to a
fiber is the second fundamental form of that fiber.  Since the
second fundamental form of a Riemannian submanifold measures the
failure of that submanifold to be totally geodesic, the invariant
$T$ measures the failure of the fibers $\pi^{-1}(b)$ to be totally
geodesic.

    The second invariant associated to a Riemannian submersion is
defined in a similar way.

\begin{defn}
For vector fields $E$ and $F$ on $M$, let $A$ be the (2,1)-tensor
field on $M$ such that
\begin{equation*}
A_{E}F=\mathcal{H}\nabla_{\mathcal{H}E}\mathcal{V}F+\mathcal{V}\nabla_{\mathcal{H}E}\mathcal{H}
F.
\end{equation*}

\end{defn}

The following proposition, found in \cite{O'Neill}, provides an
interpretation of $A$.

\begin{prop}
For $H$ and $\bar{H}$ horizontal vector fields,
$A_{H}\bar{H}=\frac{1}{2}\mathcal{V}[H,\bar{H}]$.
\end{prop}
The tensor $A$ can be interpreted as measuring the failure of the
horizontal distribution to be integrable.  Although the horizontal
distribution is pointwise isometric to the tangent space of $B$,
unless $A$ vanishes, $\mathcal{H}$ does not define a submanifold
isometric to $B$.

    The curvature equations of a Riemannian submersion are
presented in detail in \cite{Besse}.  For the sake of brevity, we
will not reproduce the general curvatures equations of a
Riemannian submersion here.

    The Einstein condition on a generic Riemannian submersion is
in general too difficult to solve.  Therefore, additional
constraints must be placed upon the submersion to make the
equations more tractable.  As mentioned above, a common means of
simplification is to look at the special case of a quotient by an
isometric action.  Since the fibers of this type of submersion are
homogeneous manifolds, the terms representing their contribution
to the curvature of the manifold are greatly simplified.

    Either as an alternative or as a complementary means of
simplification, one could set one of the two invariants $A$ or $T$
equal to zero.  Note that setting them both equal to zero would
mean that the Riemannian submersion is a simple product of the
base and a fiber endowed with the product metric.  Setting the
invariant $T$ equal to zero, since $T$ is a generalization of the
second fundamental form, is equivalent to demanding that the
fibers be totally geodesic.  For a review of submersions with
totally geodesic fibers see \cite{Besse} and the article by M.
Wang in \cite{Survey}.  Setting the invariant $A$ equal to zero is
equivalent to demanding that the horizontal distribution be
integrable.  The condition that $A=0$ is referred to as
\textit{orthogonal transitivity} by relativists.  An orthogonally
transitive action is known as a \textit{surface-orthogonal} action
elsewhere in the literature.

  We will write down only the
equations for the Ricci curvature of a quotient by an orthogonally
transitive isometric action and refer the reader to \cite{O'Neill}
for their derivation.

\begin{prop}\label{curvature}
Let $g$ be the metric on $M$ (we will sometimes set
$g(\cdot,\cdot)=\langle\cdot,\cdot\rangle$), $\hat{g}$ the induced
metric on the fibers, and $\check{g}$ the induced metric on the
quotient space $M/G$ and let $r$, $\hat{r}$ and $\check{r}$ to be
the Ricci curvatures of those three metrics.  For $X$ and $Y$
vertical vector fields and $H$ and $\bar{H}$ horizontal vector
fields, the Ricci curvature of $M$ can be expressed by the
following equations when the tensor $A$ vanishes:
\begin{equation}\label{con1}
r(X,Y)=\hat{r}(X,Y)-\langle N,T_{X}Y\rangle+(\tilde{\delta}T)(X,Y)
\end{equation}
\begin{equation}\label{con2}
r(X,H)=\langle(\hat{\delta}T)X,H\rangle+\langle\nabla_{X}N,H\rangle
\end{equation}
\begin{equation}\label{con3}
r(H,\bar{H})=\check{r}(H,\bar{H})-\langle TH,
T\bar{H}\rangle+\frac{1}{2}(\langle
\nabla_{H}N,\bar{H}\rangle+\langle\nabla_{\bar{H}}N,H\rangle)
\end{equation}
where, for $\{Y_{j}\}$ an orthonormal basis of the vertical
distribution and $\{H_{j}\}$ an orthonormal basis for the
horizontal distribution, $N=\sum_{j}T_{Y_{j}}Y_{j}$,
$\tilde{\delta}T(X,Y)=\sum_{j}\langle(\nabla_{H_{j}}T)_{X}Y,H_{j}\rangle$,
$\hat{\delta}T=-\sum_{j}(\nabla_{Y_{j}}T)_{Y_{j}}$, and $\langle
TH, T\bar{H} \rangle=\sum_{j}\langle T_{Y_{j}}H,
T_{Y_{j}}\bar{H}\rangle$.
\end{prop}

\section{Orthogonally Transitive Isometries and the Einstein Condition}
We now restrict our attention to a special kind of Riemannian
submersion: the quotient by group action.  We begin with a
definition.
\begin{defn}
Let $(M^{n},g)$ be a Riemannian manifold such that $g$ is
invariant under the action of group $G$.  The manifold, $M$, is of
\textbf{cohomogeneity-d} if and only if the difference, $d$,
between the dimension of the manifold and the dimension of the
principal orbits $G/K$ is strictly less than $n$. (Note that while
$M$ is a smooth manifold, $M/G$ is not smooth globally.)
\end{defn}
Here we have required that the dimension of the quotient space be
strictly less than the dimension of the total manifold which
implies that the dimension of the symmetry group $G$ is at least
one.

    On manifolds of cohomogeneity-one the
Riemannian submersion invariant $A$ automatically vanishes since
the one dimensional quotient manifold $M/G$ must have zero
curvature. The Einstein condition on a cohomogeneity-one manifold
therefore becomes a system of differential equations in one
variable.  Going from the cohomogeneity-one case to manifolds with
group actions of higher codimension, the complexity of the
Einstein condition increases in three different ways. Firstly, the
invariant $A$ is no longer automatically equal to zero as the
horizontal distribution orthogonal to the principal orbits is not
necessarily integrable.  Secondly, the invariant $T$ depends on
more than one variable, turning the Einstein condition into a
system of partial differential equations.  Finally, the Ricci
tensor of the base manifold $M/G$ is no longer zero in general.
There may appear to be little hope in solving the Einstein
condition in even the cohomogeneity-two case. However, if one
considers a simplification of the higher cohomogeneity case, the
Einstein condition becomes more tractable, particularly in low
dimensions.

    In light of the vanishing of the invariant $A$ in the
cohomogeneity-one case, a natural simplification of the general
higher cohomogeneity case is to set the invariant $A$ equal to
zero.  Below, after presenting certain properties of general
higher cohomogeneity manifolds, we will restrict our attention to
the case in which $A$ vanishes.

\begin{defn}
For $(M,g)$ a cohomogeneity-$d$ manifold, $M$ is an
\textbf{orthogonally transitive \linebreak cohomogeneity-d
manifold} if and only if the invariant $A$ of the Riemannanian
submersion $\pi :M\rightarrow M/G$ vanishes.
\end{defn}

\subsection{Cohomogeneity-$d$ Einstein Manifolds and Topology}
Just as in the cohomogeneity one case, we have the following
proposition.

\begin{prop}
For $M$ an Einstein manifold of cohomogeneity-$d$, if the scalar
curvature of $M$ is negative then $M$ is noncompact.
\end{prop}
\begin{proof}
If $M$ is compact with negative Ricci curvature, Bochner's theorem
states that all Killing vectors must vanish identically which
would imply that the principal orbits are of dimension zero.
\end{proof}

\begin{prop}
For $M$ a compact Einstein manifold of cohomogeneity-two or
cohomogeneity-three, the scalar curvature is positive or the
metric is flat.
\end{prop}
\begin{proof}
In light of the last proposition we need only consider the case in
which the Einstein constant is equal to zero, i.e. $r=0$.  By
Bochner's theorem, $s=0$ implies that all Killing vectors are
parallel \cite{Petersen} . It is easily seen that this implies the
vanishing of the invariants $A$ and $T$.  By equations $9.36(a-c)$
of \cite{Besse} and equations (\ref{con1})-(\ref{con3}), this
implies that $\check{r}=\hat{r}=0$.  As each principal orbit is
homogeneous, $\hat{r}=0$ implies that the Riemannian curvature of
the principal orbit is zero.  Finally, since $M/G$ is of dimension
two or three, $\check{r}=0$ implies that $\check{R}=0$.
\end{proof}
By Myer's Theorem, a compact Einstein manifold with positive
scalar curvature must have finite fundamental group.  Therefore,
we have the following rough classification of cohomogeneity-two
and -three Einstein manifolds.

  For $M$ a complete cohomogeneity-two or cohomogeneity-three
Einstein manifold, one of the following holds:
\begin{itemize}
\item $M$ is compact with positive scalar curvature and finite
fundamental group

\item $M$ is compact and flat

\item $M$ is noncompact and has nonpositive scalar curvature.

\end{itemize}

    This classification fails for manifolds of cohomogeneity
greater than three (but which still have a symmetry group of
dimension at least one), because on a manifold of dimension four
or greater Ricci flatness is not equivalent to flatness.  The
simplest counter-example to extending the above classification to
manifolds of higher cohomogeneity with symmetry group of at least
one dimension is $S^{1}\times K3$-surface. Endowed with its
standard metric, this is a cohomogeneity-four Einstein manifold
with zero scalar curvature which is not flat.
\subsection{The Einstein Equations in the Orthogonally Transitive Case}

    Let $M$ be a cohomogeneity-$d$ manifold under the orthogonally transitive action of compact
group $G$ with isotropy subgroup $K$.  The quotient by the group
action yields the Riemannian submersion $\pi:M\rightarrow M/G$
where $\textrm{dim} M/G=d$ as described above.  A principal orbit,
$\pi^{-1}(b)$, is the inverse image of a point $b$ in the interior
of the space $M/G$.  Each principal orbit is isomorphic to $G/K$
and of real codimension-$d$.  Moreover, the union of the principal
orbits forms an open dense set in $M$.  On this set, the metric
$g$, because it is invariant under an orthogonally transitive
group action, can be written as
\begin{equation*}
g=\check{g}+\hat{g}
\end{equation*}
where $\check{g}$ is the metric induced on the quotient space
$M/G$ and $\hat{g}$ is a $d$-parameter family of $G$-invariant
metrics on the principal orbit.

    To express the equations of the Einstein
condition in the orthogonally transitive case, we set
$\{X_{i}\}_{i=1,...,n-d}$ to be an invariant basis for each
principal orbits and set
\begin{equation*}
h_{ij}=g(X_{i},X_{j})=\hat{g}(X_{i},X_{j})=X^{a}_{i}X_{ja}.
\end{equation*}

For the moment, we will say no more about $\hat{g}$ other than to
note that, depending on the isotropy representation, some of the
metric elements $h_{ij}$ will be zero automatically.  We will make
this precise in the following section.

Let $\{H_{i}\}_{i=1,...,d}$ be an \textit{orthonormal} basis for
the horizontal distribution and hence for the tangent space to
$M/G$.  We are now ready to determine the equations of the
Einstein condition.
    The first equation we want to study is

\begin{equation*}
r(X_{i},X_{j})=\hat{r}(X_{i},X_{j})-\langle
N,T_{X_{i}}Y_{i}\rangle+(\tilde{\delta}T)(X_{i},X_{j})=\lambda
g(X_{i},X_{j})=\lambda h_{ij}
\end{equation*}
where $g(\cdot,\cdot)=\langle\cdot,\cdot\rangle$ as above.

By the definition of $T$
\begin{equation*}
T_{X_{i}}X_{j}=\mathcal{H}(X_{i}^{a}\nabla_{a}X_{jb}).
\end{equation*}
Note that we can choose the $X_{j}$ and the $H_{i}$ so that
$[H_{i},X_{j}]=0$.  Using this fact and the properties of tensor
$T$, we have
\begin{equation*}
\langle T_{X_{i}}X_{j},H_{k}\rangle=-\frac{1}{2}H_{a}\nabla^{a}
\langle X_{i},X_{j}\rangle.
\end{equation*}
Therefore, we have demonstrated that
\begin{equation}
T_{X_{i}}X_{j}=-\frac{1}{2}\check{\nabla}_{b}h_{ij}.
\end{equation}

    The vector $N$, referred to as the \textit{mean curvature
vector} is defined to be the trace of the tensor $T$.  That is,

\begin{equation}
N=h^{ij}(T_{X_{i}}X_{j})=-\frac{1}{2}h^{ij}\check{\nabla}_{b}h_{ij}=-\frac{1}{2}\frac{\check{
\nabla}_{b}(\textrm{det}(h))}{\textrm{det}(h)}.
\end{equation}

Immediately, we can calculate one term of the Einstein condition,
\begin{equation*}
\langle N, T_{X_{i}}X_{j}\rangle =\frac{1}{4}(\textrm{det}(h)
)^{-1}\check{\nabla}^{a}\textrm{det}(h)\check{\nabla}_{a}h_{ij}=\frac{1}{4}\frac{\langle\check{\nabla}\textrm{det}
(h),\check{\nabla}h_{ij} \rangle}{\textrm{det} (h)}.
\end{equation*}

In order to find $R_{ab}X_{i}^{a}X_{j}^{b}=r(X_{i},X_{j})$ we must
next calculate $\tilde{\delta}T$.  Recalling the definition of
$\tilde{\delta}T(X_{i}, X_{j})$ and using the fact that $T$ is a
tensor field, the action of the Levi-Civita connection on tensors
gives
\begin{equation} \label{deltaT}
\tilde{\delta}T(X_{i},
X_{j})=\sum_{k}\{\langle\nabla_{H_{k}}(T_{X_{i}}X_{j}),
H_{k}\rangle-\langle T_{\nabla_{H_{k}}X_{i}}X_{j},
H_{k}\rangle-\langle T_{\nabla_{H_{k}}X_{j}}X_{i}, H_{k}\rangle\}.
\end{equation}
Having expressed $T_{X_{i}}X_{j}$ in terms of $h_{ij}$, it is
straightforward to calculate (\ref{deltaT}).  First, we note that,
for any horizontal vector $H$,
$\nabla_{H}(T_{X_{i}}X_{j})=-\frac{1}{2}H^{b}\nabla_{b}\nabla_{a}h_{ij}$
and
\begin{equation*}
\sum_{k}\langle\nabla_{H_{k}}(T_{X_{i}}X_{j}),
H_{k}\rangle=-\frac{1}{2}\check{\nabla}^{a}\check{\nabla}_{a}h_{ij}.
\end{equation*}

The other elements of $\tilde{\delta}T$ are somewhat more
difficult to compute.  Recall that the $X_{i}$ were chosen so that
they commute with every horizontal vector field $H$, i.e.
$\nabla_{H}X_{i}=\nabla_{X_{i}}H$.  Furthermore, we note that
$\nabla_{H}X_{i}=T_{X_{i}}H$ since $A$ vanishes.  Combining these
results, we have
\begin{equation*}
-\langle T_{\nabla_{H_{k}}X_{i}}X_{j}, H_{k}\rangle-\langle
T_{\nabla_{H_{k}}X_{j}}X_{i}, H_{k}\rangle=\langle T_{X_{j}}H_{k},
\nabla_{H_{k}}X_{i}\rangle+\langle T_{X_{i}}H_{k},
\nabla_{H_{k}}X_{j}\rangle=2\langle\nabla_{H_{k}}X_{i},
\nabla_{H_{k}}X_{j}\rangle.
\end{equation*}
From \cite{Mason} we have
\begin{equation*}
\nabla_{H}X_{i}=H^{a}\nabla_{a}X_{ib}=\frac{1}{2}h^{mk}(H^{a}\check{\nabla}_{a}h_{ki})X_{mb}
\end{equation*}
and therefore
\begin{equation*}
\sum_{k=1}^{d}\langle\nabla_{H_{k}}X_{i},
\nabla_{H_{k}}X_{j}\rangle=\frac{1}{4}h^{kl}\check{\nabla}^{a}h_{ki}\check{\nabla}_{a}h_{lj}=
-\frac{1}{4}h_{ik}\check{\nabla}^{a}h^{kl}\check{\nabla}_{a}h_{lj}.
\end{equation*}
The above results prove that
\begin{equation*}
\tilde{\delta}T(X_{i},
X_{j})=-\frac{1}{2}\check{\nabla}^{a}\check{\nabla}_{a}h_{ij}+\frac{1}{2}h^{kl}\check{
\nabla}^{a}h_{ki}\check{\nabla}_{a}h_{lj}.
\end{equation*}
Finally we find that
\begin{equation}\label{RXX}
r(X_{i},X_{j})=-\frac{1}{2}\check{\nabla}^{a}\check{\nabla}_{a}h_{ij}-\frac{1}{4}(\textrm{det}
(h))^{-1}\check{\nabla}^{a}\textrm{det}
(h)\check{\nabla}_{a}h_{ij}+\frac{1}{2}h^{kl}\check{\nabla}^{a}h_{ki}\check{\nabla}_{a}h_{lj}
+\hat{r}(X_{i},X_{j})
\end{equation}
or, equivalently,
\begin{equation}\label{RXX2}
r(X_{i},X_{j})=-\frac{1}{2}h_{ik}(\textrm{det}
(h))^{-\frac{1}{2}}\check{\nabla}^{a}((\textrm{det}
(h))^{\frac{1}{2}}h^{kl}\check{\nabla}_{a}h_{ij})+\hat{r}(X_{i},X_{j}).
\end{equation}

    We now turn to the second equation of the Einstein condition
$r(X,H)=0$ where $X=X_{i}$ and $H=H_{j}$ for some $i$ and $j$.
Because the O'Neill tensor $A$ vanishes, we have that
\begin{equation*}
r(X,H)=\langle
(\hat{\delta}T)X,H\rangle+\langle\nabla_{X}N,H\rangle.
\end{equation*}
From \cite{Besse}, we see that
\begin{equation*}
\hat{\delta}T=-h^{ij}(\nabla_{X_{i}}T)_{X_{j}}.
\end{equation*}
We consider first the term $\langle\nabla_{X}N,H\rangle$.  Noting
the $N$ is $G$-invariant we see that
\begin{equation*}
\check{\nabla}_{X}\langle N,H\rangle =0\Rightarrow
\langle\nabla_{X}N,H\rangle +\langle N,\nabla_{X}H\rangle =0.
\end{equation*}
Because we have chosen $X$ to be a Killing vector field commuting
with $H$, we have that
\begin{equation*}
\mathcal{H}\nabla_{X}H=\mathcal{H}\nabla_{H}X=A_{H}X=0.
\end{equation*}
Therefore,
\begin{equation*}
\check{\nabla}_{X}\langle N,H\rangle=0
\end{equation*}
automatically.

    The term $\langle (\hat{\delta}T)X,H\rangle$ is more difficult
to analyze and we follow \cite{Bergery} closely.  In
\cite{Bergery}, $H=\frac{\partial}{\partial t}$ as they are
considering the cohomogeneity-one case.  However, we may take $H$
to be any horizontal vector field without affecting the
calculation. We have that
\begin{equation*}
r(X,H)=-h^{ij}\langle(\nabla_{X_{i}}T)_{X_{j}}X,H\rangle=h^{ij}\langle(\nabla_{X_{i}}T)_{X_{j}}
H,X\rangle
\end{equation*}
\begin{equation*}
=h^{ij}\langle\hat\nabla_{X_{i}}(T_{X_{j}}H),X\rangle-h^{ij}\langle
T_{\hat{\nabla}_{X_{i}}X_{j}}H,X\rangle.
\end{equation*}
From the symmetries of $T$, we see that
\begin{equation*}
h^{ij}\langle
T_{\hat{\nabla}_{X_{i}}X_{j}}H,X\rangle=-h^{ij}\langle
T_{\hat{\nabla}_{X_{j}}X_{i}}H,X\rangle.
\end{equation*}
    For all Killing vectors, $Z$, we know by a well-known formula that
\begin{equation}\label{Killingequation}
\langle \hat{\nabla}_{X_{i}}X_{j},Z\rangle=\frac{1}{2}(\langle
[X_{i},X_{j}],Z\rangle+\langle [X_{i},Z],X_{j}\rangle+\langle
[X_{j},Z],X_{i}\rangle).
\end{equation}

This implies (see \cite{Bergery}) that
\begin{equation*}
-h^{ij}\langle
T_{\hat{\nabla}_{X_{i}}X_{j}}H,X\rangle=-\textrm{tr}(\textrm{ad}_{\mathfrak{g}}(T_{H}X)).
\end{equation*}
This term is equal to zero if and only if $G$ is unimodular. Since
we have assumed that $G$ is compact and compact groups are
unimodular, we deduce that
\begin{equation*}
r(X,H)=h^{ij}\langle\hat\nabla_{X_{i}}(T_{X_{j}}H),X\rangle
\end{equation*}
\begin{equation*}
=h^{ij}X_{j}\langle
T_{X_{i}}H,X\rangle-h^{ij}\langle\hat{\nabla}_{X_{i}}X,T_{X_{j}}H\rangle.
\end{equation*}
The second term is equal to zero because $T$ is symmetric and
$\hat{\nabla}X$ is anti-symmetric because $X$ is a Killing vector
field \cite{Bergery}.  The first term becomes
\begin{equation*}
h^{ij}X_{j}\langle
T_{X_{i}}H,X\rangle=\frac{1}{2}h^{ij}HX_{i}\langle
X_{i},X_{j}\rangle=\frac{1}{2}(h^{ij}H\langle
\hat{\nabla}_{X_{i}}X_{j},X\rangle+h^{ij}H\langle
\hat{\nabla}_{X_{i}}X,X_{j}\rangle).
\end{equation*}
The second term vanishes because $\hat{\nabla}X$ is
anti-symmetric.  The first term can be simplified using equation
(\ref{Killingequation}) to give
\begin{equation}
r(X,H)=\frac{1}{2}h^{ij}H\langle X_{i},[X_{j},X]\rangle.
\end{equation}

Since this equation holds for all $X=X_{k}$ and $H$ horizontal,
equation $r(X,H)=0$ of the Einstein condition implies that
\begin{equation}
h^{ij}\hat{\nabla}\langle X_{i},[X_{j},X_{k}]\rangle=0
\end{equation}
for all $X_{k}$.
    Aside from the Ricci curvature of the base manifold, the only
other terms in the Ricci tensor restricted to the horizontal
distribution are $\langle TH, T\bar{H}\rangle$ and
$\langle\nabla_{H}N, \bar{H}\rangle+\langle\nabla_{\bar{H}}N,
H\rangle$ where $H$ and $\bar{H}$ are again arbitrary horizontal
vector fields. Recalling the
$N=-\frac{1}{2}h^{ij}\nabla_{b}h_{ij}$, a straightforward
calculation shows that
\begin{equation*}
\langle\nabla_{H}N,
\bar{H}\rangle=-\frac{1}{2}h^{ij}H^{a}\bar{H}^{b}\nabla_{a}\nabla_{b}h_{ij}-\frac{1}{2}H^{a}
\nabla_{a}h^{ij}\bar{H}^{b}\nabla_{b}h_{ij}.
\end{equation*}
Because both the matrix $J_{ij}$ and the Hessian
$\nabla_{a}\nabla_{b}$ are symmetric, $\langle\nabla_{H}N,
\bar{H}\rangle=\langle\nabla_{\bar{H}}N, H\rangle$.  Due to the
vanishing of $A$ and the choice of killing vectors and for
$(Y_{k})_{k=1,...,n-d}$ an orthonormal basis of the vertical
distribution,
\begin{equation*}
\langle TH, T\bar{H} \rangle=\sum_{k}\langle T_{Y_{k}}H,
T_{Y_{k}}\bar{H}\rangle=\sum_{k}\langle\nabla_{H}Y_{k},
\nabla_{\bar{H}}Y_{k}\rangle=-\frac{1}{4}H^{a}\nabla_{a}h^{ij}\bar{H}^{b}\nabla_{b}h_{ij}.
\end{equation*}
Finally, we can express the Ricci curvature restricted to the
horizontal distribution as
\begin{equation*}
r_{ab}|_{\mathcal{H}}=\check{r}_{ab}-\frac{1}{2}h^{ij}\check{\nabla}_{a}\check{\nabla}_{b}
h_{ij}-\frac{1}{4}\check{\nabla}_{a}h^{ij}\check{\nabla}_{b}h_{ij},
\end{equation*}
where $\check{r}$ is the Ricci curvature of the quotient manifold
$M/G$ endowed with its induced metric $\check{g}$.

    With the above calculations, we have proven the following
theorem.
\begin{theorem}\label{cohomdEinsteincondition}
Let $(M,g)$ be a cohomogeneity-$d$ manifold under the orthogonally
transitive action of compact connected group $G$.  Let
$\{X_{i}\}_{i=1,...,n-d}$ be an invariant basis for the vertical
distribution where $h_{ij}=g(X_{i},X_{j})$ and let $\check{g}$ be
the metric induced on the quotient space of $M/G$.  The Einstein
condition on such a manifold is equivalent to

\begin{itemize}
\item $-\frac{1}{2}\check{\triangle}h_{ij}-\frac{1}{4}\frac{\langle\check{\nabla}{\rm det} (h),
\check{\nabla}h_{ij}\rangle}{{\rm det}
(h)}-\frac{1}{2}h_{ik}\langle\check{\nabla}h^{kl},\check{
\nabla}h_{lj}\rangle+\hat{r}(X_{i},X_{j})=\lambda h_{ij}$, for all
$1\leq i,j\leq n-d$,
\item $h^{ij}\check{\nabla}\langle X_{i},[X_{j},X_{k}]\rangle=0$,
for all $1\leq k\leq n-d$, and
\item
$\check{r}_{ab}-\frac{1}{2}h^{ij}\check{\nabla}_{a}\check{\nabla}_{b}h_{ij}-\frac{1}{4}\check{
\nabla}_{a}h^{ij}\check{\nabla}_{b}h_{ij}=\lambda\check{g}_{ab}$.
\end{itemize}
\end{theorem}
For now, we will not make an effort to express $\hat{r}$ more
explicitly.  In the sequel, we will evaluate this term on a case
by case basis.  For a more general discussion in the
cohomogeneity-one case, which readily extends to the orthogonally
transitive cohomogeneity-$d$ case, we refer the reader to
\cite{WangZiller}.

\subsection{Einstein Equations in the Monotypic Case}

    In this section, we give a more thorough discussion of the form of the metric $\hat{g}$.
We then restrict our attention to the monotypic case in which the
metric and the equations of the Einstein condition take a
particularly simple form.

    We again take $(M,g)$ to be an orthogonally transitive cohomogeneity-$d$ manifold under
the action of compact connected group $G$.  Let $K$ be the
isotropy subgroup of $G$, i.e.
\begin{equation*}
K=\{g\in G:g(x)=x\}
\end{equation*}
for some $x\in\pi^{-1}(b)$ with $b$ in the interior of the
quotient space $M/G$.  As mentioned in the previous section, the
metric can be written as $g=\check{g}+\hat{g}$ where $\hat{g}$ is
a $d$-parameter family of $G$-invariant homogeneous metrics on
$G/K\cong\pi^{-1}(b)$.  To describe the space of such homogeneous
metrics, we identify the tangent space to $G/K$ with an
$\textrm{Ad}(K)$-invariant subspace of the Lie algebra of $G$.

    Let $\mathfrak{p}$ be an $Ad(K)$-invariant subspace of $\mathfrak{g}$ such that
\begin{equation*}
\mathfrak{g}=\mathfrak{k}\oplus\mathfrak{p}.
\end{equation*}
We choose once and for all a space $\mathfrak{p}$. though we note
that the subspace $\mathfrak{p}$ is not unique.  That subspace can
be identified with $T_{x}(G/K)$ for $x\in G/K$ \cite{Besse}. Under
this identification, choosing a $G$-invariant metric on $G/K$ is
equivalent to choosing an $\textrm{Ad}(K)$-invariant scalar
product on $\mathfrak{p}$.

    One can further decompose $\mathfrak{p}$ into $\textrm{Ad}(K)$-invariant summands
\begin{equation*}
\mathfrak{p}=\mathfrak{p}_{1}\oplus...\oplus\mathfrak{p}_{m}
\end{equation*}
where each $\mathfrak{p}_{i}$ is irreducible.  Set
$\textrm{dim}\mathfrak{p}_{i}=d_{i}$ for all $1\leq i\leq m$.  We
fix once and for all a bi-invariant metric $Q$ on $\mathfrak{g}$
such that $Q(\mathfrak{k},\mathfrak{p})=0$ and
$Q(\mathfrak{p}_{i},\mathfrak{p}_{j})=0$ for $i\neq j$. Any other
$G$-invariant metric on $G/K$ can be written as $Q(L\cdot,\cdot)$
where $L$ is a positive-definite symmetric
$\textrm{Ad}(K)$-equivariant linear map \cite{Kerr}.  That is,
\begin{equation*}
L:\mathfrak{p}\rightarrow\mathfrak{p}.
\end{equation*}
The space of $\textrm{Ad}(K)$-invariant scalar products is
parameterized by the matrices $L$.  Finally we let
$\{Y^{i}_{j}\}_{j=1,...,d_{i}}$ be an orthonormal basis for
$\mathfrak{p}_{i}$ relative to the metric $Q$.

    We now specialize to the monotypic case.
\begin{defn}
The isotropy representation is called \textbf{monotypic} if all of
the summands, $\mathfrak{p}_{i}$, are pairwise inequivalent.
\end{defn}

By Schur's lemma, we see that any $\textrm{Ad}(K)$-equivariant
linear map, $L$, leaves each of the $\mathfrak{p}_{i}$ invariant
in the monotypic case \cite{Bergery}.  Therefore, the
$d$-parameter family of $G$-invariant metrics $\hat{g}$ can be
written as

\begin{equation}\label{fudge}
\hat{g}(\cdot,\cdot)=\sum_{i=1}^{m}A_{i}Q(\cdot,\cdot)|_{\mathfrak{p}_{i}}
\end{equation}
where the $A_{i}$ are positive functions on $M/G$.  The resulting
cohomogeneity-$d$ metric preserves the bi-invariance in the
monotypic case. This gives the following proposition:
\begin{prop}\label{p26}
In the monotypic case, for $Y^{i}_{j}$ a vertical vector field and
$H$ a horizontal vector field on an orthogonally transitive
cohomogeneity-$d$ manifold $M$, $r(Y^{i}_{j},H)=0$.
\end{prop}
\begin{proof}
The proof follows from the proof of the analogous statement in the
cohomogeneity one case \cite{Bergery}.  To see this we note that
the second equation in (\ref{cohomdEinsteincondition}) becomes
\begin{equation*}
\sum_{i=1}^{m}\sum_{j=1}^{d_{i}}A_{i}^{-1}\check{\nabla}\langle
Y^{i}_{j},[Y^{i}_{j},Y^{k}_{l}]\rangle=0
\end{equation*}
for all $Y^{k}_{l}$.  Because $Q$ is bi-invariant, equation
(\ref{fudge}) implies that the above equation vanishes
automatically.

\end{proof}

A straightforward calculation using the equations of theorem
(\ref{cohomdEinsteincondition}) shows that for
$(Y^{1}_{1},...,Y^{m}_{d_{m}})$ a basis for the principal orbits,
the Einstein condition on an orthogonally transitive
cohomogeneity-$d$ manifold in the monotypic case becomes

\begin{equation}\label{cee1}
\hat{r}(Y^{k}_{i},
Y^{k}_{i})-\frac{1}{2}\check{\triangle}A_{i}+\frac{1}{2}\frac{\langle
\check{\nabla}A_{i},
\check{\nabla}A_{i}\rangle}{A_{i}}-\frac{1}{4}\sum_{j=1}^{m}d_{j}\frac{\langle
\check{\nabla}A_{i}, \check{\nabla}A_{j}\rangle }{A_{j}}=\lambda
A_{i}
\end{equation}
\begin{equation}\label{cee2}
\check{r}_{ab}-\frac{1}{2}\sum_{j=1}^{m}d_{j}\frac{\check{\nabla}_{a}\check{\nabla}_{b}A_{j}}
{A_{j}}+\frac{1}{4}\sum_{j=1}^{m}d_{j}\frac{\check{\nabla}_{a}A_{j}
\check{\nabla}_{b}A_{j}}{A_{j}^{2}}=\lambda\check{g}_{ab}
\end{equation}
for all $1\leq i\leq m$ and $1\leq k\leq d_{i}$ where $\hat{r}$ is
the Ricci curvature of each principal orbit, $\check{r}$ is the
Ricci curvature of the base space,
$\check{\nabla}_{a}\check{\nabla}_{b}$ is the Hessian with respect
to the base,
$\check{\triangle}=\check{\nabla}^{a}\check{\nabla}_{a}$ is the
Laplacian of $\check{g}$, and $\lambda$ is the Einstein constant
of the manifold $M$. If $(n-d)=\sum_{i=1}^{r}d_{i}$ then the trace
of (\ref{cee1}) can be written as
\begin{equation*}
-\frac{\check{\triangle}\left(A_{1}^{\frac{d_{1}}{2}}\cdot...A_{m}^{\frac{d_{m}}{2}}\right)}
{A_{1}^{\frac{d_{1}}{2}}\cdot...A_{m}^{\frac{d_{m}}{2}}}=(n-d)\lambda-\hat{s}.
\end{equation*}
This equation, or rather its non-monotypic version, will be very
useful when we consider torus principal orbits since $\hat{r}=0$
in that case.  Although torus principal orbits are not contained
in the monotypic case, it can easily be seen from the calculations
in \cite{Bergery} that $r(X,H)=0$ in that case as well.

    Finally, the scalar curvature, $S$, of metric $g$, can be written as
\begin{equation*}
S=\hat{S}+\check{S}+\sum_{j=1}^{m}d_{j}\frac{\check{\triangle}A_{j}}{A_{j}}+\frac{3}{4}
\sum_{j=1}^{m}d_{j}\frac{\langle\check{\nabla}A_{j},
\check{\nabla}A_{j}\rangle}{A_{j}^{2}}-\frac{1}{4}\sum_{j=1}^{m}\sum_{k=1}^{m}\frac{\langle
\check{\nabla}A_{j},
\check{\nabla}A_{k}\rangle}{A_{j}A_{k}}
\end{equation*}
where $\hat{S}$ denotes the scalar curvature of the fibers and
$\check{S}$ denotes the scalar curvature of the quotient metric
$\check{g}$.

\chapter{Diagonalizability of Cohomogeneity-One Einstein Metrics}

    On a manifold of cohomogeneity-one under the action of a group
$G$, we saw above that a metric is automatically diagonal in the
monotypic case.  However, when two or more of the irreducible
summands in the isotropy representation are equivalent, the metric
need not be diagonal.  The simplest example is a cohomogeneity-one
metric whose principal orbits are copies of $T^{n}$.  In this
case, the metric is 'maximally' non-diagonalizable in the sense
that all the non-diagonal entries of the fiber portion of a smooth
metric may be non-zero.  While a 'generic' cohomogeneity-one
metric is not automatically diagonalizable in the non-monotypic
case, it is known that, in certain situations, the equations of
the Einstein condition allow the metric to be diagonalized
globally.  More precisely, in these cases it is one equation in
particular of the Einstein condition which allows the metric to be
globally diagonalized.

~

\textbf{The Equation:}

~

Let $(M,g)$ be an $n$-dimensional manifold of cohomogeneity-one
under the action of a compact connected group $G$.  In the
cohomogeneity-one case, the quotient space, $M/G$, is
one-dimensional which implies that the horizontal distribution of
the Riemannian submersion $\pi:M\rightarrow M/G$ is also of
dimension one.  As a one-dimensional distribution is automatically
integrable, the O'Neill tensor $A$ vanishes. Since $A\equiv0$, the
O'Neill formulae imply that
\begin{equation*}
r\left(X,\frac{\partial}{\partial t}\right)=\left\langle(\delta
T)X,\frac{\partial}{\partial
t}\right\rangle+\left\langle\nabla_{X}N,\frac{\partial}{\partial
t}\right\rangle
\end{equation*}
where $X$ is any vertical vector field and
$\frac{\partial}{\partial t}$ is the lift of the unit length
vector field on the one-dimensional quotient space $M/G$.

Because the action of the group $G$ is orthogonally transitive the
metric can be written as
\begin{equation*}
g=dt^{2}+\hat{g}=dt^{2}+g_{t}.
\end{equation*}

    As above, let $\{X_{i}\}_{i=1,...,n-1}$ be an invariant basis of
Killing vector fields for each of the principal orbits
$\pi^{-1}(b)\cong G/K$, and let $h_{ij}=\langle
X_{i},X_{j}\rangle=\hat{g}(X_{i},X_{j})$.  Equation (\ref{fudge})
tells us that
\begin{equation*}
r\left(X,\frac{\partial}{\partial
t}\right)=\frac{1}{2}h^{ij}\frac{\partial}{\partial t}\langle
X_{i},[X_{j},X]\rangle
\end{equation*}
for all vertical vector fields $X$.  Because the vertical and
horizontal distributions are orthogonal, the Einstein condition
requires that $r\left(X,\frac{\partial}{\partial t}\right)=0$ for
all $X$.  In the Einstein case, therefore, the above equation
becomes
\begin{equation}\label{diag}
h^{ij}\frac{\partial}{\partial t}\langle X_{i},[X_{j},X]\rangle=0
\end{equation}
    Note that equation (\ref{diag}) is a necessary condition for a
metric to be Einstein but it is far from sufficient.

    Recall from proposition \ref{p26}
that equation (\ref{diag}) is automatically satisfied in the
monotypic case.  This equation is also automatically satisfied in
the abelian case in which the principal orbits are copies of
$T^{n}$.  However, in the non-monotypic case, there are many
principal orbit types in which equation (\ref{diag}) is not
automatically satisfied.

    For certain principal orbit types in the non-monotypic case, it
has been noted that the Einstein condition (specifically equation
(\ref{diag}) of the Einstein condition) forces metrics to be
diagonal. Thus far only a few individual examples of this
phenomenon are known and, to the author's knowledge, no general
study has been made of when the Einstein condition can, or cannot,
be used to diagonalize a metric.  In this chapter, we will review
in detail an example of a cohomogeneity-one Einstein four manifold
well known to relativists and will treat that example in the more
mathematical formalism due to Berard-Bergery.  Then, we will begin
a more general study of the diagonalizability of Einstein metrics
in the non-monotypic case.

\section{Einstein Four Manifolds with SU(2) Principal Orbits}

A Bianchi IX metric is a cohomogeneity-one metric on a four
manifold in which the principal orbits are copies of $SU(2)$. That
is, the group is $G=SU(2)$, the isotropy subgroup is $K=\{Id\}$,
and the Lie algebra $\mathfrak{g}$ decomposes as
\begin{equation*}
\mathfrak{g}=\mathfrak{p}_{1}\oplus\mathfrak{p}_{2}\oplus\mathfrak{p}_{3}
\end{equation*}
where each $\mathfrak{p}_{i}$ is a trivial one-dimensional
subspace of $\mathfrak{g}$. Take $X_{i}\in \mathfrak{p}_{i}$ to be
a basis for each principal orbit satisfying the standard Lie
algebra relations for $\mathfrak{su}(2)$, namely
\begin{equation*}
[X_{1},X_{2}]=X_{3},~~[X_{2},X_{3}]=X_{1},~~[X_{3},X_{1}]=X_{2}.
\end{equation*}
Defining $h_{ij}=\langle X_{i},X_{j}\rangle$ as above we see that
\begin{equation*}
{g}_{t}=(h_{ij}).
\end{equation*}
 Now, since all of the $\mathfrak{p}_{i}$ are equivalent one-dimensional trivial subspaces,
none of the off-diagonal entries in $(h_{ij})$ need be zero
automatically.

Equation (\ref{diag}) becomes the system

\begin{equation*}
h^{ij} \frac{\partial}{\partial t}\langle
X_{i},[X_{j},X_{k}]\rangle=0
\end{equation*}
for all $k$.

    We now make the important observation that given a point $x$ of
the quotient space one can, without loss of generality,
diagonalize the metric at that point.  That is,

\begin{equation*}
h_{ij}\mid_{t=x}=h^{ij}\mid_{t=x}=0
\end{equation*}
for $i\neq j$. Given this condition, equation (\ref{diag}) becomes
the following system of three equations

\begin{equation}
(h^{33}-h^{22})\dot{h}_{23}\mid_{x}=0
\end{equation}
\begin{equation}
(h^{22}-h^{11})\dot{h}_{12}\mid_{x}=0
\end{equation}
\begin{equation}
(h^{11}-h^{33})\dot{h}_{13}\mid_{x}=0
\end{equation}
where $\dot{h}_{ij}=\frac{\partial}{\partial t}h_{ij}$. If we
require that $SU(2)$ is the full isometry group (e.g. the manifold
does not have $U(2)$ symmetry) we can further assume without loss
of generality that the $h_{ii}$ are distinct. Equation
(\ref{diag}) then implies that
\begin{equation}
\dot{h}_{ij}\mid_{x}=0
\end{equation}
for $i\neq j$.

Therefore, at point $x$ the metric and the first derivative of the
metric are diagonal.  By differentiating equation (\ref{diag}) and
evaluating at the point $x$, we conclude in the same way that the
second derivatives of the non-diagonal entries must be zero at $x$
as well.  By continued differentiation of equation (\ref{diag}),
we see that all of the derivatives of the non-diagonal entries
must be zero at $x$.  It follows immediately that the non-diagonal
entries of the metric must be zero everywhere.

We have proven that equation (\ref{diag}) of the Einstein
condition allows us to globally diagonalize the metric on a
cohomogeneity-one four manifold with $SU(2)$ symmetry.

    One may ask what happens if the metric is of cohomogeneity-one
under a larger isometry group than $SU(2)$, say $U(2)$?  If the
four manifold is of cohomogeneity-one under the action of $G=U(2)$
then it has isotropy subgroup $K=U(1)$ and the Lie algebra
decomposition becomes
\begin{equation*}
\mathfrak{u}(2)=\mathfrak{u}(1)\oplus\mathfrak{p}_{1}\oplus\mathfrak{p}_{2}
\end{equation*}
where $\textrm{dim}\mathfrak{p}_{1}=1$ and
$\textrm{dim}\mathfrak{p}_{2}=2$.  These two summands are clearly
inequivalent and we are therefore in the monotypic case and the
metric is automatically diagonal.

    Having worked through a particular example in detail, we now
wish to look at the question of diagonalization of Einstein
metrics in the non-monotypic case in greater generality.

\section{The Program}

    From the example of an Einstein four manifold of cohomogeneity-one under
the action of $SU(2)$, we can derive a program for determining the
degree to which a cohomogeneity-one Einstein metric can be
diagonalized.  We here describe, in general terms, the steps of
this program.  In the following sections, we will put this program
into more detailed practice.

\begin{itemize}
\item \textbf{Step 1:} After specifying the principal orbit type
$G/K$ choose a decomposition
$\mathfrak{g}=\mathfrak{k}\oplus\mathfrak{p}$ as above.  Determine
what non-diagonal elements may exist in the metric.  Recall that
in the monotypic case there are no non-diagonal elements to the
metric.
\item \textbf{Step 2:} Assume the number of non-diagonal entries in
the metric $g_{t}$ is greater than zero.  We can diagonalize the
metric $g_{t}$ at a point $x$ in $M/G$.

\item \textbf{Step 3:}  Having diagonalized the metric at a point,
we evaluate equation (\ref{diag}) of the Einstein condition. At
the point $x$, because the metric is diagonal, we will obtain
equations involving only the derivatives of the non-diagonal
entries and not the derivatives of the diagonal elements.

\item \textbf{Step 4:} The equations obtained in Step $3$ will
allow us to deduce that certain linear combinations of the
derivatives of the non-diagonal elements are equal to zero at $x$.

\item \textbf{Step 5:} Finally, we deduce that those linear combinations
of non-diagonal elements whose derivatives vanish, must be zero on
the whole of $M$.  If sufficiently many non-diagonal elements
vanish, then we will be able to globally diagonalize the metric.
\end{itemize}

These five steps provide a broad outline of our method, the
details of each step will be made clear below.

\section{Preliminary Calculations} Let $(M^{n},g)$ be a
cohomogeneity-one Riemannian manifold under the action of a
compact connected Lie group $G$.  Each principal orbit is of the
form $G/K$ where $K$ is the isotropy subgroup.  We choose an
$Ad(K)$-invariant complement $\mathfrak{p}$ to $\mathfrak{k}$ such
that decomposition of $\mathfrak{g}$ such that
\begin{equation}\label{gdecomp}
\mathfrak{g}=\mathfrak{k}\oplus\mathfrak{p}=\mathfrak{k}\oplus\mathfrak{p}_{1}\oplus...\oplus
\mathfrak{p}_{m}
\end{equation} where  each
$\mathfrak{p_{i}}$ is an irreducible summand of dimension $d_{i}$.
As above, we fix once and for all a bi-invariant metric $Q$ on
$\mathfrak{g}$ such that $Q(\mathfrak{k},\mathfrak{p})=0$ and
$Q(\mathfrak{p}_{i},\mathfrak{p}_{j})=0$ for all $i\neq j$.

    Let $\{Y^{i}_{j}\}_{j=1,...,d_{i}}$ be an orthonormal basis
relative to the metric $Q$ for each $\mathfrak{p}_{i}$.  Because
of the bi-invariance of $Q$ we have that
\begin{equation}
Q(Y_{j}^{i},[Y_{k}^{l},Y_{r}^{s}])=-Q(Y_{r}^{s},[Y_{k}^{l},Y_{j}^{i}])
\end{equation}
for all $i,j,k,l,r,$ and $s$.

    Again, if all of the summands $\mathfrak{p}_{i}$ are distinct then the metric, $g$,
is automatically diagonal and equation
$r(X,\frac{\partial}{\partial t})=0$ is automatically satisfied.
Therefore, we will be here interested exclusively in the
non-monotypic case.

    Non-diagonal elements of the metric on a cohomogeneity-one manifold will
arise if and only if $\mathfrak{p}_{i}\simeq\mathfrak{p}_{j}$ for
two or more of the $\mathfrak{p_{i}}$.  By reordering if
necessary, we can write the decomposition of $\mathfrak{g}$ as
\begin{equation*}
\mathfrak{g}=\mathfrak{k}\oplus\mathfrak{p}_{1}\oplus...\oplus\mathfrak{p}_{r}\oplus
\mathfrak{p}_{r+1}\oplus...\oplus\mathfrak{p}_{m}
\end{equation*}
where, for all $1\leq i\leq r$,
$\mathfrak{p_{i}}\simeq\mathfrak{p}_{j}$ for some $1\leq j\leq r$
and, for $r+1\leq k\leq m$, $\mathfrak{p}_{k}$ is distinct from
all other summands.  In this case, we can write the metric on the
fibers as
\begin{equation*}
g_{t}(\cdot,\cdot)=g_{t}(\cdot,\cdot)\mid_{\mathfrak{p}_{1}\oplus...\oplus\mathfrak{p}_{r}}
\bigoplus^{\perp}_{i=r+1,...,m}A_{i}(t)Q(\cdot,\cdot)\mid_{\mathfrak{p}_{i}}.
\end{equation*}
From the bi-invariance of $Q$ we have that
\begin{equation}\label{junk}
g_{t}(Y_{j}^{i},[Y_{j}^{i},X])=0
\end{equation}
for all $r+1\leq i\leq m$, $1\leq j \leq d_{i}$, and
$X\in\mathfrak{g}$.

Therefore, using equation (\ref{junk}) and letting $\{Z_{j}\}$ be
an invariant basis for
$\mathfrak{p}_{1}\oplus...\oplus\mathfrak{p}_{r}$ with
$h_{ij}=\langle Z_{i},Z_{j}\rangle=g_{t}(Z_{i},Z_{j})$, we see
that equation (\ref{diag}) simplifies to
\begin{equation}\label{nondiag}
h^{ij}\frac{\partial}{\partial t}\langle Z_{i},[Z_{j},X]\rangle=0
\end{equation}
for $X\in\mathfrak{p}$.  (Note that in equation (\ref{diag}) the
$h_{ij}$ are the metric functions relative to basis $\{X_{i}\}$ of
$\mathfrak{p}=\mathfrak{p}_{1}\oplus...\oplus\mathfrak{p}_{m}$
while in equation (\ref{nondiag}) the $h_{ij}$ are the metric
functions relative to basis $\{Z_{i}\}$ of
$\mathfrak{p}_{1}\oplus...\oplus\mathfrak{p}_{r}$.  In the sequel,
we work with an explicit basis which will eliminate any
ambiguity.)

    Below, we make use of equation (\ref{nondiag}) extensively.
To get an idea for how it is applied, let us first consider the
simplest type of cohomogeneity-one metric which is not
automatically diagonalizable.

\section{The Case of Two Equivalent Summands} Given the
decomposition of equation (\ref{gdecomp}), consider the case in
which
\begin{itemize}
\item $\mathfrak{p}_{1}\simeq\mathfrak{p}_{2}$ and

\item all other $\mathfrak{p}_{i}$, for $i\geq 3$, are distinct.
\end{itemize}

As in the previous section, we take $Q$ to be a fixed bi-invariant
metric on $\mathfrak{g}$ and take $\{Y^{i}_{j}\}_{j=1,...,d_{i}}$
to be an orthonormal basis for each $Q\mid_{\mathfrak{p}_{i}}$. In
this case, $d_{1}=d_{2}=d$

    Now, the fiber metric $g_{t}$, following the discussion above, can be written as
\begin{equation*}
g_{t}(\cdot,\cdot)=g_{t}(\cdot,\cdot)\mid_{\mathfrak{p}_{1}\oplus\mathfrak{p}_{2}}\bigoplus^{
\perp}_{i=3,...,m}A_{i}(t)Q(\cdot,\cdot)\mid_{\mathfrak{p}_{i}}
\end{equation*}

    Before proceeding, we will need a more precise expression for
$g_{t}(\cdot,\cdot)\mid_{\mathfrak{p}_{1}\oplus\mathfrak{p}_{2}}$
in terms of the bi-invariant form
$Q\mid_{\mathfrak{p}_{1}\oplus\mathfrak{p}_{2}}$.  To do so, we
follow closely the work of Kerr in \cite{Kerr} who studied
Einstein metrics on homogeneous manifolds whose isotropy
representations contain equivalent summands.

    We want to parameterize
    $g_{t}\mid_{\mathfrak{p}_{1}\oplus\mathfrak{p}_{2}}$in terms
of $Q\mid_{\mathfrak{p}_{1}\oplus\mathfrak{p}_{2}}$ and functions
of $t$.  Put another way, we want to find the number of distinct
non-diagonal entries a general fiber metric
$g_{t}\mid_{\mathfrak{p}_{1}\oplus\mathfrak{p}_{2}}$ has.  To
answer this question, we must determine the dimension of the space
of intertwining maps
$I:\mathfrak{p}_{1}\rightarrow\mathfrak{p}_{2}$.

    Any fiber metric restricted to $\mathfrak{p}_{1}\oplus\mathfrak{p}_{2}$ can be written as
\begin{equation*}
g_{t}(\cdot,\cdot)\mid_{\mathfrak{p}_{1}\oplus\mathfrak{p}_{2}}=Q(h\cdot,\cdot)\mid_{
\mathfrak{p}_{1}\oplus\mathfrak{p}_{2}}
\end{equation*}
where $h$ is a positive definite $Ad(K)$-invariant linear map
$h:\mathfrak{p}_{1}\oplus\mathfrak{p}_{2}\rightarrow\mathfrak{p}_{1}\oplus\mathfrak{p}_{2}$.
To parameterize, the space of $Ad(K)$-equivariant maps between the
pair of equivalent representations, we first note that
$\mathfrak{p}_{1}$ and $\mathfrak{p}_{2}$ are irreducible as real
representations.  However, the complexification of, say,
$\mathfrak{p}_{1}$ is not necessarily irreducible.  After
complexifying $\mathfrak{p}_{1}$, there are three possibilities
\begin{itemize}
\item $\mathfrak{p}_{1}$ is \textbf{orthogonal} iff
$\mathfrak{p}_{1}\otimes\mathbb{C}$  is irreducible
\item $\mathfrak{p}_{1}$ is \textbf{unitary} iff
$\mathfrak{p}_{1}\otimes\mathbb{C}=\varphi\oplus\overline{\varphi}$
and $\varphi$ is not equivalent to $\overline{\varphi}$
\item $\mathfrak{p}_{1}$ is \textbf{symplectic} iff
$\mathfrak{p}_{1}\otimes\mathbb{C}=\varphi\oplus\overline{\varphi}$
and $\varphi\simeq\overline{\varphi}$.
\end{itemize}

The space of intertwining maps is $1$-dimensional in the
orthogonal case, $2$-dimensional in the unitary case, and
$4$-dimensional in the symplectic case \cite{Kerr}.  Below, we
consider each of these cases separately.

\subsection{The Orthogonal Case}

    When $\mathfrak{p}_{1}$ and $\mathfrak{p}_{2}$ are orthogonal
representations, we have that
\begin{equation*}
g_{t}|_{\mathfrak{p}_{1}\oplus\mathfrak{p}_{2}}=\begin{pmatrix} h_{11}{\rm Id}_{d} & h_{12}{\rm Id}_{d} \\
h_{12}{\rm Id}_{d} & h_{22}{\rm Id}_{d} \\
\end{pmatrix}
\end{equation*}
relative to basis
$(Y_{1}^{1},...,Y_{d}^{1},Y_{1}^{2},...,Y_{d}^{2})$ where
$h_{ij}=h_{ij}(t)$ \cite{Kerr}.  Therefore, in the orthogonal
case, there is only one non-diagonal function $h_{12}$.  When can
we conclude that $h_{12}=0$ on an Einstein metric?

First, we calculate that
\begin{equation*}
r\left(X,\frac{\partial}{\partial
t}\right)=\frac{1}{2}h^{ij}\frac{\partial}{\partial
t}\left(\sum_{k=1}^{d}\langle
Y^{i}_{k},[Y^{j}_{k},X]\rangle\right)
\end{equation*}
for all $X\in\mathfrak{p}$.  Here we have used the Einstein
summation convention and have noted that
$g_{t}(Y^{i}_{k},Y^{j}_{l})=0$ for $k\neq l$.

Having accomplished Step One of the program by determining the
number of independent non-diagonal elements of the metric, we
progress to Step Two and diagonalize the metric at point $x\in M$.
That is, $h_{12}\mid_{x}=0$.  This implies that
\begin{equation}
r\left(X,\frac{\partial}{\partial
t}\right)\mid_{x}=\frac{1}{2}h^{ii}\frac{\partial}{\partial
t}\left(\sum_{k=1}^{d}\langle
Y^{i}_{k},[Y^{i}_{k},X]\rangle\right)\mid_{x}
\end{equation}
\begin{equation*}
=\frac{1}{2}h^{11}\frac{\partial}{\partial
t}\left(\sum_{k=1}^{d}\langle
Y^{1}_{k},[Y^{1}_{k},X]\rangle\right)\mid_{x}+\frac{1}{2}h^{22}\frac{\partial}{\partial
t}\left(\sum_{k=1}^{d}\langle
Y^{2}_{k},[Y^{2}_{k},X]\rangle\right)\mid_{x}.
\end{equation*}
Recalling the form of the fiber metric determined above and the
bi-invariance of $Q$, we have that
\begin{equation*}
\langle Y^{1}_{k},[Y^{1}_{k},X]\rangle=g_{t}(
Y^{1}_{k},[Y^{1}_{k},X])=h_{11}Q(Y^{1}_{k},[Y^{1}_{k},X])+h_{12}Q(Y^{2}_{k},[Y^{1}_{k},X])
\end{equation*}
\begin{equation*}
=h_{12}Q(Y^{2}_{k},[Y^{1}_{k},X])=-h_{12}Q([Y^{1}_{k},Y^{2}_{k}],X).
\end{equation*}
Performing a similar calculation for $\langle
Y^{2}_{k},[Y^{2}_{k},X]\rangle$, we have demonstrated that
\begin{equation*}
r\left(X,\frac{\partial}{\partial
t}\right)\mid_{x}=\frac{1}{2}(h^{22}-h^{11})\dot{h}_{12}\sum_{k=1}^{d}Q([Y_{k}^{1},Y_{k}^{2}],X)
\mid_{x}.
\end{equation*}
Step Three tells us to evaluate equation
$r\left(X,\frac{\partial}{\partial t}\right)=0$ of the Einstein
condition at point $x$.  In the Einstein case, we have that
\begin{equation*}
(h^{22}-h^{11})\dot{h}_{12}\sum_{k=1}^{d}Q([Y_{k}^{1},Y_{k}^{2}],X)\mid_{x}=0
\end{equation*}
We can assume that there exists a point $x$ at which $h^{11}\neq
h^{22}$.  This could only fail if $h=h_{11}Id_{2d}$ which implies
that $g_{t}=h_{11}Q$.  In that case, diagonalizability is
automatic.  We have
\begin{equation}
\left[\sum_{k=1}^{d}Q([Y_{k}^{1},Y_{k}^{2}],X)\right]\dot{h}_{12}\mid_{x}=0.
\end{equation}
If $\sum_{k=1}^{d}Q([Y_{k}^{1},Y_{k}^{2}],X)\neq0$ for some
$X\in\mathfrak{p}$ we deduce that $\dot{h}_{12}\mid_{x}=0$.
However, since $h_{12}\mid_{x}=\dot{h}_{12}\mid_{x}=0$ we can then
conclude that $h_{12}\equiv0$ and the metric is diagonal.

Because $Q$ is nondegenerate,
$\sum_{k=1}^{d}Q([Y_{k}^{1},Y_{k}^{2}],X)$ is nonzero for some $X$
if and only if the vector $\sum_{k=1}^{d}[Y_{k}^{1},Y_{k}^{2}]$ is
nonzero.  We have proven the following theorem.

\begin{theorem}\label{twosummandtheorem}
Let $(M,g)$ be an Einstein manifold of cohomogeneity-one under the
action of a compact, connected group $G$ such that the isotropy
representation has two equivalent \textbf{orthogonal} summands of
dimension $d$, with all others distinct. The Einstein metric can
be diagonalized if
\begin{equation}
\sum_{k=1}^{d}[Y_{k}^{1},Y_{k}^{2}]_{\mathfrak{p}}\neq0
\end{equation} where the $Y_{j}^{i}$ are as defined above.
\end{theorem}
We now review two examples in which we can apply this theorem. In
the first example, equation (\ref{diag}) of the Einstein condition
can be used to diagonalize a cohomogeneity-one metric but in the
second equation (\ref{diag}) cannot be so used.

~

\textbf{Example: $G/K\cong SO(n+1)/SO(n-1)$}

~

    The Stiefel manifold is the homogeneous space $V_{2}(\mathbb{R}^{2})=SO(n+1)/SO(n-1)$ for
$n\geq3$. Let $(M,g)$ be a cohomogeneity-one manifold with
principal orbits isomorphic to Stiefel manifolds.  In this case,
the Lie algebra $\mathfrak{g}=\mathfrak{so}(n+1)$ and
$\mathfrak{k}=\mathfrak{so}(n-1)$.  Letting $E_{ij}$ the matrix
with $1$ in the $ij$th entry and $-1$ in the $ji$th entry, we
first embed $K$ in $G$ so that
\begin{equation*}
\mathfrak{k}\cong\begin{pmatrix} 0 & 0 \\ 0 & \mathfrak{so}(n-1)
\\\end{pmatrix}\subset\mathfrak{so}(n+1).
\end{equation*}
The isotropy representation is
\begin{equation*}
\mathfrak{g}=\mathfrak{k\oplus\mathfrak{p}_{1}\oplus\mathfrak{p}_{2}\oplus\mathfrak{p}_{3}}
\end{equation*}
and we can take $\mathfrak{p}_{j}=\textrm{span}\{E_{j,2+i}\mid
1\leq i\leq n-1\}$ for $j=1,2$, and
$\mathfrak{p}_{3}=\textrm{span}\{E_{12}\}$.  This decomposition is
not unique since $\mathfrak{p}_{1}\simeq\mathfrak{p}_{2}$.
Therefore, the cohomogeneity-one metric $g$ is not automatically
diagonal.

To determine whether $g$ can be diagonalized in the Einstein case,
we must calculate the vector \linebreak
$\sum_{k=1}^{n-1}[E_{1,2+k},E_{2,2+k}]_{\mathfrak{p}}$.  It is
straightforward to determine that

\begin{equation*}
\sum_{k=1}^{n-1}[E_{1,2+k},E_{2,2+k}]_{\mathfrak{p}}=-(n-1)E_{12}\neq0.
\end{equation*}
Therefore, by theorem (\ref{twosummandtheorem}), any
cohomogeneity-one Einstein metric with principal orbit isomorphic
to $SO(n+1)/SO(n-1)$ can be diagonalized globally.

    Stenzel in \cite{Stenzel} constructed Ricci-flat Kahler metrics on
the tangent bundle to the $n$-sphere, $T{S^{n}}$.  These metrics
are of cohomogeneity-one with principal orbits isomorphic to
Stiefel manifolds (the $n=2$ case is equivalent to the
Eguchi-Hanson metric and the $n=3$ case was found by Candelas and
de la Ossa in \cite{CandD}).  Dancer and Strachan in \cite{DandS}
constructed Einstein metrics of negative Einstein constant on
manifolds of this type.  In each of these cases, the metrics were
found by restricting to the diagonal case.  Up to now it was
unclear whether or not one could construct non-diagonal
cohomogeneity-one Einstein manifolds of this type which were
inequivalent to the diagonal ones.  The above proves that the
diagonal case is equivalent to the general case and requiring the
metric to be diagonal need not be considered an additional
constraint.

~

\textbf{Example: $G/K\cong SO(n_{1}+n_{2}+2)/(SO(n_{1})\times
SO(n_{2}))$}

~

    We now consider a generalization of a Stiefel manifold.  Consider $(M,g)$ a
cohomogeneity-one manifold with principal orbits isomorphic to
$SO(n_{1}+n_{2}+2)/(SO(n_{1})\times SO(n_{2}))$ where $n_{1}$ and
$n_{2}$ are positive integers.

    Letting $E_{ij}$ be the matrix with $1$ in the $ij$th entry and $-1$ in the $ji$th entry.
We can embed $SO(n_{1})\times SO(n_{2})$ in $SO(n_{1}+n_{2}+2)$ so
that on the Lie algebra level we have
\begin{equation*}
\mathfrak{k}\cong\begin{pmatrix}
0 & 0 & 0 \\
0 & \mathfrak{so}(n_{1}) & 0 \\
0 & 0 & \mathfrak{so}(n_{2}) \\
\end{pmatrix}\subset\mathfrak{so}(n_{1}+n_{2}+2).
\end{equation*}
The isotropy representation becomes
\begin{equation*}
\mathfrak{g}=\mathfrak{k}\oplus\mathfrak{p}_{1}\oplus\mathfrak{p}_{2}\oplus\mathfrak{p}_{3}
\oplus\mathfrak{p}_{4}\oplus\mathfrak{p}_{5}\oplus\mathfrak{p}_{6}
\end{equation*}
where $\mathfrak{g}=\mathfrak{so}(n_{1}+n_{2}+1)$,
$\mathfrak{k}=\mathfrak{so}(n_{1})\oplus\mathfrak{so}(n_{2})$,
$\mathfrak{p}_{1}\cong(\mathbb{R}^{n_{1}}\otimes\mathbb{R}^{n_{2}})$,
$\mathfrak{p}_{2}\simeq\mathfrak{p}_{3}\cong\mathbb{R}^{n_{1}}$,
$\mathfrak{p}_{4}\simeq\mathfrak{p}_{5}\cong\mathbb{R}^{n_{2}}$,
and $\mathfrak{p}_{6}\cong\mathbb{R}$.

We see that $\mathfrak{p}_{2}=\textrm{span}\{E_{1,2+j}\mid 1\leq
j\leq n_{1}\},$ $\mathfrak{p}_{3}=\textrm{span}\{E_{2,2+j}\mid
1\leq j\leq n_{1}\},$
$\mathfrak{p}_{4}=\textrm{span}\{E_{1,2+n_{1}+j}\mid 1\leq j \leq
n_{2}\},$ $\mathfrak{p}_{5}=\textrm{span}\{E_{2,2+n_{1}+j}\mid
1\leq j \leq n_{2}\},$ and
$\mathfrak{p}_{6}=\textrm{span}\{E_{12}\}$.  Clearly, we have two
sets of equivalent orthogonal summands.  This implies that a
general cohomogeneity-one metric of this type will have two
non-diagonal functions.
    To determine whether $g$ can be diagonalized globally in the Einstein case, we first
calculate that
\begin{equation*}
\sum_{i=1}^{n_{1}}[E_{1,2+j},E_{2,2+j}]_{\mathfrak{p}}=-n_{1}E_{12}
\end{equation*}
and
\begin{equation*}
\sum_{i=1}^{n_{2}}[E_{1,2+n_{1}+j},E_{2,2+n_{1}+j}]_{\mathfrak{p}}=-n_{2}E_{12}.
\end{equation*}

    Given this calculation, it is straightforward that the Einstein condition cannot
be diagonalized globally by the method described above.  This
follows because we have only one condition which is not sufficient
to force the two non-diagonal functions to be zero.

\subsection{The Unitary Case}

    When $\mathfrak{p}_{1}$ and $\mathfrak{p}_{2}$ are equivalent unitary
representations, we have that
$g_{t}(\cdot,\cdot)\mid_{\mathfrak{p}_{1}\oplus\mathfrak{p}_{2}}=Q((h+f)\cdot,\cdot)\mid_{
\mathfrak{p}_{1}\oplus\mathfrak{p}_{2}}$
where $h$ and $f$ are of the form
\begin{equation*}
h=\begin{pmatrix}
h_{11}\textrm{Id}_{d} & h_{12}\textrm{Id}_{d} \\
h_{12}\textrm{Id}_{d} & h_{22}\textrm{Id}_{d} \\
\end{pmatrix}
\end{equation*}
and
\begin{equation*}
f=
\begin{pmatrix}
 & & & f\textrm{Id}_{\frac{d}{2}} \\
 & & -f\textrm{Id}_{\frac{d}{2}} & \\
& -f\textrm{Id}_{\frac{d}{2}} & & \\
f\textrm{Id}_{\frac{d}{2}} & & & \\
\end{pmatrix}
\end{equation*}
relative to a basis
$(Y^{1}_{1},...,Y^{1}_{d},Y^{2}_{1},...,Y^{2}_{d})$ where
$h_{ij}=h_{ij}(t)$ and $f=f(t)$.  Therefore, in the unitary case
there are two non-diagonal functions.  When can we conclude that
$h_{12}\equiv0$ and $f\equiv0$ in the Einstein case?

    As above, we diagonalize the metric at a point $x\in M$.  That
is, $h_{12}\mid_{x}=0$ and $f\mid_{x}=0$.  At $x$, we calculate
that
\begin{equation*}
r\left(X,\frac{\partial}{\partial
t}\right)\mid_{x}=\frac{1}{2}h^{ii}\frac{\partial}{\partial
t}\left(\sum_{k=1}^{d}\langle
Y_{k}^{i},[Y_{k}^{i},X]\rangle\right)\mid_{x}
\end{equation*}
\begin{equation*}
=\frac{1}{2}h^{11}\frac{\partial}{\partial
t}\left(\sum_{k=1}^{d}\langle
Y_{k}^{1},[Y_{k}^{1},X]\rangle\right)\mid_{x}+\frac{1}{2}h^{22}\frac{\partial}{\partial
t}\left(\sum_{k=1}^{d}\langle
Y_{k}^{2},[Y_{k}^{2},X]\rangle\right)\mid_{x}
\end{equation*}
for all $X\in\mathfrak{p}$. Noting the bi-invariance of $Q$ we
calculate that

\begin{equation*}
\langle
Y_{k}^{1},[Y_{k}^{1},X]\rangle=g_{t}(Y_{k}^{1},[Y_{k}^{1},X])=-h_{12}Q([Y_{k}^{1},Y^{2}_{k}]
,X)-fQ([Y^{1}_{k},Y^{2}_{k+\frac{d}{2}}],X)
\end{equation*}
for $1\leq k\leq\frac{d}{2}$ and
\begin{equation*}
\langle
Y_{k}^{1},[Y_{k}^{1},X]\rangle=g_{t}(Y_{k}^{1},[Y_{k}^{1},X])=-h_{12}Q([Y_{k}^{1},Y^{2}_{k}]
,X)+fQ([Y^{1}_{k},Y^{2}_{k-\frac{d}{2}}],X)
\end{equation*}
for $\frac{d}{2}\leq k\leq d$ and for all $X\in\mathfrak{p}$.

    After calculating $\langle Y_{k}^{2},[Y_{k}^{2},X]\rangle$,
we have that

\begin{equation*}
r(X,\frac{\partial}{\partial
t})\mid_{x}=\frac{1}{2}(h^{22}-h^{11})\left[\dot{h}_{12}\mid_{x}\sum_{k=1}^{d}Q([Y^{1}_{k},
Y^{2}_{k}],X)+\dot{f}\mid_{x}\sum_{k=1}^{\frac{d}{2}}Q([Y_{k}^{1},Y^{2}_{k+\frac{d}{2}}]-
[Y^{1}_{k+\frac{d}{2}},Y^{2}_{k}],X)\right].
\end{equation*}
As in the orthogonal case, we can assume that $h_{11}\mid_{x}\neq
h_{22}\mid_{x}$.  Equation $r(X,\frac{\partial}{\partial
t})\mid_{x}=0$ of the Einstein condition becomes

\begin{equation}\label{unitEinst}
\dot{h}_{12}\mid_{x}\sum_{k=1}^{d}Q([Y^{1}_{k},Y^{2}_{k}],X)+\dot{f}\mid_{x}\sum_{k=1}^{
\frac{d}{2}}Q([Y_{k}^{1},Y^{2}_{k+\frac{d}{2}}]-[Y^{1}_{k+\frac{d}{2}},Y^{2}_{k}],X)=0
\end{equation}
for all $X\in\mathfrak{p}$.

We can now prove the unitary analogue of theorem
(\ref{twosummandtheorem}).

\begin{theorem}\label{2unitarysummandtheorem}  Let $(M,g)$ be an Einstein manifold of
cohomogeneity-one under the action of a compact connected group
$G$ such that the isotropy representation has two equivalent
\textbf{unitary} summands of dimension $d$, with all other
summands distinct.  The Einstein metric can be diagonalized
globally if
\begin{equation*}
{\rm dim}\left({\rm
span}\left\{\sum_{k=1}^{d}[Y^{1}_{k},Y^{2}_{k}]_{\mathfrak{p}},\sum_{k=1}^{\frac{d}{2}}
([Y_{k}^{1},Y^{2}_{k+\frac{d}{2}}]_{\mathfrak{p}}-[Y^{1}_{k+\frac{d}{2}},Y^{2}_{k}]_{
\mathfrak{p}})\right\}\right)=2
\end{equation*}
where the $Y^{i}_{j}$ are as defined above.
\end{theorem}
\begin{proof}
If the span of the above vectors is of dimension two then equation
(\ref{unitEinst}) implies that
$\dot{h}_{12}\mid_{x}=\dot{f}\mid_{x}=0$.

    Therefore, we have demonstrated that, in the Einstein case,
$h_{12}\mid_{x}=f\mid_{x}=\dot{h}_{12}\mid_{x}=\dot{f}\mid_{x}$.
This implies that $h_{12}\equiv0$ and $f\equiv0$.
\end{proof}

~

\textbf{Example: $G/K\cong SU(n+1)/U(n-1)$}

~

This is a well-known homogeneous manifold which was studied by
Calabi and others.  It can be viewed as a bundle over
$\mathbb{C}P^{n}\cong SU(n+1)/U(n)$ with fiber $S^{2n-1}\cong
U(n)/U(n-1)$.  We represent $\mathfrak{su}(n+1)$ as the set of
$(n+1)\times(n+1)$ complex matrices $A$ such that $A^{\ast}+A=0$
and $\textrm{tr}(A)=0$.  We embed $U(n-1)$ in $SU(n+1)$ so that
$\mathfrak{u}(n-1)$ sits in $\mathfrak{su}(n+1)$ as
\begin{equation*}
\mathfrak{k}\cong\mathfrak{u}(n-1)\cong\begin{pmatrix}
\mathfrak{u}(n-1)&0 \\
0& 0 \\
\end{pmatrix}
\subset\mathfrak{su}(n+1).
\end{equation*}

    Because the notation is cumbersome, we will only calculate the $n=2$ result and
leave the rest to the reader. Schematically, we can write
$\mathfrak{su}(3)$ as
\begin{equation*}
\begin{pmatrix}
ia& \alpha & \beta \\
-\bar{\alpha} & ib & \gamma \\
-\bar{\beta} & -\bar{\gamma} & i(-b-a) \\
\end{pmatrix}
\end{equation*}
where $a,b\in\mathbb{R}$ and $\alpha,\beta,$ and
$\gamma\in\mathbb{C}$.  Let
\begin{equation*}
\mathfrak{u}(1)=\textrm{span}\left\{\begin{pmatrix}ia& 0 & 0 \\
0 & 0 & 0 \\
0 & 0 & 0 \\
\end{pmatrix}\right\}.
\end{equation*}

The isotropy representation is
\begin{equation*}
\mathfrak{so}(n+1)=\mathfrak{u}(n-1)\oplus\mathfrak{p}_{1}\oplus\mathfrak{p}_{2}\oplus
\mathfrak{p}_{3}\oplus\mathfrak{p}_{4}\oplus\mathfrak{p}_{5}
\end{equation*}
where $\mathfrak{p}_{1}\simeq\mathfrak{p}_{2}$ are equivalent
unitary representations of real dimension $2n-2$ and
$\mathfrak{p}_{3}\simeq\mathfrak{p}_{4}\simeq\mathfrak{p}_{5}$ are
trivial one-dimensional representations.  Moreover,
\begin{equation*}
\mathfrak{p}_{3}\oplus\mathfrak{p}_{4}\oplus\mathfrak{p}_{5}\cong\mathfrak{su}(2).
\end{equation*}
In this case, there are two sets of equivalent summands. That is a
general metric would be of the form
\begin{equation*}
g_{t}(\cdot,\cdot)|_{\mathfrak{p}_{1}\oplus\mathfrak{p}_{2}}=Q((h+f),\cdot,\cdot)|_{
\mathfrak{p}_{1}\oplus\mathfrak{p}_{2}}
\end{equation*}
with $h$ and $f$ as above and
\begin{equation*}
g_{t}(\cdot,\cdot))|_{\mathfrak{p}_{3}\oplus\mathfrak{p}_{4}\oplus\mathfrak{p}_{5}}=
Q(\bar{h}\cdot,\cdot)|_{\mathfrak{p}_{3}\oplus\mathfrak{p}_{4}\oplus\mathfrak{p}_{5}}
\end{equation*}
where
\begin{equation*}
\bar{h}=\begin{pmatrix}\bar{h}_{11} & \bar{h}_{12}&\bar{h}_{13}\\
\bar{h}_{12}&\bar{h}_{22}&\bar{h}_{23}\\
\bar{h}_{13}&\bar{h}_{23}&\bar{h}_{33} \\
\end{pmatrix}
\end{equation*}

For our purposes we wish to study the pair of equivalent unitary
representations.  Therefore, we partially diagonalize the
cohomogeneity-one metric.

~

\textit{Assumption A:}

\begin{equation*}
\bar{h}_{12}=\bar{h}_{13}=\bar{h}_{23}=0.
\end{equation*}
~

Under this assumption, there are only two non-diagonal functions:
$h_{12}$ and $f$.  Given the above assumption, can we diagonalize
a cohomogeneity-one Einstein metric with this principal orbit
type?

In this case,
\begin{equation*}
\mathfrak{p}_{1}=\textrm{span}\left\{Y^{1}_{1}=\begin{pmatrix}0& 1 & 0 \\
-1 & 0 & 0 \\
0 & 0 & 0 \\
\end{pmatrix}, Y^{1}_{2}=\begin{pmatrix}0& i & 0 \\
i & 0 & 0 \\
0 & 0 & 0 \\
\end{pmatrix}\right\}
\end{equation*}
 and
\begin{equation*}
\mathfrak{p}_{2}=\textrm{span}\left\{Y^{2}_{1}=\begin{pmatrix}0& 0 & 1 \\
0 & 0 & 0 \\
-1 & 0 & 0 \\
\end{pmatrix}, Y^{2}_{2}=\begin{pmatrix}0& 0 & i \\
0 & 0 & 0 \\
i & 0 & 0 \\
\end{pmatrix}\right\}.
\end{equation*}  A simple calculation demonstrates that $\textrm{span}\{[Y^{1}_{1},Y^{2}_{1}]
_{\mathfrak{p}}+[Y^{1}_{2},Y^{2}_{2}]_{\mathfrak{p}},[Y^{1}_{1},Y^{2}_{2}]_{\mathfrak{p}}-
[Y^{1}_{2},Y^{2}_{1}]_{\mathfrak{p}}\}$ is of dimension two.  Therefore,
by theorem (\ref{2unitarysummandtheorem}), any cohomogeneity-one
Einstein metric with principal orbits isomorphic to $G/K\cong
SU(n+1)/U(n-1)$ such that Assumption A holds can be diagonalized
globally.

Alternatively, we could make a different assumption.

~

\textit{Assumption B:}

\begin{equation*}
h_{12}=f=0.
\end{equation*}

~

    As before, a cohomogeneity-one metric satisfying Assumption B is only partially
diagonalized as $\bar{h}$ may contain non-diagonal entries.  In
this case, however, the situation is slightly more complicated as
there are now three equivalent summands
$\mathfrak{p}_{3},\mathfrak{p}_{4},$ and $\mathfrak{p}_{5}$ which
are no longer automatically orthogonal.  We will go over the
general case of more than two equivalent summands below. However,
we do not need those results to analyze this particular case. Note
that
\begin{equation*}
\mathfrak{p}_{3}\oplus\mathfrak{p}_{4}\oplus\mathfrak{p}_{5}\cong\mathfrak{su}(2).
\end{equation*}
Therefore, given Assumption B, diagonalizing a cohomogeneity-one
Einstein metric with principal orbits isomorphic to $G/K\cong
SU(n+1)/U(n-1)$ is equivalent to diagonalizing a cohomogeneity-one
Einstein metric with $SU(2)$ principal orbits.  That, of course,
is possible.

    To summarize, a cohomogeneity-one Einstein metric with $G/K\cong SU(n+1)/U(n-1)$
can be diagonalized globally if Assumption A or Assumption B hold.
If neither hold, we cannot diagonalize the metric in the Einstein
case using this method.

\subsection{The Symplectic Case}

    When $\mathfrak{p}_{1}$ and $\mathfrak{p}_{2}$ are symplectic
representations, we have that
$g_{t}(\cdot,\cdot)\mid_{\mathfrak{p}_{1}\oplus\mathfrak{p}_{2}}=Q((h+f)\cdot,\cdot)\mid_{
\mathfrak{p}_{1}\oplus\mathfrak{p}_{2}}$
where
\begin{equation*}
h=\begin{pmatrix} h_{11}\textrm{Id}_{d} & \\
 & h_{22}\textrm{Id}_{d}
\end{pmatrix}
\end{equation*}
and
\begin{equation*}
f=\begin{pmatrix}  & & f^{1}\textrm{Id}_{\frac{d}{2}} &
f^{2}\textrm{Id}_{\frac{d}{2}} \\
& & f^{3}\textrm{Id}_{\frac{d}{2}} & f^{4}\textrm{Id}_{\frac{d}{2}} \\
f^{1}\textrm{Id}_{\frac{d}{2}} & f^{2}\textrm{Id}_{\frac{d}{2}} & & \\
f^{3}\textrm{Id}_{\frac{d}{2}} & f^{4}\textrm{Id}_{\frac{d}{2}} & & \\
\end{pmatrix}
\end{equation*}
relative to basis
$(Y^{1}_{1},...,Y^{1}_{d},Y^{2}_{1},...,Y^{2}_{d})$.  Note that we
have defined $h$ differently than we had before.  In the
symplectic case, there are four non-diagonal functions.  When can
we conclude that all the $f^{j}=0$ globally in the Einstein case?
    We can diagonalize the metric at a point $x\in M$.  That is
$f^{j}\mid_{x}=0$ for $1\leq j\leq 4$.  At $x$ we calculate that
\begin{equation*}
r\left(X,\frac{\partial}{\partial
t}\right)\mid_{x}=\frac{1}{2}h^{11}\frac{\partial}{\partial
t}\left(\sum_{k=1}^{d}\langle
Y^{1}_{k},[Y^{1}_{k},X]\rangle\right)\mid_{x}+\frac{1}{2}h^{22}\frac{\partial}{\partial
t}\left(\sum_{k=1}^{d}\langle
Y^{2}_{k},[Y^{2}_{k},X]\rangle\right)\mid_{x}
\end{equation*}
for all $X\in\mathfrak{p}$.  Making use of the bi-invariance of
$Q$ we have
\begin{equation*}
\langle
Y^{1}_{k},[Y^{1}_{k},X]\rangle=g_{t}(Y^{1}_{k},[Y^{1}_{k},X])=-f^{1}Q([Y^{1}_{k},Y^{2}_{k}],
X)-f^{3}Q([Y^{1}_{k},Y^{2}_{k+\frac{d}{2}}],X)
\end{equation*}
for $1\leq k\leq\frac{d}{2}$, and
\begin{equation*}
\langle
Y^{1}_{k},[Y^{1}_{k},X]\rangle=g_{t}(Y^{1}_{k},[Y^{1}_{k},X])=-f^{2}Q([Y^{1}_{k},Y^{2}_{k-
\frac{d}{2}}],X)-f^{4}Q([Y^{1}_{k},Y^{2}_{k}],X)
\end{equation*}
for $\frac{d}{2}+1\leq k\leq d$ for all $X\in\mathfrak{p}$.

    Similar calculations hold for $\langle
Y^{2}_{k},[Y^{2}_{k},X]\rangle$ and we have that
\begin{equation*}
r\left(X,\frac{\partial}{\partial
t}\right)\mid_{x}=\frac{1}{2}(h^{22}-h^{11})\{\dot{f}^{1}\mid_{x}\sum_{k=1}^{\frac{d}{2}}
Q([Y^{1}_{k},Y^{2}_{k}],X)+\dot{f}^{3}\mid_{x}\sum_{k=1}^{\frac{d}{2}}Q([Y^{1}_{k},Y^{2}_{k
+\frac{d}{2}}],X)
\end{equation*}
\begin{equation*}
+\dot{f}^{2}\mid_{x}\sum_{k=\frac{d}{2}+1}^{d}Q([Y^{1}_{k},Y^{2}_{k-\frac{d}{2}}],X)+
\dot{f}^{4}\mid_{x}\sum_{k=\frac{d}{2}+1}^{d}Q([Y^{1}_{k},Y^{2}_{k}],X)\}.
\end{equation*}

    We can assume, as above, that $h^{11}\mid_{x}\neq
h^{22}\mid_{x}$.  With this assumption, equation
$r(X,\frac{\partial}{\partial t})\mid_{x}=0$ of the Einstein
condition becomes

\begin{equation}\label{symp20}
\dot{f}^{1}\mid_{x}Q\left(\sum_{k=1}^{\frac{d}{2}}[Y^{1}_{k},Y^{2}_{k}],X\right)+\dot{f}^{3}
\mid_{x}Q\left(\sum_{k=1}^{\frac{d}{2}}[Y^{1}_{k},Y^{2}_{k+\frac{d}{2}}],X\right)
\end{equation}
\begin{equation*}
+\dot{f}^{2}\mid_{x}Q\left(\sum_{k=\frac{d}{2}+1}^{d}
[Y^{1}_{k},Y^{2}_{k-\frac{d}{2}}],X\right)+\dot{f}^{4}\mid_{x}Q
\left(\sum_{k=\frac{d}{2}+1}^{d}[Y^{1}_{k},Y^{2}_{k}],X\right)=0
\end{equation*}
for all $X\in\mathfrak{p}$.

    We now have the symplectic analogue to theorem
(\ref{twosummandtheorem}).
\begin{theorem}
Let $(M,g)$ be an Einstein manifold of cohomogeneity-one under the
action of a compact, connected group $G$ such that the isotropy
representation has two equivalent symplectic summands of dimension
$d$ with all others distinct.  The Einstein metric can be
diagonalized globally if
\begin{equation*}
{\rm dim}\left({\rm
span}\left\{\sum_{k=1}^{\frac{d}{2}}[Y^{1}_{k},Y^{2}_{k}]_{\mathfrak{p}},
\sum_{k=1}^{\frac{d}{2}}[Y^{1}_{k},Y^{2}_{k+\frac{d}{2}}]_{\mathfrak{p}},
\sum_{k=\frac{d}{2}+1}^{d}[Y^{1}_{k},Y^{2}_{k-\frac{d}{2}}]_{\mathfrak{p}},
\sum_{k=\frac{d}{2}+1}^{d}[Y^{1}_{k},Y^{2}_{k}]_{\mathfrak{p}}\right\}\right)=4
\end{equation*}
where the $Y^{j}_{i}$ are as defined above.
\end{theorem}
\begin{proof}
If the span of the above vectors is of dimension four then
equation (\ref{symp20}) implies that $\dot{f}^{i}\mid_{x}=0$ for
all $1\leq i\leq4$. Since all the $f^{i}\mid_{x}=0$, we deduce
that $f^{i}\equiv0$ for all $1\leq i\leq4$.
\end{proof}

\section{The Case of $r$ Equivalent Summands}

Having analyzed the orthogonal, unitary, and symplectic case for
isotropy representations with only pairs of equivalent summands,
we now consider the case in which

\begin{equation*}
\mathfrak{g}=\mathfrak{k}\oplus\mathfrak{p}_{1}\oplus...\oplus\mathfrak{p}_{r}\oplus
\mathfrak{p}_{r+1}\oplus...\oplus\mathfrak{p}_{m}
\end{equation*}
where $\mathfrak{p}_{i}\simeq\mathfrak{p}_{j}$ for $1\leq i,j \leq
r$ and all other $\mathfrak{p}_{l}$ for $l>r$ are distinct.

As above, we fix a bi-invariant metric $Q$ on $\mathfrak{g}$ where
$\{Y^{i}_{j}\}_{j=1,...,d_{i}}$ forms a $Q$-orthonormal basis for
each $\mathfrak{p}_{i}$.  The fiber metrics $g_{t}$ can be written

\begin{equation*}
g_{t}(\cdot,\cdot)=g_{t}(\cdot,\cdot)\mid_{\mathfrak{p}_{1}\oplus...\oplus\mathfrak{p}_{r}}
\bigoplus_{i=r+1,...,m}^{\perp}A_{i}(t)Q(\cdot,\cdot)\mid_{\mathfrak{p}_{i}}
\end{equation*}

Since the restriction of $g_{t}$ to
$\mathfrak{p}_{1}\oplus...\oplus\mathfrak{p}_{r}$ is
$Ad(K)$-invariant, it can be written as \linebreak
$g_{t}(\cdot,\cdot)\mid_{\mathfrak{p}_{1}\oplus...\oplus\mathfrak{p}_{r}}=Q(h\cdot,\cdot)
\mid_{\mathfrak{p}_{1}\oplus...\oplus\mathfrak{p}_{r}}$ where $h$
is a positive definite $Ad(K)$-invariant map
$h:\mathfrak{p}_{1}\oplus...\oplus\mathfrak{p}_{r}\rightarrow\mathfrak{p}_{1}\oplus...
\oplus\mathfrak{p}_{r}$. The form of $h$ depends upon whether the
representation $\mathfrak{p}_{1}$ is orthogonal, unitary, or
symplectic.  We again consider each of these case individually.

\subsection{The Orthogonal Case} When
$\mathfrak{p}_{1},...,\mathfrak{p}_{r}$ are orthogonal
representations, we have that $h$ is given by the matrix
\begin{equation*}
h=(h_{ij}Id_{d})
\end{equation*}
where $1\leq i,j\leq r$, $h_{ij}=h_{ij}(t)$ and
$d=d_{1}=...=d_{r}$.  This implies that a generic
cohomogeneity-one manifold of this type has $\frac{r(r-1)}{2}$
independent non-diagonal functions.  We want to determine when all
of these functions can be set equal to zero at every point of an
Einstein manifold without loss of generality.

We calculate that

\begin{equation*}
r\left(X,\frac{\partial}{\partial
t}\right)=\frac{1}{2}h^{ij}\frac{\partial}{\partial
t}\left(\sum_{k=1}^{d}\langle
Y^{i}_{k},[Y^{j}_{k},X]\rangle\right)
\end{equation*}
for all $X\in\mathfrak{p}$.
    We can, without loss of generality, diagonalize the metric at
a point $x\in M$ which implies that $h_{ij}\mid_{x}=0$ for $i\neq
j$.  Given this assumption, we have

\begin{equation*}
r\left(X,\frac{\partial}{\partial
t}\right)\mid_{x}=\frac{1}{2}\sum_{i=1}^{r}h^{ii}\frac{\partial}{\partial
t}\left(\sum_{k=1}^{d}\langle
Y^{i}_{k},[Y^{i}_{k},X]\rangle\right)\mid_{x}.
\end{equation*}

From the bi-invariance of $Q$, we calculate that

\begin{equation*}
\langle
Y^{i}_{k},[Y^{i}_{k},X]\rangle=g_{t}(Y^{i}_{k},[Y^{i}_{k},X])=\sum_{j=1}^{r}h_{ij}
Q(Y^{j}_{j},[Y^{i}_{k},X])
\end{equation*}
\begin{equation*}
=-\sum_{j=1}^{r}h_{ij}Q([Y^{i}_{k},Y^{j}_{k}],X).
\end{equation*}

Therefore, at $x$, the Einstein condition implies that

\begin{equation}\label{Einsteinrortho}
\sum_{1\leq i<j\leq
r}\left\{(h^{jj}-h^{ii})Q\left(\sum_{k=1}^{d}[Y^{i}_{k},Y^{j}_{k}],X\right)\right\}
\dot{h}_{ij}\mid_{x}=0
\end{equation}
for all $X\in\mathfrak{p}$.

    To simplify our presentation we define vectors $Z_{ij}$ by the
equation

\begin{equation*}
Z_{ij}=\sum_{k=1}^{d}[Y^{i}_{k},Y^{j}_{k}]_{\mathfrak{p}}
\end{equation*}
for $1\leq i<j\leq r$.

We are now ready to prove the diagonalization theorem for this
type of group action.

\begin{theorem}\label{rsummandortho}
Let $(M,g)$ be an Einstein manifold of cohomogeneity-one under the
action of a compact, connected group $G$ such that the isotropy
representation has at most $r$ equivalent \textbf{orthogonal}
summands of dimension $d$ with all others distinct.  The Einstein
metric can be diagonalized globally if
\begin{equation}
{\rm dim}({\rm span}\{Z_{ij}~|~1\leq i<j\leq
r\})=\frac{r(r-1)}{2}.
\end{equation}

\end{theorem}

\begin{proof}
As above, we can assume that $h_{ii}\mid_{x}\neq h_{jj}\mid_{x}$
for $i\neq j$.  Therefore, if the span of the $Z_{ij}$'s has
dimension $\frac{r(r-1)}{2}$, then equations
(\ref{Einsteinrortho}) become $\frac{r(r-1)}{2}$ independent
equations in the $\dot{h}_{ij}\mid_{x}$ terms.  Since there are
$\frac{r(r-1)}{2}$ such terms, we deduce that
$\dot{h}_{ij}\mid_{x}=0$ for $i\neq j$.  Since we already have
that $h_{ij}\mid_{x}=0$ for $i\neq j$, we have that
$h_{ij}\equiv0$ for all $i\neq j$.  The metric is then globally
diagonal.
\end{proof}

\subsection{The Unitary Case} When
$\mathfrak{p}_{1},...,\mathfrak{p}_{r}$ are equivalent unitary
representations of dimension $d$, we have that
$g_{t}(\cdot,\cdot)\mid_{\mathfrak{p}_{1}\oplus...\oplus\mathfrak{p}_{r}}=Q((h+f)\cdot,
\cdot)\mid_{\mathfrak{p}_{1}\oplus...\oplus\mathfrak{p}_{r}}$
relative to basis $(Y^{1}_{1},...,Y^{r}_{d})$ where $h$ is of the
form
\begin{equation*}
h=(h_{ij}Id_{d})
\end{equation*}
for $1\leq i,j\leq r$ and $f$ is defined by
\begin{equation*}
f_{ij}=g_{t}(Y^{i}_{k},Y^{j}_{k+\frac{d}{2}})
\end{equation*}
for $i<j$ and $1\leq k\leq\frac{d}{2}$ and
\begin{equation*}
f_{ij}=-g_{t}(Y^{i}_{k},Y^{j}_{k-\frac{d}{2}})
\end{equation*}
for $i<j$ and $\frac{d}{2}+1\leq k\leq d$.

    After diagonalizing the metric at point $x\in M$, that is
$h_{ij}\mid_{x}=f_{ij}\mid_{x}=0$ for $i\neq j$, equation
$r(X,\frac{\partial}{\partial x})\mid_{x}=0$ of the Einstein
condition becomes
\begin{equation}\label{unitaryrsummand}
\sum_{1\leq i<j\leq
r}(h^{jj}-h^{ii})\left\{\dot{h}_{ij}\mid_{x}\sum_{k=1}^{d}Q([Y_{k}^{i},Y_{k}^{j}],X)
+\dot{f}_{ij}\mid_{x}\sum_{k=1}^{\frac{d}{2}}Q([Y_{k}^{i},Y_{k+\frac{d}{2}}^{j}]-[Y_{k
+\frac{d}{2}}^{i},Y_{k}^{2}],X)\right\}=0
\end{equation}
for all $X\in\mathfrak{p}$.

    To simplify our presentation, we define the vectors
\begin{equation*}
Z_{ij}=\sum_{k=1}^{d}[Y_{k}^{i},Y_{k}^{j}]_{\mathfrak{p}}
\end{equation*}
and
\begin{equation*}
W_{ij}=\sum_{k=1}^{\frac{d}{2}}([Y_{k}^{i},Y^{j}_{k+\frac{d}{2}}]_{\mathfrak{p}}-[Y^{i}_{k
+\frac{d}{2}},Y^{j}_{k}]_{\mathfrak{p}})
\end{equation*}
for $1\leq i<j\leq r$.

Given the above formalism, we can now prove the following theorem.

\begin{theorem}
Let $(M,g)$ be an Einstein manifold of cohomogeneity-one under the
action of a compact, connected group $G$ such that the isotropy
representation has at most $r$ equivalent \textbf{unitary}
summands of dimension $d$ with all others distinct.  The Einstein
metric can be diagonalized globally if
\begin{equation*}
{\rm dim}({\rm span}\{Z_{ij},W_{ij}~|~1\leq i<j\leq r\})=r(r-1).
\end{equation*}

\end{theorem}

\begin{proof}
The proof of this theorem is similar to the proof of the analogous
theorem in the orthogonal case.  The only difference is that, in
the unitary case, there are $r(r-1)$ independent non-diagonal
functions.  The rest follows in an identical fashion.
\end{proof}
\subsection{The Symplectic Case}

    When $\mathfrak{p}_{1},...,\mathfrak{p}_{r}$ are equivalent
symplectic representations of dimension $d$, we have that
$g_{t}(\cdot,\cdot)\mid_{\mathfrak{p}_{1}\oplus...\oplus\mathfrak{p}_{r}}=Q((h+f)\cdot,\cdot)
\mid_{\mathfrak{p}_{1}\oplus...\oplus\mathfrak{p}_{r}}$
relative to basis $(Y^{1}_{1},...,Y^{r}_{d})$ where $h$ is of the
form
\begin{equation*}
h=\begin{pmatrix} h_{11} & & & \\
 & h_{22} & & \\
& & ... & \\
& & & h_{rr} \\
\end{pmatrix}
\end{equation*}
and $f$ is defined, for $1\leq i,j\leq r$ by
\begin{equation*}
f^{1}_{ij}=g_{t}(Y^{i}_{k},Y^{j}_{k})
\end{equation*}
for $1\leq k\leq\frac{d}{2}$,
\begin{equation*}
f^{3}_{ij}=g_{t}(Y^{i}_{k},Y^{j}_{k+\frac{d}{2}})
\end{equation*}
for $1\leq k\leq\frac{d}{2}$,
\begin{equation*}
f^{2}_{ij}=g_{t}(Y^{i}_{k},Y^{j}_{k-\frac{d}{2}})
\end{equation*}
for $\frac{d}{2}+1\leq k\leq d$, and
\begin{equation*}
f^{4}_{ij}=g_{t}(Y^{i}_{k},Y^{j}_{k})
\end{equation*}
for $\frac{d}{2}+1\leq k\leq d$.

    After diagonalizing the metric at point $x\in M$, that is
$f^{l}_{ij}\mid_{x}=0$ for all $l,i,$ and $j$, equation
$r(X,\frac{\partial}{\partial t})\mid_{x}=0$ of the Einstein
condition becomes
\begin{equation}
\sum_{1\leq i<j\leq
r}(h^{jj}-h^{ii})\{\dot{f}^{1}_{ij}\mid_{x}Q\left(\sum_{k=1}^{\frac{d}{2}}[Y^{i}_{k},Y^{j}_{k}]
,X\right)+\dot{f}^{3}_{ij}\mid_{x}Q\left(\sum_{k=1}^{\frac{d}{2}}[Y^{i}_{k},Y^{j}_{k+\frac{d}{2}}]
,X\right)
\end{equation}
\begin{equation*}
+\dot{f}^{2}_{ij}\mid_{x}Q\left(\sum_{k=\frac{d}{2}+1}^{d}[Y^{i}_{k},Y^{j}_{k-\frac{d}{2}}],
X\right)+\dot{f}^{4}_{ij}\mid_{x}Q\left(\sum_{k=\frac{d}{2}+1}^{d}[Y^{i}_{k},Y^{j}_{k}],X\right)
\}=0
\end{equation*}
for all $X\in\mathfrak{p}$.
    To simplify our presentation, we define vectors
\begin{equation*}
Z^{1}_{ij}=\sum_{k=1}^{\frac{d}{2}}[Y^{i}_{k},Y^{j}_{k}]_{\mathfrak{p}},~~~Z^{3}_{ij}=
\sum_{k=1}^{\frac{d}{2}}[Y^{i}_{k},Y^{j}_{k+\frac{d}{2}}]_{\mathfrak{p}},~~~Z^{2}_{ij}=
\sum_{k=\frac{d}{2}+1}^{d}[Y^{i}_{k},Y^{j}_{k-\frac{d}{2}}]_{\mathfrak{p}},
\end{equation*}
and
\begin{equation*}
Z^{4}_{ij}=\sum_{k=\frac{d}{2}+1}^{d}[Y^{i}_{k},Y^{j}_{k}]_{\mathfrak{p}}
\end{equation*}
for $1\leq i<j\leq r$.
    Given the above formalism, we can now prove the following
theorem.

\begin{theorem}
Let $(M,g)$ be an Einstein manifold of cohomogeneity-one under the
action of compact, connected group $G$ such that the isotropy
representation has at most $r$ equivalent \textbf{symplectic}
summands of dimension $d$ with all others distinct.  The Einstein
metric can be diagonalized globally if
\begin{equation*}
{\rm dim}({\rm span}\{Z^{l}_{ij}~|~1\leq l\leq4, 1\leq i<j\leq
r\})=2r(r-1).
\end{equation*}
\end{theorem}
\begin{proof}
The proof of this theorem is identical to the proof of the
analogous theorems in the orthogonal and unitary case, except that
in the symplectic case there are $2r(r-1)$ non-diagonal equations
which need to be eliminated.
\end{proof}

\section{The Case of a Trivial Isotropy Subgroup}
 Let $M$ be an $n$-dimensional manifold of cohomogeneity-one
under the action of group $G$ such that each principal orbit is a
copy of $G$, i.e. the isotropy subgroup $K$ is trivial.  In this
case, the Lie algebra $\mathfrak{g}$ decomposes into $n-1$
irreducible summands of dimension one
\begin{equation}
\mathfrak{g}=\mathfrak{p}_{1}\oplus...\oplus\mathfrak{p}_{n-1}.
\end{equation}
We choose a basis $\{Y^{i}\}_{i=1,...n-1}$ as described above.
(Note that since each summand is of dimension one, we have omitted
the subscript $k$.)

We have the following corollary to theorem (\ref{rsummandortho}) .

\begin{corollary}
Let $(M^{n},g)$ be an Einstein manifold of cohomogeneity-one under
the action of a compact group $G$ such that the isotropy subgroup
is trivial, i.e. $K=\{1\}$ . Equation
$r(X,\frac{\partial}{\partial t})=0$ of the Einstein condition can
be used to globally diagonalize the metric $g$ if and only if
${\rm dim} M=4$ and $G\cong SU(2)$.
\end{corollary}
\begin{proof}
As in the previous section, define
$Z_{ij}=[Y^{i},Y^{j}]\mid_{\mathfrak{p}}=[Y^{i},Y^{j}]$ for $i<j$.
Because each of the equivalent summands is one dimensional, each
is an orthogonal representation and we can apply theorem
(\ref{rsummandortho}).  In this case, $r=n-1$.  Therefore,
equation $r(X,\frac{\partial}{\partial t})=0$ of the Einstein can
be used to diagonalize the metric globally if and only if
$\textrm{dim}(\textrm{span}\{Z_{ij}\})=\frac{(n-1)(n-2)}{2}$.  But
the span of the set of $Z_{ij}$'s must be less than or equal to
the dimension of each principal orbit, $n-1$.  The inequality
$\frac{(n-1)(n-2)}{2}\leq(n-1)$ implies that $n\leq4$. Clearly,
the dimension of the cohomogeneity-one manifold is greater than
one, otherwise the principal orbits would be trivial.

    We have shown that the dimension of the manifold is two,
three, or four.  If the manifold has dimension two, then the
manifold is diagonal trivially and $r(X,\frac{\partial}{\partial t
})$ automatically.  If the manifold has dimension three, then the
group $G$ has dimension two.  The only compact Lie Group of
dimension two is abelian.  Therefore, equation
$r(X,\frac{\partial}{\partial t})$ of the Einstein condition is
satisfied automatically.  We cannot then use that equation to
diagonalize the metric.
    If the manifold has dimension four, then $G$ is of dimension
three.  The only compact Lie groups of dimension three (up to
finite quotient) are $T^{3}$ and $SU(2)$.  Since $T^{3}$ is
abelian, $r(X,\frac{\partial}{\partial t})=0$ automatically and we
cannot use this equation to diagonalize the metric.  Finally, we
have already seen how we can diagonalize an Einstein metric with
$SU(2)$ principal orbits.
\end{proof}

\chapter{Compact Einstein Four Manifolds with Torus Isometry}

    Having discussed a problem concerning Einstein metrics of cohomogeneity-one
in the previous chapter, we now consider the Einstein condition on a
manifold of cohomogeneity-two.  In particular, we specialize to the
case of a compact Einstein four manifold admitting a $T^{2}$
isometry. After discussing what is known about Einstein four
manifolds with large symmetry group, we review the topological
classification of smooth torus actions on compact, simply-connected
four manifolds due to Orlik and Raymond \cite{OandR}.  Also, we will
review the nuts and bolts formalism developed by Gibbons and Hawking
for retrieving topological data from the vanishing of Killing
vectors to the case of a compact four manifold with torus symmetry.
 Although the action of a $T^{2}$ isometry on a compact four manifold
need not be automatically orthogonally transitive, a result in
\cite{Kundt} tells us that the O'Neill tensor $A$ must vanish if the
metric is Einstein.  Using this result we are able to substantially
simplify the equations of the Einstein condition by making use of a
natural complex structure on the quotient space $M/T^{2}$.  After
this simplification, we see that the Einstein condition becomes four
second-order \textit{elliptic} equations of four functions in two
variables.
    After reviewing in detail all of the known examples of compact Einstein
four manifolds admitting a $T^{2}$-isometry we look at the special
case in which the Einstein metric is globally diagonal.  In this
special case, we derive some further consequences of the Einstein
condition.

\section{Compact Einstein Four Manifolds with Symmetry}

The Einstein condition on manifolds of dimension two and three is
equivalent to the constant curvature condition; therefore, the
Einstein condition on manifolds of dimension two and three reduces
to a topological problem.  Any manifold of dimension two admits an
Einstein metric since every two dimensional manifold admits a
metric of constant curvature \cite{Besse}.  A three dimensional
Einstein manifold has a universal covering diffeomorphic to
$\mathbb{R}^{3}$ or $S^{3}$.  At present there exists no
classification of Einstein four manifolds even in the compact
case, but some success has been achieved in classifying Einstein
four manifolds with large symmetry groups. Homogeneous four
manifolds are effectively classified by the following theorem.

\begin{theorem}\cite{Jensen}
Every homogeneous Einstein four manifold is symmetric.
\end{theorem}
  As all symmetric spaces were classified by Cartan, it is known
that the only compact homogeneous Einstein four manifolds are
$T^{4}$, $S^{4}$, $S^{2}\times S^{2}$, and $\mathbb{C}P^{2}$ as
well as certain finite quotients of these manifolds. The work of
Cartan together with the above theorem gives the following result.

\begin{theorem}
If $(M,g)$ is a compact homogeneous Einstein four manifold, then,
up to homothety, $M$ must be $T^{4}$, $S^{4}$, $\mathbb{R}P^{4}$,
$\mathbb{C}P^{2}$, $S^{2}\times S^{2}$,
$S^{2}\times\mathbb{R}P^{2}$, ($S^{2}\times S^{2}$)/$\{\pm
(1,1)\}$, or $\mathbb{R}P^{2}\times\mathbb{R}P^{2}$ endowed with
their standard metrics.
\end{theorem}

Since the Einstein condition on homogeneous four manifolds is
understood, it is logical to next consider the Einstein condition
on cohomogeneity-one four manifolds.  Cohomogeneity-one manifolds
were studied by Berard-Bergery who was able to classify all
compact Einstein cohomogeneity-one four manifolds, barring one
exception.

\begin{theorem}\cite{Madsen}
If $(M,g)$ is a compact cohomogeneity-one Einstein four manifold
which is not homogeneous and does not have principal orbits
isometric to $SU(2)$ then, up to homothety, $M$ is ($S^{2}\times
S^{2}$)/$\{\pm (\delta_{1},-1)\}$, ($S^{2}\times S^{2}$)/$\{\pm
(\delta_{2},-1)\}$, ($S^{2}\times\mathbb{R}P^{2}$)/$\{\pm
(\delta_{2},-1)\}$,
$\mathbb{C}P^{2}\sharp\overline{\mathbb{C}P^{2}}$, or
$\mathbb{C}P^{2}\sharp\mathbb{R}P^{2}$, where $\delta_{1}={\rm
diag}(-1,1,1)$ and  $\delta_{2}={\rm diag}(-1,-1,1)$.  The first
three manifolds are endowed with their standard metrics and the
last two are endowed with the D. Page metric or its
$\mathbb{Z}_{2}$ quotient
\end{theorem}
The Einstein condition on manifolds of cohomogeneity-one with
$SU(2)$ principal orbits gives rise to what are known as the
\textit{Bianchi IX} equations.  These have been studied
extensively \cite{DancerStrachan} but there exists no complete
classification of their solutions.

    Aside from the Bianchi IX case, therefore, compact cohomogeneity-one
Einstein four manifolds have been completely classified.  In order
to complete the classification of compact Einstein manifolds with
isometry group of dimension greater than or equal to one, it
remains to consider the cohomogeneity-two and cohomogeneity-three
cases. In the cohomogeneity-two case, since each principal orbit
must be a homogeneous two manifold, the only possible principal
orbit types are $S^{2}$, $\mathbb{R}P^{2}$, and $T^{2}$.  The
following theorem by Derdzinski eliminates the first two orbit
possibilities.

\begin{theorem} \cite{Derd}
If $(M,g)$ is a compact cohomogeneity-two Einstein four manifolds
with principal orbits isometric to $S^{2}$ or $\mathbb{R}P^{2}$
then $M$ must be homogeneous or cohomogeneity-one.

\end{theorem}
Therefore, we need only consider the case of $T^{2}$ principal
orbits.

Finally, a compact cohomogeneity-three manifold must have $S^{1}$
principal orbits. Below, we will study the Einstein condition on
compact cohomogeneity-two four manifolds with principal orbits
isometric to $T^{2}$.  We will also look briefly at
cohomogeneity-three four manifolds with $S^{1}$ principal orbits
by imposing the additional condition of orthogonal transitivity.

\section{Compact Simply-Connected Four Manifolds Admitting a Smooth $T^{2}$ Action}
   As in Chapter $2$, any compact Einstein four manifold with two
commuting Killing vectors must be $T^{4}$ or a manifold with
positive Einstein constant and finite fundamental group.  For the
sake of simplicity, we will, for the moment, restrict our
attention to the simply-connected case.  This restriction does not
greatly reduce the complexity of the problem but will allow for a
cleaner presentation.  Due to the work of Orlik and Raymond, there
exists a topological classification of smooth four manifolds
admitting a smooth torus action \cite{OandR}.  In this section, we
will give a brief review of their results following the summary
found in \cite{Joyce2}.

    Let $M$ be a smooth compact simply-connected four manifold admitting
a smooth, effective torus action and set $G=T^{2}=SO(2)\times
SO(2)$.
\begin{prop}\cite{OandR}
The action of $G$ on $M$ must have fixed points.
\end{prop}

At fixed points, the stabilizer is the whole of $G$.  In addition
to fixed points, $M$ admits one-dimensional subspaces upon which
the stability group is a subgroup isomorphic to $SO(2)=S^{1}$.
There are many such subgroups, each of which can be described by a
pair of integers.

    The group $G$ can be given coordinates $(\varphi,\psi)$  where
$0\leq\varphi<2\pi$ and $0\leq\psi<2\pi$.  Every subgroup of $G$
isomorphic to $SO(2)$ can be written as $G(m,n)$ where
$m\varphi+n\psi=0$ and $m$ and $n$ are coprime integers.  Any two
such subgroups $G(m,n)$ and $G(m',n')$ will generate the homology
of $G$ if and only if $mn'-m'n=\pm 1$ \cite{OandR}. We will study
not only the quotient space $M/G$ but also the quotient of $M$ by
the subgroups $G(m,n)$. Therefore, we state two propositions which
will be used below.
\begin{prop}\cite{OandR}
If $M$ is oriented, then $M/G(m,n)$ is a three space with boundary
for all coprime integer pairs $m$ and $n$.
\end{prop}
\begin{prop}\cite{OandR}
If $M$ is compact and oriented, then $M/G$ is a compact two
manifold with piecewise smooth boundary.
\end{prop}
Note that in the simply-connected case, $M/G$ will be a compact
two manifold with a piecewise smooth boundary polygon.  In
\cite{OandR}, it is demonstrated that $\partial M/G$ is a
one-dimensional manifold containing a finite number of isolated
fixed points.  Therefore, $M/G$ can be viewed as a polygon whose
interior points have $G$ orbits in $M$, whose edges away from the
vertices have $S^{1}$ orbits, and whose vertices have point orbits
in $M$.  Furthermore, the stabilizer of the orbits along the
interior of a given edge must be constant, while there is a 'jump'
in the stabilizer group at the endpoints \cite{OandR}.  Therefore
each edge of the polygon can be described by its stabilizer group
$G(m,n)$.  At the intersection of two edges, the stabilizers
$G(m,n)$ and $G(m',n')$ of the $S^{1}$ orbits of the two edges
must generate the homology of $G$. Therefore, at each vertex
$mn'-m'n=1$.

    The above results demonstrate that any smooth, effective
$T^{2}$ action on a compact, simply-connected four manifold can be
described by a finite number of coprime integer pairs,
$(m_{i},n_{i})_{1\leq i\leq k}$ where $(m_{i},n_{i})\sim
(-m_{i},-n_{i})$ and $m_{i}n_{i+1}-m_{i+1}n_{i}=\pm 1$ for $1\leq
i<k$ and $m_{k}n_{1}-m_{1}n_{1}=\pm 1$.  The following theorem
states the conditions under which two sets of integer pairs define
the same action on the same manifold.

\begin{theorem}\cite{OandR}
There exists an equivariant diffeomorphism between
$M=\{(m_{1},n_{1}),...,(m_{k},n_{k})\}$ and
$M'=\{(m'_{1},n'_{1}),...,(m'_{k'},n'_{k'})\}$ if and only if
$k=k'$ and, for some fixed integer $t$,
$(m_{i},n_{i})=(m'_{t+i},n'_{t+i})$ or
$(m_{i},n_{i})=(-m'_{t+i},-n'_{t+i})$ for $i=1,...,k$.
\end{theorem}

Given such a set of integer pairs, one can uniquely identify the
manifold corresponding to the action described.  Indeed, any such
manifold admitting a smooth $T^{2}$ action must be $S^{4}$ or
connected sums of $\mathbb{C}P^{2}$, $\overline{\mathbb{C}P^{2}}$,
and $S^{2}\times S^{2}$ \cite{OandR}.  This identification of
integer pair sets with the manifolds listed above is not a
bijection since, in general, any such manifold will admit many
inequivalent torus actions.  It is well known that $(S^{2}\times
S^{2})\sharp\mathbb{C}P^{2}$ is equivalent to
$\mathbb{C}P^{2}\sharp\overline{\mathbb{C}P^{2}}\sharp\mathbb{C}P^{2}$.
Therefore, aside from $S^{4}$ and $\sharp^{n}(S^{2}\times S^{2})$,
all compact four manifolds admitting a torus action are of the
form $k\mathbb{C}P^2\sharp l\overline{\mathbb{C}P^2}$.

    Although we will not here reproduce the entirety of the
classification, we note along with Joyce in \cite{Joyce} that
certain topological information can be readily obtained from a
collection of integer pairs defining a torus action.  Given an
action defined by $\{(m_{1},n_{1}),...,(m_{k},n_{k})\}$, because
the integer pairs are defined only up to sign, we can choose the
integer pairs so that $m_{j}\geq 0$ and $n_{j}>0$ whenever
$m_{j}=0$.  The Euler characteristic of the manifold, $M$,
admitting that action is given by
\begin{equation*}
\chi(M)=k,
\end{equation*}
where $k$ is the number of sides of the polygon forming the
boundary of the quotient manifold $M/G$.  The signature of the
manifold defined by the same torus action is
\begin{equation*}
\tau(M)=b_{+}-b_{-}=m_{1}n_{k}-m_{k}n_{1}+\sum^{k-1}_{j=1}m_{j+1}n_{j}-m_{j}n_{j+1}.
\end{equation*}

Recall that $\chi(S^{4})=2$ and $\tau(S^{4})=0$ and that
$\chi(\sharp^{n}(S^{2}\times S^{2}))=2n+2$ and
$\tau(\sharp^{n}(S^{2}\times S^{2}))=0$.  Finally, we recall that,
for $M=k\mathbb{C}P^2\sharp l\overline{\mathbb{C}P^2}$,
$\chi(M)=2+k+l$ and $\tau(M)=k-l$ \cite{Survey}.  The Euler
characteristic and the signature do not uniquely determine the
torus action or even the manifold itself.  For example,
$S^{2}\times S^{2}$ and
$\mathbb{C}P^{2}\sharp\overline{\mathbb{C}P^{2}}$ both have Euler
characteristic four and signature zero.  The difference is that
$S^{2}\times S^{2}$ is spin while
$\mathbb{C}P^{2}\sharp\overline{\mathbb{C}P^{2}}$ is not.

\begin{theorem}\cite{Donaldson2}, \cite{Freedman}
A smooth simply connected oriented four manifold is determined up to
homeomorphism by its Euler characteristic, $\chi$, its signature,
$\tau$, and whether it is spin or non-spin.
\end{theorem}

Since a given torus action can be associated to at most one
manifold, one must be able to distinguish between torus actions on
spin manifolds and those on non-spin manifolds.  By analyzing the
classification of Orlik and Raymond, it can be deduced that a
manifold is spin if and only if $m_{j+2}n_{j}-m_{j}n_{j+2}$ as
well as $m_{2}n_{k}-m_{k}n_{2}$ and $m_{1}n_{k-1}-m_{k-1}n_{1}$
are even.

    These topological invariants can be used to tell us which
manifolds admit which types of torus action.  The Euler
characteristic of a manifold $M$ is equal to the number of sides
of the polygon forming the boundary of $M/G$.  Therefore, it
immediately follows that two polygons with a different number of
sides cannot possibly bound the orbit space of the same manifold.
Below, when we endeavor to endow these manifolds with a
$T^{2}$-invariant metric, the signature of the manifold will
provide us with detailed information relating to the vanishing of
Killing vectors.

\section{The Einstein Equations}

     The classification of Orlik and Raymond was used to show that
a simply-connected compact four manifold admitting a metric with a
torus isometry must be $S^{4}$, $\sharp^{l}(S^{2}\times S^{2})$,
or $k\mathbb{C}P^{2}\sharp l\overline{\mathbb{C}P^{2}}$.  However,
the Hitchin-Thorpe inequality informs us that many of these
manifolds cannot admit Einstein metrics of any kind.

\begin{theorem}\cite{Hitchin}
If $M$ is a smooth compact oriented manifold admitting an Einstein
metric, the
\begin{equation}\label{Hitchin}
2\chi(M)\geq3|\tau(M)|
\end{equation}
with equality if and only if the manifold is locally hyper-Kahler.
\end{theorem}

Both $S^{4}$ and $\sharp^{l}(S^{2}\times S^{2})$ have zero
signature, so they satisfy the Hitchin-Thorpe inequality. However,
many of the other manifolds are excluded.  The Euler
characteristic of a manifold of type $k\mathbb{C}P^{2}\sharp
l\overline{\mathbb{C}P^{2}}$ is given by $\chi=2+k+l$ and the
signature is given by $\tau=k-l$.  Since no such manifold can be
locally hyper-Kahler, we have the following

\begin{corollary}
A manifold $k\mathbb{C}P^{2}\sharp l\overline{\mathbb{C}P^{2}}$
cannot admit an Einstein metric unless $4+5k>l>\frac{k-4}{5}$.
\end{corollary}

    When determining which of the remaining manifolds admit an Einstein metric
with a $T^{2}$ isometry, there is no reason, \textit{a priori},
that the action of the torus must be orthogonally transitive.
However, because of the commuting Killing vectors must vanish
somewhere on a compact manifold $M$, we have the following
proposition.
\begin{prop}\label{twist}\cite{Kundt}
For $M$ a compact Einstein four manifold with two commuting
Killing vectors, the invariant $A$ of the Riemannian submersion
$\pi:M\rightarrow M/G $ vanishes and, equivalently, the action of
the $G$ is orthogonally transitive.
\end{prop}
\begin{proof}
Take $X^{a}$ and $Y^{a}$ to be the two commuting Killing vectors.
Since $X$ and $Y$ commute, $\mathcal{L}_{Y}X=0$ which implies that
$\mathcal{L}_{Y}(X\wedge dX)=0$.  If we define $\phi=X\wedge dX$
and $\psi=Y\wedge dY$ ($\phi_{abc}=X_{[a}\nabla_{b}X_{c]}$ and
$\psi_{abc}=Y_{[a}\nabla_{b}Y_{c]}$) then $\mathcal{L}_{Y}\phi=0$
and $\mathcal{L}_{X}\psi=0$. The well-known \textit{twist scalars}
$\alpha$ and $\beta$ are defined to be
$\alpha=\ast(Y\wedge\phi)=Y_{[a}\phi_{bcd]}\epsilon^{abcd}$ and
$\beta=\ast(X\wedge\psi)=X_{[a}\psi_{bcd]}\epsilon^{abcd}$.

    In order to determine $\alpha$, we first calculate
$\nabla^{a}(Y_{[a}\phi_{bcd]})$.
\begin{equation*}
\nabla^{a}(Y_{[a}\phi_{bcd]})=\frac{1}{4}\nabla^{a}(Y_{a}\phi_{bcd}-Y_{b}\phi_{cda}+Y_{c}
\phi_{dab}-Y_{d}\phi_{abc}).
\end{equation*}
Knowing that $\nabla^{a}Y_{a}$ and $\mathcal{L}_{Y}\phi=0$, a
straightforward calculation shows that
\begin{equation*}
\nabla^{a}(Y_{[a}\phi_{bcd]})=\frac{1}{4}(-Y_{b}\nabla^{a}\phi_{cda}+Y_{c}\nabla^{a}\phi_{dab}
-Y_{d}\nabla^{a}\phi_{abc}).
\end{equation*}
Now,

\begin{equation*}
\nabla^{a}\phi_{abc}=\frac{1}{3}\nabla^{a}(X_{a}\nabla_{b}X_{c}+X_{b}\nabla_{c}X_{a}+X_{c}
\nabla_{a}X_{b}).
\end{equation*}
 Using the antisymmetry of $\nabla
X$ and the fact that $\nabla_{b}\nabla_{c}X_{d}=R_{abcd}X^{a}$, we
have
\begin{equation*}
\nabla^{a}\phi_{abc}=\frac{1}{3}(X^{a}R_{dabc}X^{d}+X_{b}R^{~a}_{d~
ca}X^{d}+X_{c}R^{~a}_{d~ab}X^{d}).
\end{equation*}
Using the Bianchi identities, we see that the first term vanishes
and the above becomes
\begin{equation*}
\nabla^{a}\phi_{abc}=\frac{1}{3}(-X_{b}R_{cd}X^{d}+X_{c}R_{bd}X^{d}).
\end{equation*}
Since $M$ is an Einstein manifold $R_{ab}=\lambda g_{ab}$ and we
see that $\nabla^{a}\phi_{abc}=0$ and therefore
$\nabla^{a}(Y_{[a}\phi_{bcd]})=0$.  This implies that
$\nabla_{a}\alpha=0$ and, equivalently, that $\alpha$ is constant.

    By the same argument, $\beta$ is seen to be a constant.
Indeed, since $X$ and $Y$ must vanish, we find that
$\alpha=\ast(Y\wedge\phi)=\beta=\ast(X\wedge\psi)=0$.  The
vanishing of the twist scalars is a necessary and sufficient
condition for the action of the torus to be orthogonally
transitive by \cite{Mason}.

\end{proof}
Let $X_{1}$ and $X_{2}$ be two independent commuting Killing
vectors.  As above, we set
$h=(h_{ij})=(g_{ab}X_{i}^{a}X_{b}^{j})=(X^{a}_{i}X_{ja})$ where the
entries $h_{ij}$ are functions of the coordinates of $M/G$.  The
previous theorem implies that on the open union of principal $T^{2}$
orbits, an Einstein metric can be written as
$g=g_{M/G}+h=\check{g}+h$, where $\check{g}$ is a smooth metric on
the two dimensional interior of $M/G$.  Below, we let $\nabla$
denote the Levi-Civita connection of $g$ and let $\check{\nabla}$
denote the Levi-Civita connection of $\check{g}$.  Now, because the
metric $h$ is positive definite on $M/G$, $\textrm{det}~h\geq0$.
More precisely, we know that $\textrm{det}(h)>0$ on the interior of
$M/G$ and $\textrm{det}(h)=0$ on $\partial M/G$.  More precisely,
$\textrm{det}(h)=0$ on the fixed point set of the actions of all the
$S^{1}$ subgroups of $T^{2}$.
\begin{prop}
For $(M,g)$ a compact four manifold with two commuting Killing
vectors, the Einstein condition is equivalent to the following
system of partial differential equations:
\begin{equation}\label{etorus3}
-\frac{1}{2}\check{\triangle}h_{ij}-\frac{1}{4}\frac{\langle\check{\nabla}{\rm
det}(h),\check{\nabla}h_{ij}\rangle}{{\rm
det}(h)}-\frac{1}{2}h_{ik}\langle\check{\nabla}h^{kl},\check{\nabla}h_{lj}\rangle=\lambda
h_{ij}
\end{equation}
\begin{equation}\label{etorus4}
-\frac{1}{2}h^{ij}\check{\nabla}_{a}\check{\nabla}_{b}h_{ij}-\frac{1}{4}\check{\nabla}_{a}
h^{ij}\check{\nabla}_{b}h_{ij}=-\check{r}_{ab}+\lambda\check{g}_{ab}=(\lambda-\frac{\check{s}}
{2})\check{g}_{ab}.
\end{equation}
\end{prop}
\begin{proof}
This is immediate from the calculations above giving the Einstein
equations for an orthogonally transitive action.  Because each
principal orbit is a torus, $\hat{r}=0$.  For the last line we
have used the fact that for any two-dimensional metric
$\check{g}$, the Ricci tensor can be written as
$\check{r}=(\frac{\check{s}}{2})\check{g}$.
\end{proof}
The Einstein condition on this class of manifold is a system of
partial differential equations in two variables.  This system must
satisfy a host of boundary conditions on the points of
$\partial(M/G)$. Indeed a great deal of information about the
structure of $M/G$ can be obtained by applying these boundary
conditions to the equations of the Einstein condition.  As it
stands, the Einstein condition on a compact four manifold with
$T^{2}$-symmetry is a system of six partial differential equation,
three in (\ref{etorus3}) and three in (\ref{etorus4}), in two
variables. We can show that, so long as \textit{all} of the
equations are satisfied on the boundary, then only \textit{four}
of these equations must be satisfied on the interior of $M/G$ for
the Einstein condition to hold.

    First, we write equations (\ref{etorus3}) and (\ref{etorus4}) in a
slightly different form by decomposing them into their traces and
their trace-free parts.  Let
\begin{equation*}
\rho=\sqrt{\textrm{det}(h)}.
\end{equation*}
We see that (\ref{etorus3}) can be written more compactly as
\begin{equation*}
-\frac{1}{2}h_{ik}\rho^{-1}\check{\nabla}_{a}(\rho
h^{kl}\check{\nabla}^{a}h_{lj})=\lambda h_{ij}.
\end{equation*}
Multiplying both sides of this equation by $h^{ki}$ we have
\begin{equation*}
\check{\nabla}_{a}(\rho
h^{kl}\check{\nabla}^{a}h_{jl})=-2\lambda\rho\delta^{k}_{j}.
\end{equation*}
This can be written more schematically as
\begin{equation}\label{etorus5}
\check{\nabla}_{a}(\rho h^{-1}\check{\nabla}^{a}h)=-2\lambda\rho
\textrm{Id}_{2}
\end{equation}
where $\textrm{Id}_{2}$ is the $2\times2$ identity matrix.  To
separate the trace and the trace-free parts of this equation
define matrix $K$ by
\begin{equation*}
h=\rho K.
\end{equation*}
This implies that $h^{-1}=\rho^{-1}K^{-1}$ and $\textrm{det}K=1$.
Substituting this into equation (\ref{etorus5}) implies that
\begin{equation*}
\check{\nabla}_{a}(K^{-1}\check{\nabla}^{a}(\rho K))=-2\lambda\rho
\textrm{Id}_{2}
\end{equation*}
which expands to
\begin{equation*}
(\check{\triangle}\rho)\textrm{Id}_{2}+\check{\nabla}_{a}(\rho
K^{-1}\check{\nabla}^{a}K)=-2\lambda\rho \textrm{Id}_{2}
\end{equation*}
Taking the trace of this equation gives
\begin{equation*}
2\check{\triangle}\rho+\textrm{tr}[\check{\nabla}_{a}(\rho
K^{-1}\check{\nabla}^{a}K)]=-4\lambda\rho
\end{equation*}
which is equivalent to
\begin{equation*}
2\check{\triangle}\rho+\check{\nabla}_{a}(\rho~\textrm{tr}[K^{-1}\check{\nabla}^{a}K])
=-4\lambda\rho.
\end{equation*}
Recall that
$\textrm{tr}[K^{-1}\check{\nabla}^{a}K]=(\textrm{det}K)^{-1}\check{\nabla}^{a}\textrm{det}K=0$
since $\textrm{det}K=1$.  Therefore, the trace of equation
(\ref{etorus3}) is equivalent to
\begin{equation}\label{trace1}
\check{\triangle}\rho=-2\lambda\rho.
\end{equation}
It is clear from the above calculation that the trace-free part of
equation (\ref{etorus3}) becomes
\begin{equation}\label{tracefree1}
\check{\nabla}_{a}(\rho K^{-1}\check{\nabla}^{a}K)=0.
\end{equation}
    We next perform similar calculations for equation (\ref{etorus4}).  Writing
this equation schematically gives
\begin{equation*}
\frac{1}{2}\textrm{tr}[h^{-1}\check{\nabla}_{a}\check{\nabla}_{b}h+\frac{1}{2}\check{
\nabla}_{a}h^{-1}\check{\nabla}_{b}h]=(\frac{\check{s}}{2}-\lambda)\check{g}_{ab}
\end{equation*}
which is equivalent to
\begin{equation*}
\frac{1}{2}\check{\nabla}_{a}\textrm{tr}[h^{-1}\check{\nabla}_{b}h]-\frac{1}{4}\textrm{tr}
[\check{\nabla}_{a}h^{-1}\check{\nabla}_{b}h]=(\frac{\check{s}}{2}-\lambda)\check{g}_{ab}.
\end{equation*}
Recalling that
$\textrm{tr}[h^{-1}\check{\nabla}_{b}h]=\rho^{-2}\check{\nabla}_{b}(\rho^{2})$,
this becomes
\begin{equation*}
\rho^{-1}\check{\nabla}_{a}\check{\nabla}_{b}\rho-\rho^{-2}\check{\nabla}_{a}\rho
\check{\nabla}_{b}\rho+\frac{1}{4}\textrm{tr}[h^{-1}\check{\nabla}_{a}hh^{-1}
\check{\nabla}_{b}h]=(\frac{\check{s}}{2}-\lambda)\check{g}_{ab}.
\end{equation*}
Once again we set $h=\rho K$ which implies that
\begin{equation*}
\frac{1}{4}\textrm{tr}[h^{-1}\check{\nabla}_{a}hh^{-1}\check{\nabla}_{b}h]=\frac{1}{4}
\textrm{tr}[\rho^{-2}\check{\nabla}_{a}\rho\check{\nabla}_{b}\rho
\textrm{Id}_{2}]+\frac{1}{4}\textrm{tr}[K^{-1}\check{\nabla}_{b}K]\rho^{-1}\check{\nabla}_{a}
\rho+\frac{1}{4}\textrm{tr}[K^{-1}\check{\nabla}_{a}K
K^{-1}\check{\nabla}_{b}K].
\end{equation*}
The second term vanishes, and equation (\ref{etorus4}) becomes
\begin{equation}\label{etorus6}
\rho^{-1}\check{\nabla}_{a}\check{\nabla}_{b}\rho-\frac{1}{2}\rho^{-2}\check{\nabla}_{a}\rho
\check{\nabla}_{b}\rho+\frac{1}{4}\textrm{tr}[K^{-1}\check{\nabla}_{a}K
K^{-1}\check{\nabla}_{b}K]=(\frac{\check{s}}{2}-\lambda)\check{g}_{ab}.
\end{equation}
Taking the trace of this equation implies that
\begin{equation}\label{trace2}
\rho^{-1}\check{\triangle}\rho-\frac{1}{2}\rho^{-2}\langle\check{\nabla}\rho,\check{\nabla}
\rho\rangle+\frac{1}{4}\textrm{tr}[K^{-1}\check{\nabla}_{a}K
K^{-1}\check{\nabla}^{a}K]=\check{s}-2\lambda.
\end{equation}
Combining this with equation (\ref{trace1}) gives
\begin{equation}\label{scalar}
\check{s}=-\frac{1}{2}\rho^{-2}\langle\check{\nabla}\rho,\check{\nabla}\rho\rangle+\frac{1}{4}
\textrm{tr}[K^{-1}\check{\nabla}_{a}K
K^{-1}\check{\nabla}^{a}K].
\end{equation}
For the moment, we will not write out the trace-free part of
equation (\ref{etorus4}).

    So far, we have considered $\check{g}$ in a coordinate invariant way.  For
our purposes, we wish to express $\check{g}$ in terms of
\textit{isothermal coordinates}.  As the quotient space is
two-dimensional, the induced metric $\check{g}$, like all other
metrics on a two-dimensional manifold \cite{Besse}, is conformally
flat. Therefore, we can set
\begin{equation*}
\check{g}=\Omega^{2}(dx^{2}+dy^{2})
\end{equation*}
where $\Omega=\Omega(x,y)$.  These are called \textit{isothermal
coordinates}.  A complex structure can be placed, locally, on
$M/G$ by setting
\begin{equation*}
z=x+iy.
\end{equation*}
We define the complex differential in the standard way,
\begin{equation*}
\partial=\frac{1}{2}(\frac{\partial}{\partial x}-i\frac{\partial}{\partial y}).
\end{equation*}
We wish to calculate the Einstein equations (\ref{etorus3}) and
(\ref{etorus4}), or more specifically their traces and their
trace-free parts, in terms of the function $\Omega$ and the
operators $\partial$ and $\bar{\partial}$.  These equations,
expressed in complex coordinates, are collected in the following
corollary to the above proposition.
\begin{corollary}
On a compact four-dimensional manifold $(M,g)$ admitting a
$T^{2}$-isometry, the Einstein condition is equivalent to the
following system of equations
\begin{equation}\label{comp1}
\partial\bar\partial\rho=-\frac{1}{2}\lambda\rho\Omega^{2}
\end{equation}
\begin{equation}\label{comp2}
\partial(\rho K^{-1}\bar\partial K)+\bar\partial(\rho K^{-1}\partial K)=0
\end{equation}
\begin{equation}\label{comp3}
\partial\bar\partial{\rm log}\Omega=\frac{1}{4}\frac{\partial\rho\bar\partial\rho}{\rho^{2}}
-\frac{1}{8}{\rm tr}[K^{-1}\partial K K^{-1}\bar\partial K]
\end{equation}
\begin{equation}\label{comp4}
8\frac{\partial{\rm
log}\Omega\partial\rho}{\rho}-4\frac{\partial\partial\rho}{\rho}+2\frac{\partial\rho\partial
\rho}{\rho^{2}}-{\rm
tr}[K^{-1}\partial K K^{-1}\partial K]=0
\end{equation}

where $\rho,$ $\Omega,$ and the matrix $K$ are as defined above.
\end{corollary}
\begin{proof}
The proof is a straightforward calculation; simply evaluate
equations (\ref{etorus3}) and (\ref{etorus4}) for the coordinates
given. We point out that (\ref{comp1}) is the equivalent to the
trace of (\ref{etorus3}) and (\ref{comp2}) is equivalent to the
trace-free part of (\ref{etorus3}).  Furthermore, equation
(\ref{comp3}) is obtained by subtracting (\ref{comp1}) from the
trace of equation (\ref{etorus4}).  Finally, (\ref{comp4}) is
equivalent to the trace-free part of (\ref{etorus4}).

    We note that in terms of the isothermal coordinates and the given complex
structure, $\check{s}=-8\partial\bar\partial\textrm{log}\Omega$.
\end{proof}

    For convenience we define $Q$ by
\begin{equation*}
Q=8\frac{\partial\textrm{log}\Omega\partial\rho}{\rho}-4\frac{\partial\partial\rho}{\rho}
+2\frac{\partial\rho\partial\rho}{\rho^{2}}-\textrm{tr}[K^{-1}\partial
K K^{-1}\partial K].
\end{equation*}
Equation (\ref{comp4}) is then equivalent to $Q=0$.  We have the
following proposition.

\begin{prop}
If (\ref{comp1}), (\ref{comp2}), and (\ref{comp3}) hold
(equivalently if (\ref{etorus3}) and the trace of (\ref{etorus4})
hold) then the complex quantity $\rho Q$ is holomorphic.  That is,
\begin{equation*}
\bar\partial(\rho Q)=0.
\end{equation*}
\end{prop}
\begin{proof}
From above,
\begin{equation*}
\rho
Q=8\partial\textrm{log}\Omega\partial\rho-4\partial\partial\rho+2\frac{\partial\rho\partial
\rho}{\rho}-\rho
\textrm{tr}[K^{-1}\partial K K^{-1}\partial K].
\end{equation*}
Taking the $\bar\partial$ derivative, and using the fact that
$\partial\bar\partial=\bar\partial\partial$, gives
\begin{equation*}
\bar{\partial}(\rho
Q)=8\partial\bar\partial\textrm{log}\Omega\partial\rho+8\partial\textrm{log}\Omega\partial
\bar\partial\rho-4\partial\partial\bar\partial\rho+4\frac{\partial\bar\partial\rho\partial
\rho}{\rho}-2\frac{\bar\partial\rho\partial\rho\partial\rho}{\rho^{2}}
\end{equation*}
\begin{equation*}
-\bar\partial\rho \textrm{tr}[K^{-1}\partial K K^{-1}\partial
K]-\rho \bar\partial\{\textrm{tr}[K^{-1}\partial K K^{-1}\partial
K]\}.
\end{equation*}
By equation (\ref{comp1}), we find that
\begin{equation*}
-4\partial\partial\bar\partial\rho=2\lambda\partial(\rho\Omega^{2})=2\lambda\Omega^{2}
\partial\rho+4\lambda\rho\Omega^{2}\partial
\textrm{log}\Omega,
\end{equation*}
\begin{equation*}
8\partial
\textrm{log}\Omega\partial\bar\partial\rho=-4\lambda\rho\Omega^{2}\partial
\textrm{log}\Omega,
\end{equation*}
and
\begin{equation*}
4\frac{\partial\bar\partial\rho\partial\rho}{\rho}=-2\lambda\Omega^{2}\partial\rho.
\end{equation*}
Equation (\ref{comp3}) implies that
\begin{equation*}
8\partial\bar\partial
\textrm{log}\Omega\partial\rho=2\frac{\partial\rho\partial\rho\bar\partial\rho}{\rho^{2}}-
\partial\rho
\textrm{tr}[K^{-1}\partial K K^{-1}\bar\partial K].
\end{equation*}

    Together these results imply that
\begin{equation*}
\bar\partial(\rho Q)=-\rho\bar\partial\{\textrm{tr}[K^{-1}\partial
K K^{-1}\partial K]\}-(\partial\rho)\textrm{tr}[K^{-1}\partial K
K^{-1}\bar\partial K]-(\bar\partial\rho)\textrm{tr}[K^{-1}\partial
K K^{-1}\partial K].
\end{equation*}
We calculate that
\begin{equation*}
\bar\partial\{\textrm{tr}[K^{-1}\partial K K^{-1}\partial
K]\}=2\textrm{tr}[\bar\partial (K^{-1}\partial K) K^{-1}\partial
K]
\end{equation*}
\begin{equation*}
=\textrm{tr}\{[\bar\partial (K^{-1}\partial
K)+\partial(K^{-1}\bar\partial K)] K^{-1}\partial K\}
\end{equation*}
\begin{equation*}
=\frac{1}{\rho}\textrm{tr}\{[\bar\partial (\rho K^{-1}\partial
K)+\partial(\rho K^{-1}\bar\partial K)] K^{-1}\partial
K\}-\frac{\partial\rho}{\rho}\textrm{tr}[K^{-1}\partial K
K^{-1}\bar\partial
K]-\frac{\bar\partial\rho}{\rho}\textrm{tr}[K^{-1}\partial K
K^{-1}\partial K].
\end{equation*}
    To justify the second equality we note that
\begin{equation*}
\textrm{tr}[\bar\partial (K^{-1}\partial K) K^{-1}\partial
K]=-\textrm{tr}[K^{-1}\bar\partial KK^{-1}\partial K
K^{-1}\partial K]+\textrm{tr}[K^{-1}\bar\partial\partial K
K^{-1}\partial K].
\end{equation*}
Since the operators $\partial$ and $\bar\partial$ commute and
$\textrm{tr}(ABC)=\textrm{tr}(BCA)=\textrm{tr}(CAB)$ we see that
\begin{equation*}
\textrm{tr}[\bar\partial (K^{-1}\partial K) K^{-1}\partial
K]=\textrm{tr}[\partial (K^{-1}\bar\partial K) K^{-1}\partial K].
\end{equation*}
This implies that
\begin{equation*}
\bar\partial(\rho Q)=-\textrm{tr}\{[\bar\partial (\rho
K^{-1}\partial K)+\partial(\rho K^{-1}\bar\partial K)]
K^{-1}\partial K\}.
\end{equation*}
Equation (\ref{comp2}) implies that the right-hand side vanishes.
This complete the proof.
\end{proof}

    Since the manifold $M$ is assumed to be compact, the quotient space $M/G$
is automatically compact.  By the results of Orlik and Raymond
discussed above, the quotient space $M/G$ is simply-connected.
Because $\rho Q$ is holomorphic, by Cauchy's Theorem,
\begin{equation*}
\rho Q|_{\partial(M/G)}=0\Rightarrow\rho Q|_{M/G}\equiv0.
\end{equation*}
Since $\rho>0$ on the interior of $M/G$ (this must be true for the
metric to be nondegenerate there), we see that
\begin{equation*}
\rho Q|_{M/G}\equiv0\Rightarrow Q|_{\textrm{int}(M/G)}=0
\end{equation*}
where int$(M/G)$ denotes the interior of the quotient space.  It
is clear that if $Q$ is identically zero on the interior of a
compact set, then it must be zero at the boundary as well; that
is,
\begin{equation*}
Q|_{\textrm{int}(M/G)}\Rightarrow Q|_{M/G}\equiv0.
\end{equation*}
This is precisely equation (\ref{comp4}) of the Einstein
condition.

\begin{theorem}
If equations (\ref{comp1})-(\ref{comp3}) of the Einstein condition
hold on a compact four manifold with $T^{2}$-symmetry then
equation (\ref{comp4}) need only hold at the boundary of $M/T^{2}$
for the metric to be Einstein everywhere.
\end{theorem}
In essence, we have seen that the Einstein condition on a compact
four dimensional manifold with $T^{2}$-symmetry reduces from a
system of six partial differential equations to a system of
\textit{four} partial differential equations in four variables along
with additional boundary conditions.

\section{Nuts and Bolts}
In \cite{GandH}, Gibbons and Hawking study Einstein four manifolds
admitting at least one Killing vector and reveal how topological
information of the manifold can be recovered from examining the
fixed point set a Killing vector.  On a four manifold, a given
Killing vector can vanish either on a two-surface, called a bolt,
or at a point, called a nut.  In the case of a four manifold
admitting a torus symmetry, there is a two dimensional family of
Killing vectors and any given Killing vector will vanish on a
collection of nuts and bolts.  After reviewing the theory
developed in \cite{GandH}, we will apply that theory to the
special case of a compact Einstein four manifold with two
commuting Killing vectors.

    Given an oriented four manifold $M$ admitting a one-parameter
isometry group, we denote by $X$ the Killing vector generating
that isometry group, $\mu :M\rightarrow M$.  At a fixed point $p$
of the isometric action, the Killing vector $X$ vanishes and the
map $\mu_{\ast}:T_{p}M\rightarrow T_{p}M$ becomes an isometry.  As
noted above, the matrix $\nabla X=X_{a;b}$ is an antisymmetric
$4\times 4$ matrix and at a fixed point, any such matrix must have
rank $0$, $2$, or $4$.  Were the rank of the matrix $\nabla X$ to
equal zero, then the rank would be zero everywhere and the action
of the isometry would be trivial. Therefore, we assume that the
rank is either $2$ or $4$ at a fixed point.

    If $\nabla X$ has rank two, then the Killing vector $X$
vanishes on a totally geodesic two dimensional submanifold whose
tangent space in $TM$ is invariant under the action of the map
$\mu_{\ast}$ and the submanifold of fixed points is referred to as
a \textit{bolt} . The two-dimensional distribution orthogonal to
the tangent space of this submanifold will be rotated by this
action and the period of that rotation will equal
$2\pi\kappa^{-1}$ where $\kappa$ is a rational number called the
\textit{surface gravity of the bolt}. We note along with Gibbons
and Hawking that $\kappa$ is the nonzero skew eigenvalue of
$\nabla X$ relative to an orthonormal frame.

    On the other hand, if $\nabla X$ has rank $4$ at a fixed point
$p$ of $\mu$, then $p$ must be isolated as no subspace of $T_{p}M$
is acted on trivially by the map $\mu_{\ast}$.  Any such fixed
point is called a \textit{nut} (the name refers to the fixed point
found in the Taub-NUT metric).  At a nut, the matrix $\nabla X$
will have two skew-eigenvalues $\kappa_{1}$ and $\kappa_{2}$.  The
map $\mu_{\ast}$ will rotate with periods $2\pi\kappa_{1}^{-1}$
and $2\pi\kappa_{2}^{-1}$ two orthogonal subspaces of $T_{p}$.
Nuts can be divided into sets of \textit{nuts} and
\textit{antinuts}. The point $p$ is referred to as a nut when
$\kappa_{1}\kappa_{2}$ is positive and as an antinut when
$\kappa_{1}\kappa_{2}$ is negative.

    Using fixed point theorems, Gibbons and Hawking were able to
express the Euler characteristic and the signature of the manifold
$M$ in terms of the number of nuts and antinuts and the Euler
characteristics of the bolts of a given Killing vector.  We will
here present the formulae for the Euler characteristic and the
signature without proof and refer the reader to \cite{GandH}.
Given a manifold $M$ with Killing vector $X$ vanishing on $N_{+}$
nuts, $N_{-}$ antinuts, and $n$ bolts $(B_{i})_{i=1,...,n}$ with
Euler characteristics $\chi_{i}$, the Euler characteristic of $M$
is given by
\begin{equation*}
\chi(M)=N_{+}+N_{-}+\sum_{i=1}^{n}\chi_{i}
\end{equation*}

and the signature of $M$ is given by

\begin{equation*}
\tau=N_{+}-N_{-}.
\end{equation*}

    We now want to apply this nuts and bolts formalism to the
case of a simply connected compact Einstein four manifold with two
commuting Killing vectors.  In the following discussion, a 'side'
of the polygon boundary of $M/G$ will be defined to include its
endpoints.  Given a side $A$ of the polygon, there is a Killing
vector, say $X_{1}$, which vanishes on the two sphere
$\pi^{-1}(A)$.  For simplicity, we will say that $X_{1}$ has a
bolt at $A$ when, strictly speaking, it actually has a bolt at
$\pi^{-1}(A)$.  Any other Killing vector must have a nut (or an
antinut) at each of the two vertices contained in that side. Any
two nuts located at the boundary points of a side of the polygon
can be seen to be contained in the bolt of the Killing vector
vanishing on that side.  For a fixed Killing vector $X$, a vertex
will be viewed as either an isolated nut or as part of a bolt,
i.e. a given vertex will be viewed as a nut of vector field $X$ if
and only if it is not part of a bolt of $X$.  Consequently, any
Killing vector of the torus will vanish on a certain collection of
nuts, antinuts, and bolts. From each such collection, one can
retrieve the topological information of the manifold and, of
course, any two collections associated to two different Killing
vectors on the same manifold must encode the same topological
data.

    Recalling the classification of Orlik and Raymond, any smooth torus action on a compact
simply connected smooth four manifold can be described as a set of
coprime integer pairs $(m_{i}, n_{i})_{i=1,...,k}$.  In the case
of an isometric torus action, each integer pair defines which
Killing vector of the torus is vanishing on the side of the
polygon associated to that pair.  For example, if side $A$ is
described by $(1,0)$ then it is Killing vector $X_{1}$ that
vanishes on side $A$.  Because the set of coprime integer pairs
describing a given torus action contains all of the topological
data of the manifold, the nuts and bolts data must be compatible
with the set of integer pairs.  In order to understand this
compatibility, we must be able to express the surface gravities of
the bolts in terms of the functions $F_{i}$ and their derivatives.

    Without loss of generality, we take $X_{1}$ to be the Killing
vector vanishing on side $A$ of $\partial M/G$.  Because $X_{1}$
has a bolt at $A$, the matrix $\nabla X_{1}$ has rank two.  If we
take $\partial_{1}$ and $\partial_{2}$ to be an orthonormal basis
for $M/G$ and take $\partial_{3}$ and $\partial_{4}$ to be an
orthonormal basis for the vertical distribution, then $\nabla
X_{1}$ will, on all of $M$ looks schematically like

\begin{equation*}
\nabla X_{1}=
\begin{pmatrix}
0 & 0 & -a & -b \\
0 & 0 & -c & -d \\
a & c & 0 & 0 \\
b & d & 0 & 0 \\
\end{pmatrix}.
\end{equation*}
The upper left and lower right hand two by two submatrices are
equal to zero because of orthogonal transitivity and the commuting
of the Killing vectors respectively.

    Gibbons and Hawking state that the 'surface gravity' of the
bolt at $A$ is given by the skew eigenvalue of $\nabla X_{1}$
relative to an orthonormal basis.  Put another way, the surface
gravity of the bolt at $A$, which we will label as $\kappa_{A}$,
satisfies the equation
\begin{equation}
\nabla_{a}X_{1b}\nabla^{a}X_{1}^{b}|_{A}=\kappa_{A}^{2}.
\end{equation}
    We see that $\nabla_{a}X_{1b}\nabla^{a}X_{1}^{b}|_{A}$ must be
a constant because for Killing vector $X_{1}$
\begin{equation*}
\nabla_{b}\nabla_{c}X_{1d}=R_{abcd}X_{1}^{a}
\end{equation*}
which vanishes when $X_{1}$ does.  In order to express the surface
gravity of the bolt at $A$ in terms of the $h_{ij}$ and their
derivatives we note that

\begin{equation*}
\nabla_{b}X_{1a}\nabla^{b}X_{1}^{a}|_{A}=\nabla^{b}(\frac{1}{2}\nabla_{b}h_{11})|_{A}.
\end{equation*}
This implies that
\begin{equation*}
\frac{1}{2}\triangle h_{11}|_{A}=\kappa^{2}_{A}.
\end{equation*}

On a general side of type $(m,n)$, we have the following:

\begin{prop}
The bolt corresponding to a side with stabilizer group defined by
the coprime integer pair $(m,n)$ has \textbf{surface gravity}
$\kappa$ satisfying the equation
\begin{equation}\label{T8}
m^{2}\triangle h_{11}+2mn\triangle h_{12}+n^{2}\triangle
h_{22}=2\kappa^{2}.
\end{equation}
\end{prop}

\section{Examples}

   At first blush, demanding a torus symmetry appears
extremely restrictive.  However, all known Einstein metrics on
compact four manifolds admitting a symmetry group of dimension
greater than or equal to one admit a torus isometry.  Aside from
the possible exception of the Bianchi IX case, all Einstein
metrics on compact four manifolds with symmetry group of dimension
at least two \textit{must} admit a torus isometry.  Explicit
Einstein metrics of manifolds of this type are few; in the compact
simply connected case only four such metrics are known (on
$S^{4}$, $\mathbb{C}P^{2}$, $S^{2}\times S^{2}$, and
$\mathbb{C}P^{2}\sharp\overline{\mathbb{C}P^{2}}$) and all of them
are homogeneous or cohomogeneity one.  Complex geometry has been
used to prove the existence of a $T^{2}$ invariant Kahler-Einstein
metric on $\mathbb{C}P^{2}\sharp3\overline{\mathbb{C}P^{2}}$ but
the metric is not known explicitly.  Below, we review in detail
those metrics that are known and examine the torus symmetry of
each.
\subsection{($S^{4}$, can)}
 The four sphere admits only one smooth torus action up to
equivariant diffeomorphism \cite{OandR}.  That action can be
described by the set of coprime integer pairs $\{((1,0), (0,1)\}$.
Indeed, this is the only action whose quotient space is bounded by
a polygon with two sides.  The canonical metric on the four sphere
can be written as
$d\sigma_{4}^{2}=dR^{2}+\textrm{sin}^{2}Rd\sigma_{3}^{2}$ where
$d\sigma_{3}^{2}$ is the standard metric on $S^{3}$.  If we use
the coordinates of the Hopf fibration then the metric becomes
\begin{equation*}
d\sigma^{2}_{4}=dR^{2}+\textrm{sin}^{2}Rd\theta^{2}+\textrm{sin}^{2}R\textrm{sin}^{2}\theta
d\varphi^{2}+\textrm{sin}^{2}R\textrm{cos}^{2}\theta d\psi^{2}.
\end{equation*}
$\partial_{\varphi}$ and $\partial_{\psi}$ are the two commuting
Killing vectors of the $T^{2}$ action.  Relative to these vectors
the matrices $h$ become

\begin{equation*}
h=
\begin{pmatrix}
{\rm sin}^{2}R{\rm sin}^{2}\theta & 0 \\
0 & {\rm sin}^{2}R{\rm cos}^{2}\theta \\
\end{pmatrix}.
\end{equation*}  It is easy to
check that these equations are satisfied by the above metric on
$S^{4}$ and that $\lambda=3$ and $\check{s}=2$.

    What is the structure of $M/G$?  The metric induced on the quotient space is
$dR^{2}+{\rm sin}^{2}Rd\theta^{2}$ which is the standard metric on
the two sphere.  The vertices are located at $R=0$ and $R=\pi$ as
the rank of $J$ becomes zero at those points and the sides are
located at $\theta=0$ and $\theta=\frac{\pi}{2}$.
    Note that $\partial_{\varphi}$ has a bolt at $\theta=0$ and
$\partial_{\psi}$ has a bolt at $\theta=\frac{\pi}{2}$.  Any other
Killing vector will have a nut and an antinut at $R=0$ and $R=\pi$
or vice versa.

\subsection{($\mathbb{C}P^{2}$, Fubini-Study Metric)}

Like the four sphere, two dimensional complex projective space
admits only one torus action up to equivariant diffeomorphism
\cite{OandR}. That torus action can be described by the set
$\{(0,1), (1,0), (1,1)\}$.  Note that transposing two of the
coprime integer pairs would give the torus action on
$\overline{\mathbb{C}P^{2}}$.  The Fubini-Study metric can be
written as
\begin{equation*}
ds^{2}=dR^{2}+{\rm sin}^{2}Rd\theta^{2}+{\rm sin}^{2}R{\rm
sin}^{2}\theta(1-{\rm sin}^{2}R{\rm
sin}^{2}\theta)d\varphi^{2}+{\rm sin}^{2}R{\rm
cos}^{2}\theta(1-{\rm sin}^{2}R{\rm cos}^{2}\theta)d\psi^{2}
\end{equation*}
\begin{equation*}
-2{\rm sin}^{4}R {\rm sin}^{2}\theta cos^{2}\theta d\varphi d\psi
\end{equation*}
where $\partial_{\varphi}$ and $\partial_{\psi}$ are the two
commuting Killing vectors.  The matrices $h$ and $F$ then become

\begin{equation*}
h=
\begin{pmatrix}
{\rm sin}^{2}R{\rm sin}^{2}\theta(1-{\rm sin}^{2}{\rm sin}^{2}\theta) & -{\rm sin}^{4}R
{\rm sin}^{2}\theta {\rm cos}^{2}\theta \\
-{\rm sin}^{4}R{\rm sin}^{2}\theta {\rm cos}^{2}\theta &
{\rm sin}^{2}R{\rm cos}^{2}\theta(1-{\rm sin}^{2}R{\rm cos}^{2}\theta) \\
\end{pmatrix}.
\end{equation*}
For this metric, the scalar curvature of the base manifold is
$\check{s}=2$ and the Einstein constant is $\lambda=6$.  The
quotient space of the this manifold is a triangle on the standard
two-sphere whose induced metric is the canonical sphere metric.
The sides of the triangle are located at $\theta=0$,
$\theta=\frac{\pi}{2}$, and $R=\frac{\pi}{2}$.

    Killing vector $\partial_{\phi}$ has a bolt at $\theta=0$,
Killing vector $\partial_{\psi}$ has a bolt at
$\theta=\frac{\pi}{2}$, and Killing vector
$\partial_{\phi}+\partial_{\psi}$ has a bolt at $R=\frac{\pi}{2}$.

\subsection{($S^{2}\times S^{2}$, can)}
  The standard Einstein metric on $S^{2}\times S^{2}$ can be
written as

\begin{equation*}
ds^{2}=dR^{2}+d\theta^{2}+{\rm sin}^{2}Rd\varphi^{2}+{\rm
sin}^{2}\theta d\psi^{2}
\end{equation*}
where $\partial_{\varphi}$ and $\partial_{\psi}$ are the two
commuting Killing vectors.  This metric is invariant under the torus
action defined by the set $\{(0,1), (1,0), (0,1), (1,0)\}$. Relative
to this action, the matrices $h$ become

\begin{equation*}
h=
\begin{pmatrix}
{\rm sin}^{2}R& 0 \\
0 & {\rm sin}^{2}\theta \\
\end{pmatrix}.
\end{equation*}
As in the case of the standard metric on $S^{4}$, $h_{12}\equiv 0$.
For this metric, $\check{s}=0$ and $\lambda=1$.

    The quotient space for this torus action is a rectangle whose
induced metric is the flat Euclidean metric on $\mathbb{R}^{2}$.
The sides of the rectangle are located at $R=0$, $R=\pi$,
$\theta=0$, $\theta=\pi$.

    Killing vector $\partial_{\varphi}$ has bolts at $R=0$ and $R=\pi$
and vector $\partial_{\psi}$ has bolts at $\theta=0$ and
$\theta=\pi$.
\subsection{($\mathbb{C}P^{2}\sharp\overline{\mathbb{C}P^{2}}$, D. Page Metric)}

    The Page metric on $\mathbb{C}P^{2}\sharp\overline{\mathbb{C}P^{2}}$ can be written as
\begin{equation*}
ds^{2}=\frac{3}{\lambda}(1+\nu^{2})\{\frac{1-\nu^{2}\textrm{cos}^{2}R}{3-\nu^{2}-\nu^{2}
(1+\nu^{2})\textrm{cos}^{2}R}dR^{2}+\frac{(1-\nu^{2}\textrm{cos}^{2}R)}{3+6\nu^{2}-\nu^{4}}
d\theta^{2}
\end{equation*}
\begin{equation*}
+\frac{1-\nu^{2}{\rm cos}^{2}R}{3+6\nu^{2}-\nu^{4}}{\rm
sin}^{2}\theta d\phi^{2}+\frac{3-\nu^{2}-\nu^{2}(1+\nu^{2}){\rm
cos}^{2}R}{(3+\nu^{2})^{2}(1-\nu^{2}{\rm cos}^{2}R)}{\rm
sin}^{2}R(d\psi-{\rm sin}^{2}\frac{1}{2}\theta d \phi)^{2}\}
\end{equation*}
where $\nu\approx0.2817$ is a solution to the equation
$\frac{4\nu(3+\nu^{2})}{3+6\nu^{2}-\nu^{4}}=1$ and
$\partial_{\phi}$ and $\partial_{\psi}$ are the two commuting
Killing vectors.  This metric is invariant under the torus action
defined by the set \linebreak $\{(0,1), (1,1), (0,1), (1,0)\}$.
The scalar curvature, $\check{s}$, of the quotient space, $M/G$,
is a function of $R$ only and is not constant.  It is
straightforward to calculate that
\begin{equation*}
\check{s}=\frac{2}{3}\lambda
\nu^2\frac{((\nu^4+\nu^6)\textrm{cos}^{6}R+(-3\nu^2-3\nu^4)\textrm{cos}^{4}R+(6-3\nu^2+3\nu^4)\textrm{cos}^{2}R-3
+\nu^2)}{((\nu^6+\nu^8)\textrm{cos}^{6}R+(-3\nu^4-3\nu^6)\textrm{cos}^{4}R+(3\nu^2+3\nu^4)\textrm{cos}^{2}R-1-\nu^2)}.
\end{equation*}

    The quotient space of this torus action is a four-sided
polygon with metric of nonconstant scalar curvature.  The sides of
the polygon are located at $R=0$, $R=\pi$, $\theta=0$, and
$\theta=\pi$.

    Killing vector $\partial_{\psi}$ has bolts at $R=0$ and $R=\pi$.
These geodesics are also lines of constant negative curvature of
the base.  We note that the scalar curvature on these two
geodesics is in fact equal.  Killing vector $\partial_{\phi}$ has
a bolt at $\theta=0$ and Killing vector
$\partial_{\phi}+\partial_{\psi}$ has a bolt at $\theta=\pi$.  The
curvature on these two geodesics is the same, though it is not
constant.  Rather, it symmetric about the $R=\frac{\pi}{2}$ axis
and is periodic.  It is negative at the endpoints $R=0$ and
$R=\pi$ and positive around $R=\frac{\pi}{2}$.

\section{The Diagonal Case}

    Above, we discussed the Einstein condition on compact four
manifolds with $T^{2}$ symmetry, i.e. with two commuting Killing
vectors.  Using a result in \cite{Kundt} we demonstrated that the
Einstein condition requires that the action be orthogonally
transitive.  This allowed us to write the metric as
$g=\check{g}+h$ as discussed above.  Using isothermal coordinates,
the two-dimensional metric $\check{g}$ is always diagonalizable.
However, it is not always possible to diagonalize the matrix $h$.
One could seek to classify all Einstein metrics on compact four
manifolds with two commuting Killing vectors such that
$h_{12}\equiv0$ on all of $M$.  As a first step towards this goal,
we present the Einstein equations in the diagonal case and give a
couple of preliminary results. Below, we again specialize to the
simply-connected case.

    From the topological classification of Orlik and Raymond
discussed above, any torus action can be described by a set of
coprime integer pairs $\{(m_{i}, n_{i})\}_{i=1,...,k}$.  For
$h_{12}$ to vanish on all of $M$, each $(m_{i}, n_{i})$ must equal
$(1, 0)$ or $(0, 1)$.  Furthermore, since the same Killing vector
can not vanish on adjacent sides, the number of coprime integer
pairs $k$ must be even and the pairs must alternate between $(1,
0)$ and $(0, 1)$.  Such a manifold would have Euler characteristic
$\chi=k$ and signature $\tau=0$.  Therefore, when $k=2$, $M=S^{4}$
and for $k>2$ and $l=k-2$, $M=\sharp^{l}(S^{2}\times S^{2})$.  The
standard metrics on $S^{4}$ and $S^{2}\times S^{2}$ described in
detail above are each diagonalizable.  No Einstein metrics are
known to exist on connected sums of two or more copies of
$S^{2}\times S^{2}$.

    The Einstein condition on such a manifold takes a particularly
simple form because of the vanishing of $h_{12}$.  To better
distinguish this case from the case of a general torus action we set
$h_{11}=f_{1}^{2}$ and $h_{22}=f_{2}^{2}$.  This substitution is
possible because the functions $h_{11}$ and $h_{22}$ must be
positive on the interior of $M/G$ for the metric $g$ to be positive
definite.

\begin{prop}
The metric on an Einstein four manifold admitting a
$T^{2}$-isometry such that the metric is globally diagonalizable
can be written as
$g=g_{M/G}+f_{1}^{2}d\phi^{2}+f_{2}^{2}d\psi^{2}$.  The Einstein
condition becomes the following system of differential equations
\begin{equation}\label{diag10}
\frac{\check{\triangle}f_{1}}{f_{1}}+\frac{\langle\check{\nabla}f_{1},\check{\nabla}
f_{2}\rangle}{f_{1}f_{2}}=-\lambda
\end{equation}
\begin{equation}\label{diag20}
\frac{\check{\triangle}f_{2}}{f_{2}}+\frac{\langle\check{\nabla}f_{1},\check{\nabla}
f_{2}\rangle}{f_{1}f_{2}}=-\lambda
\end{equation}
\begin{equation}\label{diag30}
\frac{\check{\triangle}f_{1}}{f_{1}}+\frac{\check{\triangle}f_{2}}{f_{2}}=\check{s}-2\lambda
\end{equation}
along with the additional boundary conditions.
\end{prop}

\section{$\sharp^{l}(S^{2}\times S^{2})$}

    The quotient space of the torus action on $S^{4}$ reviewed
above is a two-sided polygon on the two-sphere endowed with the
canonical metric.  This quotient metric has constant positive
curvature.  Similarly, the quotient space of the torus action on
$S^{2}\times S^{2}$ is a four-sided polygon on Euclidean space
endowed with its canonical metric.  This quotient metric has zero
curvature.

    Any $T^{2}$-invariant Einstein metric on $\sharp^{l}(S^{2}\times S^{2})$
which is globally diagonalizable must have a quotient space bounded
by a polygon with $2l+2$ sides. It is reasonable to ask whether the
quotient space of such an Einstein metric might, for $l\geq 2$, have
constant negative curvature. Unfortunately this supposition is too
optimistic.

\begin{prop}
The quotient space of a $T^{2}$-invariant Einstein metric on
$\sharp^{l}(S^{2}\times S^{2})$ for $l\geq 2$ which is globally
diagonalizable cannot have constant curvature.
\end{prop}
\begin{proof}
The quotient space cannot have constant nonnegative curvature
because the boundary polygon has more than four sides.  Setting
$V=f_{1}f_{2}$ we note that $V$ is nonzero on the interior of the
quotient space $M/G$ and zero on the boundary $\partial M/G$.
Therefore, $V$ must have a critical point on the interior of
$M/G$.  At that point
\begin{equation*}
\nabla V=f_{1}\nabla f_{2}+f_{2}\nabla f_{1}=0
\end{equation*}
which implies that
\begin{equation*}
f_{1}^{2}|\nabla f_{2}|^{2}+2f_{1}f_{2}\langle\nabla f_{1},\nabla
f_{2}\rangle +f_{2}^{2}|\nabla f_{1}|^{2}=0.
\end{equation*}

Equations (\ref{diag10}) and (\ref{diag20}) imply that
$\langle\nabla f_{1},\nabla
f_{2}\rangle=-\frac{\check{s}}{2}f_{1}f_{2}$.  Together with the
above we have
\begin{equation*}
f_{1}^{2}|\nabla f_{2}|^{2}-\check{s}f_{1}^{2}f_{2}^{2}
+f_{2}^{2}|\nabla f_{1}|^{2}=0.
\end{equation*}

On the interior of $M/G$, the functions $f_{1}$ and $f_{2}$ must
be nonzero in order to ensure smoothness at interior points.
Therefore, $\check{s}$ must be nonnegative.  This implies that
$\check{s}$ is nonconstant.
\end{proof}

    Though our prediction concerning the quotient space structure
of an Einstein metric failed, we are able to say something about
the bolts of such an Einstein metric should it exist.

    Let $A_{i}$ denote the $(l+1)$-sides of type $(1,0)$ and let $B_{i}$ denote
the $(l+1)$-sides of type $(0,1)$.  Let $\textrm{Area}(A_{i})$ and
$\textrm{Area}(B_{i})$ denote the areas of the bolts corresponding
to side $A_{i}$ and $B_{i}$.  If we set the period of both the
$S^{1}$-actions to be $2\pi$, then the surface gravity becomes
\begin{equation*}
\kappa_{A_{i}}^{2}=\kappa_{B_{i}}^{2}=1
\end{equation*}
for all $1\leq i\leq l+1$ \cite{GandH}.

    If we let $\eta_{A_{i}}$ and $\eta_{B_{i}}$ be the \textit{outward} pointing
unit normal vector fields to $A_{i}$ and $B_{i}$ respectively.
The surface gravity equation implies that
\begin{equation*}
\eta_{A_{i}}f_{1}=\eta_{B_{i}}f_{2}=-1.
\end{equation*}
Furthermore, along a side $B_{i}$, the Einstein condition requires
that $\check{\nabla}f_{1}$ be tangent to all $B_{i}$ and the
$\check{\nabla}f_{2}$ be tangent to all $A_{i}$.  We deduce
immediately that
\begin{equation*}
\eta_{A_{i}}f_{2}=\eta_{B_{i}}f_{1}=0.
\end{equation*}

    Consider equation (\ref{diag10}) of the Einstein condition.  This equation can be
written as
\begin{equation*}
f_{2}\check{\triangle}f_{1}+\langle\check{\nabla}f_{1},\check{\nabla}f_{2}\rangle=-\lambda
f_{1}f_{2}.
\end{equation*}
Taking the integral of both sides gives
\begin{equation}\label{Green1}
\int_{M/G}f_{2}\check{\triangle}f_{1}+\langle\check{\nabla}f_{1},\check{\nabla}f_{2}\rangle
d\textrm{vol}_{\check{g}}=-\lambda\int_{M/G}f_{1}f_{2}d\textrm{vol}_{\check{g}}.
\end{equation}
The right-hand side of this equation becomes
\begin{equation*}
-\lambda\int_{M/G}f_{1}f_{2}d\textrm{vol}_{\check{g}}=-\frac{\lambda}{4\pi^{2}}\textrm{Vol}(M)
\end{equation*}
where $\textrm{Vol}(M)$ is the total volume of the manifold.

    To analyze the left-hand side of (\ref{Green1}), we make use of Green's first identity.
Let $u$ and $v$ be smooth functions on $\Omega$ a compact region
with piecewise smooth boundary.  Let $(\partial\Omega)_{i}$ be the
$i$th codimension-one element of the boundary and let $\eta_{i}$
be the outward pointing unit normal vector field to the face
$(\partial\Omega)_{i}$. Green' first identity states that
\begin{equation*}
\int_{\Omega}v\triangle u+\langle\nabla u,\nabla v\rangle
d\textrm{vol}_{\Omega}=\sum_{i}\int_{(\partial\Omega)_{i}}v(\eta_{i}u)d\omega_{i}
\end{equation*}
where $d\omega_{i}$ is the area form of the face
$(\partial\Omega_{i})$.

    Using this identity, the left-hand side of equation (\ref{Green1}) becomes
\begin{equation*}
\int_{M/G}f_{2}\check{\triangle}f_{1}+\langle\check{\nabla}f_{1},\check{\nabla}f_{2}\rangle
d\textrm{vol}_{\check{g}}=\sum_{i=1}^{l+1}\int_{A_{i}}f_{1}(\eta_{A_{i}}f_{2})dA_{i}+
\sum_{i=1}^{l+1}\int_{B_{i}}f_{1}(\eta_{B_{i}}f_{2})dB_{i}
\end{equation*}
where $dA_{i}$ and $dB_{i}$ are the elements of lengths of sides
$A_{i}$ and $B_{i}$.  Recalling that $f_{1}|_{A_{i}}=0$ and
$\eta_{B_{i}}f_{2}=-1$.  We see
\begin{equation*}
\int_{M/G}f_{2}\check{\triangle}f_{1}+\langle\check{\nabla}f_{1},\check{\nabla}f_{2}\rangle
d\textrm{vol}_{\check{g}}=-\sum_{i=1}^{l+1}\int_{B_{i}}f_{1}dB_{i}=-\frac{1}{2\pi}
\sum_{i=1}^{l+1}\textrm{Area}(B_{i}).
\end{equation*}
Performing a similar calculation for equation (\ref{diag20}), we
see that
\begin{equation*}
2\pi\sum_{i=1}^{l+1}\textrm{Area}(A_{i})=2\pi\sum_{i=1}^{l+1}\textrm{Area}(B_{i})=\lambda
\textrm{Vol}(M)
\end{equation*}
in the Einstein case.  What we have shown is that the Einstein
condition requires the sum of the areas of the bolts of
$\partial_{\phi}$ to equal the sum of areas of the bolts of
$\partial_{\psi}$.  Moreover, each of these sums is a multiple of
the total volume of the manifold.  It would be interesting to see
if we could devise an analogue of this result in the non-diagonal
case.

\section{Einstein Four Manifolds with an Orthogonally Transitive
$S^{1}$ Action}

  Note that in the following discussion, we could consider the more general case of a manifold
with one Killing vector whose associated action is orthogonally
transitive.  Only when we restrict ourselves to the compact case
will it be necessary for the action to be that of an $S^{1}$
isometry.  A manifold with an orthogonally transitive $S^{1}=G$
action is a cohomogeneity $(n-1)$ manifold with Riemannian
submersion invariant $A\equiv 0$.  On such a manifold the metric
on an open dense submanifold can be written as
$g=\check{g}+f^{2}d\varphi^{2}$ where $\check{g}$ is the induced
metric on the quotient space, $\phi$ is the coordinate on $S^{1}$,
and $f$ is a function of the base manifold.  The quotient space
$M/G$ is an $(n-1)$ dimensional manifold and the Einstein
condition on $M$ becomes a system of partial differential
equations in the variables of quotient space.  As this case is
monotypic and $\hat{r}\equiv 0$, the Einstein condition is given
by the following two equations

\begin{equation}\label{ch3-1}
\frac{\check{\triangle}f}{f}=-\lambda
\end{equation}
\begin{equation}\label{ch3-2}
\check{r}-\frac{\check{\textrm{D}}\textrm{d}f}{f}=\lambda\check{g}.
\end{equation}
where $\check{\textrm{D}}\textrm{d}f$ is the Hessian of the metric
$\check{g}$.
    Taking the trace of (\ref{ch3-2}), we have
\begin{equation*}
\frac{\check{\triangle}f}{f}=\check{s}+(1-n)\lambda.
\end{equation*}
This implies that $\check{s}=(n-2)\lambda$; therefore, the scalar
curvature of the base manifold is a constant of the same sign as
the Einstein constant.  Given the above observations, one can
rephrase the Einstein condition in the following way

\begin{prop}
The Einstein condition on a manifold $M^{n}$ with an orthogonally
transitive $S^{1}$ isometry is equivalent to the equation
\begin{equation}\label{ch3-3}
\check{{\rm D}}{\rm d}f-\check{\triangle}f\check{g}-f\check{r}=0
\end{equation}

subject to the constraint that $\check{s}=(n-2)\lambda$ is the
Einstein constant.
\end{prop}

Equation (\ref{ch3-3}) also arises in a very different context.
For $M$ a compact manifold, the scalar curvature map is a
quasi-linear map of second order
\begin{equation*}
s:\mathcal{M}\rightarrow C^{\infty}M
\end{equation*}
\begin{equation*}
g\mapsto s_{g}
\end{equation*}
where $\mathcal{M}$ is the space of all Riemannian metrics which
is open in $\Gamma^{2}M$ the space of all symmetric bilinear
differential forms \cite{Besse}.

    The scalar curvature map has at every point $g$ of
$\mathcal{M}$ a differential

\begin{equation*}
s'_{g}:\Gamma^{2}M\rightarrow C^{\infty}M
\end{equation*}
such that
\begin{equation*}
s'_{g}h=-\triangle_{g}({\rm
tr}_{g}h)+\delta_{g}(\delta_{g}h)-g(r_{g},h)
\end{equation*}
where $\delta$ is the formal adjoint of $d$.

    The formal adjoint of the map $s'_{g}$, which is denoted by
$(s'_{g})^{\ast}$, is an overdetermined elliptic operator
\begin{equation*}
(s'_{g})^{\ast}:C^{\infty}M\rightarrow\Gamma^{2}M
\end{equation*}
defined by
\begin{equation*}
(s'_{g})^{\ast}f=\textrm{Dd}_{g}f-(\triangle_{g}f)g-fr_{g}.
\end{equation*}
Returning to equation (\ref{ch3-3}), we see that, if $M$ is an
Einstein manifold, then $f$ must lie in the kernel of
$(\check{s}'_{\check{g}})^{\ast}$.  We have the following
proposition of J. P. Bourguignon

\begin{prop}\cite{Bour}
If ${\rm ker}(\check{s}'_{\check{g}})^{\ast}\neq 0$ then either
$(M/G, \check{g})$, which is a manifold of dimension $(n-1)$ is
Ricci-flat and ${\rm
ker}(\check{s}'_{\check{g}})^{\ast}=\mathbb{R}\cdot \mathbb{I}$ or
the scalar curvature is a strictly positive constant
$\frac{\check{s}_{\check{g}}}{n-2}$
\end{prop}
\chapter{Kahler Toric Manifolds}

    Above, in our study of Einstein four manifolds with torus symmetry
we noted that the known $T^{2}$-invariant Einstein metrics on
$S^{2}\times S^{2},$ $\mathbb{C}P^{2}$, and
$\mathbb{C}P^{2}\sharp3\overline{\mathbb{C}P^{2}}$ are
Kahler-Einstein.  Compact Kahler manifolds of dimension $2n$
admitting a $T^{n}$ isometry, of which these Kahler-Einstein
metrics are examples, have been studied in great detail by
algebraic geometers, differential geometers, and physicists
because they have many novel geometric properties.  Their geometry
can be described using complex coordinates, by cones and fans as
in \cite{Fulton2}, or using symplectic Darboux coordinates, by
polytopes in $\mathbb{R}^{n}$ as in \cite{Guilleman},
\cite{Delzant}, and \cite{Abreu}.  Much work has been done to
classify all the Kahler-Einstein \cite{Mabuchi} and constant
scalar curvature Kahler metrics \cite{Donaldson} on these
manifolds but many open questions remain. In this chapter, we will
review the basic definitions of Kahler toric geometry following
the work of Delzant, Guillemin, and Abreu who discuss Kahler toric
geometry in symplectic coordinates and will study the Einstein and
extremal Kahler conditions on these manifolds.  While this section
contains mostly background material for later work, we also
present a simplification of the scalar curvature equation written
by Abreu in \cite{Abreu2}.

\section{Kahler Toric Geometry}

\begin{defn}
A \textbf{Kahler Toric Manifold} is a Kahler manifold
$(M^{2n},J,\omega)$ which admits a smooth Hamiltonian $T^{n}$
action preserving the Kahler structure.
\end{defn}

On the level of Riemannian geometry, such a manifold gives rise to
the following submersion (over an open dense set)
\begin{equation}
T^{n}\rightarrow M\rightarrow M/T^{n}
\end{equation}
with torus principal orbits.  Note that $M/T^{n}$ is not a smooth
manifold but is a stratified space.  In the general Riemannian
case, the O'Neill invariants $A$ and $T$ would be nontrivial.
However, in the Kahler toric case, we see that the number of
invariants of the submersion is greatly reduced.
\begin{theorem}
On a Kahler toric manifold, the O'Neill invariant $A$ of the
associated Riemannian submersion is equal to zero.
\end{theorem}

\begin{proof}
Take $\{\partial_{\phi_{i}}\}_{i=1,..,n}$ to be a basis for
$T^{n}$.  Because the torus action is Hamiltonian, there exists a
moment map
\begin{equation*}
\mu:M\rightarrow (\mathfrak{t}^{n})^{\ast}\cong
(\mathbb{R}^{n})^{\ast}
\end{equation*}
such that

\begin{equation*}
d\mu(v)(X)=\omega(X,v)=g(JX,v)
\end{equation*}
for all $v\in TM$ and $X\in\mathfrak{t}^{n}$.  Because $J$ is
one-to-one , a vector field $v$ is vertical if and only if
$d\mu({v})=0$.  Let  $v=\partial_{\phi_{i}}$ be such a vertical
vector field.  Therefore,
$g(JX,\partial_{\phi_{i}})=-g(X,J\partial_{\phi_{i}})=0$ for all
$X\in\mathfrak{t}^{n}$.  This implies that $J\partial_{\phi_{i}}$
is a horizontal vector field.  Therefore we may take
$\{-J\partial_{\phi_{i}}\}_{i=1..n}$ to be a basis for the
horizontal distribution.  To prove that $A=0$ we will demonstrate
that this is in fact a commuting basis implying that the
horizontal distribution is integrable.

    The vector fields $\partial_{\phi_{i}}$ are by assumption
automorphisms of the complex structure.  Such a vector field
satisfies the following

\begin{equation*}
\mathcal{L}_{\partial_{\phi_{i}}}J=0\Leftrightarrow
[\partial_{\phi_{i}},JX]=J[\partial_{\phi_{i}},X]
\end{equation*}
for all $X\in TM$.  From this equation we deduce immediately that
\begin{equation*}
[\partial_{\phi_{i}},J\partial_{\phi_{j}}]=[\partial_{\phi_{i}},\partial_{\phi_{j}}]=0
\end{equation*}
for all $i$ and $j$.  As $J$ is an integrable complex structure,
its Nijenhuis tensor, $N_{J}$, vanishes.  This implies that
\begin{equation*}
N_{J}(\partial_{\phi_{i}},\partial_{\phi_{j}})=[J\partial_{\phi_{i}},J\partial_{\phi_{j}}]-
J[J\partial_{\phi_{i}},\partial_{\phi_{j}}]-J[\partial_{\phi_{i}},J\partial_{\phi_{j}}]-
[\partial_{\phi_{i}},\partial_{\phi_{j}}]=0.
\end{equation*}
It follows immediately that
\begin{equation*}
[J\partial_{\phi_{i}},J\partial_{\phi_{j}}]=0
\end{equation*}
for all $i$ and $j$.

Therefore, the basis for the horizontal distribution is integrable
and the O'Neill tensor $A$, which is a measure of the failure of
the horizontal distribution to be integrable, is equal to zero.
\end{proof}
The above proof demonstrated that the basis
$\{-J\partial_{\phi_{i}}\}_{i=1,..n}$ is a coordinate basis for
the horizontal distribution; we can set
$\partial_{u_{i}}=-J\partial_{\phi_{i}}$.  Because the vector
fields $\partial_{\phi_{i}}$ and the complex structure $J$ are
invariant, we find that $\{\partial_{u_{i}}\}_{i=1,..n}$ is a
basis for the quotient space $M/T^{n}$.

    We are now in a position to describe the metric on a Kahler
toric manifold.  First, we set
$h_{ij}=g(\partial_{\phi_{i}},\partial_{\phi_{j}})$.  Because the
metric is invariant under the action of $J$, we see that
$h_{ij}=g(\partial_{u_{i}},\partial_{u_{j}})$.  We have shown that
a Kahler toric metric, $g$, is of the form
\begin{equation*}
g=
\begin{pmatrix}
h & 0 \\
0 & h \\
\end{pmatrix}
\end{equation*} where $h=(h_{ij})$ and the manifold $M$ has been
given the real coordinates
$(u_{1},...,u_{n},\phi_{1},...,\phi_{n})$.

The complex structure $J$ relative to this basis can be written as
\begin{equation*}
J=
\begin{pmatrix}
0 & -\textrm{Id}_{n} \\
\textrm{Id}_{n} & 0 \\
\end{pmatrix}
\end{equation*}
Complex coordinates on this manifold can then be defined as
$z_{i}=u_{i}+i\phi_{i}$ where the torus acts by
$t\cdot(u+i\phi)=u+i(\phi+t)$.

The Kahler form $\omega$ is given by
$\omega(\cdot,\cdot)=g(J\cdot,\cdot)$ which implies that
\begin{equation*}
\omega=
\begin{pmatrix}
0 & -h \\
h & 0 \\
\end{pmatrix}.
\end{equation*}
Because the metric is Kahler, there exists, locally, a potential
$\eta$ for the Kahler form.  The function $\eta$ can be chosen to
be a function of the $u_{i}$'s only since the Kahler form is
invariant under the action of the torus.  It is straightforward to
demonstrate that
\begin{equation*}
h_{kl}=\frac{\partial^{2}\eta}{\partial u_{k}\partial u_{l}}.
\end{equation*}

    Above, we have described the Kahler toric metric in terms of
standard complex coordinates.  To do this, we first fixed the
complex structure $J$ to be the standard skew-symmetric matrix and
arrived at different Kahler metrics by altering the Kahler
potential $\eta$.  Next, we will look at Kahler toric geometry
relative to a different set of coordinates by fixing the
symplectic structure and varying the complex structure $J$.  We
will see that the different complex structures can be
parameterized by a 'symplectic' potential.  This feature is unique
to Kahler toric geometry and the existence of these 'symplectic'
coordinates will greatly simplify curvature calculations.

    To effect this change in coordinates, we note that
$\check{g}=(h_{ij})$ in the complex coordinate system.  Now, the
moment map $\mu$ of the Hamiltonian torus action is given by
$\mu(u,\phi)=\frac{\partial\eta}{\partial u}$ \cite{Abreu}.  We
construct a new coordinate system
$(x_{1},...,x_{n},\phi_{1},...,\phi_{n})$ on $M$ by defining
$x_{i}=\frac{\partial\eta}{\partial u_{i}}$ and fixing the
coordinate $\phi$ (effectively we are mapping $\mathbb{R}^n\times
T^{n}$ to $\mathbb{R}^n\times T^{n}\cong P^{o}\times T^{n}$ where
$P^{o}$ is the interior of the image of the moment map).  Under
the action of this map the Kahler form $\omega$ becomes
$\omega=dx\wedge d\phi$.  The Kahler toric metric has been
transformed into Darboux coordinates.
    In these coordinates, the complex structure becomes

\begin{equation*}
J=
\begin{pmatrix}
0 & -h \\
h^{-1} & 0 \\
\end{pmatrix}
\end{equation*}
and the metric becomes
\begin{equation*}
g=
\begin{pmatrix}
h^{-1} & 0 \\
0 & h \\
\end{pmatrix}.
\end{equation*}
In effect, all we have done is to change the coordinate on the
quotient space $M/T^{n}$ so that
$\check{g}=(h_{ij}(x_{1},...,x_{n}))$.

    Guillemin, in \cite{Guilleman}, demonstrated that $h^{-1}$
itself derives from a potential in symplectic coordinates just as
$h$ does in complex coordinates.  That is
\begin{equation*}
h^{kl}=\frac{\partial^{2}\Phi}{\partial x_{k}\partial x_{l}}
\end{equation*}
where $\Phi=\Phi(x_{1},...,x_{n})$ does not depend on the
$\phi_{i}$ because the metric is invariant under the torus action.
Furthermore, the potentials $\eta$ and $\Phi$ are Legendre dual in
the sense that
$\eta(u)+\Phi(x)=\sum_{j}\frac{\partial\eta}{\partial
u_{j}}\frac{\partial\Phi}{\partial x_{j}}$.  We will refer to the
coordinate system $(x,\phi)$ as the symplectic coordinates and to
the potential $\Phi$ as the symplectic potential to distinguish it
from the Kahler potential $\eta$.  This notation is somewhat
misleading in that $\Phi$ determines the complex structure while
the symplectic structure is fixed in standard Darboux form.

    Thus far, all of our calculations have been local.  We now wish to study the global
structure of a Kahler toric manifold.  To do so, we ignore for the
moment the complex structure and concentrate only on the
symplectic structure.

    Let $(M^{2n},\omega)$ be a compact symplectic manifold with an effective Hamiltonian
action of $T^{n}$. Because the action is Hamiltonian there exists
a $T^{n}$ equivariant moment map
$\mu:M\rightarrow(\mathfrak{t}^{n})^{\ast}\cong(\mathbb{R}^{n})^{\ast}$.
It is well-known that the image of such a moment map is a convex
polytope $P=\mu(M)\subset\mathbb{R}^{n}$.  Delzant gave the
following characterization of the polytopes that are the moment
polytopes of symplectic toric manifolds.
\begin{defn}
A convex polytope $P$ in $\mathbb{R}^{n}$ is \textbf{Delzant} if
\item 1. $n$ edges meet at every vertex $p$
\item 2. the edges meeting at $p$ are rational, i.e. each edge is of the
form $p+tv_{i}$ where $v_{i}\in\mathbb{Z}^{n}$
\item 3. the $v_{1},...,v_{n}$ form a basis for $\mathbb{Z}^{n}$.
\end{defn}

    Delzant associated to every polytope satisfying the above condition a symplectic toric
manifold $(M_{P},\omega_{P})$.  That manifold is the symplectic
quotient of $\mathbb{R}^{2d}$ (where $d$ is the number of sides of
the polytope) endowed with its standard symplectic structure.
This procedure is described in \cite{Delzant} and we will not
review it here.   In this procedure, every Delzant polytope is in
one-to-one correspondence with a symplectic toric manifold.

    Viewed in this way, a polytope does not designate a compatible complex structure and
therefore does not itself determine a Kahler toric metric.  That
said, one can associate to each polytope a 'canonical' compatible
complex structure to the toric symplectic manifold.  This is done
by taking $\mathbb{R}^{2d}\cong\mathbb{C}^{d}$ with its standard
linear complex structure.  This complex structure descends via the
symplectic quotient to a compatible complex structure on the
symplectic toric manifold $(M_{P},\omega_{P})$. One can therefore
assign a 'canonical' Kahler toric metric to each polytope.  It is
important to emphasize that this 'canonical' metric is not the
only Kahler toric metric whose image under the moment map $\mu$ is
the given polytope $P$.  A different compatible $T^{n}$-invariant
complex structure would give a different Kahler toric metric.

    Delzant was able to determine the 'symplectic' potential $\Phi_{P}$ of the 'canonical'
Kahler toric metric associated to a polytope $P$.  After reviewing
the construction of this potential, we discuss the work of Abreu
which characterizes all other Kahler toric metrics that can be
associated to polytope $P$.

    One can describe a Delzant polytope by a set of inequalities
$\langle x, \mu_{m}\rangle\geq\lambda_{m}$ for $m=1,...,d$ where
$d$ is the number of faces of the polytope, $\lambda_{m}$ is a
constant, and $\mu_{m}$ is a primitive element of the lattice
$\mathbb{Z}^{n}\subset\mathbb{R}^{n}$.  Define the function
$l_{m}(x)$ by
\begin{equation*}
l_{m}(x)=\langle x,\mu_{m}\rangle-\lambda_{m}.
\end{equation*}
Delzant showed that the 'symplectic' potential
\begin{equation*}
\Phi_{P}=\frac{1}{2}\sum_{m=1}^{d}l_{m}\textrm{log}(l_{m})
\end{equation*}
determines a smooth Kahler toric metric on the manifold $M$ with
symplectic structure defined by the polytope $P$.  Note that there
are many other Kahler metrics associated to a given polytope just
as there are many complex structures compatible with the
symplectic structure defined by the polytope.  It is possible to
choose different potentials which will give rise to different
complex structures.

    It is worthwhile to pause to determine the significance of the
constants $\mu_{m}$ and $\lambda_{m}$.  Given a polytope defined
by a set of inequalities, what happens if we change the values of
a particular $\mu_{m}$ or a particular $\lambda_{m}$?  If one
chooses a different primitive element, $\mu_{m}$, of the integer
lattice, there are two possibilities.  First, the resulting
polytope may no longer be Delzant and may not give rise to a
smooth manifold. Second, if the resulting polytope is Delzant then
the corresponding symplectic structure is altered.  The associated
differentiable manifold may be the same or different but the
symplectic structure will have been modified.  On the other hand,
if we perturb slightly the value of $\lambda_{m}$ the
corresponding differentiable manifold is unchanged while the
cohomology class of the associated symplectic form $\omega$ is
perturbed slightly. Therefore, fixing a polytope
$\{\mu_{m},\lambda_{m}\}$ is equivalent to fixing a manifold, a
symplectic structure on the manifold, and the Kahler class of any
metric on the manifold.  Varying the potential $\Phi$ while
leaving fixed the polytope $P$ changes the complex structure of
the Kahler metric \cite{Guilleman}.

    As mentioned above, the 'canonical' Kahler metric constructed
by Delzant is not the only Kahler metric which can be associated
to a given polytope.  Abreu found a means of describing all the
compatible complex structures.

\begin{theorem}\label{Abreutheorem}\cite{Abreu}
Let $(M_{P},\omega_{P})$ be a toric symplectic manifold associated
to a Delzant polytope $P\subset\mathbb{R}^{n}$ and $J$ any
compatible complex structure.  Then $J$ is determined by a
potential $\Phi$ of the form
\begin{equation*}
\Phi=\Phi_{P}+\Psi
\end{equation*}
where $\Psi$ is smooth on $P$, the matrix
$h^{-1}=\left(\frac{\partial^{2}\Phi}{\partial x^{k}\partial
x^{l}}\right)$ is positive definite on $P^{0}$ and has determinant
\begin{equation*}
{\rm
det}(h^{-1})=\left[\delta(x)\prod_{m=1}^{d}l_{m}(x)\right]^{-1}
\end{equation*}
where $\delta$, which depends on $\Psi$, is smooth and strictly
positive on the whole of $P$.
\end{theorem}

To summarize, any Kahler toric metric can be described by a
symplectic potential which is the sum of the 'canonical' potential
associated to the polytope $P$ (which is determined by the
symplectic structure) and a function satisfying the conditions of
theorem (\ref{Abreutheorem}).

\section{Curvature Equations and
Conditions}

In this section, we will review the expression of the scalar
curvature of a Kahler toric manifold in symplectic coordinates
discovered by Abreu \cite{Abreu}.  We will also present a slightly
simplified version of the scalar curvature equation. Next, we will
study the Extremal Kahler condition and will derive the full
Einstein equation on a Kahler toric manifold

\subsection{The Scalar Curvature Equation}

    Before we calculate the scalar curvature in symplectic
coordinates, recall that the matrix $h_{kl}$, when expressed in
complex coordinates $(u,\phi)$, is derived from a potential $\eta$
such that  $h_{kl}=\frac{\partial^{2}\eta}{\partial u^{k}\partial
u^{l}}$.  Conversely, when expressed in symplectic coordinates,
$(x,\phi)$, the inverse matrix $h^{kl}$ is derived from a
potential $\Phi$ such that
$h^{kl}=\frac{\partial^{2}\Phi}{\partial x^{k}\partial x^{l}}$.
Given the definitions of $x$ and $u$ of the previous section we
readily derive the following relations
\begin{equation}\label{toricvector1}
h^{kl}\frac{\partial}{\partial u^{l}}=\frac{\partial}{\partial
x^{k}}
\end{equation}
and
\begin{equation}\label{toricvector2}
h_{kl}\frac{\partial}{\partial x^{k}}=\frac{\partial}{\partial
u^{l}}.
\end{equation}
Note: Here we are using a modified Einstein summation convention.
\textit{Whenever} an index is repeated twice on the same side of
an equality, we sum over that index.  We will use this modified
convention consistently throughout so there should be no chance
for confusion with the typical convention in which one only sums
over repeated indices when one of indices is raised and one is
lowered.

    It is well-known that in complex coordinates $(u,\phi)$, the
scalar curvature of a Kahler toric manifold can be written as
\begin{equation}\label{complexscalar}
S=-h^{ij}\frac{\partial^{2}\textrm{log}(\textrm{det}(h))}{\partial
u^{i}\partial u^{j}}.
\end{equation}
The following theorem gives an expression for the scalar curvature
in terms of the symplectic coordinates.

\begin{theorem}\cite{Abreu2}
The scalar curvature of a Kahler toric manifold can be expressed,
relative to the coordinates $(x,\phi)$ as
\begin{equation}\label{symplecticscalar}
S=-\frac{\partial^{2}h_{ij}}{\partial x^{i}\partial x^{j}}
\end{equation}
where $h_{ij}=h_{ij}(x_{1},...,x_{n})$ as above.
\end{theorem}
\begin{proof}
Using equation (\ref{toricvector2}) we see that
\begin{equation*}
S=-h^{ij}\frac{\partial^{2}\textrm{log}(\textrm{det}(h))}{\partial
u^{i}\partial u^{j}}=-\frac{\partial}{\partial
x^{i}}\left(\frac{\partial \textrm{log}(\textrm{det}(h))}{\partial
u^{i}}\right)
\end{equation*}
\begin{equation*}
=\frac{\partial}{\partial x^{i}}\left(\frac{\partial
\textrm{log}(\textrm{det}(h^{-1}))}{\partial
u^{i}}\right)=\frac{\partial}{\partial
x^{j}}\left(h_{ij}\frac{\partial
\textrm{log}(\textrm{det}(h^{-1}))}{\partial x^{i}}\right).
\end{equation*}
Following Abreu, we make use of the following identity: For
$V=(v^{ab})$ and $Y=(y^{ab})$ where $V$ and $Y$ are symmetric and
$Y$ is positive definite then
\begin{equation*}
v^{ab}\frac{\partial}{\partial
y^{ab}}(\textrm{log}(\textrm{det}Y))=\textrm{trace}(VY^{-1}).
\end{equation*}
Using this identity for $Y=h^{ab}$ and $V=h_{ij}$ we have
\begin{equation*}
S=\frac{\partial}{\partial x^{j}}\left(h_{ij}\frac{\partial
\textrm{log}(\textrm{det}(h^{-1}))}{\partial
x^{i}}\right)=\frac{\partial}{\partial
x^{j}}\left(h_{ij}\frac{\partial
\textrm{log}(\textrm{det}(h^{-1}))}{\partial h^{ab}}\frac{\partial
h^{ab}}{\partial x^{i}}\right)
\end{equation*}
\begin{equation*}
=\frac{\partial}{\partial x^{j}}\left(h_{ij}\frac{\partial h^{ab}
}{\partial x ^{i}}h_{ab}\right)=\frac{\partial}{\partial
x^{j}}\left(h_{ij}\frac{\partial h^{ib}}{\partial x
^{a}}h_{ab}\right)=-\frac{\partial}{\partial
x^{j}}\left(\frac{\partial h_{ij}}{\partial x
^{a}}h^{ib}h_{ab}\right)
\end{equation*}
this implies immediately that
\begin{equation*}
S=-\frac{\partial^{2}h_{ij}}{\partial x^{i}\partial x^{j}}.
\end{equation*}
\end{proof}

 We are able to give a modified version of the Abreu formula which
is somewhat simpler and will prove useful later.  First, we will
require the following lemma.

\begin{lemma} On a Kahler toric manifold, setting $h_{ij}=\rm{det}(h)\mathcal{M}_{ij}$,
\begin{equation}\label{adjugateid}
\frac{\partial}{\partial x^{i}}\left(\frac{1}{{\rm
det}(h)}h_{ij}\right)=\frac{\partial\mathcal{M}_{ij}}{\partial
x^{i}}=0
\end{equation}
for all $j$.
\end{lemma}
\begin{proof}
Note that $\mathcal{M}$ is the adjugate matrix of $h^{-1}$. Recall
that in symplectic coordinates $h^{kl}$ derives from a potential
$\phi$.  After expressing, $h^{-1}$ in terms of this potential,
and calculating the adjugate matrix $\mathcal{M}$ it is easy to
check that equation (\ref{adjugateid}) holds for $n=1$ or $2$.  In
the $n=1$ case, $h^{11}$ is the only term and
$\mathcal{M}_{11}=1$.  Equation (\ref{adjugateid}) holds easily.
In the $n=2$ case,
\begin{equation*}
h^{-1}=\begin{pmatrix} \Phi_{x_{1}x_{1}} & \Phi_{x_{1}x_{2}} \\
\Phi_{x_{1}x_{2}} & \Phi_{x_{2}x_{2}} \\
\end{pmatrix}
\end{equation*} and
\begin{equation*}
\mathcal{M}=\begin{pmatrix} \Phi_{x_{2}x_{2}} & -\Phi_{x_{1}x_{2}} \\
-\Phi_{x_{1}x_{2}} & \Phi_{x_{1}x_{1}} \\
\end{pmatrix}.
\end{equation*}
Equation (\ref{adjugateid}) again follows immediately.  A
straightforward but cumbersome induction argument shows that
equation (\ref{adjugateid}).  We omit the proof here.

\end{proof}

\begin{theorem}
The scalar curvature of a Kahler toric manifold can be written in
coordinate $(x,\phi)$ as
\begin{equation}\label{KTVScalar}
S=-\frac{1}{{\rm det}(h)}h_{ij}\frac{\partial^{2}{\rm
det}(h)}{\partial x^{i}\partial x^{j}}.
\end{equation}
\end{theorem}
\begin{proof}
 The
scalar curvature can be written as
\begin{equation*}
S=-\frac{\partial^{2}h_{ij}}{\partial x^{i}\partial
x^{j}}=-\frac{\partial^{2}((\textrm{det}(h))\mathcal{M}_{ij})}{\partial
x^{i}\partial x^{j}}=-\frac{\partial}{\partial
x^{i}}\left(\mathcal{M}_{ij}\frac{\partial
\textrm{det}(h)}{\partial x^{j}}\right)-\frac{\partial}{\partial
x^{i}}\left(\textrm{det}(h)\frac{\partial
\mathcal{M}_{ij}}{\partial x^{j}}\right).
\end{equation*}
From the proof of Abreu's formula we see that

\begin{equation*}
S=-\frac{\partial}{\partial x^{j}}\left(h_{ij}\frac{\partial
\textrm{log}(\textrm{det}(h))}{\partial
x^{i}}\right)=-\frac{\partial}{\partial
x^{i}}\left(\mathcal{M}_{ij}\frac{\partial
\textrm{det}(h)}{\partial x^{j}}\right).
\end{equation*}

    Using equation(\ref{adjugateid}), the scalar curvature becomes
\begin{equation*}
S=-\frac{\partial}{\partial
x^{i}}\left(\mathcal{M}_{ij}\frac{\partial
\textrm{det}(h)}{\partial
x^{j}}\right)=-\frac{\partial\mathcal{M}_{ij}}{\partial
x^{i}}\frac{\partial \textrm{det}(h)}{\partial
x^{j}}-\mathcal{M}_{ij}\frac{\partial^{2}\textrm{det}(h)}{\partial
x^{i}\partial
x^{j}}=-\mathcal{M}_{ij}\frac{\partial^{2}\textrm{det}(h)}{\partial
x^{i}\partial x^{j}}.
\end{equation*}
\end{proof}
Equation (\ref{KTVScalar}) and equation(\ref{complexscalar}) are
very similar in form; this similarity is obscured by the form of
Abreu because the equation he derived contains elements which
cancel upon inspection.  Therefore, it appears that the form of
the scalar curvature function given by (\ref{KTVScalar}) will be
more appropriate for performing calculations.

    Below, we will make frequent use of equation
(\ref{adjugateid}) while analyzing both Kahler toric and fiberwise
Kahler toric manifolds.

    Before beginning our discussion of the extremal Kahler
condition, we first calculate, for use later, the Laplacians of
the Kahler toric metric, $g$, and of the quotient metric,
$\check{g}$.  The general formula for the Laplacian of a metric
$g$ is

\begin{equation*}
\Delta f=\frac{1}{\sqrt{\textrm{det}(g)}}\frac{\partial}{\partial
y^{j}}\left(\sqrt{\textrm{det}(g)}g^{ij}\frac{\partial f}{\partial
y^{i}}\right)
\end{equation*}
where $(y_{1},...,y_{m})$ is a basis for an $m$-dimensional
manifold.  In symplectic coordinates, a Kahler toric metric $g$
can be written as $g=\begin{pmatrix}h^{-1} & 0 \\ 0 & h
\end{pmatrix}$.  This implies that $\textrm{det}(g)=1$.  Restricted to
invariant functions, the Laplacian becomes
\begin{equation}
\Delta f=\frac{\partial}{\partial x^{j}}\left(h_{ij}\frac{\partial
f}{\partial x^{i}}\right)
\end{equation}
for $f=f(x_{1},...,x_{n})$.  Defining the matric $\mathcal{M}$ as
above and recalling equation (\ref{adjugateid}), we have
\begin{equation*}
\Delta f=\frac{\partial}{\partial
x^{j}}\left(\textrm{det}(h)\mathcal{M}_{ij}\frac{\partial
f}{\partial
x^{i}}\right)=\textrm{det}(h)\mathcal{M}_{ij}\frac{\partial^{2}f}{\partial
x^{i}\partial x^{j}}+\mathcal{M}_{ij}\frac{\partial
\textrm{det}(h) }{\partial x^{i}}\frac{\partial f}{\partial
x^{j}}+\textrm{det}(h)\frac{\partial\mathcal{M}_{ij}}{\partial
x^{i}}\frac{\partial f }{\partial x^{j}}.
\end{equation*}
From this we immediately have that
\begin{equation}\label{Laplacian}
\Delta f=h_{ij}\frac{\partial^{2}f}{\partial x^{i}\partial
x^{j}}+\frac{1}{\textrm{det}(h)}h_{ij}\frac{\partial
\textrm{det}(h)}{\partial x^{i}}\frac{\partial f}{\partial x^{j}}.
\end{equation}
Next, we calculate $\check{\Delta}$, the Laplacian of the quotient
metric $\check{g}=(h^{-1})$. Using the general formula for the
Laplace operator, we see that
\begin{equation*}
\check{\Delta}
f=\frac{1}{\sqrt{\textrm{det}(h^{-1})}}\frac{\partial}{\partial
x^{j}}\left(\sqrt{\textrm{det}(h^{-1})}h_{ij}\frac{\partial
f}{\partial
x^{i}}\right)=\sqrt{\textrm{det}(h)}\frac{\partial}{\partial
x^{j}}\left(\frac{1}{\sqrt{\textrm{det}(h)}}\textrm{det}(h)\mathcal{M}_{ij}\frac{\partial
f}{\partial x^{i}}\right)
\end{equation*}Using equation (\ref{adjugateid})
\begin{equation*}
\check{\Delta}
f=\textrm{det}(h)\mathcal{M}_{ij}\frac{\partial^{2}f}{\partial
x^{i}\partial x^{j}}+\frac{1}{2}\mathcal{M}_{ij}\frac{\partial
\textrm{det}(h)}{\partial x^{j}}\frac{\partial f}{\partial x^{i}}.
\end{equation*}

Therefore,
\begin{equation}\label{quotLaplace}
\check{\Delta} f=h_{ij}\frac{\partial^{2}f}{\partial x^{i}\partial
x^{j}}+\frac{1}{2}\frac{1}{\textrm{det}(h)}h_{ij}\frac{\partial
\textrm{det}(h)}{\partial x^{j}}\frac{\partial f}{\partial x^{i}}.
\end{equation}
This formula allows us to prove the following proposition
\begin{prop}
On a Kahler toric manifold with the metric defined in symplectic
coordinates as above, we have the following identity
\begin{equation}
\check{\Delta}\left(\sqrt{{\rm
det}(h)}\right)=-\left(\frac{S}{2}\right)\sqrt{{\rm det}(h)}.
\end{equation}
Therefore, on a manifold of constant scalar curvature $S$ the
function $\sqrt{{\rm det}(h)}$ is an eigenfunction of the quotient
Laplacian with eigenvalue $\frac{S}{2}$.
\end{prop}
\begin{proof} Equation (\ref{quotLaplace}) acting on function
$\sqrt{\textrm{det}(h)}$ implies that
\begin{equation*}
\check{\Delta}\sqrt{{\rm
det}(h)}=h_{ij}\frac{\partial^{2}\sqrt{\textrm{det}(h)}}{\partial
x^{i}\partial
x^{j}}+\frac{1}{2}\frac{1}{\textrm{det}(h)}h_{ij}\frac{\partial
\textrm{det}(h)}{\partial x^{j}}\frac{\partial
\sqrt{\textrm{det}(h)}}{\partial x^{i}}.
\end{equation*}
Expanding the derivatives and noting the equation
(\ref{KTVScalar}) for the scalar curvature finishes the proof.

\end{proof}
\subsection{The Extremal Kahler Condition}

    As mentioned in the introduction, Einstein metrics are critical points of
the scalar curvature functional
$\textbf{S}:g\rightarrow\int_{M}S_{g}\textrm{vol}_{g}$. Calabi
introduced another functional on metrics restricted to a fixed
Kahler class.  The functional is defined to be
\begin{equation}
\textbf{Ss}:g\rightarrow\int_{M}s_{g}^{2}\textrm{vol}_{g}.
\end{equation}
\begin{defn}
Critical points of the functional \textbf{Ss} are called
\textbf{Extremal Kahler Metrics}.
\end{defn}
By deriving the Euler-Lagrange equation for this functional,
Calabi proved the following theorem.
\begin{theorem}\cite{Calabi}
A Kahler metric is extremal if and only if the gradient $\nabla S$
is a holomorphic vector field.
\end{theorem}
It follows immediately from this theorem that any constant scalar
curvature metric is in fact extremal Kahler.  Also, Calabi found
that, if one exists, an extremal Kahler metric is unique within a
given Kahler class \cite{Calabi}.

What is the extremal Kahler condition on a Kahler toric manifold?
In complex coordinates, the metric can be written as
$g=\begin{pmatrix} h & 0 \\ 0 & h\end{pmatrix}$ relative to
coordinates $(u,\phi)$. The complex coordinates themselves can be
written as $z_{i}=u_{i}+i\phi_{i}$.  Since the metric is invariant
under the action of $T^{n}$, the scalar curvature, $S$, is a
function of the $u_{i}$ only.  The gradient of $S$ can be written
\begin{equation*}
\nabla S=\left(h^{kl}\frac{\partial S}{\partial
u^{l}}\right)\frac{\partial}{\partial z_{k}}.
\end{equation*}
Since the coefficients, $h^{kl}\frac{\partial S}{\partial u^{l}}$,
of the above equation are real, the vector field $\nabla S$ is
holomorphic if and only if these coefficients are constant.  The
extremal Kahler condition is equivalent to the following system of
equations:
\begin{equation*}
h^{kl}\frac{\partial S}{\partial u^{l}}=\frac{\partial S}{\partial
x^{k}}=\alpha_{k}
\end{equation*}
where the $\alpha_{i}$ are constants.

The above proves the following theorem due to Abreu:

\begin{theorem}\cite{Abreu2}
A Kahler Toric metric is Extremal Kahler if and only if the scalar
curvature satisfies
\begin{equation}
S=\sum_{k=1}^{n}\alpha_{k}x_{k}+\beta
\end{equation}
where $\beta$ and the $\alpha_{k}$ are constants and the $x_{i}$
are as defined above.
\end{theorem}
Later, we will return to the question of finding extremal Kahler
toric metrics and we discuss the intimate connection between
extremal Kahler and Kahler-Einstein four manifolds discovered by
Derdzinski.

\subsection{The Einstein Condition}

To calculate the Einstein condition on a Kahler toric manifold in
symplectic coordinates, we begin by looking at the Einstein
condition in complex coordinates $(z_{1},...,z_{n})$ where
$z_{i}=u_{i}+i\phi_{i}$.  On a Kahler manifold, the Ricci form,
$\rho(\cdot,\cdot)=r(J\cdot,\cdot)$, is given by the expression
\begin{equation*}
\rho=-i\partial\overline{\partial}\textrm{log}(\textrm{det}(g_{a\overline{\beta}})).
\end{equation*}
The Kahler-Einstein condition becomes equivalent to the equality
\begin{equation*}
\rho=\lambda\omega.
\end{equation*}
It is straightforward to evaluate this equation on a Kahler toric
manifold.   A Kahler-Einstein toric manifold must satisfy the
equation
\begin{equation}\label{KTMComEin}
\frac{\partial^{2}\textrm{log}(\textrm{det}(h))}{\partial
u^{i}\partial u^{j}}=-2\lambda h_{ij}
\end{equation} for all $i$ and $j$, where $\lambda$ is the Einstein constant.

Converting this equation into symplectic coordinates gives the
following theorem.
\begin{theorem}
On a Kahler toric manifold, the Einstein condition, when expressed
in symplectic coordinates $(x,\phi)$ becomes the equation
\begin{equation}
h_{ik}\frac{\partial}{\partial x^{k}}\left(h_{lj}\frac{\partial
{\rm log}({\rm det}(h))}{\partial x^{l}}\right)=-2\lambda h_{ij}.
\end{equation} for all $i$ and $j$.
\begin{proof}
From equations (\ref{KTMComEin}) and (\ref{toricvector2}), this
formula follows immediately.
\end{proof}
\end{theorem}

This form of the Einstein condition is somewhat cumbersome and we
will simplify it by making use of the adjugate matrix
$\mathcal{M}$ defined above.

\begin{corollary}
The Einstein condition on a Kahler toric manifold is equivalent to
the equation
\begin{equation}\label{KTMEinSymp}
\mathcal{M}_{ik}\frac{\partial}{\partial
x^{k}}\left(\mathcal{M}_{lj}\frac{\partial {\rm det}(h)}{\partial
x^{l}}\right)=-2\lambda\mathcal{M}_{ij}
\end{equation} for all $i$ and $j$.
\end{corollary}
\begin{proof}
Recall that the matrix $\mathcal{M}$ is defined by the expression
$h_{kl}=\textrm{det}(h)\mathcal{M}_{kl}$.  Using the equation from
the previous theorem, we have
\begin{equation*}
\textrm{det}(h)\mathcal{M}_{ik}\frac{\partial}{\partial
x^{k}}\left(\textrm{det}(h)\mathcal{M}_{lj}\frac{\partial
\textrm{log}(\textrm{det}(h))}{\partial x^{l}}\right)=-2\lambda
\textrm{det}(h)\mathcal{M}_{ij}.
\end{equation*}
The Corollary follows immediately.

\end{proof}

\chapter{Einstein and Extremal Kahler Metrics on Kahler Toric Manifolds}

    In this chapter, we study the Einstein and extremal Kahler conditions
on Kahler toric manifolds, with particular emphasis on real
dimensions four and six. First, we demonstrate how to integrate
the equations of the Einstein condition to obtain a Legendre
transform of the Monge-Ampere equations.  Next, we study the
Kahler-Einstein condition on
$\mathbb{C}P^{2}\sharp3\overline{\mathbb{C}P^{2}}$ and write down
explicitly the Kahler-Einstein metric on
$\mathbb{P}(\mathcal{O}_{\mathbb{C}P^{1}\times\mathbb{C}P^{1}}\oplus\mathcal{O}_{\mathbb{C}
P^{1}\times\mathbb{C}P^{1}}(1,-1))$
in symplectic coordinates.  Finally, after reviewing the
construction of extremal Kahler metrics on
$\mathbb{C}P^{2}\sharp\overline{\mathbb{C}P^{2}}$, we construct a
new two parameter family of extremal Kahler metrics on
$\mathbb{P}(\mathcal{O}_{\mathbb{C}P^{1}\times\mathbb{C}P^{1}}\oplus\mathcal{O}_{\mathbb{C}
P^{1}\times\mathbb{C}P^{1}}(1,-1))$
which contains two one-parameter families of constant scalar
curvature metrics.
\section{Kahler-Einstein Toric Manifolds}

    Having calculated the Einstein condition on a Kahler toric manifold, we know
demonstrate how to integrate those equations, and reduce the
Einstein condition to a single equation involving only the
potential and its first and second derivatives.

\subsection{Integration of the Einstein Equations}

    As calculated above, the Einstein condition on a Kahler toric metric in arbitrary
dimension can be written as
\begin{equation*}
\mathcal{M}_{ik}\frac{\partial}{\partial
x^{k}}\left(\mathcal{M}_{lj}\frac{\partial
\textrm{det}(h)}{\partial x^{l}}\right)=-2\lambda\mathcal{M}_{ij}
\end{equation*}
where $\mathcal{M}=\frac{1}{\textrm{det}(h)}h$.  Note that on the
interior of the polytope $\mathcal{M}$ is a positive definite
matrix and therefore invertible.  Multiplying both sides of this
equations by $\mathcal{M}^{mi}$ gives
\begin{equation*}
\mathcal{M}^{mi}\mathcal{M}_{ik}\frac{\partial}{\partial
x^{k}}\left(\mathcal{M}_{lj}\frac{\partial
\textrm{det}(h)}{\partial
x^{l}}\right)=-2\lambda\mathcal{M}^{mi}\mathcal{M}_{ij}.
\end{equation*}
This implies that
\begin{equation*}
\delta^{m}_{k}\frac{\partial}{\partial
x^{k}}\left(\mathcal{M}_{lj}\frac{\partial
\textrm{det}(h)}{\partial x^{l}}\right)=-2\lambda\delta^{m}_{j}.
\end{equation*}
Integrating these equations gives
\begin{equation*}
\mathcal{M}_{lj}\frac{\partial \textrm{det}(h)}{\partial
x^{l}}=-2\lambda x_{j}+C_{j}
\end{equation*}
for all $1\leq j\leq n$ where $2n$ is the dimension of the Kahler
toric manifold and the $c_{j}$ are constants. With the first
integration accomplished, we attempt the second.

    We begin, as before, by multiplying both sides of the equation by $\mathcal{M}^{mj}$
giving
\begin{equation*}
\mathcal{M}^{mj}\mathcal{M}_{lj}\frac{\partial
\textrm{det}(h)}{\partial x^{l}}=\mathcal{M}^{mj}(-2\lambda
x_{j}+C_{j}).
\end{equation*}
Since $\mathcal{M}=\frac{1}{\textrm{det}(h)}h$,
$\mathcal{M}^{ml}=(\textrm{det}(h))h^{ml}$.  Substituting this
into the above equation gives
\begin{equation*}
\delta^{m}_{l}\frac{\partial
\textrm{log}(\textrm{det}(h))}{\partial x^{l}}=(-2\lambda
x_{j}+C_{j})h^{mj}.
\end{equation*}
We recall that $h^{mj}=\frac{\partial^{2}\Phi}{\partial
x^{m}\partial x^{j}}$ where $\Phi$ is the symplectic potential.
This implies that
\begin{equation*}
\frac{\partial \textrm{log}(\textrm{det}(h))}{\partial
x^{l}}=(-2\lambda x_{j}+C_{j})\frac{\partial^{2}\Phi}{\partial
x^{l}\partial x^{j}}
\end{equation*}
Integrating this equation, we have the following theorem.

\begin{theorem}
Let $(M^{2n},g,J)$ be a Kahler toric manifold with symplectic
potential $\Phi=\Phi(x_{1},...,x_{n})$ as above. The Einstein
condition is equivalent to the equation
\begin{equation}\label{KEToricSymp}
{\rm log}({\rm det}(h))=\sum_{j=1}^{n}(-2\lambda
x_{j}+C_{j})\frac{\partial \Phi}{\partial x^{j}}+2\lambda\Phi+E
\end{equation}
where $E$ and the $C_{j}$ are constants and $h^{-1}$ is the matrix
whose $ij$th entry is $\frac{\partial^{2}\Phi}{\partial
x^{i}\partial x^{j}}$.
\end{theorem}

    As in the four dimensional case, we can perform these integrations using the complex
coordinates $(u_{1},...,u_{n})$ as well. In these coordinates, the
Einstein condition can be written as
\begin{equation*}
\frac{\partial^{2}\textrm{log}(\textrm{det}(h))}{\partial
u^{i}\partial u^{j}}=-2\lambda h_{ij}
\end{equation*}
for all $1\leq i,j\leq n$.  Noting that
$\frac{\partial^{2}\eta}{\partial x^{i}\partial x^{j}}=h_{ij}$ and
integrating the above equation gives the following theorem.
\begin{theorem} Let $(M^{2n},g,J)$ be a Kahler toric manifold with complex potential
$\eta=\eta(u_{1},...,u_{n})$ as above.  The Einstein condition is equivalent to the equation
\begin{equation}\label{KEToricComp}
{\rm log}({\rm det}(h))=-2\lambda \eta+\sum_{j=1}^{n}C_{j}u_{j}+E
\end{equation}
where  $E$ and the $C_{j}$ are constants.
\end{theorem}

    This is a well-known expression for the Einstein condition in complex coordinates
on a Kahler toric manifold.  This reduction of the Kahler-Einstein
equations to the real Monge-Ampere equations in the Kahler toric
case was studied by Calabi in \cite{C}.

    Returning to equation (\ref{KEToricSymp}), we want to determine the significance
of the constants $C_{i}$.  As a first step, we fix a scale for the
metric by setting the Einstein constant, $\lambda$, equal to one.
A polytope is defined in the symplectic picture by a collections
of hyperplanes.  Of course, translating in the coordinates $x_{i}$
does not alter the metric.  Neither does transforming the
coordinates by an action of $SL(n,\mathbb{Z})$.  We can use these
two symmetries to fix the location of $n$ sides of the polytope.
Pick any vertex $p$ of the polytope.  By translation in the
$x_{i}$ we can set $p=(-1,...,-1)$.  Because the polytope must be
Delzant to define a Kahler toric manifold we know that there are
exactly $n$ faces meeting at $p$.  By an action of
$SL(n,\mathbb{Z})$, we can set $p$ to be the intersection of the
planes $1+x_{i}=0$ for all $i$.

    Abreu demonstrated that the potential of any Kahler metric whose polytope
contains a vertex of this type can be written as
\begin{equation*}
\Phi=\frac{1}{2}\sum_{i=1}^{n}(1+x_{i})\textrm{log}(1+x_{i})+\Psi
\end{equation*}
where $\Psi$ and its derivatives are smooth in a small open
neighborhood of the point $p=(-1,...,-1)$. Consider the a point
$x$ such that $-1<x_{i}<-1+\epsilon$ for $\epsilon>0$.  Leave
$x_{i}$ fixed for $i>1$.  As $x_{1}\longrightarrow-1$,
$\frac{\partial^{2}\Phi}{\partial x^{1}\partial
x^{1}}\longrightarrow\frac{1}{2}\frac{1}{1+x_{1}}$ and
$\frac{\partial^{2}\Phi}{\partial x^{i}\partial x^{j}}$ approaches
a constant if $i$ or $j$ does not equal one.  Furthermore, a
straightforward calculation shows that
$\frac{1}{\textrm{det}(h)}\frac{\partial \textrm{det}(h)}{\partial
x_{1}}\longrightarrow\frac{1}{1+x_{1}}$ as $x_{1}\rightarrow-1$.

Consider the following consequence of the Einstein condition when
$\lambda=1$
\begin{equation*}
\frac{\partial \textrm{log}(\textrm{det}(h))}{\partial x^{1}}=(-2
x_{j}+C_{j})\frac{\partial^{2}\Phi}{\partial x^{1}\partial x^{j}}.
\end{equation*}
As $x_{1}\longrightarrow-1$, the left-hand side of this equation
approaches $\frac{1}{1+x_{1}}$ and the right-hand side approaches
$(2+C_{1})\frac{1}{2}\frac{1}{1+x_{1}}$.  Equating these implies
that $C_{1}=0$.  Performing a similar calculation as
$x_{i}\longrightarrow-1$ for $i\neq 1$, we deduce immediately that
$C_{i}=0$ for all $i$.  The above proves the following
proposition.

\begin{prop}
Let $(M^{2n},g,J)$ be a Kahler-Einstein toric manifold with
Einstein constant $\lambda=1$. Through translation and an
$SL(n,\mathbb{Z})$ transformation, we can without loss of
generality assume that the polytope has faces $x_{i}=-1$.  Given
these assumptions, the constants $C_{i}$ of the Einstein condition
(\ref{KEToricSymp}) vanish and the Einstein condition becomes
\begin{equation}\label{KET}
-{\rm log}\left({\rm
det}(h^{-1})\right)=\sum_{j=1}^{n}-2x_{j}\frac{\partial
\Phi}{\partial x^{j}}+2\Phi+E.
\end{equation}
\end{prop}
Below, we will assume that our Kahler Einstein metrics satisfy the
conditions of the preceding proposition.

Before specializing to the compact case, we pause to mention work
done in the case in which the Kahler toric manifold is also
hyperkahler.  Because the Einstein constant is zero, the manifold
must be noncompact.  In this, the real Monge-Ampere equations
reduce further to a linear system.  For a more detailed account,
see \cite{Hitchin2}.

\subsection{Holomorphic Obstructions}

    Any Kahler-Einstein toric metric on a compact manifold $(M,g,J)$ must have
a positive Einstein constant.  That is, $M$ must be a Fano
manifold with $c_{1}(M)>0$.  There are two known holomorphic
obstructions to the existence of Einstein metrics on Fano
manifolds.  The first is due to Matsushima and involves the Lie
algebra of holomorphic vector fields which we denote by
$\mathfrak{h}(M)$.
\begin{theorem}\cite{Matsushima}
If a compact Fano manifold $N$ admits a Kahler-Einstein metric
with positive Einstein constant, then $\mathfrak{h}(N)$ is
reductive.
\end{theorem}

    The second obstruction is due to Futaki.  Let $N^{2n}$ be a Fano manifold.  One
can find a Kahler metric whose Kahler form $\omega$ is contained
in the first Chern class: $\omega\in c_{1}(N)$.  Let
$\rho_{\omega}$ be the Ricci form of the Kahler metric.  It is
well-known that $\rho_{\omega}\in c_{1}(N)$ also.  Therefore,
there exists a function $F_{\omega}$ such that
\begin{equation*}
\rho_{\omega}-\omega=\partial\bar{\partial}F_{\omega}.
\end{equation*}
Futaki defined the following map
$\mathcal{F}:\mathfrak{h}(N)\rightarrow\mathbb{C}$ by the equation
\begin{equation*}
\mathcal{F}(X)=\int_{N}XF_{\omega}\omega^{n}
\end{equation*}
for all $X\in\mathfrak{h}(N)$.  The map $\mathcal{F}$ is
independent of the choice of $\omega\in c_{1}(N)$ \cite{Futaki}.
The map $\mathcal{F}$ is known as the Futaki invariant.

The following theorem demonstrates that the Futaki invariant is an
obstruction to the existence of Kahler-Einstein metrics on Fano
manifolds.
\begin{theorem}\cite{Futaki}
If $N$ admits a Kahler-Einstein metric then $\mathcal{F}\equiv0$.
\end{theorem}
This condition can also be expressed in terms of the moment map
$\mu:N\rightarrow\mathfrak{h}^{\ast}$. Let $\mu_{X}$ be defined by
the relation $\langle\mu(p),X\rangle=\mu_{X}(p)$ for all $p\in N$
and $X\in\mathfrak{h}(N)$.
\begin{corollary}\cite{Futaki2}
The Futaki invariant satisfies
\begin{equation*}
\mathcal{F}(X)=-\int_{N}\mu_{X}\omega^{n}.
\end{equation*}
\end{corollary}
That is, the Futaki invariant vanishes when the barycenter of the
moment map lies at zero. The above corollary allows us to
calculate the Futaki invariant for a Kahler toric manifold
$(M^{2n},g,J)$. Because the Kahler structure is $T^{n}$ invariant,
we can, without loss of generality, choose $F_{\omega}$ to be a
function of the $u_{i}$ only. (Here is the $u_{i}$ are the real
part of the complex coordinates given above.)  Consider the (real)
holomorphic vector field $\frac{\partial}{\partial u^{i}}$.  As
discussed above, we know that $\mu_{\frac{\partial}{\partial
u^{i}}}=\frac{\partial\eta}{\partial u^{i}}=x_{i}$.  Furthermore,
we recall that the 'symplectic' coordinates $(x,\phi)$ are Darboux
coordinates so $\omega^{n}=dx_{1}\wedge...\wedge dx_{n}\wedge
d\phi_{1}\wedge...\wedge d\phi_{n}$. We calculate that
\begin{equation*}
\mathcal{F}(\frac{\partial}{\partial
u^{i}})=-\int_{M}x_{i}\omega^{n}=-(2\pi)^{n}\int_{\Delta}x_{i}dx
\end{equation*}
for all $i$, where $\Delta$ is the polytope defined in symplectic
coordinates $x_{i}$. We have the following proposition
\begin{prop}\cite{Mabuchi}
If $(M,g,J)$ is a compact Kahler-Einstein toric manifold then
\begin{equation}
\int_{\Delta}x_{i}dx=0
\end{equation}
for all $i$.
\end{prop}

\subsection{Kahler-Einstein Toric Manifolds in Dimensions Two and Four}

    In dimension two, the only compact Kahler toric manifold is $\mathbb{C}P^{1}$.  The
canonical Kahler metric is defined by the potential
\begin{equation*}
\Phi=\frac{1}{2}[(1+x_{1})\textrm{log}(1+x_{1})+(1-x_{1})\textrm{log}(1-x_{1})].
\end{equation*}
It is easily seen that this potential is Einstein.  Therefore, the
Kahler-Einstein condition on two-dimensional Kahler toric
manifolds (and indeed all two dimensional manifolds) is completely
understood.

    However, there remains one open problem concerning compact
Kahler-Einstein toric manifolds in dimension four.  It is
well-known \cite{YandT} that the only Kahler toric Fano manifolds
whose Lie algebra of holomorphic vector fields is reductive and
whose Futaki invariant vanishes are $\mathbb{C}P^{2}$,
$\mathbb{C}P^{1}\times\mathbb{C}P^{1}$, and
$\mathbb{C}P^{2}\sharp3\overline{\mathbb{C}P^{2}}$.  The last is
the blow-up of $\mathbb{C}P^{2}$ at three points in general
position.  The potentials of the Kahler-Einstein metrics on
$\mathbb{C}P^{2}$ and $\mathbb{C}P^{1}\times\mathbb{C}P^{1}$ are
both known explicitly.  They are
\begin{equation*}
\Phi=\frac{1}{2}[(1+x_{1})\textrm{log}(1+x_{1})+(1+x_{2})\textrm{log}(1+x_{2})+(1-x_{1}
-x_{2})\textrm{log}(1-x_{1}-x_{2})]
\end{equation*}
and
\begin{equation*}
\Phi=\frac{1}{2}[(1+x_{1})\textrm{log}(1+x_{1})+(1+x_{2})\textrm{log}(1+x_{2})+(1-x_{1})
\textrm{log}(1-x_{1})+(1-x_{2})\textrm{log}(1-x_{2})]
\end{equation*}
respectively.  The Kahler-Einstein potential are simply the
potentials of the canonical Kahler metric on the polytopes defined
by these potentials.

    While $\mathbb{C}P^{2}\sharp3\overline{\mathbb{C}P^{2}}$ is known to admit a
Kahler-Einstein toric metric (cf. \cite{Siu} and \cite{YandT}) the
potential is not known explicitly.  Much, however, is known about
this metric.  From the theorem of Abreu (\ref{Abreutheorem}) given
above, we know that the potential of the Kahler-Einstein metric
can be written as $\Phi=\Phi_{0}+f(x_{1},x_{2})$ where $\Phi_{0}$
is the potential of the canonical Kahler metric and $f$ is a
function whose second derivatives are smooth on the polytope.  For
ease of presentation, set $x_{1}=x$ and $x_{2}=y$.  We can write
$\Phi_{0}$ as
\begin{equation*}
\Phi_{0}=\frac{1}{2}[(1+x)\textrm{log}(1+x)+(1+y)\textrm{log}(1+y)+(1-x-y)\textrm{log}(1-x-y)
\end{equation*}
\begin{equation*}
+(1-x)\textrm{log}(1-x)+(1-y)\textrm{log}(1-y)+(1+x+y)\textrm{log}(1+x+y)].
\end{equation*}
This potential $\Phi_{0}$ is associated to polytope $P_{0}$
pictured in the figure below.

\pagebreak

\begin{center}
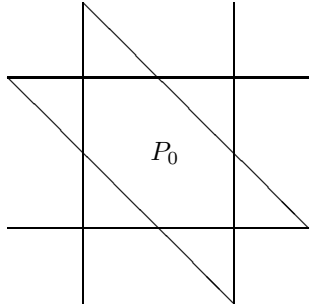
\textbf{Figure $1$}:~~\textit{Polytope associated to $\Phi_{0}$ on
$\mathbb{C}P^{2}\sharp3\overline{\mathbb{C}P^{2}}$}
\end{center}

\setlength{\unitlength}{1cm}
\begin{picture}(1,1)
\put(6,-4){\line(0,1){4}} \put(5,-3){\line(1,0){4}}
\put(8,-4){\line(0,1){4}} \put(6,0){\line(1,-1){3}}
\put(5,-1){\line(1,0){4}} \put(5,-1){\line(1,-1){3}}
\put(6.9,-2.1){$P_{0}$}
\end{picture}

~~

~

~

~

~

~

~

~

Substituting this expression for $\Phi$ into equation (\ref{KET}),
the Einstein condition becomes
\begin{equation}\label{newEinstein}
\textrm{log}(1-x^{2}+1-y^{2}+1-(x+y)^{2}+[1-x^{2}+1-(x+y)^{2}](1-y^{2})f_{yy}+[1-y^{2}+1-
(x+y)^{2}](1-x^{2})f_{xx}
\end{equation}
\begin{equation*}
-2(1-x^{2})(1-y^{2})f_{xy}+(1-x^{2})(1-y^{2})(1-(x+y)^{2})(f_{xx}f_{yy}-f_{xy}^{2}))=2xf_{x}
-2f+2yf_{y}-E.
\end{equation*}

    This equation, being equivalent to the Monge-Ampere equation,
is highly non-linear and cannot be integrated by any known
techniques.  However, one could try different types of functions
in an attempt to find the solution to (\ref{newEinstein}) that
corresponds to the Kahler-Einstein metric on
$\mathbb{C}P^{2}\sharp3\overline{\mathbb{C}P^{2}}$.

    Let $f$ be a function of only one coordinate say $x$.  We could ask whether there
are any solutions to equation (\ref{newEinstein}) such that
$f=f(x)$.  If $f$ does not depend on $y$ then (\ref{newEinstein})
becomes
\begin{equation*}
\textrm{log}(1-x^{2}+1-y^{2}+1-(x+y)^{2}+[1-y^{2}+1-(x+y)^{2}](1-x^{2})f_{xx})=2xf_{x}-2f-E.
\end{equation*}
The right-hand side of this equation depends only on $x$,
therefore the same must be true of the left-hand side.  The
left-hand side of the equations is a function of $x$ only if and
only if
\begin{equation*}
f_{xx}=\frac{1}{x^{2}-1}.
\end{equation*}
The function $f$ can then be written as
\begin{equation*}
f=\frac{1}{2}[(x-1)\textrm{log}(1-x)-(1+x)\textrm{log}(1+x)].
\end{equation*}
The resulting potential $\Phi$ becomes
\begin{equation*}
\Phi=\frac{1}{2}[(1+y)\textrm{log}(1+y)+(1-x-y)\textrm{log}(1-x-y)+(1-y)\textrm{log}(1-y)+
(1+x+y)\textrm{log}(1+x+y)]
\end{equation*}
which is the potential for the Kahler-Einstein metric on
$S^{2}\times S^{2}$.  Therefore, the Siu metric does not have
potential the form $\Phi=\Phi_{0}+f(x).$

    Alternatively, one could look for solutions to equation (\ref{newEinstein}) in which
$f=p(x)+q(y)$.  Substituting this expression for $f$ into equation
(\ref{newEinstein}) gives
\begin{equation*}
\textrm{log}(1-x^{2}+1-y^{2}+1-(x+y)^{2}+[1-x^{2}+1-(x+y)^{2}](1-y^{2})q_{yy}+[1-y^{2}+1
-(x+y)^{2}](1-x^{2})p_{xx}
\end{equation*}
\begin{equation*}
+(1-x^{2})(1-y^{2})(1-(x+y)^{2})p_{xx}q_{yy}=2xp_{x}-2p-2q+2yq_{y}-E.
\end{equation*}
Evaluating this equation at $y=1$ and $y=-1$ gives
\begin{equation*}
\textrm{log}(1-x^{2}+1-(x+1)^{2}+[1-(x+1)^{2}](1-x^{2})p_{xx})=2xp_{x}-2p+F
\end{equation*}
and
\begin{equation*}
\textrm{log}(1-x^{2}+1-(x-1)^{2}+[1-(x-1)^{2}](1-x^{2})p_{xx})=2xp_{x}-2p+G
\end{equation*}
where $F$ and $G$ are constants.  Equating these two expressions
implies that
\begin{equation*}
p=\frac{(2H-2)x^{2}+(2H-1)x+(1-H)}{x(x^{2}-1)((H-1)x+(2H+2))}
\end{equation*}
for some constant $H$. However, to define a smooth function on
$\mathbb{C}P^{2}\sharp3\overline{\mathbb{C}P^{2}}$ all of the
second derivatives must be smooth on the polytope.  This cannot be
the case for any value of $H$. Therefore, the Siu metric does not
have potential the form $\Phi=\Phi_{0}+p(x)+q(y)$.

    In a sense, it is not surprising that the function $f$ is not a function of
one variable or a sum of functions in $x$ and $y$.  Potentials of
this form do not respect the 'hexagonal' symmetry of the polytope
defining the symplectic structure on
$\mathbb{C}P^{2}\sharp3\overline{\mathbb{C}P^{2}}$.  One could
suppose the potential $\Phi$ corresponding to the Siu metric would
respect that six-fold symmetry.  One could propose that the
function $f$ which will give rise to the Siu metric will be of the
form
\begin{equation*}
f=\frac{1}{2}\sum_{j=1}^{m}[(A_{i}+B_{i}x)\textrm{log}(A_{i}+B_{i}x)+(A_{i}+B_{i}y)
\textrm{log}(A_{i}+B_{i}y)+(A_{i}-B_{i}(x+y))\textrm{log}(A_{i}-B_{i}(x+y))
\end{equation*}
\begin{equation*}
+(A_{i}-B_{i}x)\textrm{log}(A_{i}-B_{i}x)+(A_{i}-B_{i}y)\textrm{log}(A_{i}-B_{i}y)+(A_{i}
+B_{i}(x+y))\textrm{log}(A_{i}+B_{i}(x+y))]
\end{equation*}
where $A_{i}$ and $B_{i}$ are constants on the function $f$ is
smooth on the polytope.  Unfortunately, testing such an ansatz
directly proves impossible using MAPLE as the calculations are too
lengthy.  We were able to show that the Siu metric potential does
not have this form for $m=1$.

    Below, we write down explicitly a six-dimensional Kahler-Einstein metric in
Kahler toric coordinates.  We note that the potential of this
metric can be written in the form given in the above ansatz.

\subsection{A Kahler-Einstein Toric Metric in Dimension Six}

    In \cite{Sakane}, Sakane proved the existence of a Kahler-Einstein
metric on
$\mathbb{P}(\mathcal{O}_{\mathbb{C}P^{1}\times\mathbb{C}P^{1}}\oplus\mathcal{O}_{\mathbb{C}
P^{1}\times\mathbb{C}P^{1}}(1,-1))$ though he did not write the
metric down explicitly.  This is a six-dimensional Kahler toric
manifold. In fact, it is a cohomogeneity-one manifold.  Sakane's
proof of the existence of this metric makes use of the full
symmetry group. (His construction will be discussed below in the
section on fiberwise Kahler toric manifolds.) Here, we will
express the Kahler-Einstein metric on this manifold in the
'symplectic' coordinates of Kahler toric geometry.
    For ease of expression, we set $x_{1}=x$, $x_{2}=y$, and $x_{3}=z$.  Following
the theorem of Abreu, the symplectic potential, $\Phi$, of the
Kahler-Einstein metric can be written as
\begin{equation*}
\Phi=\Phi_{0}+f
\end{equation*}
where
\begin{equation*}
\Phi_{0}=\frac{1}{2}[(1+x)\textrm{log}(1+x)+(1-x)\textrm{log}(1-x)+(1+y)\textrm{log}(1+y)
+(1+z)\textrm{log}(1+z)
\end{equation*}
\begin{equation*}
+(1-x-y)\textrm{log}(1-x-y)+(1+x-z)\textrm{log}(1+x-z)]
\end{equation*}
and $f$ and its derivatives are smooth on the polytope whose
canonical potential is $\Phi_{0}$. Because the metric is of
cohomogeneity-one, $f=f(x)$.  Substituting this expression into
equation (\ref{KET}), the Einstein condition becomes
\begin{equation*}
\textrm{log}(f_{xx}(4-x^{2})(1-x^{2})+3(2-x^{2}))-\textrm{log}(4)=2xf_{x}-2f-E.
\end{equation*}
Differentiating both sides with respect to $x$ we have
\begin{equation*}
(x^{2}-4)(x^2-1)f_{xxx}-2x(x^{2}-4)(x^{2}-1)f_{xx}^{2}+2x(5x^{2}-11)f_{xx}-6x=0.
\end{equation*}
Through trial and error inputing rational functions into MAPLE for
$f_{xx}$, we find that
\begin{equation*}
f_{xx}=\frac{(x^{2}-10)}{(x^{2}-4)(x^{2}-7)}.
\end{equation*}
Integrating this twice, we can write $f$ in the form
\begin{equation*}
f=\frac{1}{2}[(-2+x)\textrm{log}(-2+x)+(-2-x)\textrm{log}(-2-x)+(1+\frac{\sqrt{7}}{7}x)
\textrm{log}(1+\frac{\sqrt{7}}{7}x)+(1-\frac{\sqrt{7}}{7}x)\textrm{log}(1-\frac{\sqrt{7}}{7}x)].
\end{equation*}
Note that $\Phi$ is equivalent to the canonical potential obtained
by adding four lines outside of the polytope $\Delta_{0}$.
\section{Extremal Kahler Metrics on Kahler Toric Manifolds}

    As proved above, a Kahler toric metric is extremal if and only if
\begin{equation}\label{KTExtKah}
-\frac{\partial^{2}h_{ij}}{\partial x^{i}\partial
x^{j}}=-\frac{1}{\textrm{det}(h)}h_{ij}\frac{\partial^{2}\textrm{det}(h)}{\partial
x^{i}\partial x^{j}}=\sum_{k=1}^{n}\alpha_{k}x_{k}+\beta
\end{equation}
where $h^{ij}=\frac{\partial^{2}\Phi}{\partial x^{i}\partial
x^{j}}$.  Clearly, any metric which is of constant scalar
curvature is extremal.  In this section we want to study the
existence of extremal Kahler metrics on Kahler toric manifolds.
After reviewing the construction of a one-parameter family of
extremal Kahler metrics on
$\mathbb{C}P^{2}\sharp\overline{\mathbb{C}P^{2}}$ we will
construct a two-parameter (up to scale) family of extremal Kahler
metrics on a six-dimensional Kahler toric manifold:
$\mathbb{P}(\mathcal{O}_{\mathbb{C}P^{1}\times\mathbb{C}P^{1}}\oplus\mathcal{O}_{\mathbb{C}
P^{1}\times\mathbb{C}P^{1}}(1,-1))$.
We shall see that within this two-parameter family of extremal
Kahler metrics there exists a one-parameter subfamily of constant
scalar curvature metrics, one of which is the Einstein metric
described above.

\subsection{Extremal Kahler Metrics on $\mathbb{C}P^{2}\sharp\overline{\mathbb{C}P^{2}}$}

    Extremal Kahler metrics on the blow-up of $\mathbb{C}P^{2}$ were originally found
by Calabi in \cite{Calabi} but were not written down explicitly.
These metrics were written in Kahler toric coordinates by Abreu in
\cite{Abreu2}.  Abreu wrote the metrics explicitly by converting
Calabi's non-explicit solutions into Kahler toric coordinates.
Our presentation will provide the same results as Abreu; however,
we will arrive at the solutions without reference to Calabi's
method.  Rather, we will obtain the solutions by solving equation
(\ref{KTExtKah}) directly.  This will be our method for obtaining
new extremal Kahler metrics below.

    Up to scale, we can write any canonical Kahler toric potential on $\mathbb{C}P^{2}\sharp
    \overline{\mathbb{C}P^{2}}$ as
\begin{equation*}
\Phi_{0}=\frac{1}{2}[(1+x)\textrm{log}(1+x)+(1+y)\textrm{log}(1+y)+(a-x-y)\textrm{log}
(a-x-y)+(1-x)\textrm{log}(1-x)]
\end{equation*}
where $a>0$ and we have set $x_{1}=x$ and $x_{2}=y$ \cite{Abreu}.
Varying $a$ is equivalent to changing the Kahler class of the
metric.  Any Kahler toric metric on
$\mathbb{C}P^{2}\sharp\overline{\mathbb{C}P^{2}}$ can be written
as $\Phi=\Phi_{0}+f$ where $f$ satisfies theorem
(\ref{Abreutheorem}).  This manifold has more than the
$T^{2}$-symmetry evident in these coordinates.  There is in fact a
$U(2)$-symmetry which implies that the manifold is actually of
cohomogeneity-one.  This translates into the conditions $f=f(x)$
and $S=\alpha x+\beta$ in the extremal Kahler case \cite{Abreu}.

    A straightforward calculation shows that
\begin{equation*}
h^{-1}=\begin{pmatrix}
\Phi_{xx} & \Phi_{xy} \\
\Phi_{xy} & \Phi_{yy} \\
\end{pmatrix}=\begin{pmatrix}
\frac{1}{1-x^{2}}+\frac{1}{2}\frac{1}{a-x-y}+f_{xx} & \frac{1}{2}\frac{1}{a-x-y} \\
\frac{1}{2}\frac{1}{a-x-y} & \frac{1}{2}\frac{1}{1+y}+\frac{1}{2}\frac{1}{a-x-y} \\
\end{pmatrix}.
\end{equation*}
Setting $P(x)=f_{xx}(x)+\frac{1}{1-x^{2}}$, and calculating the
inverse of the above matrix gives
\begin{equation*}
h=\begin{pmatrix}
\frac{2(a+1-x)}{2P(a+1-x)+1} & \frac{-2(1+y)}{2P(a+1-x)+1} \\
\frac{-2(1+y)}{2P(a+1-x)+1} & \frac{2(1+y)(2P(a-x-y)+1)}{2P(a+1-x)+1} \\
\end{pmatrix}.
\end{equation*}
To simplify this further, we set $Q=P+\frac{1}{2}\frac{1}{a+1-x}$.
After making this substitution, equation (\ref{KTExtKah}) becomes
\begin{equation}\label{one-pointblow-up}
-\frac{\partial^{2}(\frac{1}{Q})}{\partial
x^{2}}+\frac{2}{a+1-x}\frac{\partial(\frac{1}{Q})}{\partial
x}+\frac{4}{a+1-x}-\alpha x-\beta=0.
\end{equation}
We are able to solve this equation by setting
\begin{equation*}
\frac{1}{Q}=\frac{1}{(a+1)-x}\int\left[\int\left(\frac{4}{a+1-x}-\alpha
x-\beta\right)((a+1)-x)dx\right]dx.
\end{equation*}
A brief calculation shows that this indeed gives a solution for
(\ref{one-pointblow-up}).  Solving for $f_{xx}$ we find that
\begin{equation*}
f_{xx}=\frac{1}{2}\frac{1}{x-a-1}+\frac{a}{ax^{2}-(3a^{2}+6a+2)x+3a^{3}+9a^{2}+7a+2}.
\end{equation*}
This metric has scalar curvature
\begin{equation*}
S=\frac{6}{3a^{2}+6a+2}(-2ax+a^{2}+4a+2).
\end{equation*}
For $a>0$ we see that the metric is never of constant scalar
curvature.  It is intriguing to note that when $a=0$ we retrieve
the Fubini-Study metric on $\mathbb{C}P^{2}$.  In a sense, the
extremal Kahler metrics on
$\mathbb{C}P^{2}\sharp\overline{\mathbb{C}P^{2}}$ and the
Fubini-Study metric on $\mathbb{C}P^{2}$ are part of a single one
parameter family of metrics.

    Before proceeding to our construction of a new family of extremal Kahler metrics,
we want to determine the form of the potential $f$ which involves
performing two more integrations. When writing down the canonical
Kahler metric for a given polytope the potentials $\Phi_{0}$ can
be written schematically as
\begin{equation}\label{potentialform}
\Phi_{0}=\sum_{k=1}^{l}[(a_{k}+\sum_{i=1}^{d}b_{ik}x_{i})\textrm{log}(a_{k}+\sum_{i=1}^{d}
b_{ik}x_{i})]
\end{equation}
where $l$ is the number of $(d-1)$-dimensional faces of the
polytope and $2d$ is the real dimension of the Kahler toric
manifold.  While it is possible to construct Kahler toric metrics
whose symplectic potentials cannot be written in this form one
could ask whether all extremal Kahler metrics might have
potentials of this form.  Let us test this conjecture on the one
parameter family of extremal Kahler metrics we have just
constructed.

    By inspection, we see that the potential can be written in the form
of (\ref{potentialform}) if and only if the polynomial
$ax^{2}-(3a^{2}+6a+2)x+3a^{3}+9a^{2}+7a+2$ has two real roots.
This quadratic equation has two real roots when the discriminant
is positive.  That is $-3a^{4}+20a^{2}+16a+4>0$. This is true for
small values of $a>0$.  However, for $a$ large enough the
discriminant of this equation is negative.  Therefore, not all
extremal Kahler metrics on Kahler toric manifolds have symplectic
potentials of the form (\ref{potentialform}).

    While not every extremal Kahler metric on $\mathbb{C}P^{2}\sharp\overline{\mathbb{C}P^{2}}$
has symplectic potential of the form (\ref{potentialform}), there
is one privileged metric which does.  When $a=1$, the resulting
extremal Kahler metric, $g$, is conformal to an Hermitian-Einstein
metric, $h$, found by Page. That is,
\begin{equation*}
h=S^{-2}g
\end{equation*}
where $S$ is the scalar curvature of the extremal Kahler metric
$g$.

    When $a=1$, we have
\begin{equation*}
f_{xx}=\frac{1}{2}\frac{1}{x-2}+\frac{1}{x^{2}-11x+21}.
\end{equation*}
This form of the metric was first found by Abreu in \cite{Abreu}.
However, he did not write down the function $f$ which determines
the alteration in the canonical potential necessary to obtain this
metric.  Performing the simple double integral (which of course
does not determine $f$ uniquely) we can write $f$ as
\begin{equation*}
f=\frac{1}{2}(-2+x)\textrm{log}(-2+x)+\frac{1}{2}\left(-1-\frac{11}{\sqrt{37}}+\frac{2}{
\sqrt{37}}x\right)\textrm{log}\left(-1-\frac{11}{\sqrt{37}}+\frac{2}{\sqrt{37}}x\right)
\end{equation*}
\begin{equation*}
+\frac{1}{2}\left(-1+\frac{11}{\sqrt{37}}-\frac{2}{\sqrt{37}}x\right)\textrm{log}\left(-1+
\frac{11}{\sqrt{37}}-\frac{2}{\sqrt{37}}x\right).
\end{equation*}

    We make a point of writing the symplectic potential down in this from
because it shows that this privileged metric has potential which
can be written in the form (\ref{potentialform}).  We believe this
to be interesting for the following reasoning.  If one wants to
describe a particular symplectic toric four manifold, one can
exhibit a collection of lines in the plane which define a
particular polytope.  This manifold, if Kahler toric, can be
endowed with a canonical Kahler toric metric by adding together
functions each of which is determined by one of the lines defining
the polytope.  Abreu, in theorem (\ref{Abreutheorem}), then gives
a method of modifying this potential to obtain all other possible
Kahler toric metrics with the symplectic structure defined by the
polytope. This is precisely the method we used to determine the
extremal Kahler metric on
$\mathbb{C}P^{2}\sharp\overline{\mathbb{C}P^{2}}$.

    We could change the Kahler structure by an inverted method.  We could
begin with a collection of lines in the plane that define a
polytope.  We could first add lines which do not intersect the
polytope at any point and then write down the canonical Kahler
toric metric determined by the lines defining the polytope and by
the 'extra' lines lying outside of the polytope. For example, to
define a symplectic structure on the blow-up of $\mathbb{C}P^{2}$
at one-point we draw the following four lines in the $xy$-plane:
\begin{equation*}
l_{1}=1+x=0,~~~l_{2}=1-x=0,~~~l_{3}=1+y=0,~~~l_{4}=1-x-y=0.
\end{equation*}
The canonical Kahler metric is then defined by the potential
$\Phi_{0}=\frac{1}{2}\sum_{j=1}^{4}l_{j}\textrm{log}l_{j}$.  To
obtain the extremal Kahler metric conformal to the Page metric one
could, instead of directly modifying the potential, first add
three lines to those defining the polytope.  That is, set
\begin{equation*}
{l}_{5}=-2+x=0,~~~{l}_{6}=-1-\frac{11}{\sqrt{37}}+\frac{2}{\sqrt{37}}x=0,~~~{l}_{7}=-1+
\frac{11}{\sqrt{37}}-\frac{2}{\sqrt{37}}x=0.
\end{equation*}

    The extremal Kahler metric conformal to the Page metric can then
be viewed as the 'canonical' Kahler metric defined by the
\textit{seven} lines $l_{1}$ to $l_{7}$:
$\Phi=\frac{1}{2}\sum_{j=1}^{7}l_{j}\textrm{log}l_{j}$.   A
diagram of these lines appears in figure $2$ below.  Note that $P$
corresponds to the original polytope on which the metric is
defined.

~

\begin{center}
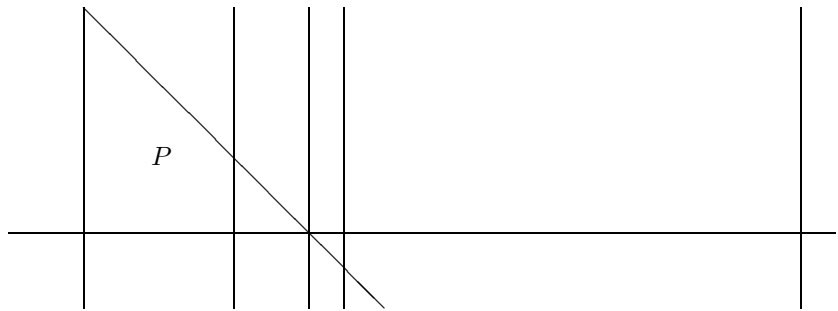
\textbf{Figure $2$}:~~\textit{Lines defining the extremal Kahler
metric conformal to the Page metric on
$\mathbb{C}P^{2}\sharp\overline{\mathbb{C}P^{2}}$}
\end{center}

\setlength{\unitlength}{1cm}
\begin{picture}(1,1)
\put(3,-4){\line(0,1){4}} \put(2,-3){\line(1,0){11}}
\put(5,-4){\line(0,1){4}} \put(3,0){\line(1,-1){4}}
\put(6,-4){\line(0,1){4}} \put(12.541,-4){\line(0,1){4}}
\put(6.459,-4){\line(0,1){4}}
 \put(3.9,-2.1){$P$}
\end{picture}

~~

~

~

~

~

~

~

~

    As demonstrated above, not all of the extremal Kahler manifolds
can be constructed in this way.  It would be interesting to know
whether there is any significance to the metric which can be
described by a collection of lines (or hypersurfaces in higher
dimensions.)  One could ask whether all constant scalar curvature
metrics or all Kahler-Einstein metrics on Kahler toric manifolds
have potentials which can be written in this form.  Finally, one
could ask whether there is a deeper combinatorial significance to
the location of the lines added to give this special metric.
Unfortunately, we are unable at this time to provide an answer to
any of these questions.  Below, when discussing Hermitian-Einstein
metrics in dimension four, we will revisit the Page metric.

\subsection{New Extremal Kahler Metrics on $\mathbb{P}(\mathcal{O}_
{\mathbb{C}P^{1}\times\mathbb{C}P^{1}}\oplus\mathcal{O}_{\mathbb{C}P^{1}\times\mathbb{C}P^{1}}
(1,-1))$}

    The complex manifold $\mathbb{P}(\mathcal{O}_{\mathbb{C}P^{1}\times\mathbb{C}P^{1}}
    \oplus\mathcal{O}_
{\mathbb{C}P^{1}\times\mathbb{C}P^{1}}(1,-1))$ is a manifold of
six real dimensions which is invariant under the action of
$T^{3}$.  In fact, like the blow-up of $\mathbb{C}P^{2}$ at one
point, this manifold is of cohomogeneity-one.  In the previous
section, we have written down the Kahler-Einstein metric first
discovered by Sakane on this manifold in Kahler toric coordinates.
In this section, we will find a two-parameter family of extremal
Kahler metrics.  Within this two-parameter family lies two
one-parameter families of constant scalar curvature metrics and,
of course, the Kahler-Einstein metric we identified earlier.

    For ease of presentation, let $x=x_{1}, y=x_{2},$ and $z=x_{3}$.  The potential of
any extremal Kahler metric on this manifold can be written as
$\Phi=\Phi_{0}+f$ where
\begin{equation*}
\Phi_{0}=\frac{1}{2}[(1+x)\textrm{log}(1+x)+(1-x)\textrm{log}(1-x)+(1+y)\textrm{log}(1+y)
+(1+z)\textrm{log}(1+z)
\end{equation*}
\begin{equation*}
+(a-x-y)\textrm{log}(a-x-y)+(c+x-z)\textrm{log}(c+x-z)]
\end{equation*}
where $a>0$ and $c>0$ and the second derivatives of $f$ are smooth
on the polytope.  A diagram of the polytope is given in figure
$3$.  Note that the labels $a$ and $c$ give the lengths of the
corresponding sides.

\begin{center}
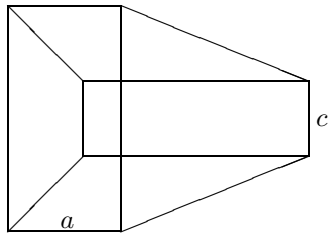
\textbf{Figure $3$}:~~\textit{Polytope associated to $\Phi_{0}$ on
$\mathbb{P}(\mathcal{O}_{\mathbb{C}P^{1}\times\mathbb{C}P^{1}}\oplus\mathcal{O}_{\mathbb{C}
P^{1}\times\mathbb{C}P^{1}}(1,-1))$}
\end{center}

~

\setlength{\unitlength}{1cm}
\begin{picture}(1,1)
\put(6,-2){\line(1,0){3}} \put(6,-2){\line(0,1){1}}
\put(9,-2){\line(0,1){1}} \put(6,-1){\line(1,0){3}}
\put(6,-2){\line(-1,-1){1}} \put(5,-3){\line(1,0){1.5}}
\put(6.5,-3){\line(5,2){2.5}} \put(6.5,-3){\line(0,1){3}}
\put(5,-3){\line(0,1){3}} \put(6.5,0){\line(5,-2){2.5}}
\put(6,-1){\line(-1,1){1}} \put(5,-0){\line(1,0){1.5}}
\put(5.7,-2.95){$a$}\put(9.1,-1.6){$c$}
\end{picture}

~~

~

~

~

~

~

~

~
 Because this manifold is of cohomogeneity-one, we can search for
extremal Kahler metrics which are of cohomogeneity-one.  That is,
we can set $f=f(x)$. The presence of these extra symmetries tells
us that, if the metric is extremal Kahler, then $S=\alpha
x+\beta$, i.e. it is dependant on $x$ only.

    Setting $P(x)=f_{xx}+\frac{1}{1-x^{2}}$, a straightforward calculation shows that
\begin{equation*}
h^{-1}=\begin{pmatrix}
\Phi_{xx} & \Phi_{xy} & \Phi_{xz}\\
\Phi_{xy} & \Phi_{yy} & \Phi_{yz} \\
\Phi_{xz} & \Phi_{yz} & \Phi_{zz}\\
\end{pmatrix}=\begin{pmatrix}
P+\frac{1}{2}\frac{1}{a-x-y}+\frac{1}{2}\frac{1}{c+x-z} & \frac{1}{2}\frac{1}{a-x-y} &
\frac{1}{2}\frac{1}{c+x-z} \\
\frac{1}{2}\frac{1}{a-x-y} & \frac{1}{2}\frac{1}{1+y}+\frac{1}{2}\frac{1}{a-x-y} & 0 \\
\frac{1}{2}\frac{1}{c+x-z} & 0 & \frac{1}{2}\frac{1}{1+z}+\frac{1}{2}\frac{1}{c+x-z} \\
\end{pmatrix}.
\end{equation*}
As we did when constructing the extremal Kahler metrics on
$\mathbb{C}P^{2}\sharp\overline{\mathbb{C}P^{2}}$, we make one
final change of variable by setting
$Q=P+\frac{1}{2}\frac{1}{x+c+1}+\frac{1}{2}\frac{1}{a+1-x}$. Given
this substitution a straightforward but cumbersome calculation
demonstrates that equation (\ref{KTExtKah}) becomes
\begin{equation}\label{sixdimextkah}
-\frac{\partial^{2}}{\partial
x^{2}}\left(\frac{1}{Q}\right)+\left(\frac{2}{a+1-x}+\frac{-2}{c+1+x}\right)\frac{\partial}
{\partial
x}\left(\frac{1}{Q}\right)+\frac{2}{(x+c+1)(a+1-x)}\left(\frac{1}{Q}\right)
\end{equation}
\begin{equation*}
+\frac{4}{a+1-x}+\frac{4}{x+c+1}-\alpha x-\beta=0.
\end{equation*}

    We are able to solve this equation by, for the moment, inspection.  Equation
(\ref{sixdimextkah}) has solution
\begin{equation*}
\frac{1}{Q}=\frac{1}{(a+1-x)(x+c+1)}\int\left[\int\left(\frac{4}{a+1-x}+\frac{4}{c+1+x}-\alpha
x-\beta\right)(a+1-x)(x+c+1)dx\right]dx
\end{equation*}
Performing this double integration and solving for $f_{xx}$ we
find that
\begin{equation*}
f_{xx}=\frac{1}{2}\frac{1}{x-a-1}+\frac{1}{2}\frac{-1}{x+c+1}+
\end{equation*}
\begin{equation*}
((-6c-6a^2c+6ac^2+6a)x+6c+10a^3c+70ac+40a^2c+40ac^2+10c^3a+5a^2c^2+6a)\{6(a-c)(1-ac)x^3
\end{equation*}
\begin{equation*}
+(10a^3c-26c-42ac-8-20ac^2-20c^2-26a-20a^2-20a^2c+10c^3a-25a^2c^2)x^2
\end{equation*}
\begin{equation*}
+(98a^2c+32a^2-20c^3-32c^2+20a^3-2c-30c^3a^2-60c^3a-98ac^2+30a^3c^2+2a+60a^3c)x
\end{equation*}
\begin{equation*}
+8+30a^3c^3+34a+34c+114ac+52c^2+52a^2+184a^2c+184ac^2+257a^2c^2
\end{equation*}
\begin{equation*}
+20c^3+20a^3+90a^3c^2+70c^3a+90c^3a^2+70a^3c\}^{-1}
\end{equation*}
While the calculations are somewhat unwieldy, it is not too
difficult to check that this defines a positive definite metric
with the appropriate smoothness conditions for all values of $a>0$
and $c>0$.
\begin{prop}
There exists, up to scale, a two-parameter family of extremal
Kahler metrics on \linebreak
$\mathbb{P}(\mathcal{O}_{\mathbb{C}P^{1}\times\mathbb{C}P^{1}}\oplus\mathcal{O}_{\mathbb{C}
P^{1}\times\mathbb{C}P^{1}}(1,-1))$:
one for every value of $a>0$ and $c>0$.
\end{prop}
    We further calculate that
\begin{equation*}
\alpha=\frac{60(a-c)(ac-1)}{(30a^2c+15a^2c^2+10a^2+16a+56ac+30ac^2+4+16c+10c^2)}
\end{equation*}
and
\begin{equation*}
\beta=\frac{6(10a^2+26a+20ac^2+10c^2+20a^2c+62ac+8+26c+5a^2c^2)}{(30a^2c+15a^2c^2+10a^2
+16a+56ac+30ac^2+4+16c+10c^2)}.
\end{equation*}
When $\alpha=0$, the metric is of constant scalar curvature.  This
occurs when $a=c$ or $a=\frac{1}{c}$. This proves the following
proposition.
\begin{prop}
There exists, up to scale, two one-parameter families of constant
scalar curvature Kahler metrics on
$\mathbb{P}(\mathcal{O}_{\mathbb{C}P^{1}\times\mathbb{C}P^{1}}\oplus\mathcal{O}_{
\mathbb{C}P^{1}\times\mathbb{C}P^{1}}(1,-1))$:
one for every value of $a>0$ and $c>0$ such that $a=c$ or
$a=\frac{1}{c}$.
\end{prop}
    These two families of constant scalar curvature Kahler metrics intersect at the point
$(a,c)=(1,1)$.  This point of intersection corresponds to the
Kahler-Einstein metric found by Sakane and presented in the
previous section.  Figure $4$ below illustrate how the two
families of constant scalar curvature Kahler metrics sits inside
the two-parameter family of extremal Kahler metrics.

~

\begin{center}
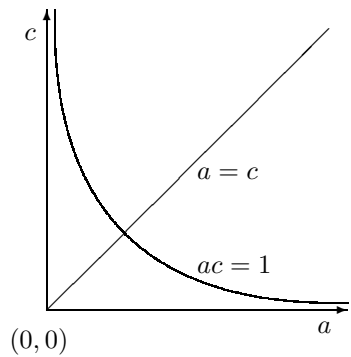
\textbf{Figure $4$}:~~\textit{The space of extremal Kahler and
constant scalar curvature Kahler metrics on
$\mathbb{P}(\mathcal{O}_{\mathbb{C}P^{1}\times\mathbb{C}P^{1}}\oplus\mathcal{O}_{\mathbb{C}P^{1}
\times\mathbb{C}P^{1}}(1,-1))$}
\end{center}

\setlength{\unitlength}{1cm}
\begin{picture}(1,1)
\put(5,-4){\vector(1,0){4}} \put(5,-4){\vector(0,1){4}}
\put(5,-4){\line(1,1){3.75}} \put(8.6,-4.3){$a$}
\put(4.7,-.4){$c$} \qbezier(5.1,-0)(5,-4)(9,-3.9)
\put(7,-3.5){$ac=1$} \put(7,-2.25){$a=c$} \put(4.5,-4.5){$(0,0)$}
\end{picture}

~

~

~

~

~

~

~

~

~

~
    We now ask what happens if we allow $a$ and $c$ to take zero values.  We see
immediately that if $a=c=0$ then we recover the Fubini-Study
metric on $\mathbb{C}P^{3}$.  If $a=0$ and $c>0$ (or equivalently
$c=0$ and $a>0$) , then we recover a family of extremal Kahler
metrics (of non-constant scalar curvature) on another
six-dimensional Kahler toric manifold which is a
$\mathbb{C}P^{2}$-bundle over $\mathbb{C}P^{1}$.

\subsection{Hermitian-Einstein Manifolds in Dimension Four}

    While the question of the \textit{existence} of Kahler-Einstein metrics with
positive Einstein constant is completely understood in dimension
four, less is known about the existence of Hermitian-Einstein
metrics. LeBrun has used a theorem of Derdzinski to determine a
list of compact manifolds which may admit positive
Hermitian-Einstein metrics.

\begin{theorem}\cite{LeBrun}
If $M$ is a four-dimensional compact manifold endowed with a
non-Kahler Hermitian-Einstein metric $g$, then $M$ is a blow-up of
$\mathbb{C}P^{2}$ at one, two, or three points in general
position.  Any such metric would have a $T^{2}$-symmetry.
\end{theorem}

    Derdzinski showed that in dimension four, extremal Kahler metrics satisfying
an additional second-order equation are conformally Einstein.
This result is summarized in the following theorem.

\begin{theorem}\cite{Derd}
$(M,g,J)$ is a four dimensional extremal Kahler manifold with
nonconstant scalar curvature $S$ such that
\begin{equation}\label{extkahler}
S^{3}+6S\Delta S-12\langle\nabla S,\nabla S\rangle=\beta,
\end{equation}
where $\beta$ is a constant, if and only if the metric $S^{-2}g$
is Hermitian-Einstein where it is defined.
\end{theorem}

Note that the resulting Einstein metric is not Kahler but is
Hermitian.

In this section, we want to apply this theorem to the Kahler toric
case in dimension four.  That is, we want to express equation
(\ref{extkahler}) in symplectic coordinates.  Recall that the
scalar curvature of an extremal Kahler toric manifold can be
written as $S=\sum_{i=1}^{n}\alpha_{i}x_{i}+\beta$ where $\beta$
and the $\alpha_{i}$ are constants.  We now have the following
corollary to the previous theorem.

\begin{corollary}
Let $(M,g,J)$ be a four-dimensional extremal Kahler toric manifold
with scalar curvature $S=\sum_{i=1}^{2}\alpha_{i}x_{i}+\beta$. The
metric $S^{-2}g$ is locally Einstein if and only if
\begin{equation}
\left(\sum_{i=1}^{2}\alpha_{i}x_{i}+\beta\right)^{3}+6\left(\sum_{i=1}^{2}\alpha_{i}x_{i}+
\beta\right)\left(\frac{1}{{\rm det}(h)}h_{ij}\frac{\partial {\rm
det}(h)}{\partial
x^{j}}\alpha_{i}\right)-12h_{ij}\alpha_{i}\alpha_{j}=\gamma
\end{equation}
where $\gamma$ is a constant.
\end{corollary}
\begin{proof}
To prove this theorem, we need only evaluate equation
(\ref{extkahler}) on an extremal Kahler toric manifold. To do this
let us look at the terms $\Delta S$ and $\langle\nabla S,\nabla
S\rangle$.  From the discussion above, we see that
\begin{equation*}
\Delta S= h_{ij}\frac{\partial^{2}S}{\partial x^{i}\partial
x^{j}}+\frac{1}{\textrm{det}(h)}h_{ij}\frac{\partial
\textrm{det}(h)}{\partial x^{i}}\frac{\partial S}{\partial
x^{j}}=\frac{1}{\textrm{det}(h)}h_{ij}\frac{\partial
\textrm{det}(h)}{\partial x^{i}}\alpha_{j}
\end{equation*}
in the extremal Kahler case.  Furthermore,
\begin{equation*}
\langle\nabla S,\nabla S\rangle=h_{ij}\frac{\partial S}{\partial
x^{i}}\frac{\partial S}{\partial
x^{j}}=h_{ij}\alpha_{i}\alpha_{j}.
\end{equation*}
The corollary follows immediately.
\end{proof}

    This theorem and its corollary provide a strategy for finding the Hermitian-Einstein
metrics (if they exist) on
$\mathbb{C}P^{2}\sharp2\overline{\mathbb{C}P^{2}}$ and
$\mathbb{C}P^{2}\sharp3\overline{\mathbb{C}P^{2}}$.  First we
could look for all (if they exist) extremal Kahler metrics on
these manifolds, then check which (if any) satisfy the Derdzinski
condition.

    Alternatively, one could take a slightly different but ultimately equivalent
approach.  Namely, one could try to solve the four-dimensional
Hermitian-Einstein condition directly.  Knowing that any such
manifold would be conformal (with a known conformal factor) to an
extremal Kahler toric metric, we can write the Hermitian-Einstein
equation in terms of the Kahler toric coordinates of the extremal
Kahler metric to which it is conformal.

    Let $\bar{g}$ be a four-dimensional Hermitian-Einstein toric metric on a
compact manifold.  The above shows that such a metric is conformal
to an extremal Kahler toric metric $g$, such that
\begin{equation*}
\bar{g}=S^{-2}g.
\end{equation*}

    Now, because $g$ is Kahler toric we can write it in symplectic
coordinates.  Recalling that in the extremal case one can express
the scalar curvature in symplectic coordinates as
$S=\alpha_{1}x_{1}+\alpha_{2}x_{2}+\beta$, the Hermitian-Einstein
metric can be written as
\begin{equation*}
\bar{g}=\begin{pmatrix}
\frac{1}{(\alpha_{1}x_{1}+\alpha_{2}x_{2}+\beta)^{2}}h^{-1} & \\
 & \frac{1}{(\alpha_{1}x_{1}+\alpha_{2}x_{2}+\beta)^{2}}h \\
\end{pmatrix}
\end{equation*} relative to coordinates $(x_{1},x_{2},\phi_{1},\phi_{2})$
where $h^{-1}=\frac{\partial^{2}\Phi}{\partial x^{i}\partial
x^{j}}$.  We want to calculate the Einstein condition on $\bar{g}$
in terms of the $h$ and the constants $\alpha_{i}$ and $\beta$. To
do so we need to calculate the Ricci tensor $\bar{r}$.  Once that
is calculated the Einstein condition is, as usual,
$\bar{r}=\lambda\bar{g}$.  In \cite{Besse}, there is a formula for
the change in Ricci tensor under a conformal transformation. Using
that formula we have
\begin{equation*}
\bar{r}=r+2\frac{\textrm{\textrm{Dd}}S}{S}+\left(\frac{\triangle
S}{S}-3\frac{\langle \nabla S,\nabla S\rangle}{S^{2}}\right)g
\end{equation*}
where $r$ is the Ricci tensor, $\textrm{Dd}$ is the Hessian, and
$\triangle$ is the Laplacian of $g$.  Because both $g$ and
$\bar{g}$ are complex, they and their Ricci tensor are by
construction invariant under the action of the complex structure
$J$.  Therefore, to calculate the Ricci tensor $\bar{r}$ we need
only calculate the values of
$\bar{r}(\partial_{\phi_{i}},\partial_{\phi_{j}})$.

    Two of terms,  $\frac{\triangle S}{S}$ and $\frac{\langle \nabla S,\nabla S\rangle}
{S^{2}}$, we have already calculated in our proof of the above
corollary.  We need to calculate the terms
$\textrm{Dd}S(\partial_{\phi_{i}},\partial_{\phi_{j}})$.  By
definition,
\begin{equation*}
\textrm{Dd}S(\partial_{\phi_{i}},\partial_{\phi_{j}})=\frac{\partial^{2}S}{\partial\phi^{i}
\partial\phi^{j}}-\left(\nabla_{\frac{\partial}{\partial\phi^{i}}}\frac{\partial}{\partial
\phi^{j}}\right)S.
\end{equation*}
The first term is automatically because $S$ does not depend
explicitly on the coordinates $\phi_{i}$.  By calculating the
Christoffel symbols directly, we have that
\begin{equation*}
-\left(\nabla_{\frac{\partial}{\partial\phi^{i}}}\frac{\partial}{\partial\phi^{j}}\right)S=
\frac{1}{2}h_{ik}\frac{\partial
h_{lj}}{\partial x^{k}\alpha_{l}}
\end{equation*}
where we have here used the summation convention for repeated
indices.  Finally, we recall that
\begin{equation*}
{r}(\partial_{\phi_{i}},\partial_{\phi_{j}})=-\frac{1}{2}h_{ik}\frac{\partial}{\partial
x^{k}}\left(h_{lj}\frac{\partial
\textrm{log}(\textrm{det}(h))}{\partial x^{l}}\right).
\end{equation*}
    The above calculations prove the following.

\begin{prop}
The Hermitian metric $\bar{g}=S^{-2}g$ defined above is Einstein
if and only if
\begin{equation}\label{HermEinstToric}
-\frac{1}{2}h_{ik}\frac{\partial}{\partial
x^{k}}\left(h_{lj}\frac{\partial {\rm log}({\rm det}(h))}{\partial
x^{l}}\right)+\frac{1}{S}h_{ik}\frac{\partial h_{lj}}{\partial
x^{k}}\alpha_{l}
\end{equation}
\begin{equation*}
+\left(\frac{1}{S}h_{kl}\frac{\partial {\rm log}({\rm
det}(h))}{\partial
x^{k}}\alpha_{l}-\frac{3}{S^{2}}h_{kl}\alpha_{k}\alpha_{l}\right)h_{ij}=\frac{\lambda}{S^{2}}
h_{ij}
\end{equation*}
for $1\leq i,j\leq 2$ where
$S=\alpha_{1}x_{1}+\alpha_{2}x_{2}+\beta$.
\end{prop}

    Note that when all of the $\alpha_{i}=0$ then equation (\ref{HermEinstToric}) is
equivalent to the Kahler-Einstein condition on Kahler toric four
manifolds.  There is only one known solution to equation
(\ref{HermEinstToric}) when $S$ is not a constant, i.e. when the
metric $\bar{g}$ is non-Kahler Hermitian-Einstein.  That is the
Page metric on $\mathbb{C}P^{2}\sharp\overline{\mathbb{C}P^{2}}$.
Above, we have already written down the extremal Kahler metric
conformal to the Page metric.  Note that in that case
\begin{equation*}
S=-\frac{12}{11}x+\frac{42}{11}=-\frac{12}{11}x_{1}+\frac{42}{11}.
\end{equation*}

Although the calculations are tedious and lengthy, this metric is
a solution of equation (\ref{HermEinstToric}).  The only other
possible compact solutions to equation (\ref{HermEinstToric})
would be on the blow-up of $\mathbb{C}P^{2}$ at two or three
points.

    The Kahler-Einstein condition can be integrated to the Monge-Ampere equation
as described above.  It would be interesting to determine whether
equation (\ref{HermEinstToric}) can be integrated to give a
modified version of that equation.  Unfortunately, we have been
unable to perform such an integration.

\chapter{Fiberwise Kahler Toric Manifolds}

    Kahler toric geometry enjoys many special features which are absent in Kahler
geometry in general.  In particular, the presence of a Legendre
transformation from complex to symplectic coordinates allows one
to fix a symplectic structure and vary the compatible complex
structure via the potential $\Phi$.  This is not always possible
in Kahler geometry. Furthermore, one can explicitly write down a
metric in each Kahler class in terms of rational functions in the
symplectic coordinates. In this section, we develop a
generalization of Kahler toric geometry, which we refer to as
fiberwise Kahler toric geometry, which shares many of these
features.

    We construct fiberwise Kahler toric manifolds in two very different ways.  Construction
one begins by considering a Kahler manifold of cohomogeneity-$d$
under a semisimple group action.  We impose two conditions: we
assume that each principal orbit is a $T^{d}$-bundle over a
coadjoint orbit $G/L$ and we demand that the isotropy
representation of $\mathfrak{g}$ is composed of \textit{distinct}
$\textrm{Ad}(L)$-invariant summands (this is the monotypic case).
We will see that the metric of such a manifold can be viewed as a
Kahler toric metric along with a collection of linear functions in
the symplectic coordinates of the Kahler toric metric.

    Alternatively, construction two begins by first considering a
$d$-parameter family of principal $T^{d}$-bundles over a product
of $m$ Kahler-Einstein manifolds, $(M_{i},g^{\ast}_{i})$.  Let
$(u_{1},..,u_{d})$ denote the $d$ parameters.  We then construct a
metric
\begin{equation*}
g=h_{ij}(u)du_{i}\otimes
du_{j}+h_{ij}(u)\theta_{i}\otimes\theta_{j}+\sum_{k=1}^{m}A_{k}(u)\pi^{\ast}g^{\ast}_{k}
\end{equation*}
where the $\theta_{i}$ are the connection forms of the
$T^{d}$-bundle, the $\pi^{\ast}g^{\ast}_{k}$ are the lifts of the
Kahler-Einstein metrics on the $M_{k}$, and the $A_{k}$ are
nonnegative functions.  By defining a complex structure which
satisfies $J\frac{\partial}{\partial
u_{i}}=\frac{\partial}{\partial\phi_{i}}$ where
$\frac{\partial}{\partial\phi_{i}}$ is the vector field dual to
$\theta_{i}$, we construct a Kahler metric which, in form, is
identical to the one obtained by the previous construction.

    Therefore, given completely different geometric assumptions,
we construct metrics with identical properties.  That is, while
the manifolds constructed may be different, their metrics and
curvature equations are of the same form. We will refer to both
kinds as fiberwise Kahler toric manifolds. This chapter will be
devoted to making the above discussion more precise and to
analyzing the Kahler-Einstein and extremal Kahler conditions on
these manifolds.

\section{Fiberwise Kahler Toric Geometry}

    We now present our two constructions of fiberwise Kahler toric manifolds.

\subsection{Construction One}

    Take $M$ to be a Kahler manifold of cohomogeneity-$d$ under
the action of a compact connected semisimple Lie group $G$ where
\begin{equation}\label{FKToricSub}
G\rightarrow(M^{n},g,J)\rightarrow(M/G,\check{g}).
\end{equation}
Note that $(M/G,\check{g})$ is not a smooth manifold in general.
However, on an open dense set in $M$ the submersion onto the open
interior of $M/G$ is smooth.  Furthermore, we assume that the
metric is irreducible and is not hyperkahler so that the action of
$G$ preserves $J$ and $\omega$.

    Because $G$ is semisimple, there exists a $G$-equivariant moment
map, $\mu$, for the action of $G$ on $M$.  That is
\begin{equation*}
\mu:M\rightarrow\mathfrak{g}^{\ast}
\end{equation*} such that
\begin{equation*}
d\mu(v)(X)=\omega(X,v)=g(JX,v)
\end{equation*}
where $v\in TM$ and $X\in\mathfrak{g}$.  Note that we have
suppressed the distinction between the element of $\mathfrak{g}$
and the Killing field on $M$ induced by that element.  We have the
following lemma:

\begin{lemma}\label{FKTlemma1}
Each principal orbit $G/K$ is a fiber bundle over a coadjoint
orbit $G/L$
\begin{equation*}
L/K\rightarrow G/K\rightarrow G/L
\end{equation*}
such that ${\rm dim}(L/K)=d$.
\end{lemma}
\begin{proof}
The fiber bundle above is induced by the moment map, $\mu$, which,
when restricted to a principal orbit $G/K$, maps

\begin{equation*}
\mu:G/K\rightarrow G/L
\end{equation*}
where $G/L$ is a coadjoint orbit in $\mathfrak{g}^{\ast}$.  We
assumed that $G$ was connected and semisimple, this implies that
$G/L$ is simply-connected.  Therefore, $L$ and $L/K$ are
connected.

    For $v$ a vector tangent to the fiber $L/K$, $d\mu(v)=0$.
This is true if and only if
\begin{equation*}
d\mu(v)(X)=g(JX,v)=-g(X,Jv)=0
\end{equation*}
for all $X\in\mathfrak{g}$.

Therefore, $Jv$ is orthogonal to the principal orbits; $Jv$
belongs to the horizontal distribution.  The vector $Jv$ is
orthogonal to the principal orbits if and only if $d\mu(v)=0$.  As
$J$ is one-to-one we deduce immediately that
$\textrm{dim}(L/K)=\textrm{dim}(M/G)=d$.
\end{proof}

We have seen how the complex structure, $J$, sends horizontal
vector fields to vertical vector fields tangent to the fiber
$L/K$.  In the cohomogeneity-one case, $\textrm{dim}(L/K)=1$ and
$L/K\cong S^{1}$.  By contrast, $L/K$ in the general
cohomogeneity-$d$ case need not be abelian.  It is natural to
first consider the case in which $L/K$ is abelian.

\textbf{ASSUMPTION} \textbf{ONE}:
\begin{equation*}
L/K\cong T^{d}.
\end{equation*}

    As described above, the Lie
algebra, $\mathfrak{g}$ can be decomposed as
$\mathfrak{g}=\mathfrak{k}\oplus\mathfrak{p}$.  Given assumption
one, we have that

\begin{equation*}
\mathfrak{p}=\mathfrak{t}^{d}\oplus\mathfrak{p}_{1}\oplus...\oplus\mathfrak{p}_{m}
\end{equation*}
where $\mathfrak{t}^{d}$ is the space on which $Ad(K)$ acts
trivially and the $\mathfrak{p_{i}}$ are non-trivial irreducible
$K$-modules. Let $d_{i}$ denote the real dimension of the summand
$\mathfrak{p}_{i}$.  Note that
$\mathfrak{k}\oplus\mathfrak{t}^{d}$ is the Lie algebra of the
group $L$.  We make one more assumption.

\textbf{ASSUMPTION} \textbf{TWO}:

~

~~~~~~~~~~~~~~~~~~~~~~~~~~~~~~~~~~~~~~~~~~~~~~~~~~~~~\textit{The
$\mathfrak{p}_{i}$ are distinct.}

~

    On a Kahler toric manifold, the horizontal distribution of the
quotient by group action
\begin{equation*}
T^{n}\rightarrow (M^{2n},g)\rightarrow (M/G,\check{h})
\end{equation*}
is integrable.  Is this the case when $G$ is semisimple and
Assumption One holds?

\begin{theorem}
The O'Neill invariant $A$ of the submersion (\ref{FKToricSub}) is
zero when assumption one holds.  Equivalently, the action of $G$
on $M$ is \textbf{orthogonally transitive}.
\end{theorem}
\begin{proof}
Take $\{\partial_{\phi_{i}}\}_{i=1,..,d}$ to be a basis for the
Lie algebra of $T^{d}\cong L/K$.  Clearly,
$[\partial_{\phi_{i}},\partial_{\phi_{j}}]=0$ for all $i$ and $j$.
We demonstrated in the proof of the last lemma that
$\{J\partial_{\phi_{i}}\}_{i=1,..,d}$ will form a basis for the
horizontal distribution.  $\partial_{\phi_{i}}$ is a Killing
vector field and an automorphism of the complex structure $J$;
therefore, Besse tells us that
\begin{equation*}
\mathcal{L}_{\partial_{\phi_{i}}}J=0\Leftrightarrow[\partial_{\phi_{i}},JX]=J[\partial_{
\phi_{i}},X]
\end{equation*}
for all $X\in TM$. From this equation we deduce immediately that
\begin{equation*}
[\partial_{\phi_{i}},J\partial_{\phi_{j}}]=[\partial_{\phi_{i}},\partial_{\phi_{j}}]=0
\end{equation*}
for all $i$ and $j$.  As $J$ is an integrable complex structure,
its Nijenhuis tensor, $N_{J}$ vanishes.  This implies that
\begin{equation*}
N_{J}(\partial_{\phi_{i}},\partial_{\phi_{j}})=[J\partial_{\phi_{i}},J\partial_{\phi_{j}}]
-J[J\partial_{\phi_{i}},\partial_{\phi_{j}}]-J[\partial_{\phi_{i}},J\partial_{\phi_{j}}]-[
\partial_{\phi_{i}},\partial_{\phi_{j}}]=0.
\end{equation*}
It follows immediately that
\begin{equation*}
[J\partial_{\phi_{i}},J\partial_{\phi_{j}}]=0
\end{equation*}
for all $i$ and $j$.

Therefore, $\{J\partial_{\phi_{i}}\}_{i=1,..,d}$ forms a commuting
basis for the horizontal distribution.  The existence of such a
basis implies that the horizontal distribution is integrable.
\end{proof}

The above proof demonstrates that
$\{-J\partial_{\phi_{i}}\}_{i=1,..,d}$ forms a coordinate basis
for the horizontal distribution and for the quotient space $M/G$.
Set $\partial_{u_{i}}=-J\partial_{\phi_{i}}$ for all $i$.

    Dancer and Wang in \cite{DandW} analyzed the structure of
Kahler manifolds of cohomogeneity-one under the action of a
semisimple Lie group such that assumption two holds, the $d=1$
case.  Our results in the rest of this section are a
generalization of their work to the cohomogeneity-$d$ case.

    So far we have made no mention of the Kahler structure of the
metric.  Before we do, we can deduce the general form of the
metric.  Let $Q$ be a bi-invariant metric on $\mathfrak{g}$ such
that $Q(\mathfrak{k},\mathfrak{p})=0$.  Furthermore, we have that
$Q(\mathfrak{p}_{i},\mathfrak{p}_{j})$ for all $i$ and $j$. Choose
an orthonormal basis $\{Y^{i}_{j}\}_{j=1,...,d_{i}}$ for
$\mathfrak{p}_{i}$ relative to the background metric $Q$.  For any
metric $g$ on $M$ we have, by Assumption Two and Schur's Lemma,
that
\begin{equation*}
g(\cdot,\cdot)|_{\mathfrak{p}_{1}\oplus...\oplus\mathfrak{p}_{m}}=\bigoplus^{\bot}_{i}
A_{i}Q(\cdot,\cdot)|_{\mathfrak{p}_{i}}
\end{equation*}
where the $A_{i}$ are positive functions on the interior of the
quotient space $M/G$.  That is $g(Y^{i}_{j},Y^{i}_{j})=A_{i}$ for
all $i,$ $j$, and $k$.  Furthermore, set
$h_{ij}=g(\partial_{\phi_{i}},\partial_{\phi_{j}})$.  Because of
the invariance of the metric under the complex structure $J$, we
deduce immediately that,
$h_{ij}=g(\partial_{u_{i}},\partial_{u_{j}})$.  Therefore, the
metric, $g$, on the interior of the quotient space can be written
as
\begin{equation*}
g=\begin{pmatrix}
h_{ij}&&&& \\
&h_{ij}&&&  \\
&& A_{1}&&  \\
&&& ... &  \\
&&&& A_{m}  \\
\end{pmatrix}
\end{equation*}
relative to the basis
$(\partial_{u_{1}},...,\partial_{u_{d}},\partial_{\phi_{1}},...,\partial_{\phi_{d}},
Y^{1}_{1},...,Y^{m}_{d_{m}})$.
The metric is determined by $\frac{d(d+1)}{2}+m$ functions of the
coordinates $(u_{1},...,u_{d})$.  In addition to being pairwise
inequivalent $K$-modules, the $\mathfrak{p_{i}}$ are also pairwise
inequivalent $L$-modules.  For each $p$ in the interior of $M/G$,
let $g_{p}$ be the restriction of the metric to the principal
orbit that is the inverse image of the point $p$.  Therefore
$g_{p}$ induces an $\textrm{Ad(}L)$-invariant inner product on
$\mathfrak{q}=\mathfrak{p_{1}}\oplus...\oplus\mathfrak{p_{m}}$
which can be identified with the tangent space to $G/L$.  The
induced $G$-invariant metric on $G/L$ will be denoted by
$g_{p}^{\ast}$.

~

\textbf{The Kahler Condition}:

~

In the Kahler case, $d\omega=0$ where
$\omega(\cdot,\cdot)=g(J\cdot,\cdot)$.

\begin{theorem}\label{FKTTheorem1}
Let $(M,g,J)$ be a Fiberwise Kahler Toric manifold.  The Kahler
condition is equivalent to

\item (i)
$\omega_{p}([X,Y]_{\mathfrak{p}},Z)+\omega_{p}([Y,Z]_{\mathfrak{p}},X)+\omega_{p}([Z,X]_{
\mathfrak{p}},Y)=0$

\item
(ii)
$\partial_{u_{i}}\omega(X,Y)+\omega(\partial_{u_{i}},[X,Y]_{\mathfrak{p}})=0$

where $p$ is a point in the interior of $M/G$ and
$X,Y,Z\in\mathfrak{p}$ and
\item
(iii) $\bar{g}=\begin{pmatrix} h_{ij} & \\ & h_{ij}
\\ \end{pmatrix}$ is locally Kahler toric
relative to basis
$(\partial_{u_{1}},...,\partial{u_{d}},\partial_{\phi_{1}},...,\partial_{\phi_{d}})$,
i.e. $h_{kl}=\frac{\partial^{2}\eta}{\partial u^{k}\partial
u^{l}}$ for some function $\eta=\eta(u_{1},...,u_{d})$.
\end{theorem}

\begin{proof}
First, consider arbitrary vector fields $X$, $Y$, and $Z$ in $TM$.
The following equality always holds for any two-form $\omega$
\begin{equation*}
d\omega(X,Y,Z)=X(\omega(Z,Y))+Y(\omega(X,Z))+Z(\omega(Y,X))+\omega([X,Y]_{\mathfrak{p}},Z)
+\omega([Y,Z]_{\mathfrak{p}},X)+\omega([Z,X]_{\mathfrak{p}},Y)
\end{equation*}
In the Kahler case, $d\omega=0$.  If we choose
$X,Y,Z\in\mathfrak{p}$, the invariance of the Kahler form under
the action of the group $G$ immediately gives $(i)$.  Setting
$Z=\partial_{u_{i}}$, we obtain $(ii)$ by noting the
$[\partial_{u_{i}},X]=0$ for all $X\in\mathfrak{p}$.  Part $(iii)$
can be obtained by noting that the metric $\bar{g}$ is
well-defined locally relative to the basis given and all the
vector fields commute.  Therefore, locally, the complex structure
and metric are that of a toric manifold and $\omega$ restricted to
the basis given is precisely $\bar{\omega}$.  If $\omega$ is
closed then $\bar{\omega}$ is necessarily closed and the metric
$\bar{g}$ is Kahler toric.  The only remaining case to check is
when $X=\partial_{u_{i}},$ $Z=\partial_{u_{j}},$ and
$Y=Y^{k}_{m}\in\mathfrak{p}_{k}$.  The Kahler condition is
automatically satisfied since each $\partial_{u_{i}}$ commutes
with $Y^{k}_{m}$.  Therefore, the equations listed in the theorem
are the only ones not automatically satisfied.
\end{proof}

    For every such manifold, there is a natural projection from $(M,g)$ onto $(G/L,g^{\ast})$.
The fiber of this projection is a Kahler toric manifold
$(N,\bar{g})$ such that
\begin{equation*}
T^{d}\rightarrow(N,\bar{g})\rightarrow(N/T^{d},\check{\bar{g}})\cong(M/G,\check{g}).
\end{equation*}
Moreover, the existence of this Kahler toric manifold will allow
us to later make use of the 'symplectic' coordinates obtained by
changing the coordinates of the base manifold $M/G$.  Before we do
that, however, let us first analyze in greater detail the action
of the complex structure $J$.

\begin{prop}
Let $(M,g,J)$ be a Kahler manifold of cohomogeneity-$d$ under the
action of a semisimple Lie group such that assumptions one and two
hold.  There is a coadjoint orbit $G/L$ with $K\subset L$ such
that
\item (1) for each point $p$ in the interior of the quotient space
$M/G$, $g_{p}=g\mid_{\pi^{-1}(p)}$ induces a $G$-invariant metric
$g_{p}^{\ast}$ on $G/L$ and the projection
$\mu:(G/K,g_{p})\rightarrow(G/L,g_{p}^{\ast})$ is a Riemannian
submersion with totally geodesic $T_{d}$ fibers.
\item (2) there is a complex structure $J^{\ast}$, which does not
depend on the choice of $p$, on $G/L$ and relative to this complex
structure $g_{p}^{\ast}$ is Kahler.
\end{prop}
\begin{proof}
Part $(1)$ follows directly from the discussion above.  To obtain
$(2)$, note that $J_{p}$ is $\textrm{Ad}(K)$-invariant and
skew-symmetric with respect to the metric $g_{p}$.  Therefore, we
deduce that $J_{p}$ leaves each $\mathfrak{p}_{i}$ invariant
because of assumption two.  The above discussion also demonstrates
that $J_{p}$ is independent of the choice of $p$.  Because each
summand $\mathfrak{p}_{i}$ is irreducible, all invariant metrics
differ only by a constant.  This implies that $J\mid_{q}$ is
skew-symmetric with respect to the Killing form.

    From Lemma (\ref{FKTlemma1}) we see that
$\textrm{ad}(\partial_{\phi_{i}})$ is $\textrm{Ad}(K)$-invariant
for every $i=1,...,d$ and if, of course, skew with respect to the
Killing form. On each $\mathfrak{p}_{i}$, $J$ is proportional to
$\textrm{ad}(\partial_{\phi_{i}})$ and also commutes with that
operator. This implies that $J\mid_{q}$ is in fact
$\textrm{Ad}(L)$-invariant and gives rise to an almost-complex
structure $J^{\ast}$ on $G/L$.  To complete the proof of this
proposition, we need to show that this complex structure is
integrable and that the induced invariant metric on $G/L$ is
Kahler with respect to this complex structure.

~

\textit{Integrability}:

~

    The integrability of $J^{\ast}$ is equivalent to the vanishing
of the Nijenhuis tensor which is equivalent to
\begin{equation*}
[J^{\ast}X,J^{\ast}Y]_{\mathfrak{q}}-[X,Y]_{\mathfrak{q}}-J^{\ast}[X,J^{\ast}Y]
_{\mathfrak{q}}-J^{\ast}[J^{\ast}X,Y]_{\mathfrak{q}}=0
\end{equation*}
for all $X,Y\in\mathfrak{q}$.  Now, the submersion
$\mu:G/K\times\{p\}\rightarrow G/L$ is independent of the point
$p$ as discussed above.  Take $\check{X}$ and $\check{Y}$ to be
basic vector fields of this submersion.  We have that
$d\mu[\check{X},\check{Y}]=[d\mu(\check{X}),d\mu(\check{Y})]$.
Also, we see that $d\mu(J\check{X})=J^{\ast}(d\mu(\check{X}))$
because of the $\textrm{Ad}(L)$ invariance of $J$.  This together
with the integrability of $J$ implies the integrability of
$J^{\ast}$.

~

\textit{Kahler}:

~

    We now must demonstrate that the metric $(g_{p}^{\ast},J^{\ast})$
is Kahler.  The Kahler form, $\omega_{p}^{\ast}$, is closed by
part $(i)$ of (\ref{FKTTheorem1}).
\end{proof}

    Instead of the construction above, we could construct manifolds of this type by the
following process.  Let $G/L$ be a coadjoint orbit, whose isotropy
representation is composed of inequivalent summands, with complex
structure $J^{\ast}$.  Take $g^{\ast}_{p}$ to be a $d$-parameter
family of Kahler metrics on $G/L$.  The first Chern class of
$(G/L,J^{\ast})$ is positive and therefore all class in
$H^{2}(G/L;\mathbb{Z})$ are of type $(1,1)$.  Every $T^{d}$-bundle
over $G/L$ is of the form $G/K$ where $L/K\cong T^{d}$.

    Any principal $T^{d}$-bundle, $G/K$, over $G/L$ ($\mu:G/K\rightarrow
G/L$) is characterized by $d$ invariant $S^{1}$ connections
$\theta_{i}$, $1\leq i\leq d$.  We can construct a metric
\begin{equation*}
g=h_{ij}du_{i}\otimes
du_{j}+h_{ij}\theta_{i}\otimes\theta_{j}+\mu^{\ast}g^{\ast}_{p}
\end{equation*}
locally on $G/K\times\mathbb{R}^{d}$ where
$h_{ij}=h_{ij}(u_{1},...,u_{d})$.  We define the complex structure
$J$ to be the $J^{\ast}$ when restricted to vector field tangent
to $G/L$ and to otherwise satisfy the equations
$J\partial_{u_{i}}=\partial_{\phi_{i}}$ where
$\partial_{\phi_{i}}$ is the dual vector field to the connection
$\theta_{i}$.  Furthermore, we demand that
$h_{ij}=\frac{\partial^{2}\eta}{\partial_{u_{i}}\partial_{u_{j}}}$.
It follows that this defines a Kahler metric which is invariant
under the action of $G$.

    We return now to the Kahler condition.  For $X,Y\in\mathfrak{q}$, the Kahler condition
implies that
\begin{equation*}
\partial_{u_{i}}(\omega(X,Y))+\omega(\partial_{u_{i}},[X,Y]_{\mathfrak{p}})=0
\end{equation*}
for all $1\leq i\leq d$.  As mentioned above, the moment map
induces a fibration $\mu:G/K\rightarrow G/L$ on each principal
orbit.  This is a principal $T^{d}$-bundle over a coadjoint orbit.
Since, as proved above, each fiber is totally geodesic, the
O'Neill invariant $T^{\mu}$ of this submersion vanishes.  However,
the invariant $A^{\mu}$ of the submersion does not vanish.  From
the definition of
\begin{equation*}
A^{\mu}_{X}Y=\frac{1}{2}\mathcal{V}[X,Y]_{\mathfrak{p}}
\end{equation*}
for $X,Y\in\mathfrak{q}$ where $\mathcal{V}$ is the vertical
projection onto the $T^{d}$ fibers.  Letting
$\{\partial_{\phi_{i}}\}$ be a basis from $\mathfrak{t}^{d}$
letting $\omega^{\ast}_{j}=\omega^{\ast}_{\mathfrak{p}_{j}}$, we
have that
\begin{equation*}
A^{\mu}_{X}Y=-\frac{1}{2}\sum_{k=1}^{d}\sum_{j=1}^{m}b_{kj}\omega^{\ast}_{j}(X,Y)\partial_{
\phi_{k}}
\end{equation*}
where the $b_{jk}$ are the eigenvalues of the curvature form of
the connection $\theta_{j}$.

    Take $X_{j},Y_{j}\in\mathfrak{p}_{j}$, the Kahler condition
implies that
\begin{equation*}
\partial_{u_{i}}(\omega(X_{i},Y_{i}))+\omega(\partial_{u_{i}},[X_{j},Y_{j}]_{\mathfrak{p}_{j}})
=0
\end{equation*}
which implies that
\begin{equation*}
\partial_{u_{i}}(g(X_{i},Y_{i}))=g(\partial_{\phi_{i}},\sum_{k=1}^{d}b_{kj}\omega^{\ast}_{j}
(X_{j},Y_{j})\partial_{\phi_{k}}).
\end{equation*}
Using the form of the metric $g$, it is straightforward to
calculate that the Kahler condition implies that
\begin{equation*}
\partial_{u_{i}}A_{j}=\sum_{k=1}^{d}h_{ik}b_{kj}.
\end{equation*}
multiplying both sides by $h^{-1}$ and making use of the summation
convention, we have that
\begin{equation}\label{FKTKahlercomp}
h^{ik}\partial_{u_{i}}A_{j}=b_{kj}
\end{equation}
for all $1\leq k\leq d$ and all $1\leq j\leq m$.  Therefore, up to
constants, the Kahler condition completely determines the
functions $A_{j}=A_{j}(u_{1},...,u_{d})$ in terms of the $h_{ij}$
and the $b_{ij}$.  Although the $A_{j}$ satisfy
(\ref{FKTKahlercomp}), this equation is somewhat difficult to use
directly as determining an explicit formula for the $A_{j}$ would
involve integrating the functions $h_{ij}$.  If, however, we
change coordinates on the quotient space, we can place equation
(\ref{FKTKahlercomp}) in a form which can be explicitly
integrated.

    We have demonstrated that within every manifold of this type is a Riemannian submersion
over coadjoint orbit $G/L$ such that the fiber of this submersion
is a Kahler toric manifold with metric
\begin{equation*}
\bar{g}=\begin{pmatrix}h_{ij} & \\
 & h_{ij}
\end{pmatrix}
\end{equation*}
relative to vector fields
$(\partial_{u_{1}},...,\partial_{\phi_{d}})$ such that
$h_{ij}=\frac{\partial^{2}\eta}{\partial_{u_{i}}\partial_{u_{j}}}$.
Instead, of these complex coordinates, we could express the Kahler
toric metric in symplectic coordinates by setting
$x_{i}=\frac{\partial\eta}{\partial u_{i}}$.  With this change, we
have
\begin{equation*}
\bar{g}=\begin{pmatrix}h^{ij} & \\
 & h_{ij}
\end{pmatrix}
\end{equation*}
relative to basis
$(\partial_{x_{1}},...,\partial_{x_{d}},\partial_{\phi_{1}},...,\partial_{\phi_{d}})$.

    We can now view the functions $A_{j}$ as functions of the
$x_{i}$.  After expressing equation (\ref{FKTKahlercomp}) in these
coordinates, we obtain the following lemma.

\begin{lemma}
The Kahler condition implies that the functions
$A_{j}(x_{1},...,x_{d})$ are of the form
\begin{equation}
A_{j}=\sum_{k=1}^{d}b_{kj}x_{k}+a_{j}
\end{equation}
where the $a_{j}$ are constants.
\end{lemma}
\begin{proof}
From calculations given in the section on Kahler toric manifolds,
$h^{kl}\frac{\partial}{\partial u^{l}}=\frac{\partial}{\partial
x_{k}}$.  Equation (\ref{FKTKahlercomp}) implies that
\begin{equation*}
h^{ik}\partial_{u_{i}}A_{j}=\frac{\partial}{\partial
x^{k}}A_{j}=b_{kj}
\end{equation*}
for all $1\leq j\leq m$ and all $1\leq k\leq d$.  The lemma
follows immediately.
\end{proof}

    We pause here to discuss the significance of the constants $b_{kj}$ and $a_{j}$.  Each
principal orbit $G/K$ is a principal $T^{d}$-bundle over coadjoint
orbit $G/L$.  Given a particular $G/L$, a principal $T^{d}$-bundle
over that orbit is determined by the constants $b_{kj}$.  While
changing the constants $a_{j}$ alters the metric is does not
change the principal orbit type.

These above results are summarized in the following theorem:

\begin{theorem}\label{FKTth1}
The metric on a Kahler manifold $(M,g,J)$ of cohomogeneity-$d$
under the action of a semisimple Lie group such that assumptions
one and two hold can be written locally as
\begin{equation*}
g=\begin{pmatrix} h^{ij} & & & & \\
& h_{ij} & & & \\
& & (\sum_{k=1}^{d}b_{k1}x_{k}+a_{1})Id_{d_{1}} & & \\
& & & ... & \\
& & & & (\sum_{k=1}^{d}b_{km}x_{k}+a_{m})Id_{d_{m}} \\
\end{pmatrix}
\end{equation*}
relative to basis
$(\partial_{x_{1}},...,\partial_{x_{d}},\partial_{\phi_{1}},...,\partial_{\phi_{d}},
Y^{1}_{1},...,Y^{m}_{d_{m}})$
and
\begin{equation*}
h^{kl}=\frac{\partial^{2}\Phi}{\partial x^{k}\partial x^{l}}
\end{equation*}
for $\Phi=\Phi(x_{1},...,x_{d})$.
\end{theorem}
\subsection{Construction Two}

    We now construct Kahler manifolds with metrics of the same form using different geometric
assumptions.  Let $(M_{k},J_{k}^{\ast},g_{k}^{\ast})$ be compact
Kahler-Einstein manifolds of real dimension $d_{k}$ with positive
Einstein constant for all $1\leq k\leq m$.  We normalize each of
these manifolds so that the Einstein constant is equal to one.
Furthermore, let $P$ be a principal $T^{d}$-bundle over
$M_{1}\times...\times M_{m}$. Manifolds of this type were studied
by Wang and Ziller in \cite{WangZiller} and we rely on their
construction here.

    A principal $T^{d}$-bundle $\pi:P\rightarrow M_{1}\times...\times M_{m}$ is characterized
by $d$ connections $\theta_{i}$ which satisfy
$d\theta_{i}=\pi^{\ast}(\sum_{k=1}^{m}b_{ik}\omega_{k}^{\ast})$
for constants $b_{ij}$.  We can construct a metric on such a torus
bundle by writing
\begin{equation*}
h_{ij}\theta_{i}\otimes\theta_{j}+\sum_{k=1}^{m}A_{k}\pi^{\ast}g_{k}^{\ast}
\end{equation*}
where the $A_{k}$ are, for the moment, constants and $h_{ij}$ is a
positive-definite matrix of constants.  (Note that here we are
summing over $i$ and $j$.)  It is important to note the product
manifold $(M_{1}\times...\times
M_{m},\sum_{k=1}^{m}A_{k}g_{k}^{\ast})$, is itself not necessarily
Einstein, though it is of constant scalar curvature. The map $\pi$
is a Riemannian submersion with totally geodesic fibers.  Let
$A^{\pi}$ denote the O'Neill invariant of this submersion.  By
definition,
\begin{equation*}
A^{\pi}_{X}Y=\frac{1}{2}\mathcal{V}[X,Y]
\end{equation*}
for $X$ and $Y$ vector fields in the horizontal distribution of
the submersion $\pi$.   If we let $\partial_{\phi_{i}}$ be the
dual vector field to the connection $\theta_{i}$, it is
straightforward to calculate (see \cite{WangZiller} for details)
that
\begin{equation*}
A^{\pi}_{X}Y=-\frac{1}{2}\sum_{j=1}^{d}\sum_{k=1}^{m}b_{jk}\omega_{k}^{\ast}(X,Y)\partial_{
\phi_{j}}.
\end{equation*}

    Next, we consider a $d$-parameter family of such metrics, i.e. we set
$h_{ij}=h_{ij}(u_{1},...,u_{d})$ and
$A_{k}=A_{k}(u_{1},...,u_{d})$.  We now construct what is to be
our fiberwise Kahler metric by setting
\begin{equation*}
g=h_{ij}du_{i}\otimes
du_{j}+h_{ij}\theta_{i}\otimes\theta_{j}+\sum_{k=1}^{m}A_{k}\pi^{\ast}g_{k}^{\ast}.
\end{equation*}
Locally, we now have a Riemannian manifold $(M,g)$ of dimension
$n=2d+\sum_{k=1}^{m}d_{k}$.  To turn this into a Kahler manifold,
we first need to define a complex structure $J$ which is
compatible with this metric. If $X$ is the lift of a vector field
on the product of the $M_{i}$ then set
$JX=\sum_{k=1}^{m}\pi^{\ast}J_{k}^{\ast}X$. Also, let
$J\partial_{u_{i}}=\partial_{\phi_{i}}$ for all $i$.  These
requirements completely determine a compatible almost complex
structure and it is straightforward to see that this almost
complex structure is in fact integrable.  Thus, locally, $(M,g,J)$
is a complex manifold.  We note that the Einstein condition on
metrics of this type were studied in both the Hermitian and the
Kahler case for $d=1$ by Wang and Wang in \cite{WangWang}.

    We define the two form $\omega$ in the usual way by the equation $\omega(\cdot,\cdot)
    =g(J\cdot,\cdot)$.  The Kahler
condition is equivalent to $d\omega=0$.  Using a well-known
identity, this condition is equivalent to
\begin{equation*}
X(\omega(Z,Y))+Y(\omega(X,Z))+Z(\omega(Y,X))+\omega([X,Y],Z)+\omega([Y,Z],X)+\omega([Z,X],Y)=0
\end{equation*}
for all $X,Y,$ and $Z$ in $TM$.  Consider the induced metric
\begin{equation}\label{Kahlercondition}
\bar{g}=h_{ij}du_{i}\otimes
du_{j}+h_{ij}\theta_{i}\otimes\theta_{j}
\end{equation}
with induced complex structure
$\bar{J}=J|_{\textrm{span}\{\partial_{u_{1}},...,\partial_{\phi_{d}}\}}$.
Let $\bar{\omega}(\cdot,\cdot)=\bar{g}(\bar{J},\cdot,\cdot)$. From
equation (\ref{Kahlercondition}), it is clear that $d\omega=0$
implies that $d\bar{\omega}=0$.  Put another way, if $g$ is Kahler
then $\bar{g}$ is Kahler toric, i.e.
$h_{kl}=\frac{\partial^{2}\eta}{\partial u^{k}\partial u^{l}}$ for
some function $\eta=\eta(u_{1},...,u_{d})$.

    Next, let $X, Y$ and $Z$ be the lifts of vector fields on the product of the
Kahler-Einstein manifolds $(M_{k},g^{\ast}_{k}$. For such vector
fields $d\omega(X,Y,Z)=0$ because each of the $\omega^{\ast}_{k}$
are closed.  Lastly, we see that, given $X$ and $Y$ as above,
$d\omega$ implies that
\begin{equation*}
\partial_{u_{i}}\omega(X,Y)+\omega(\partial_{u_{i}},[X,Y])=0
\end{equation*}
for all $1\leq i\leq d$.  This equation implies that
\begin{equation*}
\partial_{u_{i}}(g(X,Y))=g(\partial_{\phi_{i}},\sum_{j=1}^{d}b_{jk}\sum_{k=1}^{m}\omega^{
\ast}_{k}(X,Y)\partial_{\phi_{j}}).
\end{equation*}
From the form of the metric determined above, this becomes
\begin{equation}\label{Kahlercondition2}
\partial_{u_{i}}A_{k}=\sum_{j=1}^{d}h_{ij}b_{jk}
\end{equation}
for all $1\leq j\leq d$ and all $1\leq k\leq m$.  Up to a
constant, the Kahler condition completely determines the function
$A_{i}(u)$. We can express this equation in a simplified form if
we express the Kahler toric metric $\bar{g}$ in terms of
symplectic coordinates. That is, let
$x_{i}=\frac{\partial\eta}{\partial_{u_{i}}}$ for all $\leq i\leq
d$.  Given this coordinate change
\begin{equation*}
\bar{g}=
\begin{pmatrix}
h^{ij} & \\
  & h_{ij}\\
\end{pmatrix}
\end{equation*}
relative to basis $(\partial_{x_{1}},...,\partial_{\phi_{d}})$.
Given this coordinate change, equation (\ref{Kahlercondition2})
implies that
\begin{equation*}
A_{k}=\sum_{j=1}^{d}b_{jk}x_{j}+a_{j}.
\end{equation*}
This completes the calculation of the Kahler condition.

    As in construction one, we wish to describe more precisely the significance of
the constants $b_{jk}$ and $a_{j}$.  Each point in the interior of
$M/G$ is a principal $T^{d}$-bundle over the product of
Kahler-Einstein $M_{1}\times...\times M_{m}$.  Given a product
manifold, the various principal $T^{d}$-bundles are described by
the constants $b_{jk}$.  Varying the constants $a_{j}$ alters the
metric but does not alter these principal torus bundles.  However,
the $a_{j}$ in part determine the zero sets of the functions
$A_{j}$.  Whether or not some of the $A_{j}$ vanish at special
orbits will modify the manifold on which the metric is defined.
While the $a_{j}$ do not affect the principal orbit type, they do
affect the topology of the total manifold $M$.

The above results are summarized in the following theorem:
\begin{theorem}\label{FKTth2}
The metric on Kahler manifold $(M,g,J)$ obtained via construction
two can be written locally as
\begin{equation*}
g=\begin{pmatrix} h^{ij} & & & & \\
& h_{ij} & & & \\
& & (\sum_{k=1}^{d}b_{k1}x_{k}+a_{1})Id_{d_{1}} & & \\
& & & ... & \\
& & & & (\sum_{k=1}^{d}b_{km}x_{k}+a_{m})Id_{d_{m}} \\
\end{pmatrix}
\end{equation*}
relative to basis
$(\partial_{x_{1}},...,\partial_{x_{d}},\partial_{\phi_{1}},...,\partial_{\phi_{d}},
Y^{1}_{1},...,Y^{m}_{d_{m}})$
and
\begin{equation*}
h^{kl}=\frac{\partial^{2}\Phi}{\partial x^{k}\partial x^{l}}
\end{equation*}
for $\Phi=\Phi(x_{1},...,x_{d})$ where
$(Y^{k}_{1},...,Y^{k}_{d_{k}})$ is an orthonormal basis for the
Kahler-Einstein manifold $(M_{k},g^{\ast}_{k})$.
\end{theorem}
    Note that this metric has exactly the same form as the metric determined by construction
one.  This justifies the following definition:
\begin{defn}
We refer to a Kahler manifold $(M,g,J)$ obtained via either
construction one or construction two as a \textbf{fiberwise Kahler
toric manifold}.
\end{defn}

When calculating the curvature equations on a fiberwise Kahler
toric manifold, we will not distinguish between to two types of
manifolds determined by the above constructions.  Since the
curvature equations are identical, it is not necessary to
distinguish between them for the purpose of calculations. Finally,
we state the following definition:
\begin{defn}
A fiberwise Kahler toric manifold is of \textbf{cohomogeneity-$d$}
if $h_{ij}$ is a $d\times d$ matrix.
\end{defn}
In the case of manifolds obtained via construction one, the
appellation is strictly true as the manifolds are of
cohomogeneity-$d$.  In the case of manifolds obtained via
construction two, the metrics need only have a $T^{d}$ symmetry
and will therefore not necessarily be of cohomogeneity-$d$.  As
the form of the metric on the two types of manifolds are the same,
we believe this abuse will not cause confusion.

    Before discussing various curvature conditions on fiberwise Kahler toric manifolds, we
make a few comments about the coordinates we are using.  The
coordinates on the induced Kahler toric metric, $\bar{g}$, are the
symplectic Darboux coordinates.  However, the metric on the
fiberwise Kahler manifold presented here is not written in Darboux
coordinates at all.  Furthermore, these metrics are not written in
complex coordinates.  We believe that it is best to think of
fiberwise Kahler toric manifolds, and these coordinates, as Kahler
toric manifolds with metrics written in symplectic coordinates
along with additional structure.  The connection with Kahler toric
geometry will be central in our work below.

\section{Curvature Equations and Conditions}

    In this section, we calculate the scalar curvature equation
and the equations of the Einstein condition on a fiberwise Kahler
toric manifold.  After presenting these equations, we will examine
the extremal Kahler condition on manifolds of this type.

\subsection{The Einstein Equations}

    We demonstrated above that the action of the group $G$ is
orthogonally transitive in the case of a fiberwise Kahler toric
manifold.  The full Einstein equations on a manifold of
cohomogeneity-$d$ under an orthogonally transitive group action
are given in the section on Riemannian submersions.  Not
surprisingly, the equations given in theorem
(\ref{cohomdEinsteincondition}) simply greatly in the fiberwise
Kahler toric case. We begin with the following lemma.

\begin{lemma}
If $(M,g,J)$ is a fiberwise Kahler toric manifold, then
\begin{equation*}
r(\partial_{\phi_{i}},\partial_{u_{j}})=0
\end{equation*}
for all $i,j$.  This implies that the second equation of the
Einstein condition in theorem (\ref{cohomdEinsteincondition}) is
automatically satisfied.
\end{lemma}
\begin{proof}
Recall the second equation of the Einstein condition given in
theorem (\ref{cohomdEinsteincondition}).  The form of a fiberwise
Kahler toric metric implies that this equation becomes
\begin{equation*}
h^{ij}\check{\nabla}\langle\partial_{\phi_{i}},[\partial_{\phi_{j}},X]\rangle=0
\end{equation*}
for all $X\in\mathfrak{p}$.  However, the vector fields
$\partial_{\phi{i}}$ commute with all other vector fields in
$\mathfrak{p}$.  Therefore, the left-hand side of this equation is
automatically equal to zero.
\end{proof}
Next we recall that $\partial_{u_{i}}=-J\partial_{\phi_{i}}$.
Because the metric is Kahler it is invariant under the action of
the complex structure $J$.  Furthermore, on a Kahler manifold, the
Ricci tensor is also invariant under the action of the complex
structure.  That is,
$r(\partial_{u_{i}},\partial_{u_{j}})=r(J\partial_{\phi_{i}},J\partial_{\phi_{j}})=r(
\partial_{\phi_{i}},\partial_{\phi_{j}})$.
The invariance of the metric and the Ricci tensor under the action
of the complex structure proves the following lemma.

\begin{lemma}
On a fiberwise Kahler toric manifold, if the first equation of the
Einstein condition in (\ref{cohomdEinsteincondition}) is satisfied
then the third equation is satisfied automatically.
\end{lemma}

    We further note that $r(\partial_{\phi_{i}},Y_{i}^{j})=0$ for all $i$ and $j$ and
$r(Y_{k}^{l},Y_{i}^{j})=0$ for all $k\neq i$ and $l\neq j$.
Therefore, the Einstein condition on a fiberwise Kahler toric
manifold contains two types of equations:
$r(\partial_{\phi_{i}},\partial_{\phi_{j}})=\lambda
g(\partial_{\phi_{i}},\partial_{\phi_{j}})$ and
$r(Y^{i}_{k},Y^{i}_{k})=\lambda g(Y^{i}_{k},Y^{i}_{k})$.  To
calculate these equations we will need to determine the quantities
$\hat{r}(\partial_{\phi_{i}},\partial_{\phi_{j}})$ and
$\hat{r}(Y^{i}_{k},Y^{i}_{k})$.

    Noting that each principal orbit is a principal torus bundle
over a coadjoint orbit and using Wang and Ziller's work in
\cite{WangZiller} we calculate that

\begin{equation*}
\hat{r}(Y^{i}_{k},Y^{i}_{k})=r^{\ast}(Y^{i}_{k},Y^{i}_{k})-\frac{1}{2}\frac{h_{kl}b_{ki}
b_{li}}{A_{i}}
\end{equation*}
where we have summed over the repeated indices $k$ and $l$. Dancer
and Wang proved that $r^{\ast}(Y^{i}_{k},Y^{i}_{k})=1$ in the
cohomogeneity-one case.  However, their proof works without
alteration in the cohomogeneity-$d$ case.  Therefore,
\begin{equation*}
\hat{r}(Y^{i}_{k},Y^{i}_{k})=1-\frac{1}{2}\frac{h_{kl}b_{ki}b_{li}}{A_{i}}.
\end{equation*}
    Similarly, we calculate that
\begin{equation*}
\hat{r}(\partial_{\phi_{i}},\partial_{\phi_{j}})=\frac{1}{4}h_{ik}h_{jl}\sum_{r=1}^{m}
\frac{d_{r}b_{kr}b_{lr}}{A_{r}^{2}}
\end{equation*}
Using the first equation in (\ref{cohomdEinsteincondition}) we
calculate that
\begin{equation*}
r(Y^{i}_{k},Y^{i}_{k})=-\frac{1}{2}\check{\triangle}A_{i}-\frac{1}{4}\frac{\langle
\check{\nabla}\textrm{det}(h),\check{\nabla}A_{i}\rangle}{\textrm{det}(h)}-\frac{1}{4}\frac{
\langle\check{\nabla}V,\check{\nabla}A_{i}\rangle}{V}
-\frac{1}{2}A_{i}\langle\check{\nabla}(\frac{1}{A_{i}}),\check{\nabla}A_{i}\rangle+1
-\frac{1}{2}\frac{h_{kl}b_{ki}b_{li}}{A_{i}}
\end{equation*}
where $V=\prod_{r=1}^{m}A_{r}^{d_{r}}$.  Noting that, just as in
the Kahler toric case, \linebreak
$\check{\triangle}(\cdot)=h_{kl}\frac{\partial^{2}(\cdot)}{\partial
x^{k}\partial
x^{l}}+\frac{1}{2}\frac{1}{\textrm{det}(h)}h_{kl}\frac{\partial
\textrm{det}(h)}{\partial x^{k}}\frac{\partial(\cdot)}{\partial
x^{l}}$. Since the $A_{j}$ are linear in the $x_{i}$, we see that
$\frac{\partial^{2}A_{j}}{\partial x^{k}\partial x^{l}}=0$ for all
$j,k,l$.  We calculate that
\begin{equation*}
-\frac{1}{2}A_{i}\langle\check{\nabla}(\frac{1}{A_{i}}),\check{\nabla}A_{i}\rangle=\frac{1}{2}
\frac{1}{A_{i}}h_{kl}\frac{\partial
A_{i}}{\partial x^{k}}\frac{\partial A_{i}}{\partial
x^{l}}=\frac{1}{2}\frac{h_{kl}b_{ki}b_{li}}{A_{i}}.
\end{equation*}

Equation $r(Y^{i}_{k},Y^{i}_{k})=\lambda g(Y^{i}_{k},Y^{i}_{k})$
of the Einstein condition becomes
\begin{equation*}
-\frac{1}{2}\frac{\langle\check{\nabla}\textrm{det}(h),\check{\nabla}A_{i}\rangle}{\textrm{det}(h)}
-\frac{1}{4}\frac{\langle
\check{\nabla}V,\check{\nabla}A_{i}\rangle}{V}+1=\lambda A_{i}
\end{equation*}
for all $1\leq i\leq m$ where
$A_{i}=\sum_{r=1}^{d}b_{ri}x_{r}+a_{i}$.  Expanding this equation
gives
\begin{equation*}
\frac{1}{\textrm{det}(h)}h_{kl}\frac{\partial
\textrm{det}(h)}{\partial x^{l}}\frac{\partial A_{i}}{\partial
x^{k}}+\frac{1}{2}h_{kl}\frac{\partial\textrm{log}V}{\partial
x^{l}}\frac{\partial A_{i}}{\partial x^{k}}-2=-2\lambda
\sum_{k=1}^{d}b_{ki}x_{k}+a_{i}.
\end{equation*}
This becomes
\begin{equation*}
(h_{kl}\frac{\partial \textrm{log}(\textrm{det}(h)
V^{\frac{1}{2}})}{\partial x^{l}}+2\lambda
x_{k})b_{ki}=2(1-\lambda a_{i})
\end{equation*}
for all $1\leq i\leq m$.

    Next, we need to express equation $r(\partial_{\phi_{i}},\partial_{\phi_{j}})=\lambda
g(\partial_{\phi_{i}},\partial_{\phi_{j}})$ in terms of the
$h_{kl}$ and the $A_{i}$.  The first equation in
(\ref{cohomdEinsteincondition}) becomes
\begin{equation*}
-\frac{1}{2}\check{\triangle}h_{ij}-\frac{1}{4}\frac{\langle\check{\nabla}\textrm{det}(h),
\check{\nabla}h_{ij}\rangle}{\textrm{det}(h)}-\frac{1}{4}\frac{\langle\check{\nabla}V,\check{
\nabla}h_{ij}\rangle}{V}-\frac{1}{2}h_{ik}\langle\check{\nabla}h^{kl},\check{\nabla}h_{lj}
\rangle+\frac{1}{4}h_{ik}h_{jl}\sum_{r=1}^{m}\frac{d_{r}b_{kr}b_{lr}}{A_{r}^{2}}=\lambda
h_{ij}.
\end{equation*}
Using the calculations for the Einstein condition of a Kahler
toric manifold performed in the previous section, we see that the
above equation becomes
\begin{equation*}
-\frac{1}{2}h_{ik}\frac{\partial}{\partial
x^{k}}(\frac{1}{\textrm{det}(h)}h_{lj}\frac{\partial
\textrm{det}(h)}{\partial
x^{l}})-\frac{1}{4}\frac{\langle\check{\nabla}V,\check{\nabla}h_{ij}\rangle}{V}+\frac{1}
{4}h_{ik}h_{jl}\sum_{r=1}^{m}\frac{d_{r}b_{kr}b_{lr}}{A_{r}^{2}}=\lambda
h_{ij}
\end{equation*}
for all $1\leq i,j\leq d$.

    Expanding this equation we see that
\begin{equation*}
h_{ik}\frac{\partial}{\partial
x^{k}}(\frac{1}{\textrm{det}(h)}h_{lj}\frac{\partial
\textrm{det}(h)}{\partial x^{l}})+h_{ik}\frac{\partial
h_{lj}}{\partial x^{k}}\frac{\partial
\textrm{log}(V^{\frac{1}{2}})}{\partial
x^{l}}+h_{ik}h_{jl}\frac{\partial^{2}\textrm{log}(V^{\frac{1}{2}})}{\partial
x^{k}\partial x^{l}}=-2\lambda h_{ij}
\end{equation*}
which is equivalent to
\begin{equation*}
h_{ik}\frac{\partial}{\partial x^{k}}(h_{lj}\frac{\partial
\textrm{log}(\textrm{det}(h)V^{\frac{1}{2}})}{\partial
x^{l}})=-2\lambda h_{ij}.
\end{equation*}

    We collect the above calculations in
the following theorem
\begin{theorem}
The Einstein condition on a fiberwise Kahler toric manifold is
given by
\begin{equation}\label{FKTEinsteinA}
h_{ik}\frac{\partial}{\partial x^{k}}\left(h_{lj}\frac{\partial
{\rm log}({\rm det}(h)V^{\frac{1}{2}})}{\partial
x^{l}}\right)=-2\lambda h_{ij}
\end{equation}
for all $1\leq i,j\leq d$ and
\begin{equation}\label{FKTEinsteinB}
\left(h_{kl}\frac{\partial {\rm log}({\rm det}(h)
V^{\frac{1}{2}})}{\partial x^{l}}+2\lambda
x_{k}\right)b_{ki}=2(1-\lambda a_{i})
\end{equation}
for all $1\leq i\leq m$, where $\lambda$ is the Einstein constant,
$h^{kl}=\frac{\partial^{2}\Phi}{\partial x^{k}\partial x^{l}}$,
$A_{j}=\sum_{r=1}^{d}b_{ri}x_{r}+a_{j}$ for all $1\leq j\leq m$,
and $V=\prod_{r=1}^{m}A_{r}^{d_{r}}$.
\end{theorem}
\subsection{The Scalar Curvature Equation}

    Having calculated the equations of the Einstein condition, it
is a simple matter to calculate the scalar curvature equation (the
scalar curvature being the trace of the Ricci tensor).

\begin{theorem}
The scalar curvature, $S$, of a fiberwise Kahler toric manifold
can be written as
\begin{equation}\label{FKTScalar}
S=-\frac{1}{{\rm det}(h)}h_{kl}\frac{\partial^{2}{\rm
det}(h)}{\partial x^{k}\partial
x^{l}}-\sum_{r=1}^{m}d_{r}\frac{h_{kl}\frac{\partial {\rm
det}(h)}{\partial x^{k}}b_{lr}}{{\rm
det}(h)A_{r}}+\frac{1}{2}\sum_{r=1}^{m}d_{r}\frac{h_{kl}b_{kr}b_{lr}}{A_{r}^{2}}
\end{equation}
\begin{equation*}
-\frac{1}{4}\sum_{r=1}^{m}\sum_{s=1}^{m}d_{r}d_{s}\frac{h_{kl}b_{kr}b_{ls}}{A_{r}A_{s}}+
\sum_{r=1}^{m}\frac{d_{r}}{A_{r}}
\end{equation*}
\end{theorem}
\begin{proof}
From the equations of the Einstein condition, we calculated that
\begin{equation*}
S=\bar{S}-\sum_{r=1}^{m}d_{r}\frac{\langle\check{\nabla}{\rm
det}(h),\check{\nabla}A_{r}\rangle}{{\rm
det}(h)A_{r}}+\frac{1}{2}\sum_{r=1}^{m}d_{r}\frac{\langle
\check{\nabla}A_{r},\check{\nabla}A_{r}\rangle}{A_{r}^{2}}
-\frac{1}{4}\sum_{r=1}^{m}\sum_{s=1}^{m}d_{r}d_{s}\frac{\langle
\check{\nabla}A_{r},\check{\nabla}A_{s}\rangle}{A_{r}A_{s}}+\sum_{r=1}^{m}\frac{d_{r}}{A_{r}}
\end{equation*}
where $\bar{S}$ is the scalar curvature of the associated Kahler
toric metric $\bar{g}$.  The theorem follows immediately from the
form of the $A_{r}$ and the scalar curvature of a Kahler toric
metric given in the previous section.
\end{proof}
\subsection{The Extremal Kahler Condition}
    As discussed in the previous section, a Kahler metric is
extremal if and only if the gradient of the scalar curvature,
$\nabla S$ is a holomorphic vector field.  Expressed in the
complex coordinates
$\partial_{z_{i}}=\partial_{u_{i}}+i\partial_{\phi_{i}}$, the
gradient of the scalar curvature of a fiberwise Kahler toric
metric is
\begin{equation*}
\nabla S=\left(h^{kl}\frac{\partial S}{\partial
u^{l}}\right)\frac{\partial}{\partial z^{k}}.
\end{equation*}
Because the coefficients, $h^{kl}\frac{\partial S}{\partial
u^{l}}$, are real, the vector field $\nabla S$ is holomorphic if
and only if these coefficients are constant.  The extremal Kahler
condition is equivalent to
\begin{equation*}
h^{kl}\frac{\partial S}{\partial u^{k}}=\frac{\partial S}{\partial
x^{k}}=\alpha_{k}
\end{equation*}
where the $\alpha_{k}$ are constants.  This gives the following
theorem.
\begin{theorem}
A fiberwise Kahler toric metric is extremal Kahler if and only if
its scalar curvature satisfies
\begin{equation*}
S=\sum_{k=1}^{d}\alpha_{k}x_{k}+\beta
\end{equation*}
where $\beta$ and the $\alpha_{k}$ are constants.
\end{theorem}
    Note the striking parallel between the extremal Kahler
condition on a fiberwise Kahler toric metric and the extremal
Kahler condition on a Kahler toric metric given in the previous
section. This provides further support for the idea that fiberwise
Kahler toric metrics can be viewed as Kahler toric metrics with
additional structure.

\chapter{Cohomogeneity-One Fiberwise Kahler Toric Manifolds}

    The simplest class of fiberwise Kahler toric manifolds are those
of cohomogeneity-one, i.e. $d=1$.  If $(M,g,J)$ is a
cohomogeneity-one manifold under the action of a compact connected
semisimple Lie group $G$, then $\textrm{dim}(L/K)=1$.  Therefore,
$L/K\cong S^{1}$ and assumption one is automatically satisfied.
Assumption two, however, is not automatically satisfied, although
it is satisfied for a generic $S^{1}$-bundle over a coadjoint
orbit $G/L$.  The explanation for this can be found in
\cite{DandW}.  Assumption two holds by definition in the fiberwise
Kahler toric case and we will from now on assume that the summands
in the isotropy representation are distinct.

    From the calculations preformed above, the metric on a
cohomogeneity-one fiberwise Kahler toric metric is of the form
\begin{equation*}
g=\begin{pmatrix} \frac{1}{h} & & & & \\
 & h & & & \\
 & & A_{1}Id_{d_{1}} & & \\
 & & & ... & \\
 & & & & A_{m}Id_{d_{m}} \\
\end{pmatrix}
\end{equation*}
relative to basis
$(\partial_{x},\partial_{\phi},Y^{1}_{1},...,Y^{m}_{d_{m}})$ where
we have set $h_{11}=h$, $x=x_{1}$, and
$\partial_{\phi}=\partial_{\phi_{1}}$.  Recall that $h=h(x)$ and
$A_{i}=A_{i}(x)$ for all $1\leq i\leq m$.

    Furthermore, from the Kahler condition, we know that
$A_{i}=b_{i}x+a_{i}$ for all $i$ where $b_{i}$ and $a_{i}$ are
constants.  The metric can be written
\begin{equation*}
g=\begin{pmatrix} \frac{1}{h} & & & & \\
 & h & & & \\
 & & (b_{1}x+a_{1})Id_{d_{1}} & & \\
 & & & ... & \\
 & & & & (b_{m}x+a_{m})Id_{d_{m}} \\
\end{pmatrix}.
\end{equation*}
The metric on a manifold of this type is determined by one
function, $h$, in one variable, $x$, and $2m$ constants.  Note
that $m$ of the $2m$ constants, the $b_{i}$ determine the
principal orbit types.

    We now want to examine the constant scalar curvature, extremal
Kahler, and Einstein conditions on a cohomogeneity-one fiberwise
Kahler toric manifold  After deriving explicit local solutions for
each of these conditions, we will look at the global structure of
such metrics.  The Einstein case was already solved by Dancer and
Wang in \cite{DandW}.  However, they did not, nor to the author's
knowledge has anyone else, derive explicit solutions in the
constant scalar curvature or the extremal Kahler case. Before,
proceeding to new results, we first review the work of Dancer and
Wang.
\section{Explicit Integration of the Einstein Equations}
Let $\lambda$ be the Einstein constant.  Following the
calculations of the previous section, there are two types of
equations in the Einstein condition.  The first such equation
arise from the condition
$r(\partial_{\phi},\partial_{\phi})=\lambda
g(\partial_{\phi},\partial_{\phi}))$.  Because the Ricci tensor,
like the metric, is invariant under the action of complex
structure, this equation contains the same information as the
equation $r(\partial_{u},\partial_{u})=\lambda
g(\partial_{u},\partial_{u}))$ and therefore
$r(\partial_{x},\partial_{x})=\lambda
g(\partial_{x},\partial_{x}))$.  The second type of equation comes
from $r(Y^{i}_{k},Y^{i}_{k})=\lambda g(Y^{i}_{k},Y^{i}_{k})$.  The
Einstein condition is equivalent to the following system of $m+1$
equations
\begin{equation}\label{FKTEinstein1}
-\frac{1}{2}h''-\frac{1}{4}h'\left(\sum_{j=1}^{m}\frac{d_{j}b_{j}}{b_{j}x+a_{j}}\right)
+\frac{1}{4}h\left(\sum_{j=1}^{m}\frac{d_{j}b^{2}_{j}}{(b_{j}x+a_{j})^{2}}\right)=\lambda
\end{equation}
and
\begin{equation}\label{FKTEinstein2}
-\frac{1}{2}h'\frac{b_{i}}{b_{i}x+a_{i}}-\frac{1}{4}\frac{b_{i}}{b_{i}x+a_{i}}h\left(
\sum_{j=1}^{m}\frac{d_{j}b_{j}}{b_{j}x+a_{j}}\right)+\frac{1}{b_{i}x+a_{i}}=\lambda
\end{equation}
for all $1\leq i\leq m$, where $h'=\frac{dh}{dx}$.  Equation
(\ref{FKTEinstein2}) is equivalent to $1=\lambda a_{i}$ if
$b_{i}=0$ and
\begin{equation}\label{FKTEinstein3}
-\frac{1}{2}h'-\frac{1}{4}h\left(\sum_{j=1}^{m}\frac{d_{j}b_{j}}{b_{j}x+a_{j}}\right)=\lambda
x+\frac{a_{i}\lambda-1}{b_{i}}
\end{equation}
when $b\neq0$.  We noted above that not all of the $b_{i}$ are
equal to zero.

    We see immediately that differentiating equation
(\ref{FKTEinstein3}) gives equation (\ref{FKTEinstein1}).
Therefore, equation (\ref{FKTEinstein1}) of the Einstein condition
is automatically satisfied if equation (\ref{FKTEinstein2}) is. We
now have the following theorem.

\begin{theorem}\cite{DandW}
A cohomogeneity-one fiberwise Kahler toric metric is Einstein if
and only if
\begin{itemize}
\item $1=\lambda a_{i}$, if $b_{i}=0$
\item $\frac{1-\lambda a_{i}}{b_{i}}=\frac{1-\lambda a_{k}}{b_{k}}=D$, if
$b_{i}b_{k}\neq0$
\item
$h'+\frac{1}{2}h(\sum_{j=1}^{m}\frac{d_{j}b_{j}}{b_{j}x+a_{j}})+2\lambda
x-2D=0$.
\end{itemize}
\end{theorem}
\begin{proof}
The first equation follows immediately from equation
(\ref{FKTEinstein3}) when $b_{i}=0$.  The second equation follows
from equation (\ref{FKTEinstein3}) when $i=k$.  The third equation
follows from the second.
\end{proof}

Dancer and Wang were able to integrate the Einstein equation to
obtain an explicit solution for $h$
\begin{theorem}\cite{DandW}
If $(M,g,J)$ is a cohomogeneity-one fiberwise Kahler toric metric,
then
\begin{equation}
h=\frac{2}{\prod_{j=1}^{m}(b_{j}x+a_{j})^{\frac{d_{j}}{2}}}\int(D-\lambda
x)\prod_{j=1}^{m}(b_{j}x+a_{j})^{\frac{d_{j}}{2}}dx.
\end{equation}
\end{theorem}

Note that the integrand of the expression for $h$ is always a
polynomial as each $d_{j}$ is even and greater than or equal to
two.  Therefore, one can always, in principle, perform this
integration and express $h$ as a rational function in $x$.
\section{Explicit Integration of the Scalar Curvature Equation}

From the equations of the Einstein condition, it is
straightforward to calculate the scalar curvature of $(M,g,J)$
\begin{theorem}
If $(M,g,J)$ is a cohomogeneity-one fiberwise Kahler toric
manifold, the scalar curvature, $S=S(x)$, is given by
\begin{equation}\label{FKTscalODE}
S=-h''-h'\left(\sum_{j=1}^{m}\frac{d_{j}b_{j}}{b_{j}x+a_{j}}\right)+h\left(\frac{1}{2}
\sum_{j=1}^{m}\frac{d_{j}b_{j}^{2}}{(b_{j}x+a_{j})^{2}}-\frac{1}{4}\left(\sum_{j=1}^{m}
\frac{d_{j}b_{j}}{b_{j}x+a_{j}}\right)^{2}\right)+\sum_{j=1}^{m}\frac{d_{j}}{b_{j}x+a_{j}}.
\end{equation}
\end{theorem}
Note that $-h''$ is the scalar curvature $\bar{S}$ of the embedded
two-dimensional Kahler toric metric.

    We have been able to integrate this equation for \textit{any}
scalar curvature $S$.  Moreover, we can obtain \textit{explicit}
local solutions in the constant scalar curvature and extremal
Kahler cases.
\begin{theorem}
If $(M,g)$ is a cohomogeneity-one fiberwise Kahler toric manifold
the function $h$ satisfies
\begin{equation}\label{FKThsolution}
h=\frac{1}{\prod_{j=1}^{m}(b_{j}x+a_{j})^{\frac{d_{j}}{2}}}\int\left[\int\left(
\sum_{j=1}^{m}\frac{d_{j}}{b_{j}x+a_{j}}-S\right)\prod_{j=1}^{m}(b_{j}x+a_{j})^{
\frac{d_{j}}{2}}dx\right]dx
\end{equation}
where $S$ is the scalar curvature of the manifold.
\end{theorem}
\begin{proof}
The explicit solution in the Einstein case suggests that $h$ is of
the form
\begin{equation*}
h=\frac{1}{\prod_{j=1}^{m}(b_{j}x+a_{j})^{\frac{d_{j}}{2}}}\int\left(\prod_{j=1}^{m}
(b_{j}x+a_{j})^{\frac{d_{j}}{2}}\right)\Theta
dx
\end{equation*}
for some function $\Theta=\Theta(x)$.  With this assumption the
scalar curvature equation (\ref{FKTscalODE}) becomes a first-order
equation in $\Theta$.
\begin{equation*}
S=-\Theta'-\frac{1}{2}\Theta\left(\sum_{j=1}^{m}d_{j}\frac{b_{j}}{b_{j}x+a_{j}}\right)+
\sum_{j=1}^{m}\frac{d_{j}}{b_{j}x+a_{j}}.
\end{equation*}

    This equation is a simple first-order linear ordinary differential equation which we
can integrate explicitly.  We see that
\begin{equation*}
\Theta=\frac{1}{\prod_{j=1}^{m}(b_{j}x+a_{j})^{\frac{d_{j}}{2}}}\int\left(\sum_{j=1}^{m}
\frac{d_{j}}{b_{j}x+a_{j}}-S\right)\prod_{j=1}^{m}(b_{j}x+a_{j})^{\frac{d_{j}}{2}}dx.
\end{equation*}
Substituting this expression for $\Theta$ into our expression for
$h$ in terms of theta above completes the proof.
\end{proof}
Note that if we are looking for metrics of constant scalar
curvature $S=\beta$, the we can perform the double integration in
(\ref{FKThsolution}) explicitly and thereby obtain an explicit
description of \textit{any} cohomogeneity-one fiberwise Kahler
toric metric of constant scalar curvature in terms of
\textit{rational} functions in $x$.

    Furthermore, if we are looking instead for extremal Kahler
metrics then, as proved above, $S=\alpha x+\beta$ where $\alpha$
and $\beta$ are constants.  Once again we can perform explicitly
the double integration in (\ref{FKThsolution}), and express any
extremal Kahler cohomogeneity-one fiberwise Kahler toric metric in
terms of rational functions in $x$.

\section{Special Orbits and Completeness}
    A general cohomogeneity-one metric must have zero, one, or two
special orbits.  If the manifold is compact, then any extremal
Kahler metric will have two special orbits.   This follows from
the fact that the quotient space of any compact cohomogeneity-one
manifold is either a closed interval or a copy of $S^{1}$
\cite{Bergery}.  If the quotient space is $S^{1}$ then there are
no special orbits and the metric must be periodic.  However, this
implies that all of the functions $A_{i}=b_{i}+a_{i}$ are periodic
which can only occur if all of the $b_{i}$ vanish.  In this case
the metric is a product metric.  We have assumed above that not
all of the $b_{i}$ are zero so the metric must have two special
orbits.  (Below we will present some examples in which all of the
$b_{i}$ are zero but in general we assume that at least one is
non-zero.)  We will see below that a complete cohomogeneity-one
fiberwise Kahler toric metric must have at least one.  After
determining the possible special orbit types, we will determine
the smoothness conditions which $h$ and the constants $a_{i}$ and
$b_{i}$ must satisfy for the metric to be extended smoothly over
the principal orbits.

    Let $G/H$ be a special orbit with $K\subset H$.  Via
translation, we can assume that the special orbit is located at
$x=0$.

\begin{lemma}\cite{DandW}
The moment map $\mu$ gives us a fiber bundle over
\begin{equation*}
L^{\ast}/H\rightarrow G/H\rightarrow G/L^{\ast}
\end{equation*}
where $G/L^{\ast}$ is a coadjoint orbit.  Furthermore, if $V$
denotes the normal slice to the special orbit $G/H$ then either
$H=L^{\ast}$ or the complex structure, $J$, induces an
$H$-equivariant isomorphism $ker{d\mu}\cong V$.

\end{lemma}
\begin{proof}
We know that $V$ is an irreducible $H$-module and that if $v\in
\textrm{ker}(d\mu)$ then $Jv\in V$.  These imply that $J$ acts as
an $H$-equivariant monomorphism from $\textrm{ker}(d\mu)$ to $V$.
If $\textrm{ker}(d\mu)$ is trivial at the special orbit then
$H=L^{\ast}$.  If not, then the irreducibility of $V$ requires
that the map is an isomorphism.  Finally, we note that
$\pi(G/L^{\ast})=\{1\}$ and therefore $L^{\ast}/H$ is not a
discrete set of more than one point.
\end{proof}

    Dancer and Wang were able to demonstrate that, at a special
orbit in the cohomogeneity-one fiberwise Kahler toric case,
$H=L^{\ast}$.  As we will attempt to extend this argument to the
cohomogeneity-$d$ case below, we will review their proof.

    We begin by noting the following important facts: $L\subset
L^{\ast}$ by semicontinuity and $H/K\cong S^{k}$ for some
non-negative odd integer $k$ for smoothness.

    At the special orbits, we have the following Lie algebra
decomposition of $K$-modules
\begin{equation*}
\mathfrak{g}=\mathfrak{k}\oplus\mathfrak{m}_{1}\oplus\mathfrak{m}_{2}\oplus\mathfrak{m}_{3}
\end{equation*}
where $\mathfrak{h}=\mathfrak{k}\oplus\mathfrak{m}_{1}$ and
$\mathfrak{m}_{3}=\textrm{ker}(d\mu)$.  If $H=L^{\ast}$ then
$\textrm{ker}(d\mu)=0$. If not, the preceding lemma proves that
$J$ maps $\mathfrak{m}_{3}$ isomorphically onto
$1\oplus\mathfrak{m}_{1}\approx V$.  This would mean that
$\mathfrak{m}_{3}\cong 1\oplus\mathfrak{m}_{1}$ as $K$-modules.
However, this contradicts Assumption Two as $\mathfrak{m}_{1}$
would then appear twice in the isotropy representation.  We
conclude that $\mathfrak{m}_{1}=0$.  This is equivalent to saying
the collapsing sphere $H/K$ has dimension zero.  For the metric to
be extended smoothly over this special orbit we require that
$h(0)>0$ while $A_{i}>0$ and $A'_{i}=0$.  This contradicts the
equations of the Kahler condition.  Since $H=L^{\ast}$ and
$\textrm{ker}(d\mu)=0$, the $S^{1}$ fiber must collapse at any
special orbit and $h(0)=0$.

    The above discussion can be summarized as follows.  Each
special orbit occurs at a boundary point of the one-dimensional
quotient space.  At that point the $S^{1}$ fiber of the principal
orbit \textit{must} collapse.  A portion of the coadjoint orbit
$G/L$ may or may not collapse along with it.  That is, $H/K\cong
S^{1}$ or $H/K\cong S^{k}$ for $k>1$.  This sphere admits the
fibration
\begin{equation*}
S^{1}\cong L/K\rightarrow S^{k}\rightarrow H/L.
\end{equation*}
It remains to determine the possible $H/L$.  It is straightforward
to see that $H/L\cong\mathbb{C}P^{\frac{k-1}{2}}$.

    Next, we look at the necessary conditions on the functions $h$
$A_{i}$ for the metric to be extended smoothly over a special
orbit.  Because of the invariance of the metric under translation
in $x$, we can fix the location of any one special orbit to be
$x=0$ where the principal orbits occur at $x>0$.  We can write
$\mathfrak{m}_{1}$ as
\begin{equation*}
\mathfrak{m}_{1}=\mathfrak{p}_{1}\oplus...\oplus\mathfrak{p}_{l}.
\end{equation*}
It is, of course, possible that $l=0$; in fact, by the
classification of coadjoint orbits, $l\leq 2$ \cite{DandW}.

    For the metric to be smooth at $x=0$, we must have
\begin{equation*}
h(0)=0,~~~\frac{dh}{dx}(0)=2.
\end{equation*}

Before determining what the smoothness conditions on the $A_{i}$
are, we note that the smoothness conditions on $h$ at the special
orbit are identical to the smoothness conditions at a special
orbit of the toric Kahler manifold (see \cite{Abreu}) for
smoothness conditions on special orbits of Kahler toric
manifolds).  The implicit manifold $(N,\bar{g})$ mentioned above
is in fact a global Kahler toric manifold of cohomogeneity-one.
That is, every fiberwise Kahler toric manifold $(M,g,J)$ can be
viewed as a toric Kahler manifold $(N,\bar{g},\bar{J})$ along with
a collection of $3m$ constants (the $a_{i}$, $b_{i}$, and $d_{i}$)
along with conditions on these constants to insure that the metric
is positive definite and smooth at the special orbits.

    We return now to the conditions on the functions $A_{i}$ at the special orbit $x=0$ necessary for smoothness.  Clearly,
\begin{equation*}
A_{i}(0)=0
\end{equation*}
for $1\leq i\leq l$.  Since $A_{j}=b_{j}x+a_{j}$, we see that
$a_{i}=0$ for $1\leq i\leq l$.  Since the principal orbits are
located at positive values of $x$, $b_{i}>0$ for $1\leq i\leq l$.

    There is one more condition on the $b_{i}$ for $1\leq i\leq l$ necessary for the metric to
be extended smoothly over the special orbit at $x=0$.  To obtain
this condition, we evaluate the scalar curvature equation
(\ref{FKTscalODE}) at $x=0$.  The scalar curvature is constant at
$x=0$ if the metric is smooth there so the left-hand side of this
equation is smooth.  However, the right-hand side appears to
contain poles of order one at $x=0$.  By inspection, the
right-hand side will be smooth at $x=0$ if and only if
\begin{equation*}
\sum_{j=1}^{l}\frac{d_{j}}{b_{j}}=\frac{k^{2}-1}{2}
\end{equation*}
where $k-1=\sum_{j=1}^{l}d_{j}$.

    Let us consider in more detail the case in which each principal orbit is an $S^{1}$-bundle
over a product of Kahler-Einstein manifolds $M_{1}\times...\times
M_{m}$ as in construction two.  Again, two things can happen at
the special orbit $x=0$.  The $S^{1}$-fiber can collapse alone or
the $S^{1}$-fiber can collapse along with, say,
$M_{1}\cong\mathbb{C}P^{\frac{d_{1}}{2}}$.  In the second case,
$l=1$ and $k-1=d_{1}$ and we deduce that
\begin{equation}\label{smooth1}
b_{1}=\frac{2}{d_{1}+2}
\end{equation}
for the metric to be extended smoothly over the special orbit.  We
take this opportunity to be more precise about the meaning of the
$b_{i}$ from construction two.  Take, say, $M_{i}$ to be a
Kahler-Einstein manifold with positive Einstein constant.  The
first Chern class of $M_{i}$ can be written as
\begin{equation*}
c_{1}(M_{i})=p_{i}\alpha_{i}
\end{equation*}
where $p_{i}$ is a positive integer and $\alpha_{i}$ is an
indivisible integral cohomology class.  Furthermore, if $d_{i}$ is
the real dimension of $M_{i}$, then $p\leq\frac{d_{i}+2}{2}$ with
equality if and only if $M_{i}\cong\mathbb{C}P^{\frac{d_{i}}{2}}$
\cite{Ko-No}.  If $\omega_{i}$ is the Kahler form of the
Kahler-Einstein metric on $M_{i}$ then we have that
$[\omega_{i}]=2\pi p_{i}\alpha_{i}$.  Finally, the curvature form,
$\Omega_{i}$, of the connection $\theta_{i}$ of the $S^{1}$-bundle
over $M_{i}$ can be written as $\Omega_{i}=-2\pi q_{i}\alpha_{i}$
where $q_{i}$ is some integer.  From the definition of $b_{i}$
given above we deduce that
\begin{equation*}
b_{i}=-\frac{q_{i}}{p_{i}}.
\end{equation*}
In the case being considered
$M_{1}\cong\mathbb{C}P^{\frac{d_{1}}{2}}$ so
$p_{1}=\frac{d_{1}+2}{2}$.  Therefore, equation (\ref{smooth1})
implies that $q=-1$.

    Having determined the smoothness conditions for the metric to be extended over a
special orbit, we now address the problem of finding
\textit{complete} constant scalar curvature and extremal Kahler
metrics. Cohomogeneity-one metrics are defined on intervals in
$\mathbb{R}$.  Endpoints can, a priori, occur at poles of $h$,
zeros of $h$, zeros of the $A_{i}$, or at $\pm\infty$.

    Let $p$ be an endpoint of the interval on which the metric is defined (of course $p$
could be $\pm\infty$). Given the form of the metric $g$, the
geodesic distance to an endpoint $p$ can be written as

\begin{equation*}
\int_{\ast}^{p}h^{-\frac{1}{2}}dx.
\end{equation*}

If $p$ is a pole of $h$, then the geodesic distance is finite and
there are no complete metrics.  If $h$ is nonzero and $A_{i}=0$
for some $i$ then again the geodesic distance is finite.  We must
have a special orbit at the endpoint which is impossible because
$h=0$ at all special orbits.  Therefore, if the metric is to be
complete, each endpoint must be either a zero of $h$ (possibly
also a zero of one or more of the $A_{i}$) or $\pm\infty$.
Furthermore, if the manifold is noncompact, the above discussion
proves that one of the endpoints must be $\pm\infty$.

    For $g$ an extremal Kahler metric, the scalar curvature is linear in $x$:
$S=\alpha x+\beta$.  In equation (\ref{FKThsolution}), we express
$h$ as a double integral involving $S$.  If $\alpha=0$, then
$h=O(x^{2})$ or $O(x)$ as $x\rightarrow\pm\infty$ and the geodesic
distance is infinite if $p=\pm\infty$.  A constant scalar
curvature metric will be complete as we approach infinity.  On the
other hand, if $\alpha\neq0$, then $h=O(x^{3})$ as
$x\rightarrow\pm\infty$ and the geodesic distance is
\textit{finite} and the metric is not complete at $\pm\infty$.
This proves the following theorem.

\begin{theorem}
Let $(M,g,J)$ be a non-compact cohomogeneity-one fiberwise Kahler
toric manifold.  There are no complete extremal Kahler metrics on
$M$ which are not of constant scalar curvature
\end{theorem}

\section{Non-Compact Metrics of Constant Scalar Curvature}

    As proved above, any complete extremal Kahler metric on a non-compact cohomogeneity-one
fiberwise Kahler toric manifold must have constant scalar
curvature. Therefore, we can restrict our attention to the search
for constant scalar curvature metrics. From the form of the metric
we deduce that $g$ must have at least one special orbit.  This can
be seen by noting that each $A_{i}$ vanishes at some value of $x$.
Moreover, if the manifold is non-compact, then there are no other
special orbits. By translating in $x$ if necessary, we take the
location of the special orbit to be $x=0$.  We construct the
metric so that the principal orbits occur at $x>0$ and that the
metric extends to $+\infty$.

    Let $S=\beta$, a constant.  For the metric to be smooth and positive definite, the
following conditions must be satisfied
\begin{itemize}
\item $h(0)=0$ and $h'(0)=2$,
\item $A_{i}=b_{i}x+a_{i}>0$ for $x\in(0,\infty)$ which implies that $b_{i}\geq0$ and
$a_{i}\geq0$ for all $i$, and
\item $h(x)>0$ for $x\in(0,\infty)$.
\end{itemize}
    If $A_{i}(0)=0$ for some $i$ (i.e. if $a_{i}=0$ for some $i$) then an additional smoothness
condition must be satisfied.  Let us first consider the case in
which $A_{i}=0$ for $x\in(0,\infty)$ (that is $a_{i}>0$).  In this
case the manifold is a complex line bundle over a coadjoint orbit
(or a product of Kahler-Einstein manifolds).

    The function $h$ must satisfy equation (\ref{FKThsolution}).  We set $h=\frac{P(x)}{Q(x)}$
where
\begin{equation*}
P(x)=\int\left[\int\left(\sum_{j=1}^{m}\frac{d_{j}}{b_{j}x+a_{j}}-\beta\right)\prod_{j=1}^{m}
(b_{j}x+a_{j})^{\frac{d_{j}}{2}}dx\right]dx
\end{equation*}
and
\begin{equation*}
Q(x)=\prod_{j=1}^{m}(b_{j}x+a_{j})^{\frac{d_{j}}{2}}.
\end{equation*}

    Note that because we have assumed that $b_{j}x+a_{j}>0$ for $x\in[0,\infty)$, we have
that $Q(x)>0$ for $x\in[0,\infty)$.  For the metric to be
positive, $P>0$ for $x\in(0,\infty)$.  Both $P$ and $Q$ are
polynomials in $x$ because each of the $d_{i}$ is even and is
greater than or equal to two. After integration, $P$ can be
written as
\begin{equation*}
P(x)=f+ex+...-\beta\frac{\left(\prod_{j=1}^{m}b_{j}^{\frac{d_{j}}{2}}\right)}{(l+2)(l+1)}
x^{l+2}
\end{equation*}
where $e$ and $f$ are the constants of integration and
$l=\sum_{j=1}^{m}\frac{d_{j}}{2}$.  The smoothness conditions at
$x=0$ determine the constants.  To see this we note that $h(0)=0$
implies that $P(0)=0$ as $Q(0)>0$ by assumption.  Therefore
\begin{equation*}
h(0)=0 \Rightarrow P(0)=0\Rightarrow f=0.
\end{equation*}
Next, we note that $h'(0)=\frac{P'(0)}{Q(0)}$ since $P(0)=0$.  We
see that
\begin{equation*}
h'(0)=2\Rightarrow\frac{P'(0)}{Q(0)}=2\Rightarrow\frac{e}{\prod_{j=1}^{m}a_{j}^{
\frac{d_{j}}{2}}}=2.
\end{equation*}
This implies that
\begin{equation*}
P(x)=2\prod_{j=1}^{m}a_{j}^{\frac{d_{j}}{2}}x+...-\beta\frac{\left(\prod_{j=1}^{m}b_{j}^{
\frac{d_{j}}{2}}\right)}{(l+2)(l+1)}x^{l+2}.
\end{equation*}
    It remains only to check that $h(x)>0$ for $x\in(0,\infty)$.  This is equivalent to
demanding that $P(x)>0$ for $x\in(0,\infty)$.  Now, as
$x\rightarrow\infty$,
\begin{equation*}
P(x)\sim-\beta\frac{\left(\prod_{j=1}^{m}b_{j}^{\frac{d_{j}}{2}}\right)}{(l+2)(l+1)}x^{l+2}.
\end{equation*}
If $\beta>0$, then $P(x)\rightarrow-\infty$ as
$x\rightarrow\infty$ and the metric will fail to be positive for
large values of $x$.  Therefore, $\beta\leq0$.

    If $\beta\leq0$, then we see that in the polynomial $P(x)$, all of the coefficients
of the $x^{k}$ are positive as $b_{i}\geq0$ and $a_{i}>0$ for all
$i$.  Therefore, $P(x)>0$ for $x\in(0,\infty)$.  All of the
smoothness conditions and the positivity condition are satisfied
for all values of $a_{i}>0$.  As discussed above, it is the values
of the $b_{i}$'s which determine the principal orbit type.  So
given a collection of $b_{i}>0$, we see that there exists a scalar
flat metric for all values of $a_{i}$ as well as a metric with
scalar curvature, say, $-1$.  We have the following theorem.

\begin{theorem}
Let $(M,g)$ be a noncompact cohomogeneity-one fiberwise Kahler
toric manifold defined on the interval $x\in[0,\infty)$ such that
$b_{i}\geq0$ for all $i$ (i.e. $M$ is a complex line bundle) as
described above.  Every such manifold admits a scalar-flat metric
($\beta=0$) as well as a negative scalar curvature metric
($\beta=-1$) for every value of $a_{i}>0$.  That is, each such
manifold admits an $m$-parameter family of scalar-flat metrics and
an $m$-parameter family of metrics with scalar curvature equal to
$-1$.
\end{theorem}
On the other hand, if $a_{i}=0$ when $1\leq i\leq l$ for some $l$
then there is the additional smoothness condition that
$\sum_{i=1}^{l}\frac{d_{i}}{b_{i}}=\frac{k^{2}-1}{2}$.

    Consider the case in which each principal orbit is an $S^{1}$-bundle over the product
of Kahler-Einstein manifolds $M_{1}\times...\times M_{m}$.
Consider a special orbit in which $A_{1}(0)=0$ and $A_{i}(0)>0$
for all  $2\leq i\leq m$.  We demonstrated above that
$M_{1}\cong\mathbb{C}P^{\frac{d_{1}}{2}}$, $a_{1}=0$, and
$b_{1}=\frac{2}{d_{1}+2}$. In the constant scalar curvature case
the function $h$ can be written as $h=\frac{P(x)}{Q(x)}$ where
\begin{equation*}
P=\int\left[\int\left(\frac{d_{1}(d_{1}+2)}{2x}+\sum_{j=2}^{m}\frac{d_{j}}{b_{j}x+a_{j}}-
\beta\right)x^{\frac{d_{1}}{2}}\prod_{j=2}^{m}(b_{j}x+a_{j})^{\frac{d_{j}}{2}}dx\right]dx
\end{equation*}
and
\begin{equation*}
Q=x^{\frac{d_{1}}{2}}\prod_{j=2}^{m}(b_{j}x+a_{j})^{\frac{d_{j}}{2}}.
\end{equation*}
Both $P$ and $Q$ are polynomials and we can write $P$
schematically as
\begin{equation*}
P=-\beta\frac{\prod_{j=2}^{m}b_{j}^{\frac{d_{j}}{2}}}{(l+2)(l+1)}x^{l+2}+...+2\left(
\prod_{j=2}^{m}a_{j}^{\frac{d_{j}}{2}}\right)x^{\frac{d_{1}+2}{2}}+ex+f
\end{equation*}
where $l=\sum_{j=1}^{m}\frac{d_{j}}{2}=\frac{n-2}{2}$ and $e$ and
$f$ are the constants of integration.  The polynomial $Q$ can be
written as
\begin{equation*}
Q=\left(\prod_{j=2}^{m}b_{j}^{\frac{d_{j}}{2}}\right)x^{l}+...+\left(\prod_{j=2}^{m}a_{j}^{
\frac{d_{j}}{2}}\right)x^{\frac{d_{1}}{2}}.
\end{equation*}
It is straightforward to check that $h(0)=0$ and $h'(0)=2$ if and
only if $e=f=0$.

    All that remains is to determine when $h(x)>0$ for all $x\in(0,\infty)$.  Recalling
that $b_{i}>0$ and $a_{i}>0$ for all $2\leq i\leq m$, we see that
we obtain a smooth positive definite metric if $\beta\leq0$.  We
see, as above, that $h$ will not be positive for large values of
$x$ if $\beta$ is positive.  Finally, we note that the constant
scalar curvature metrics we have here constructed are on
$\mathbb{C}^{\frac{d_{1}+2}{2}}$-bundles over
$M_{2}\times...\times M_{m}$.

\begin{theorem}
Let $(M,g)$ be a fiberwise Kahler toric manifold in which each
principal orbit is an $S^{1}$-bundle over $M_{1}\times...\times
M_{m}$ such that $M_{1}\cong\mathbb{C}P^{\frac{d_{1}}{2}}$ and
$b_{1}=\frac{2}{d_{1}+2}$.  Let $g$ be defined on the interval
$x\in[x,\infty)$ such that $b_{i}\geq0$ for all $1\leq i\geq m$.
In this case $M$ is a $\mathbb{C}^{\frac{d_{1}+2}{2}}$-bundle over
$M_{2}\times...\times M_{m}$.  Every such manifold admits a
scalar-flat metric and a metric with scalar curvature equation to
$-1$ for all values of $a_{i}>0$ for $2\leq i\leq m$.
\end{theorem}
\section{Compact Extremal Kahler Metrics}

    Having found all the non-positive constant scalar curvature metrics on non-compact
cohomogeneity-one fiberwise Kahler toric metrics, we now turn our
attention to the case of extremal Kahler metrics on compact
manifolds. In the compact case, the metric must have two special
orbits at which $h$ and possibly some $A_{i}$'s vanish.  The
presence of a second special orbit makes the problem of finding
constant scalar curvature Kahler metrics substantially more
difficult in the compact case than it was in the non-compact case.
In fact, there are known examples of compact cohomogeneity-one
fiberwise Kahler toric manifolds which do not admit any constant
scalar curvature metrics (e.g. any Hirzebruch surface other that
$S^{2}\times S^{2}$ with its product metric \cite{Besse}).
However, it is reasonable to suppose that every compact
cohomogeneity-one fiberwise Kahler toric metric could admit an
extremal Kahler metric.

    We begin by looking at the case of $\mathbb{C}P^{1}$-bundles.

\subsection{Extremal Kahler Metrics on $\mathbb{C}P^{1}$-Bundles}

    Let $(M,g,J)$ be a compact fiberwise Kahler toric manifold of cohomogeneity-one.  If
both of the special orbits occur when the $S^{1}$-fiber collapses
and $A_{i}\neq0$ for all $i$, then the total manifold $M$ is a
$\mathbb{C}P^{1}$-bundle over a coadjoint orbit $G/L$ or the
product of Kahler-Einstein manifolds $M_{1}\times...\times M_{m}$
each with positive Einstein constant (though of course the product
metric need not be Einstein).

    By translation, we can fix the location of the first special orbit to be at $x=0$.  By
rescaling if necessary, we can fix the location of the second
special orbit to be located at, say, $x=1$.  We have assumed that
at both of the special orbits $A_{i}\neq0$ for all $i$.  For the
metric to be smooth on the special orbits and positive for
$x\in(0,1)$ the following conditions must be satisfied:
\begin{itemize}
\item $h(0)=0$ and $h'(0)=2$
\item $h(1)=0$ and $h'(1)=-2$
\item $A_{i}=b_{i}x+a_{i}>0$ for $x\in[0,1]$ which implies that $a_{i}>0$ for all $i$ and
$-\frac{a_{i}}{b_{i}}<0$ or $-\frac{a_{i}}{b_{i}}>1$ for all $i$
\item $h(x)>0$ for $x\in(0,1)$.
\end{itemize}

    We ask which manifolds of this type admit extremal Kahler metrics.  In the extremal
Kahler case, $S=\alpha x+\beta$ and
\begin{equation*}
h=\frac{1}{\prod_{j=1}^{m}(b_{j}x+a_{j})^{\frac{d_{j}}{2}}}\int\left[\int\left(\sum_{j=1}^{m}
\frac{d_{j}}{b_{j}x+a_{j}}-\alpha
x-\beta\right)\prod_{j=1}^{m}(b_{j}x+a_{j})^{\frac{d_{j}}{2}}dx\right]dx.
\end{equation*}
We write $h$ as $h(x)=\frac{P(x)}{Q(x)}$ where
$P(x)=\int\left[\int(\sum_{j=1}^{m}\frac{d_{j}}{b_{j}x+a_{j}}-\alpha
x-\beta)\prod_{j=1}^{m}(b_{j}x+a_{j})^{\frac{d_{j}}{2}}dx\right]dx$
and $Q(x)=\prod_{j=1}^{m}(b_{j}x+a_{j})^{\frac{d_{j}}{2}}$. By
construction $Q(x)>0$ for $x\in[0,1]$.  Schematically,
\begin{equation*}
P(x)=f+ex+...-\alpha\frac{\prod_{j=1}^{m}b_{j}^{\frac{d_{j}}{2}}}{(l+3)(l+2)}x^{l+3}
\end{equation*}
where $e$ and $f$ are the constants of integration and
$l=\sum_{j=1}^{m}\frac{d_{j}}{2}$.  The condition $h(0)=0$ holds
if and only if $P(0)=0$.  This implies that $f=0$.  Next, consider
the condition $h'(0)=2$.  Since $P(0)=0$, we have that
$h'(0)=\frac{P'(0)}{Q(0)}$ which implies that
$e=2\prod_{j=1}^{m}a_{j}^{\frac{d_{j}}{2}}$.  The two constants of
integration are then completely determined by the smoothness
conditions at the special orbit located at $x=0$.  Moreover,

\begin{equation*}
P(x)=2\prod_{j=1}^{m}a_{j}^{\frac{d_{j}}{2}}x+...-\alpha\frac{\prod_{j=1}^{m}b_{j}^{\frac{d_{j}}
{2}}}{(l+3)(l+2)}x^{l+3}.
\end{equation*}

There are two remaining smoothness conditions at the second
special orbit.  The condition $h(1)=0$ is equivalent to $P(1)=0$.
This can be satisfied by fixing the value of, say, the free
variable $\beta$.  The condition $h'(1)=-2$ is equivalent to
\begin{equation*}
P'(1)=-2Q(1)=-2\prod_{j=1}^{m}(b_{j}+a_{j})^{\frac{d_{j}}{2}}
\end{equation*}
since $P(1)=0$.  We can always choose $\alpha$ so that this
equation is satisfied.  The resulting function $h$ has zeroes of
order one at the points $x=0$ and $x=1$ and its first derivative
satisfies the smoothness conditions at those points.  However, it
is not immediately obvious that the resulting function $h$ will be
greater than zero on the open interval $(0,1)$.  That is, the
function could have zeroes in the interval $(0,1)$.

    Before proving this sections main theorem, we pause to recall the significance of the
constants $b_{i}$ and $a_{i}$.  The $b_{i}$ determine the
principal orbit type; changing the values of the $b_{i}$ changes
the connection of the $S^{1}$-bundle.  It can therefore be said
that it is the $b_{i}$ that determine the complex manifold.  On
the other hand, changing the $a_{i}$ (so long as one does not
violate smoothness) does not change the manifold.  It does however
change the Kahler class of the Kahler form.

    Therefore, we consider the $b_{i}$ to be fixed, and we take the $a_{i}$ to be free
parameters.  Above, we have shown that for every value of $a_{i}$
we can construct a function $h_{a_{1}...a_{m}}$ which satisfies
the extremal Kahler condition at each point and also satisfies all
of the smoothness conditions.  If $h_{a_{1}...a_{m}}$ is positive
for $x\in(0,1)$, then we have a global extremal Kahler metric.
What remains is to determine at what values of $a_{i}$ the
function $h_{a_{1}...a_{m}}>0$ on that open interval.

    There are two possible questions that we could ask. First, we could ask which manifolds
admit at least one extremal Kahler metric.  Second, we could ask
which Kahler classes admit an extremal Kahler metric.  In the
following theorem, we give a complete answer to the first question
for manifolds of the type being considered.
\begin{theorem}
Let $(M,g,J)$ be a fiberwise Kahler toric manifold of
cohomogeneity-one such that $M$ is a $\mathbb{C}P^{1}$-bundle over
a coadjoint orbit $G/L$ whose isotropy representation has $m$
distinct summands or the product of $m$ Kahler- Einstein manifolds
$M_{k}$.  \textbf{Every} such manifold admits, up to scaling, an
$m$-parameter family of extremal Kahler metrics.
\end{theorem}
\begin{proof}
The $a_{i}$ are free parameters and for each value of the $a_{i}$
we can construct a function $h_{a_{1}...a_{m}}$ with zeros of
order one at $x=0,1$ which satisfies the smoothness conditions.
The function $h_{a_{1}...a_{m}}$ defines an extremal Kahler metric
on $M$ if and only if
\begin{equation*}
h_{a_{1}...a_{m}}>0~~\textrm{for}~~x\in(0,1).
\end{equation*}
In this proof, by looking at the limit of the functions
$h_{a_{1}...a_{m}}$ as $a_{i}\rightarrow\infty$ for \textit{all}
$i$, we will demonstrate that the function $h_{a_{1}...a_{m}}$ is
positive on $x\in(0,1)$ for the $a_{i}$ sufficiently large.

    Let $h_{a_{1}...a_{m}}=\frac{P_{a_{1}...a_{m}}}{Q_{a_{1}...a_{m}}}$ as above.  Grouping
the terms by the order of the $a_{i}$ we write $P_{a_{1}...a_{m}}$
as
\begin{equation*}
P_{a_{1}...a_{m}}=\prod_{j=1}^{m}a_{j}^{\frac{d_{j}}{2}}(2x-\frac{\alpha
x^{3}}{6}-\frac{\beta x^{2}}{2})+H
\end{equation*}
where $H$ consists of lower order terms in the $a_{i}$.  Also,
\begin{equation*}
Q_{a_{1}...a_{m}}=\prod_{j=1}^{m}a_{j}^{\frac{d_{j}}{2}}+G
\end{equation*}
where $G$ consists of lower order terms in the $a_{i}$.

    For $a_{i}$ sufficiently large for all $i$,
\begin{equation*}
P_{a_{1}...a_{m}}\sim\prod_{j=1}^{m}a_{j}^{\frac{d_{j}}{2}}(2x-\frac{\alpha
x^{3}}{6}-\frac{\beta x^{2}}{2})
\end{equation*}
and
\begin{equation*}
Q_{a_{1}...a_{m}}\sim\prod_{j=1}^{m}a_{j}^{\frac{d_{j}}{2}}.
\end{equation*}
    Therefore,
\begin{equation*}
h_{a_{1}...a_{m}}\sim 2x-\frac{\alpha x^{3}}{6}-\frac{\beta
x^{2}}{2}.
\end{equation*}

The conditions $h(1)=0$ and $h'(1)=-2$ imply that for the $a_{i}$
sufficiently large
\begin{equation*}
\alpha\sim0,
\end{equation*}
\begin{equation*}
\beta\sim4,
\end{equation*}
and
\begin{equation*}
h_{a_{1}...a_{m}}\sim2x(1-x)\equiv h_{\infty}.
\end{equation*}
Recall that each of the $h_{a_{1}...a_{m}}$ is a smooth function
on $[0,1]$ with zeros of order one at $x=0,1$.  As the $a_{i}$
tend to infinity for all $i$, these functions approach a smooth
function, $h_{\infty}$, on $[0,1]$ with zeros of order one at
$x=0,1$.  Since $h_{\infty}$ is positive on the interval $(0,1)$,
we deduce that, for the $a_{i}$ sufficiently large for all $i$,
the function $h_{a_{1}...a_{m}}$ must be positive on the interval
$(0,1)$ also. Therefore, for the $a_{i}$ sufficiently large,
$h_{a_{1}...a_{m}}$ defines an extremal Kahler metric on $M$.  We
therefore have the desired $m$ parameter family of extremal Kahler
metrics.

Note that while $\alpha\rightarrow0$ as the $a_{i}$ approach
infinity, $\alpha$ may be non-zero and we cannot deduce that the
metrics are of constant scalar curvature for the $a_{i}$
sufficiently large.  In fact, there are known examples of
manifolds of this type which admit no constant scalar curvature
Kahler metrics.
\end{proof}

This theorem provides a wealth of new examples of extremal Kahler
metrics.  Furthermore, these metrics can be written down
\textit{explicitly} in terms of rational functions.  Previously, a
one-parameter family of extremal Kahler metrics were found on
$\mathbb{C}P^{1}$-bundles over $\mathbb{C}P^{l}$ \cite{Besse}. The
extremal Kahler metrics we have found can be viewed as a
generalization of those metrics.

\subsection{Extremal Kahler Metrics on $\mathbb{C}P^{\frac{d_{1}+2}{2}}$-Bundles}

    In the previous section, we assumed that only the $S^{1}$-fiber collapsed at the special
orbits.  In this section, we set
\begin{equation*}
A_{1}(0)=0
\end{equation*}
while keeping $A_{i}(x)>0$ when $x\in[0,1]$ for all $2\leq i\leq
m$.  This condition implies that $a_{1}=0$ and $b_{1}>0$.

    For simplicity of presentation, we restrict ourselves to the case in which each principal
orbit is an $S^{1}$-bundle over $M_{1}\times...\times M_{m}$.  As
described above, for the metric to be extended smoothly over the
special orbit, we must have
$M_{1}\cong\mathbb{C}P^{\frac{d_{1}}{2}}$ and
$b_{1}=\frac{2}{d_{1}+2}$.  In the extremal Kahler case, the
function $h$ can be written as $h(x)=\frac{P(x)}{Q(x)}$ where
\begin{equation*}
P=\int\left[\int\left(\frac{d_{1}(d_{1}+2)}{2x}+\sum_{j=2}^{m}\frac{d_{j}}{b_{j}x+a_{j}}-\alpha
x-\beta\right)x^{\frac{d_{1}}{2}}\prod_{j=2}^{m}(b_{j}x+a_{j})^{\frac{d_{j}}{2}}dx\right]dx
\end{equation*}
and
\begin{equation*}
Q=x^{\frac{d_{1}}{2}}\prod_{j=2}^{m}(b_{j}x+a_{j})^{\frac{d_{j}}{2}}.
\end{equation*}
Both $P$ and $Q$ are polynomials in $x$; they can be written
schematically as
\begin{equation*}
P=-\alpha\frac{\prod_{j=2}^{m}b_{j}^{\frac{d_{j}}{2}}}{(l+3)(l+2)}x^{l+3}+...+2\left(\prod_{
j=2}^{m}a_{j}^{\frac{d_{j}}{2}}\right)x^{\frac{d_{1}+2}{2}}+ex+f
\end{equation*}
where $l=\sum_{j=1}^{m}\frac{d_{j}}{2}=\frac{n-2}{2}$ and $e$ and
$f$ are the constants of integration.  The polynomial $Q$ can be
written as
\begin{equation*}
Q=\left(\prod_{j=2}^{m}b_{j}^{\frac{d_{j}}{2}}\right)x^{l}+...+\left(\prod_{j=2}^{m}a_{j}^{
\frac{d_{j}}{2}}\right)x^{\frac{d_{1}}{2}}.
\end{equation*}
It is straightforward to check that $h(0)=0$ and $h'(0)=2$ if and
only if $e=f=0$.  There are two remaining smoothness conditions
are $h(1)=0$ and $h'(1)=-2$.  As in the previous section, these
two conditions will determine the values of the constants $\alpha$
and $\beta$.  After selecting $\alpha$ and $\beta$ to satisfy
those conditions, we need only check that the resulting function
$h$ is positive for $x\in(0,1)$.

    We distinguish between two cases: $m=1$ and $m>1$.  If $m=1$, then it is straightforward
to calculate that
\begin{equation*}
h=-\frac{4\alpha}{(d_{1}+6)(d_{1}+4)}x^{3}-\frac{4\beta}{(d_{1}+4)(d_{1}+2)}x^{2}+2x.
\end{equation*}
Solving the equations $h(1)=0$ and $h'(0)=-2$, we see that
$\alpha=0$, $\beta=\frac{1}{2}(d_{1}+4)(d_{1}+2)$, and
\begin{equation*}
h=2x(1-x).
\end{equation*}
This gives the Fubini-Study metric on
$M\cong\mathbb{C}P^{\frac{d_{1}+2}{2}}$.

    The case $m>1$ is the more interesting.  Since the conditions $h(1)=0$ and $h'(1)=-2$
determine the values of $\alpha$ and $\beta$, we have shown that
for every value of $a_{i}$ we can construct a function
$h_{a_{2}...a_{m}}$ which satisfies the extremal Kahler condition
at each point and also satisfies all of the smoothness conditions.
If $h_{a_{2}...a_{m}}>0$ for $x\in(0,1)$, then we have constructed
a global extremal Kahler metric on a
$\mathbb{C}P^{\frac{d_{1}+2}{2}}$-bundle over
$M_{2}\times...\times M_{m}$.
\begin{theorem}
Let $(M,g,J)$ be a fiberwise Kahler toric metric of
cohomogeneity-one such that $M$ is a
$\mathbb{C}P^{\frac{d_{1}+2}{2}}$-bundle over a product of
Kahler-Einstein manifolds $M_{2}\times...\times M_{m}$ as
described above.  \textbf{Every} such manifold admits, up to
scaling, an $(m-1)$ parameter family of extremal Kahler metrics.
\end{theorem}
\begin{proof}
We have already considered the case in which $m=1$.  When $m>1$,
the function $h_{a_{2}...a_{m}}$ defines an extremal Kahler metric
on $M$ if and only if
\begin{equation*}
h_{a_{2}...a_{m}}>0~\textrm{for}~x\in(0,1).
\end{equation*}
Letting
$h_{a_{2}...a_{m}}=\frac{P_{a_{2}...a_{m}}}{Q_{a_{2}...a_{m}}}$ as
above, we see that for $a_{i}$ sufficiently large for all $i$,
\begin{equation*}
P_{a_{2}...a_{m}}\sim\left(\prod_{j=2}^{m}a_{j}^{\frac{d_{j}}{2}}\right)\left(-\frac{4\alpha}
{(d_{1}+6)(d_{1}+4)}x^{\frac{d_{1}}{2}+3}-\frac{4\beta}{(d_{1}+4)(d_{1}+2)}x^{\frac{d_{1}}{2}+2}
+2x^{\frac{d_{1}}{2}}\right)
\end{equation*}
and
\begin{equation*}
Q_{a_{2}...a_{m}}\sim\left(\prod_{j=2}^{m}a_{j}^{\frac{d_{j}}{2}}\right)x^{\frac{d_{1}}{2}}.
\end{equation*}
Therefore,
\begin{equation*}
h_{a_{2}...a_{m}}\sim-\frac{4\alpha}{(d_{1}+6)(d_{1}+4)}x^{3}-\frac{4\beta}{(d_{1}+4)(d_{1}+2)}
x^{2}+2x.
\end{equation*}

    The conditions $h(1)=0$ and $h'(0)=-2$ imply that for the $a_{i}$ sufficiently large
\begin{equation*}
\alpha\sim0,
\end{equation*}
\begin{equation*}
\beta\sim\frac{1}{2}(d_{1}+4)(d_{1}+2),
\end{equation*}
and
\begin{equation*}
h_{a_{2}...a_{m}}\sim2x(1-x)\equiv h_{\infty}.
\end{equation*}
The function $h_{\infty}$ is positive on the interval $x\in(0,1)$.
Therefore, by the same reasoning used in the proof of the theorem
of the previous section, we deduce that $h_{a_{2}...a_{m}}>0$ for
$x\in(0,1)$ for the $a_{i}$ sufficiently large.  This gives the
$(m-1)$-parameter family of extremal Kahler metrics.
\end{proof}

\subsection{Constant Scalar Curvature Kahler Metrics on $\mathbb{C}P^{1}$-Bundles}

    We have just demonstrated the existence of extremal Kahler metrics on compact
codimension-one fiberwise Kahler toric manifolds in which only the
$S^{1}$ fiber collapses on the special orbits.  Globally, these
manifolds are $\mathbb{C}P^{1}$-bundles.  We found that all such
manifolds admit at least one family of extremal Kahler metrics. In
this section, we ask which if these $\mathbb{C}P^{1}$-bundles
admit at least one constant scalar curvature metrics.  Unlike in
the general extremal Kahler case, we will see below that there are
manifolds of this type which admit no constant scalar curvature
metric.

    In the compact constant scalar curvature case,

\begin{equation*}
h=\frac{1}{\prod_{j=1}^{m}(b_{j}x+a_{j})^{\frac{d_{j}}{2}}}\int\left[\int\left(\sum_{j=1}^{m}
\frac{d_{j}}{b_{j}x+a_{j}}-\beta\right)\prod_{j=1}^{m}(b_{j}x+a_{j})^{\frac{d_{j}}{2}}dx\right]dx
\end{equation*}
where $\beta$ is a positive constant.  By translation and
rescaling, we again set the location of the special orbits to be
$x=0$ and $x=1$. We again assume that at both special orbits
$A_{i}\neq0$ for all $1\leq i\leq m$.  The same four smoothness
conditions hold in the constant scalar curvature case as did in
the extremal Kahler case.

    If one of the $b_{i}$ vanishes, say $b_{1}$, then the function $h$ becomes
\begin{equation*}
h=\frac{1}{\prod_{j=2}^{m}(b_{j}x+a_{j})^{\frac{d_{j}}{2}}}\int\left[\int\left(\sum_{j=2}^{m}
\frac{d_{j}}{b_{j}x+a_{j}}-(\beta-\frac{d_{1}}{a_{1}})\right)\prod_{j=2}^{m}(b_{j}x+a_{j})^{
\frac{d_{j}}{2}}dx\right]dx.
\end{equation*}
Setting $\bar\beta=\beta-\frac{d_{1}}{a_{1}}$, we see that the
problem of finding a constant scalar curvature metric is unaltered
by the presence of absence of terms with $b_{i}=0$.  Therefore, we
can without loss of generality assume that $b_{i}\neq0$ for all
$i$.  As it is the $b_{i}$'s which determine the manifold, the
question we ask is: Given a collection of $b_{i}$, when is it
possible to choose the $a_{i}$ so that the metric is of constant
scalar curvature and all of the smoothness conditions can be
satisfied?

    There are two general cases: either $b_{i}>0$ (equivalently $b_{i}<0$) for all
$1\leq i \leq m$ or some of the $b_{i}$ are positive and some are
negative.  We consider the second case in the following theorem:

\begin{theorem}
Let $(M,g,J)$ be a fiberwise Kahler toric manifold of
cohomogeneity-one such that $M$ is a $\mathbb{C}P^{1}$-bundle as
above. If there exists $b_{i}>0$ for some $1\leq i\leq m$ and
$b_{j}<0$ for some $1\leq j\leq m$, then $M$ admits at least a
$1$-parameter family of constant positive scalar curvature Kahler
metric.
\end{theorem}
\begin{proof}
By reordering if necessary, let $b_{1},...,b_{l}>0$ and
$b_{l+1},...,b_{m}<0$ for some $2\leq l\leq m$.  Let
$c_{i}=\frac{a_{i}}{|b_{i}|}$.  We have
\begin{equation*}
V^{\frac{1}{2}}=\prod_{j=1}^{m}(b_{j}x+a_{j})^{\frac{d_{j}}{2}}=\prod_{j=1}^{m}|b_{j}|^{
\frac{d_{j}}{2}}\prod_{j=1}^{l}(x+c_{j})^{\frac{d_{j}}{2}}\prod_{j=l+1}^{m}(c_{j}-x)^{
\frac{d_{j}}{2}}
\end{equation*}
and
\begin{equation*}
\sum_{j=1}^{m}\frac{d_{j}}{b_{j}x+a_{j}}=\sum_{j=1}^{l}\frac{\frac{d_{j}}{|b_{j}|}}{c_{j}+x}
+\sum_{j=l+1}^{m}\frac{\frac{d_{j}}{|b_{j}|}}{c_{j}-x}.
\end{equation*}
To prove this theorem, we need only find constant scalar curvature
metrics for some $1$-parameter family of the $a_{i}$.  To this
end, we make the following assumption which simplifies the
calculations.

Set
\begin{equation*}
c_{j}=c>0
\end{equation*}
for all $1\leq j\leq l$ and set
\begin{equation*}
c_{i}=e>1
\end{equation*} for all $l+1\leq i\leq m$.

This implies that
\begin{equation*}
h=\frac{1}{(x+c)^{p}(e-x)^{q}}\int\left[\int\left(\frac{r}{x+c}+\frac{t}{e-x}-\beta\right)
(x+c)^{p}(e-x)^{q}dx\right]dx
\end{equation*}
where $p=\sum_{j=1}^{l}\frac{d_{j}}{2}$,
$q=\sum_{j=l+1}^{m}\frac{d_{j}}{2}$,
$r=\sum_{j=1}^{l}\frac{d_{j}}{|b_{j}|}$, and
$t=\sum_{j=l+1}^{m}\frac{d_{j}}{|b_{j}|}$ which are all greater
than zero and independent of the $c_{j}$.

Set $h(x)=\frac{T(x)}{L(x)}$ where
\begin{equation*}
T(x)=\int\left[\int\left(\frac{r}{x+c}+\frac{t}{e-x}-\beta\right)(x+c)^{p}(e-x)^{q}dx\right]dx
\end{equation*}
and
\begin{equation*}
L(x)=(x+c)^{p}(e-x)^{q}.
\end{equation*}

    Schematically, we can write $T$ as
\begin{equation*}
T=F_{p+q+2}x^{p+q+2}+...+F_{w}x^{w}+...+Ex+F
\end{equation*}
where $E$ and $F$ are the constants of integration and the $F_{w}$
are polynomials in $c,$ $e,$ $r,$ $t,$ $p,$ $q,$ and $\beta$.  For
the metric to be smooth at the special orbit $x=0$, we must have
$h(0)=0$ and $h'(0)=2$.  These conditions determine the constants
$E$ and $F$.  Because $L(x)>0$ for $x\in[0,1]$, $h(0)=0$ implies
that $T(0)=0$ which is equivalent to saying that $F=0$.
Furthermore, we see that $h'(0)=2$ implies that
$\frac{T'(0)}{L(0)}=2$.  This is equivalent to $E=2c^{p}e^{q}$.

    The two remaining conditions $h(1)=0$ and $h'(1)=-2$ determine $\beta$ and place a
condition on $c$ and $e$.

    To obtain a constant scalar curvature metric, it will be helpful to group the terms of
$T$ and $L$ by the order of $c$ and $e$.  It is straightforward to
see that
\begin{equation*}
T(x)=[c^{p-1}e^{q-1}(re+tc)\frac{x^{2}}{2}+...]-\beta[c^{p}e^{q}\frac{x^{2}}{2}+c^{p-1}e^{q-1}
(pe-qc)\frac{x^{3}}{6}+...]+2c^{p}e^{q}x.
\end{equation*}
Here the $'...'$ refers to lower order terms in $c$ and $e$.  The
condition $h(1)=0$ implies that $T(1)=0$.  This condition
determines the value of $S=\beta$,
\begin{equation}\label{csckbeta}
\beta=\frac{2c^{p}e^{q}+\frac{1}{2}c^{p-1}e^{q-1}(re+tc)+...}{\frac{1}{2}c^{p}e^{q}
+\frac{1}{6}c^{p-1}e^{q-1}(pe-qc)+...}.
\end{equation}
As $c\rightarrow\infty$ and $e\rightarrow\infty$ we see that
$\beta\rightarrow4$.  The next condition at the special orbit
$x=1$ is $h'(1)=-2$ which is equivalent to
$\frac{T'(1)}{L(1)}=-2$.  This condition can be interpreted as a
condition on $c$ and $e$.  To calculate this condition, we first
note that
\begin{equation*}
L(1)=c^{p}e^{q}+c^{p-1}e^{q-1}(pe-qc)+...
\end{equation*}
and
\begin{equation*}
T'(1)=[2c^{p}e^{q}+c^{p-1}e^{q-1}(re+tc)+...]-\beta[c^{p}e^{q}+\frac{1}{2}c^{p-1}e^{q-1}
(pe-qc)+...].
\end{equation*}
The condition $h'(1)=-2$ becomes
\begin{equation*}
\beta[c^{p}e^{q}+\frac{1}{2}c^{p-1}e^{q-1}(pe-qc)+...]=4c^{p}e^{q}+c^{p-1}e^{q-1}(re+tc)
+2c^{p-1}e^{q-1}(pe-qc)+...
\end{equation*}
Substituting equation (\ref{csckbeta}) into this equation we have
the condition
\begin{equation*}
\frac{2}{3}c^{2p-1}e^{2q-1}(pe-qc)-\frac{1}{12}c^{2p-2}e^{2q-2}(pe-qc)(re+tc)+\frac{1}{3}
c^{2p-2}e^{2q-2}(pe-qc)^{2}+...=0
\end{equation*}
This condition is equivalent to
\begin{equation*}
\frac{2}{3}c^{2p-1}e^{2q-1}(pe-qc)+K(c,e)=0
\end{equation*}
where $K=K(c,e)$ is a polynomial in $c$ and $e$ in which each of
the terms in $K$ is of the form $K_{\eta\xi}c^{\eta}e^{\xi}$ such
that $\eta\leq2p$, $\xi\leq2q$, $\eta+\xi<2p+2q-1,$ and the
$K_{\eta\xi}$ are constants depending on the fixed values of $r,$
$t,$ $p,$ and $q$. Dividing by $c^{2p-1}e^{2q-1}$ implies that
\begin{equation}\label{eccondition}
B(c,e)\equiv pe-qc+\frac{3}{2}\frac{K(c,e)}{c^{2p-1}e^{2q-1}}=0.
\end{equation}
As $e$ and $c$ go to infinity,
\begin{equation*}
\frac{3}{2}\frac{K(c,e)}{c^{2p-1}e^{2q-1}}\rightarrow\bar{K}
\end{equation*}
for some constant $\bar{K}$. For all $0<\epsilon\ll 1$, there
exists $N_{\epsilon}\gg0$ such that
$|\frac{3}{2}\frac{K(c,e)}{c^{2p-1}e^{2q-1}}-\bar{K}|<\epsilon$
whenever $pe,qc>N_{\epsilon}$.  Set $pe=2N_{\epsilon}+2\epsilon$.
If $qc=2N_{\epsilon}+\bar{K}+\epsilon$ then the quantity
$B=B(c,e)>0$.  If, on the other hand,
$qc=2N_{\epsilon}+\bar{K}+3\epsilon$, then $B<0$.  By the
intermediate value theorem, there exists
$qc\in(2N_{\epsilon}+\bar{K}+\epsilon,2N_{\epsilon}+\bar{K}+3\epsilon)$
such that $B=0$ and equation (\ref{eccondition}) is satisfied.

Therefore as $qc\rightarrow\infty$ the condition $h'(1)=-2$ can be
solved by setting $pe$ equal to a value increasingly close to
$qc-\bar{K}$.  Since $q$ and $p$ are positive constants, we have a
family of solutions to the smoothness conditions in which $c$ and
$e$ increase together.  These solutions have zeros of order one at
$x=0,1$.  We must show that the solutions are positive for
$x\in(0,1)$ for $e$ and $c$ sufficiently large.  It is
straightforward to see that as $c$ and $e$ approach infinity
\begin{equation*}
h\sim 2x(1-x).
\end{equation*}
The family of solutions to the smoothness conditions tend to a
function which has zeros of order one at $x=0,1$ and is positive
for $x\in(0,1)$.  By the same argument used in the proof of the
existence of extremal Kahler metrics above, we see that for $c$
and $e$ sufficiently large, the function $h$ will be positive for
$x\in(0,1)$ and will define a constant scalar curvature metric.

To summarize, we have shown the existence of a one-parameter
family of constant scalar curvature metrics for $c$ and $e$ and
hence the $a_{i}$ sufficiently large.

\end{proof}

This existence proof fails when $b_{i}>0$ (or $b_{i}<0$) for all
$1\leq i\leq m$.  Indeed, there are known examples where there is
no constant scalar curvature metric in which all of the $b_{i}$'s
have the same sign.  The best known examples are the Hirzebruch
surfaces which will be discussed below.

\subsection{Constant Scalar Curvature Metrics on $\mathbb{C}P^{\frac{d_{1}+2}{2}}$-Bundles}

    In this section, we set
\begin{equation*}
A_{1}(0)=0
\end{equation*}
while continuing to set $A_{i}(x)>0$ when $x\in[0,1]$ for all
$2\leq i\leq m$.  This condition is equivalent to setting
$a_{1}=0$.  In this case, the manifold $M$ is a
$\mathbb{C}P^{\frac{d_{1}+2}{2}}$-bundle over
$M_{2}\times...\times M_{m}$.

    As in the extremal Kahler situation discussed above, we restrict to the case to the
fiberwise Kahler manifolds of Construction Two in which each
principal orbit is an $S^{1}$-bundle over $M_{1}\times...\times
M_{m}$.  For the metric to be extended smoothly over the special
orbit at $x=0$, we must have
$M_{1}\cong\mathbb{C}P^{\frac{d_{1}}{2}}$ and
$b_{1}=\frac{2}{d_{1}+2}$.

    We have already seen that we retrieve the Fubini-Study metric on $M\cong\mathbb{C}
    P^{\frac{d_{1}+2}{2}}$
when $m=1$.  Therefore, we need only consider the case in which
$m>1$.

    As above, we can without loss of generality assume that none of the $b_{i}$ are
zero.  We distinguish between three general cases: $1.$ $b_{i}>0$
for all $2\leq i\leq m$; $2.$ $b_{i}<0$ for all $2\leq i\leq m$;
and $3.$ for $2\leq i\leq m$ some of the $b_{i}$ are positive and
some are negative.  Obviously, in the third case $m\geq3$.  We
consider the third case in the following theorem:
\begin{theorem}
Let $(M,g,J)$ be a fiberwise Kahler toric manifold obtained from
Construction Two such that each principal orbit is an
$S^{1}$-bundle over $M_{1}\times...\times M_{m}$ such that
$M_{1}\cong\mathbb{C}P^{\frac{d_{1}}{2}}$ and $M$ is a
$\mathbb{C}P^{\frac{d_{1}+2}{2}}$-bundle over
$M_{2}\times...\times M_{m}$.  If there exist $b_{i}>0$ for some
$2\leq i \leq m$ and $b_{j}<0$ for some $2\leq j\leq m$, then $M$
admits at least a $1$-parameter family of constant positive scalar
curvature Kahler metrics.
\end{theorem}
\begin{proof}
By reordering if necessary, let $b_{2},...,b_{l}>0$ and
$b_{l+1},...,b_{m}>0$ for some $3\leq l\leq m$.  Let
$c_{i}=\frac{a_{i}}{|b_{i}|}$.  We see that
\begin{equation*}
\prod_{j=1}^{m}(b_{j}x+a_{j})^{\frac{d_{j}}{2}}=\frac{2}{d_{1}+2}x^{\frac{d_{1}}{2}}
\prod_{j=2}^{m}|b_{j}|^{\frac{d_{j}}{2}}\prod_{j=2}^{l}(x+c_{j})^{\frac{d_{j}}{2}}
\prod_{j=l+1}^{m}(c_{j}-x)^{\frac{d_{j}}{2}}
\end{equation*}
and
\begin{equation*}
\sum_{j=1}^{m}\frac{d_{j}}{b_{j}x+a_{j}}=\frac{d_{1}(d_{1}+2)}{2x}+\sum_{j=2}^{l}\frac{
\frac{d_{j}}{|b_{j}|}}{c_{j}+x}+\sum_{j=l+1}^{m}\frac{\frac{d_{j}}{|b_{j}|}}{c_{j}-x}.
\end{equation*}
To find our $1$-parameter family of constant scalar curvature
metrics, we make the following assumption to simplify
calculations.  Set
\begin{equation*}
c_{j}=c>0
\end{equation*}
for all $2\leq j\leq l$ and set
\begin{equation*}
c_{i}=e>1
\end{equation*}
for all $l+1\leq i\leq m$.

    This implies that
\begin{equation*}
h=\frac{1}{x^{k}(x+c)^{p}(e-x)^{q}}\int\left[\int\left(\frac{2k(k+1)}{x}+\frac{r}{x+c}+\frac{t}
{e-x}-\beta\right)x^{k}(x+c)^{p}(e-x)^{q}dx\right]dx
\end{equation*}
where $k=\frac{d_{1}}{2},$ $p=\sum_{j=2}^{l},$
$q=\sum_{j=l+1}^{m},$ $r=\sum_{j=2}^{l}\frac{d_{j}}{|b_{j}|},$ and
$t=\sum_{j=l+1}^{m}\frac{d_{j}}{|b_{j}|}$ which are all greater
than zero and independent of the $a_{i}$.

    Set $h(x)=\frac{T(x)}{L(x)}$ where
\begin{equation*}
T(x)=\int\left[\int\left(\frac{2k(k+1)}{x}+\frac{r}{x+c}+\frac{t}{e-x}-\beta\right)x^{k}(x+c)^{p}
(e-x)^{q}dx\right]dx
\end{equation*}
and
\begin{equation*}
L(x)=x^{k}(x+c)^{p}(e-x)^{q}.
\end{equation*}

    Schematically, we can write $T$ as
\begin{equation*}
T=F_{k+p+q+2}x^{k+p+q+2}+...+F_{w}x^{w}+...+F_{k+2}x^{k+2}+2c^{p}e^{q}x^{k+1}+Ex+F
\end{equation*}
where $E$ and $F$ are the two constant of integration and the
$F_{w}$ are polynomials in $c,$ $e,$ $r,$ $t,$ $p,$ $q,$ and
$\beta$.  For the metric to be smooth at the special orbit $x=0$,
we must have $h(0)=0$ and $h'(0)=2$.  It is straightforward to see
that these two conditions hold if and only if $E=F=0$.  The two
remaining conditions $h(1)=0$ and $h'(1)=-2$ determine $\beta$ and
place a condition on $c$ and $e$.

    It will be helpful to group the terms of $T$ and $L$ by the order of $c$ and $e$.  We see,
after factoring out $x^{k}$ from $T$ and $L$, that
\begin{equation*}
W=\frac{T(x)}{x^{k}}=[2c^{p}e^{q}x+\frac{c^{p-1}e^{q-1}}{(k+2)(k+1)}(re+tc)x^{2}+...]
\end{equation*}
\begin{equation*}
-\beta[\frac{c^{p}e^{q}}{(k+2)(k+1)}x^{2}+\frac{c^{p-1}e^{q-1}}{(k+3)(k+2)}(pe-qc)x^{3}+...]
\end{equation*}
and
\begin{equation*}
Z=\frac{L(x)}{x^{k}}=c^{p}e^{q}+c^{p-1}e^{q-1}(pe-qc)x+...
\end{equation*}
The condition $h(1)=0$ holds if and only if $W(1)=T(1)=0$.  This
conditions determines the value of $\beta$,
\begin{equation*}
\beta=\frac{2c^{p}e^{q}+\frac{c^{p-1}e^{q-1}}{(k+2)(k+1)}(re+tc)+...}{\frac{c^{p}e^{q}}
{(k+2)(k+1)}+\frac{c^{p-1}e^{q-1}}{(k+3)(k+2)}(pe-qc)+...}.
\end{equation*}
As $c\rightarrow\infty$ and $e\rightarrow\infty$, we see that
$\beta\rightarrow2(k+2)(k+1)$.  The second condition as the
special orbit $x=1$, $h'(1)=-2$, is equivalent to
$\frac{W'(1)}{Z(1)}=-2$.  Having found an expression for $\beta$,
this condition becomes an equation in $c$ and $e$.  We note that
\begin{equation*}
W'(1)=[2c^{p}e^{q}+\frac{2c^{p-1}e^{q-1}}{(k+2)(k+1)}(re+tc)+...]-\beta[\frac{2c^{p}e^{q}}
{(k+2)(k+1)}+\frac{3c^{p-1}e^{q-1}}{(k+3)(k+2)}(pe-qc)+...]
\end{equation*}
and
\begin{equation*}
Z(1)=c^{p}e^{q}+c^{p-1}e^{q-1}(pe-qc)+...
\end{equation*}
The condition $h'(1)=-2$ becomes
\begin{equation*}
\beta[\frac{2c^{p}e^{q}}{(k+2)(k+1)}+\frac{3c^{p-1}e^{q-1}}{(k+3)(k+2)}(pe-qc)+...]=4c^{p}e^{q}
+2c^{p-1}e^{q-1}(pe-qc)+\frac{2c^{p-1}e^{q-1}}{(k+2)(k+1)}(re+tc).
\end{equation*}
Combining this with our expression for $\beta$ above implies that
\begin{equation}\label{X}
\frac{4c^{2p-1}e^{2q-1}}{(k+3)(k+2)(k+1)}(pe-qc)+K(c,e)=0
\end{equation}
where $K=K(c,e)$ is a polynomial is $c$ and $e$ in which each of
the terms is of the form $K_{\eta\xi}c^{\eta}e^{\xi}$ such that
$\eta\leq2p$, $\xi\leq2q,$ $\eta+\xi<2p+2q-1,$ and the
$K_{\eta\xi}$ are constants depending on the fixed values of $r,$
$t,$ $p,$ and $q$.

    Just as in our construction of constant scalar curvature metrics
on $\mathbb{C}P^{1}$-bundles this equation can always be solved
with $e$ and $c$ both arbitrarily large.  We do not write the
details here as they are almost identical to those of the proof of
the last theorem in the previous section. We note that, as $c$ and
$e$ approach infinity
\begin{equation*}
h\sim2x(1-x).
\end{equation*}
Therefore, for $c$ and $e$ arbitrarily large and satisfying
(\ref{X}), $h(x)$ will be positive for $x\in(0,1)$.  This results
gives a one-parameter family of constant scalar curvature metrics
and completes the proof.
\end{proof}

\subsection{Four-Dimensional Extremal Kahler Manifolds}

    On a four dimensional Kahler toric manifold, the number of summands, $m$, is
equal to one and the dimension of that summand is $d_{1}=d=2$.
Therefore, each principal orbit must be an $S^{1}$-bundle over
$\mathbb{C}P^{1}$.    To simplify the presentation, we set
$b_{1}=b$ and $a_{1}=a$.  From the definition of $b$, we see that
\begin{equation*}
b=-\frac{q}{2}
\end{equation*}
where $q$ is an integer and each value of $q$ determines a
distinct principal orbit type.  In dimension four, the fiberwise
Kahler toric metric can be written as
\begin{equation*}
g=\begin{pmatrix}
\frac{1}{h} & & \\
& h & \\
& & (bx+a)Id_{2} \\
\end{pmatrix}.
\end{equation*}
    Because the manifold is compact, there must be two special orbits.  Via
translation in $x$, fix the first special orbit at $x=-1$.  By
rescaling, we can set the second special orbit at $x=1$.  For the
metric to be positive on $[-1, 1]$, $a>0$ and $bx+a>0$ when
$x\in(-1,1)$.

    From the discussion above, we notice two possibilities for the special orbits.  Either
both special orbits are copies of $\mathbb{C}P^{1}$ or one of the
special orbits is $\mathbb{C}P^{1}$ and the other is a point. In
the first case, $bx+a>0$ on $x\in[-1,1]$ and in the second case,
$a=b$ or $a=-b$. We will address these two case separately.

~

\textbf{Case 1: One Special Orbit is a Point}:

~

If the special orbit at $x=-1$ is a point, then $a=b$.  The
dimension of the sphere vanishing as that point is three: $k=3$.
The smoothness condition demands that
$\frac{d}{b}=\frac{k^{2}-1}{2}=\frac{3^{2}-1}{2}$ which implies
that $b=a=\frac{1}{2}$.  Plugging this into the equation and
integrating, we retrieve the Fubini-Study metric on
$\mathbb{C}P^{2}$.

~

\textbf{Case 2: Both Special Orbits are Copies of
$\mathbb{C}P^{1}$}:

~

    When both the special orbits are two-spheres, the manifold $M$ will be a Hirzebruch
surface indexed by the integer $q$. Calabi found that each
Hirzebruch surface admits, up to scaling, a one-parameter family
of extremal Kahler metrics \cite{Besse}. However, these metrics
were not written down explicitly.  In fact, Hwang and Simanca
prove that within each of these one-parameter families is an
extremal Kahler metric locally equivalent to an Hermitian-Einstein
metric \cite{Hwang}.  This correspondence, however, defines
globally a smooth Hermitian-Einstein four manifold when $q=0,1$.

For this metric to be extremal Kahler, we have shown that the
function $h$ must satisfy the equation
\begin{equation}\label{FKTHirzExtremal}
h=\frac{1}{bx+a}\int\left[\int\left(\frac{2}{bx+a}-(\alpha
x+\beta)\right)(bx+a)dx\right]dx
\end{equation}
for $x\in(-1,1)$.  For the metric to be smooth at the special
orbits $h(-1)=h(1)=0$, $h'(-1)=2$, and $h'(1)=-2$. Performing this
simple double integral and imposing the boundary conditions we
find that
\begin{equation}\label{alphaFKT}
\alpha=\frac{6b(2a-1)}{3a^{2}-b^{2}},
\end{equation}
\begin{equation*}
\beta=\frac{6(a-b^{2}+a^{2})}{3a^{2}-b^{2}},
\end{equation*}
and
\begin{equation}
h=\frac{(1-x)(1+x)((2a-1)b^{2}x^{2}+2(3a^{2}-b^{2})bx+6a^{3}+b^{2}-4b^{2}a)}{2(3a^{3}-b^{2})
(bx+a)}.
\end{equation}

    For $a$ sufficiently large, this function is positive for $x\in(-1,1)$.  Therefore,
for each value of $b$ (i.e. on each Hirzebruch surface), there
exists a one-parameter family of extremal Kahler metrics.

    It remains to determine when these extremal Kahler metrics are of constant scalar
curvature.  That is, when do we get smooth metrics with
$\alpha=0$?  From the form of $\alpha$ in the extremal Kahler
case, $\alpha=0$ if and only if $b=0$ or $a=\frac{1}{2}$ by
equation (\ref{alphaFKT}).  If $b=0$ then the $S^{1}$ bundle over
$\mathbb{C}P^{1}$ is trivial and the manifold
$\mathbb{C}P^{1}\times\mathbb{C}P^{1}$ with its product metric.
Moreover, the only extremal Kahler metrics on
$\mathbb{C}P^{1}\times\mathbb{C}P^{1}$ are of constant scalar
curvature. If $a=\frac{1}{2}$, then $|q|<\frac{1}{2}$.  This
follows from the fact that $a+b>0$ and $a-b>0$ for the metric to
be positive at the special orbits.  However, $q$ is an integer and
this is satisfied if and only if $q=0$.  Therefore, $b$ must equal
zero and the manifold is a again a product.

    Note: Because each $\mathbb{C}P^{1}$ in the fibers is $S^{1}$-invariant, these
manifolds are Kahler toric manifolds whose polytopes are
four-sided.  In fact, these are the only Kahler toric four
manifolds whose polytopes have four sides.  Transforming the
metric from fiberwise Kahler toric coordinates to Kahler toric
coordinates is straightforward.  The extremal Kahler metrics on
the blow-up of $\mathbb{C}P^{2}$ at one point (the Hirzebruch
surface with $q=1$ or $-1$) was written by Abreu in Kahler toric
coordinates \cite{Abreu}.

\subsection{Six-Dimensional Extremal Kahler Manifolds}

    While the classification of four-dimensional fiberwise Kahler toric manifolds admitting was
already known, there remain open questions in dimension six.  Here
we will classify all six-dimensional fiber Kahler toric manifolds
admitting extremal Kahler metrics and write those metrics down
explicitly in terms of rational functions. In dimension six, each
principal orbit is an $S^{1}$-bundle over a four-dimensional
Kahler-Einstein manifold.  The only four manifolds admitting a
Kahler-Einstein metric with positive Einstein constant are
$\mathbb{C}P^{2}$, $\mathbb{C}P^{1}\times\mathbb{C}P^{1}$, or the
blow-up of $\mathbb{C}P^{2}$ at $k$ distinct points where $3\leq
k\leq 8$ such that no three points lie on a line and no six points
lie on a conic \cite{YandT}.

    When the base manifold of the $S^{1}$-bundle is any manifold other
than $\mathbb{C}P^{1}\times\mathbb{C}P^{1}$, the metric can be
written as
\begin{equation*}
g=\begin{pmatrix}
\frac{1}{h} & & \\
& h & \\
& & (bx+a)Id_{4} \\
\end{pmatrix}
\end{equation*}
since the Kahler-Einstein metrics on $\mathbb{C}P^{2}$ and the
blow-ups are irreducible.

On the other hand, if each principal orbit is an $S^{1}$-bundle
over $\mathbb{C}P^{1}\times\mathbb{C}P^{1}$, the metric looks like
\begin{equation*}
g=\begin{pmatrix}
\frac{1}{h} & & & \\
& h & & \\
& & (b_{1}x+a_{1})Id_{2} \\
& & & (b_{2}x+a_{2})Id_{2}
\\
\end{pmatrix}.
\end{equation*}

    To the author's knowledge, it had not been determined before whether or
not manifolds of this type admit extremal Kahler metrics.

    For ease of presentation, we set $b_{1}=b$, $a_{1}=a$, $b_{2}=d$, and $a_{2}=c$.
The values of $b$ and $d$ determine the principal orbit type.  As
each principal orbit is an $S^{1}$ bundle over
$\mathbb{C}P^{1}\times\mathbb{C}P^{1}$, $b=-\frac{p}{2}$ and
$d=-\frac{q}{2}$ where $p$ and $q$ are integers \cite{DandW}. Let
$P(p,q)$ denote the principal orbit type.  Note that $P(p,q)\cong
P(-p,-q)$.

    Again, by translation in $x$ and rescaling of the metric, we can set the location
of the two special orbits to be $x=-1$ and $x=1$.  First, we
assume that both the special orbits are $S^{2}\times S^{2}$, i.e.
\begin{equation*}
M=\mathbb{P}(\mathcal{O}_{\mathbb{C}P^{1}\times\mathbb{C}P^{1}}\oplus\mathcal{O}_{
\mathbb{C}P^{1}\times\mathbb{C}P^{1}}(p,q))
\end{equation*}

a $\mathbb{C}P^{1}$-bundle over
$\mathbb{C}P^{1}\times\mathbb{C}P^{1}$.

 The metric is positive if and only if $a>0$,
$c>0$, $|\frac{a}{b}|>1$, $|\frac{c}{d}|>1$, and $h>0$ when
$x\in(-1,1)$.  Finally, smoothness at the special orbits requires
that $h(-1)=h(1)=0$, $h'(-1)=2$, and $h'(1)=-2$.

    In the extremal Kahler case, $h$ satisfies (\ref{FKThsolution}) with $S=\alpha x+\beta$.
Integrating this equation and imposing smoothness conditions we
find that

~

$h=\frac{1}{2}(1-x)(x+1)((6a^2d^2bc+6db^2ac^2-3ba^2d^2-3b^2c^2d+b^3d^2+b^2d^3)x^3
+(3cb^2d^2+3ab^2d^2+25dba^2c^2-5a^3d^2+10c^3b^2a+10a^3d^2c-5a^2d^3b-5c^2b^3d-5c^3b^2+3d^3b^3
-2d^2b^2ac)x^2+(-8a^2d^2bc-8db^2ac^2+30a^2c^3b-b^3d^2+6cd^2b^3-10c^3b^3-10a^3d^3-b^2d^3+3ba^2d^2
+3b^2c^2d+6ad^3b^2+30a^3c^2d)x-5c^2b^3d+13dba^2c^2-5a^2d^3b+5a^3d^2+5c^3b^2+30a^3c^3+3d^3b^3
-20a^3d^2c-20c^3b^2a+16d^2b^2ac-3ab^2d^2-3cb^2d^2)/((-5c^2b^2+15a^2c^2+4dbac-5a^2d^2+3d^2b^2)
(bx+a)(dx+c))$.

~

    The constants $\alpha$ and $\beta$ in the scalar curvature satisfy
\begin{equation}
\alpha=\frac{10(6a^2cd+6abc^2+bd^2+db^2-3c^2b-3a^2d)}{-5c^2b^2+15a^2c^2+4dbac-5a^2d^2+3d^2b^2}
\end{equation}
and
\begin{equation*}
\beta=\frac{2(15c^2a-5cb^2+4dbc+15a^2c-5ad^2+4bad-6dbac+15a^2c^2-15c^2b^2-15a^2d^2+9d^2b^2)}
{-5c^2b^2+15a^2c^2+4dbac-5a^2d^2+3d^2b^2}.
\end{equation*}

    All of the smoothness conditions are satisfied for all values of $a$, $b$, $c$, and $d$.
However, the function $h$ is not necessarily positive for all
values of these variable.  It is positive for $a$ and $c$
sufficiently large. Therefore, for every choice of $b$ and $d$
(that is for every principal orbit type) there exists, up to
scale, a two parameter family of extremal Kahler metrics.

\begin{theorem}
Let $M$ be
$\mathbb{P}(\mathcal{O}_{\mathbb{C}P^{1}\times\mathbb{C}P^{1}}\oplus\mathcal{O}_{\mathbb{C}
P^{1}\times\mathbb{C}P^{1}}(p,q))$ endowed with a
cohomogeneity-one fiberwise Kahler toric metric as described
above.  Every such manifold admits, up to scaling, a two parameter
family of extremal Kahler metrics.
\end{theorem}

    It remains to determine which of these manifolds admit metrics of constant scalar curvature.  To determine this, we need only
determine when $\alpha=0$ while the smoothness and positiveness
conditions are satisfied.

\begin{theorem}
Let
$M=\mathbb{P}(\mathcal{O}_{\mathbb{C}P^{1}\times\mathbb{C}P^{1}}\oplus\mathcal{O}_{\mathbb{C}
P^{1}\times\mathbb{C}P^{1}}(p,q))$,
with $p$ and $q$ as above.  For $(p,q)=(0,0)$, the manifold is the
product of three copies of $\mathbb{C}P^{1}$ and admits, up to
scale, a two parameter family of constant scalar curvature
metrics.  If $p>0$ and $q<0$, then the manifold admits at least a
one-parameter family of constant scalar curvature metrics within
its two-parameter family of extremal Kahler metrics.  For $p\geq0$
and $q>0$, the manifold admits no constant scalar curvature Kahler
metrics.
\end{theorem}
\begin{proof}
    Recall that $b=-\frac{p}{2}$ and $d=-\frac{q}{2}$ for $p,q\in\mathbb{Z}$ and that $a>\frac{
|p|}{2}$ and $c>\frac{|q|}{2}$.  The metric has constant scalar
curvature if $\alpha=0$ which is equivalent to
\begin{equation*}
3a^{2}d(2c-1)+3c^{2}b(2a-1)+bd(b+d)=0.
\end{equation*}
If $b=d=0$, then $\alpha=0$ automatically and there is a
two-parameter family of extremal Kahler metrics. On the other hand
if $b\geq0$ and $d>0$ then $\alpha\neq0$ because $a$ and $c$ are
both greater then $\frac{1}{2}$.

Finally, if $b>0$ and $d<0$, then we can solve for $c$ in terms of
$a$. The metric has constant scalar curvature if
\begin{equation}\label{csckcond}
c=\frac{-6a^{2}d\pm\sqrt{9a^{4}d^{2}+(3b^{3}d+3b^{2}d^{2}-9bda^{2})(1-2a)}}{3b(2a-1)}.
\end{equation}
Consider the first of these solutions (taking '$+$') as
$a\rightarrow\infty$, $c\rightarrow\infty$. Therefore, for
sufficiently large $a$ we are certain to have a constant scalar
curvature metric with $h$ positive on $(-1,1)$.  This gives the
one-parameter family of constant scalar curvature metrics.

\end{proof}

    Note that this theorem determines which of these \textit{manifolds} admit a
constant scalar curvature metric.  It further demonstrates that
those manifolds which do admit a constant scalar curvature Kahler
metric in fact admit at least a one-parameter family.  They may
admit more.  In fact,
$\mathbb{C}P^{1}\times\mathbb{C}P^{1}\times\mathbb{C}P^{1}$
trivially admits a full two-parameter family of such metrics. What
is more interesting is that
$\mathbb{P}(\mathcal{O}_{\mathbb{C}P^{1}\times\mathbb{C}P^{1}}\oplus\mathcal{O}_{\mathbb{C}
P^{1}\times\mathbb{C}P^{1}}(1,-1))$,
as demonstrated in a previous section, admits \textit{two}
one-parameter families of constant scalar curvature metrics.  This
occurs because both solutions for the constant $c$ in terms of the
constant $a$ can be solved without violating the smoothness and
positive definiteness conditions.

\begin{prop}
For
$M=\mathbb{P}(\mathcal{O}_{\mathbb{C}P^{1}\times\mathbb{C}P^{1}}\oplus\mathcal{O}_{\mathbb{C}
P^{1}\times\mathbb{C}P^{1}}(p,q))$, $M$ admits, up to scale, more
than a one-parameter family of constant scalar curvature Kahler
metrics if and only if $(p,q)=(0,0)$ or $(1,-1)$.
\end{prop}
\begin{proof}
The proof consists of demonstrating that the '$-$' solution of
(\ref{csckcond}) only gives a valid solution for those two values
of $p$ and $q$.  This is trivially done but awkward to present. We
omit the details.
\end{proof}

The existence of Kahler-Einstein metrics on manifolds of this type
was already answered by Sakane.

\begin{theorem}\cite{Sakane}
For
$M=\mathbb{P}(\mathcal{O}_{\mathbb{C}P^{1}\times\mathbb{C}P^{1}}\oplus\mathcal{O}_{
\mathbb{C}P^{1}\times\mathbb{C}P^{1}}(p,q))$,
$M$ admits a Kahler-Einstein metric if and only if $(p,q)=(0,0)$
or $(1,-1)$.
\end{theorem}

    We note that the manifolds admitting Kahler-Einstein metrics are precisely those that
admit more than a one-parameter of constant scalar curvature
metrics. It would be interesting to study whether there is a
connection between the existence of Kahler-Einstein metrics and
the 'number' of constant scalar curvature Kahler metrics.

\chapter{Cohomogeneity-$d$ Fiberwise Kahler Toric Manifolds}

    In the cohomogeneity-one case, the Einstein and extremal Kahler conditions were seen
to be second-order ordinary differential equations of one function
in one variable, $h$, which we could integrate explicitly. The
Einstein and extremal Kahler conditions in the cohomogeneity-$d$
case are partial differential equations and we must therefore use
a slightly different strategy of integration.  Our integration
procedure is similar to that employed to integrate the Einstein
equations in the Kahler toric case.
\section{Explicit Integration of the Einstein Condition}

    Above, we calculate that the Einstein condition on a fiberwise Kahler toric metric
is equivalent to equations (\ref{FKTEinsteinA}) and
(\ref{FKTEinsteinB}). We begin our integration of the Einstein
equations by equation (\ref{FKTEinsteinA}):
\begin{equation*}
h_{ik}\frac{\partial}{\partial x^{k}}\left(h_{jl}\frac{\partial
\textrm{log}(\textrm{det}(h)V^{\frac{1}{2}})}{\partial
x^{l}}\right)=-2\lambda h_{ij}.
\end{equation*}
Multiplying both sides of this equation by $h^{im}$ we find that
\begin{equation*}
\frac{\partial}{\partial x^{k}}\left(h_{lj}\frac{\partial
\textrm{log}(\textrm{det}(h)V^{\frac{1}{2}})}{\partial
x^{l}}\right)=-2\lambda\delta^{k}_{j}.
\end{equation*}
Integrating these equations we find that
\begin{equation}\label{FKTEinsteinC}
h_{lj}\frac{\partial
\textrm{log}(\textrm{det}(h)V^{\frac{1}{2}})}{\partial
x^{l}}=-2\lambda x_{j}+C_{j}
\end{equation}
for all $1\leq j\leq d$ where the $C_{j}$ are constants.  We can
integrate these equations by first multiplying both sides of
equation (\ref{FKTEinsteinC}) by $h^{mj}$ to obtain
\begin{equation*}
\frac{\partial
\textrm{log}(\textrm{det}(h)V^{\frac{1}{2}})}{\partial
x^{l}}=h^{lj}(-2\lambda x_{j}+C_{j}).
\end{equation*}
Recall that $h^{lj}=\frac{\partial^{2}\Phi}{\partial x^{l}\partial
x^{j}}$ and therefore
\begin{equation*}
\frac{\partial
\textrm{log}(\textrm{det}(h)V^{\frac{1}{2}})}{\partial
x^{l}}=(-2\lambda x_{j}+C_{j})\frac{\partial^{2}\Phi}{\partial
x^{l}\partial x^{j}}.
\end{equation*}
These equations integrate to
\begin{equation*}
\textrm{log}(\textrm{det}(h)V^{\frac{1}{2}})=\sum_{j=1}^{d}(-2\lambda
x_{j}+C_{j})\frac{\partial\Phi}{\partial x^{j}}+2\lambda\Phi+E
\end{equation*}
where $E$ is a constant.

    In light of equation (\ref{FKTEinsteinC}), equation (\ref{FKTEinsteinB}) of the
Einstein condition becomes

\begin{equation*}
\left(h_{kl}\frac{\partial
\textrm{log}(\textrm{det}(h)V^{\frac{1}{2}})}{\partial
x^{l}}+2\lambda
x_{k}\right)b_{ki}=\sum_{k=1}^{d}C_{k}b_{ki}=2(1-\lambda a_{i})
\end{equation*}
for all $1\leq i\leq m$.  These results prove the following
theorem.
\begin{theorem}
Let $(M^{n},g,J)$ be a fiberwise Kahler toric manifold of
cohomogeneity-$d$.  The Einstein condition is equivalent to
equation
\begin{equation}\label{FKTEinsteinD}
{\rm log}({\rm det}(h)V^{\frac{1}{2}})=\sum_{j=1}^{d}(-2\lambda
x_{j}+C_{j})\frac{\partial\Phi}{\partial x^{j}}+2\lambda\Phi+E
\end{equation}
and the relation
\begin{equation}\label{FKTEinsteinE}
\begin{pmatrix}
b_{11} & \cdot & \cdot & \cdot & b_{d1} \\
\cdot & \cdot & & & \cdot \\
\cdot & & \cdot & & \cdot \\
\cdot & & & \cdot & \cdot \\
b_{1m} & \cdot & \cdot & \cdot & b_{dm} \\
\end{pmatrix}\begin{pmatrix}
C_{1} \\
\cdot \\
\cdot \\
\cdot \\
C_{d} \\
\end{pmatrix}=\begin{pmatrix}
2(1-\lambda a_{1}) \\
\cdot \\
\cdot \\
\cdot \\
2(1-\lambda a_{m})
\end{pmatrix}
\end{equation}
where the $a_{r}$, the $b_{jr}$, the $C_{j}$ and $E$ are constants
and $V$, $\Phi$, and $\rm{det}h$ are as defined above.
\end{theorem}
Expressed in complex coordinates, this theorem becomes the
following.

\begin{theorem}
Let $(M,g,J)$ be a fiberwise Kahler toric manifold with complex
coordinates $\eta=\eta_{u_{1},...,u_{d}}$ as above.  The Einstein
condition is equivalent to the equation
\begin{equation*}
{\rm log}({\rm
det}(h)V^{\frac{1}{2}})=-2\lambda\eta+\sum_{j=1}^{d}C_{j}u_{j}+E
\end{equation*}
and the relation
\begin{equation}
\begin{pmatrix}
b_{11} & \cdot & \cdot & \cdot & b_{d1} \\
\cdot & \cdot & & & \cdot \\
\cdot & & \cdot & & \cdot \\
\cdot & & & \cdot & \cdot \\
b_{1m} & \cdot & \cdot & \cdot & b_{dm} \\
\end{pmatrix}\begin{pmatrix}
C_{1} \\
\cdot \\
\cdot \\
\cdot \\
C_{d} \\
\end{pmatrix}=\begin{pmatrix}
2(1-\lambda a_{1}) \\
\cdot \\
\cdot \\
\cdot \\
2(1-\lambda a_{m})
\end{pmatrix}.
\end{equation}
In complex coordinates, the functions $A_{i}$ satisfy
\begin{equation*}
A_{i}=\sum_{j=1}^{d}b_{ji}\frac{\partial\eta}{\partial
u_{j}}+a_{i}.
\end{equation*}
\end{theorem}
From the Einstein condition, we see that much of the structure of
the induced Kahler toric metric $\bar{g}$ determines much of the
structure of the total fiberwise Kahler toric metric in the
Einstein case.  For example, we have the following theorem.
\begin{theorem}
Let $(M,g)$ be a cohomogeneity-$d$ fiberwise Kahler toric metric
such that $g$ is Kahler-Einstein.  If, locally, $(N,\bar{g})$ can
be written as a product $(N_{1}\times
N_{2},\bar{g}_{1}\oplus\bar{g}_{2})$ of two Kahler toric metrics
where $\rm{dim} N_{1}=2l_{1}$ and $\rm{dim} N_{2}=2l_{2}$ then
$(M,g)$ splits as $(M_{1}\times M_{2},g_{1}\oplus g_{2})$ where
$M_{1}$ and $M_{2}$ are cohomogeneity-$l_{1}$ and
codimension-$l_{2}$ fiberwise Kahler toric manifolds.
\end{theorem}
\begin{proof}
If $(N,\bar{g})$ can be written as a product metric then the
potential $\Phi(x_{1},...,x_{d})$ can be written as
$\Phi=\Phi_{1}(x_{1},...,x_{l_{1}})+\Phi_{2}(x_{l_{1}+1},...,x_{d})$
(here we might have transformed by $SL(d,\mathbb{Z})$ and
reordered the coordinates).  The right-hand side of equation
(\ref{FKTEinsteinD}) splits as the sum of two functions: one in
the variables $x_{1},...,x_{l_{1}}$ and one in the coordinates
$x_{l_{1}},...,x_{d}$.  On the left-hand side, the
$\textrm{log}(\textrm{det}(h))$ term similarly splits.  Therefore,
the $\textrm{log}(V^{\frac{1}{2}})$ term must split in the same
way which occurs if and only if the total manifold is a product.

\end{proof}
\section{The Futaki Invariant on Fiberwise Kahler Toric Manifolds}

    Above, we reviewed the known holomorphic obstructions to the existence
of Kahler-Einstein metrics with positive Einstein constant.  More
specifically, we demonstrated how to calculate the Futaki
invariant on Kahler toric manifolds with positive first Chern
class in terms of the symplectic coordinates.  In this section we
want to calculate the Futaki invariant on fiberwise Kahler toric
manifolds in terms of the coordinates introduced above . We will
see that the formulae we obtain are very similar to those obtained
in the Kahler toric case.

    Let $(M^{2n},g,J)$ be a compact fiberwise Kahler toric manifold with
positive first Chern class: $c_{1}(M)>0$.  If a manifold of this
type admits a Kahler-Einstein metric then the Futaki invariant
must be zero.  Recall that for $X\in\mathfrak{h}(M)$ a holomorphic
vector field, the Futaki invariant satisfies
\begin{equation*}
\mathcal{F}(X)=-\int_{N}\mu_{X}\omega^{n}.
\end{equation*}

    To calculate the Futaki invariant on a fiberwise Kahler toric manifold, consider
the (real) holomorphic vector fields $\frac{\partial}{\partial
u^{i}}$.  As shown above, we know that
$\mu_{\frac{\partial}{\partial
u^{i}}}=\frac{\partial\eta}{\partial u^{i}}=x_{i}$.  Furthermore,
given the form of the metric of a fiberwise Kahler toric manifold,
we see that
\begin{equation*}
\omega^{n}=\prod_{j=1}^{m}A_{j}^{\frac{d_{j}}{2}}dx_{1}\wedge..\wedge
dx^{d}\wedge...\wedge d\phi_{1}\wedge...\wedge d\phi_{d}\wedge
d\textrm{vol}^{\ast}
\end{equation*}
where $d\textrm{vol}^{\ast}$ is the volume form of the metric
$g^{\ast}$ associated to either the coadjoint orbit $G/L$ or the
product of Kahler-Einstein manifolds $M_{1}\times...\times M_{m}$.
If $\textrm{vol}^{\ast}$ is the volume associated
$d\textrm{vol}^{\ast}$ then we calculate that
\begin{equation*}
\mathcal{F}\left(\frac{\partial}{\partial
u^{i}}\right)=-\int_{M}x_{i}\omega^{n}=-(\textrm{vol}^{\ast})(2\pi)^{d}\int_{\triangle}x_{i}
\prod_{j=1}^{m}A_{j}^{\frac{d_{j}}{2}}dx
\end{equation*}
for all $i$ where $\triangle$ is the polytope defined in
symplectic coordinates.

We have the following proposition.

\begin{prop}
If $(M,g,J)$ is a compact fiberwise Kahler toric manifold with a
Kahler-Einstein metric then
\begin{equation*}
\int_{\triangle}x_{i}V^{\frac{1}{2}}dx=0
\end{equation*}
where
$V^{\frac{1}{2}}=\prod_{j=1}^{m}A_{j}^{\frac{d_{j}}{2}}=\prod_{j=1}^{m}(\sum_{r=1}^{d}
b_{rj}x_{r}+a_{j})^{\frac{d_{j}}{2}}$.
\end{prop}
In the case when $d=1$, Koiso and Sakane proved that this
condition was in fact sufficient.  It would be interesting to
determine whether the vanishing of the Futaki invariant is
sufficient for the existence of a Kahler-Einstein metric for
higher values of $d$.

\chapter{Hermitian Metrics of Constant Scalar Curvature}

    In our study of fiberwise Kahler toric manifolds with $d=1$, we were able to integrate the
constant scalar curvature equation and express the resulting
metrics explicitly in terms of rational functions.  We will now
perform this analysis on manifolds of a very similar structure
which are Hermitian and non-Kahler.  Manifolds of this type were
first studied by Wang and Wang in \cite{WangWang} who provide
solutions to the Einstein condition on such manifolds.  After
reviewing their construction of these manifolds, we will calculate
the scalar curvature equation.

    We begin, as we did in constructing fiberwise Kahler toric manifolds, by letting
$(M_{j},g_{j}^{\ast})$ be Kahler-Einstein manifolds of real
dimension $d_{j}$ with positive Einstein constant for $1\leq j\leq
m$.  Let $P$ be a principal $S^{1}$-bundle over the product
manifold $M_{1}\times...\times M_{m}$.  The first Chern class can
be written as $c_{1}(M_{j})=p_{j}\alpha_{j}$ where $p_{j}$ is a
positive integer and $\alpha_{j}$ is an indivisible class in
$H^{2}(M_{j},\mathbb{Z})$ \cite{Ko-No}.  Let
$\chi=\sum_{j=1}^{m}b_{j}\pi^{\ast}_{j}\alpha_{j}$, where
$\pi_{j}$ is the projection from $P$ onto the manifold $M_{j}$, be
the Euler class which determines the principal $S^{1}$-bundle as
in \cite{WangZiller}.  Finally, let $\theta$ be the connection of
this $S^{1}$-bundle whose curvature represents the Euler class.

    We define a metric $g$ on $P\times I$ where $I$ is some open interval in $\mathbb{R}$ by
the equation
\begin{equation*}
g=dt^{2}+h(t)\theta\otimes\theta+\sum_{j=1}^{m}A_{j}(t)\pi_{j}^{\ast}g^{\ast}_{j}
\end{equation*}
where $h$ and the $A_{j}$ are positive and smooth functions on
$I$.  Let $M$ be the manifold that we will construct by extending
$g$ over special orbits.

    We perform a change of variable by letting $dx=\sqrt{h}dt$.  Under
this transformation, the metric becomes
\begin{equation*}
g=\frac{1}{h}dx^{2}+h(x)\theta\otimes\theta+\sum_{j=1}^{m}A_{j}(x)\pi_{j}^{\ast}g^{\ast}_{j}.
\end{equation*}

This metric contains, trivially, a two-dimensional Kahler toric
metric $\bar{g}=\frac{1}{h}dx^{2}+h(t)\theta\otimes\theta$ with
complex structure $\bar{J}$ which satisfies
$\bar{J}(\partial_{\phi})=h^{2}\partial_{x}$ where
$\partial_{\phi}$ is the dual vector field to the connection
$\theta$.  By combining $\bar{J}$ and the lifts of the complex
operators $J_{i}$ in the natural way, the manifold $(M,g)$ can be
endowed with a complex structure $J$.

    In these variables the Einstein condition becomes (see M. Wang's report in \cite{Survey})
the following system of equations
\begin{equation}\label{WangEinstOne}
-\frac{1}{2}h''-h'\left(\sum_{j=1}^{m}\frac{d_{j}}{4}\frac{A_{j}'}{A_{j}}\right)+h
\left[\sum_{j=1}^{m}\frac{d_{j}}{4}\left(\frac{A_{j}'}{A_{j}}\right)^{2}-\sum_{j=1}^{m}
\frac{d_{j}}{2}\frac{A_{j}''}{A_{j}}\right]=\lambda
\end{equation}
\begin{equation}\label{WangEinstTwo}
-\frac{1}{2}h''-h'\left(\sum_{j=1}^{m}\frac{d_{j}}{4}\frac{A_{j}'}{A_{j}}\right)+h
\left(\sum_{j=1}^{m}\frac{d_{j}}{4}\frac{b_{j}^{2}}{A_{j}^{2}}\right)=\lambda
\end{equation}
and
\begin{equation}\label{WangEinsteinThree}
-\frac{1}{2}h'\frac{A_{i}'}{A_{i}}+h\left[-\frac{1}{2}\frac{A_{i}''}{A_{i}}+\frac{1}{2}
\left(\frac{A_{i}'}{A_{i}}\right)^{2}-\frac{A_{i}'}{A_{i}}\left(\sum_{j=1}^{m}\frac{d_{j}}{4}
\frac{A_{j}'}{A_{j}}\right)-\frac{1}{2}\frac{b_{i}^{2}}{A_{i}^{2}}\right]+\frac{1}{A_{i}}
=\lambda
\end{equation}
for all $1\leq i\leq m$.  (Note that $f'=\frac{df}{dx}$ for all
functions $f(x)$.)

    Before proceeding to the question of finding constant scalar curvature metrics we
first review some of the results found in \cite{WangWang}.  By
subtracting equation (\ref{WangEinstOne}) from
(\ref{WangEinstTwo}) we find that
\begin{equation}\label{sumofmus}
\sum_{j=1}^{m}d_{j}\left[2\frac{A_{j}''}{A_{j}}-\left(\frac{A_{j}'}{A_{j}}\right)^{2}+
\frac{b_{j}^{2}}{A_{j}^{2}}\right]=0.
\end{equation}

    Following Wang and Wang we set
\begin{equation*}
\mu_{j}=\left[2\frac{A_{j}''}{A_{j}}-\left(\frac{A_{j}'}{A_{j}}\right)^{2}+\frac{b_{j}^{2}}
{A_{j}^{2}}\right]
\end{equation*}

    Although unable to solve the Einstein condition in full generality, Wang and
Wang were able to obtain explicit solutions when $\mu_{j}=0$ for
all $1\leq i\leq m$.  They were able to assign geometric meaning
to this condition.   They say that demanding that the $\mu_{j}=0$
for all $j$ was equivalent to demanding that the Riemannian
curvature tensor, $R$, exhibit a particular kind of $J$
invariance.

~

\textbf{Assumption:}

\begin{equation*}
R(X,Y,Z,W)=R(JX,JY,JZ,JW)
\end{equation*}
for all $X,Y,Z,$ and $W$ in $TM$ where $R$ is the curvature tensor
of $g$.  Hermitian metrics of this type were studied in \cite{gr}
and \cite{Falcitelli}.

Let $h'$ denote differentiation with respect to the new coordinate
$x$.  We have the following proposition.

\begin{prop}\cite{WangWang}
The above assumption is equivalent to the condition that
\begin{equation}\label{Hermcodim1cond}
\mu_{j}=0
\end{equation}
for all $1\leq j\leq m$.
\end{prop}
\begin{proof}
The proof is a straightforward calculation.  We refer the reader
to \cite{WangWang} which gives this equation in slightly different
coordinates.
\end{proof}

Equation (\ref{Hermcodim1cond}) admits two types of solutions.
Namely,
\begin{equation}\label{solution1}
A_{j}=b_{j}x+a_{j}
\end{equation}
or
\begin{equation}\label{solution2}
A_{j}=e_{j}(x+c_{j})^{2}-\frac{1}{4}\frac{b_{j}^{2}}{e_{j}}
\end{equation}
where $a_{j}, c_{j},$ and $e_{j}$ are the constants of
integration.

    Note that solution (\ref{solution1}) is equivalent to the Kahler condition if
it holds for all $j$.  If, however, solution (\ref{solution1})
does not hold for some $j$, the manifold will be Hermitian and
non-Kahler.  We will not discuss their work further except to note
that in the Kahler case, the resulting metric is fiberwise Kahler
toric.

    We now return to the question of integrating the constant
scalar curvature equation on manifolds of this type.  From
equations (\ref{WangEinstOne})-(\ref{WangEinsteinThree}), it is
straightforward to retrieve the equation for the scalar curvature,
$S$, of this kind of metric.  We have
\begin{equation}\label{WangScalar}
S=-h''-h'\left(\sum_{j=1}^{m}d_{j}\frac{A_{j}'}{A_{j}}\right)
\end{equation}
\begin{equation*}
+h\left[-\sum_{j=1}^{m}d_{j}\frac{A_{j}''}{A_{j}}+\frac{3}{4}\sum_{j=1}^{m}d_{j}\left(
\frac{A_{j}'}{A_{j}}\right)^{2}-\frac{1}{4}\left(\sum_{j=1}^{m}d_{j}\frac{A_{j}'}{A_{j}}
\right)^{2}-\sum_{j=1}^{m}\frac{d_{j}}{4}\frac{b_{j}^{2}}{A_{j}^{2}}\right]+\sum_{j=1}^{m}
\frac{d_{j}}{A_{j}}.
\end{equation*}
    Although this equation is only second-order in $h$ we are unable to integrate it in the
general case.  However, we can integrate this equation explicitly
for a large subclass.  To see how this is done, let us recall our
method of integration in the Kahler case. The Kahler condition
requires that $A_{j}=b_{j}x+a_{j}$ for all $1\leq j\leq m$.  We
were able to integrate this equation to find that
\begin{equation*}
h=\frac{1}{\prod_{j=1}^{m}(b_{j}x+a_{j})^{\frac{d_{j}}{2}}}\int\left[\int\left(\sum_{j=1}^{m}
\frac{d_{j}}{b_{j}x+a_{j}}-S\right){\prod_{j=1}^{m}(b_{j}x+a_{j})^{\frac{d_{j}}{2}}}dx\right]dx.
\end{equation*}

    This suggests a possible form for $h$ in the general case.  We ask when $h$ is
of this form with $A_{i}$ a general function instead of
$b_{i}x+a_{i}$.

    The answer is provided by the following theorem.

\begin{theorem}
The scalar curvature equation (\ref{WangScalar}) is satisfied for
\begin{equation}\label{GenConstScal}
h=\frac{1}{\prod_{j=1}^{m}A_{j}^{\frac{d_{j}}{2}}}\int\left[\int\left(\sum_{j=1}^{m}\frac
{d_{j}}{A_{j}}-S\right){\prod_{j=1}^{m}A_{j}^{\frac{d_{j}}{2}}}dx\right]dx.
\end{equation}
if and only if equation (\ref{sumofmus}),
\begin{equation*}
\sum_{j=1}^{m}d_{j}\mu_{j}\equiv\sum_{j=1}^{m}d_{j}\left[2\frac{A_{j}''}{A_{j}}-\left(\frac{
A_{j}'}{A_{j}}\right)^{2}+\frac{b_{j}^{2}}{A_{j}^{2}}\right]=0,
\end{equation*}
holds.
\end{theorem}
\begin{proof}
First we calculate the first and second derivatives of $h$.  They
are
\begin{equation*}
h'=\left(-\sum_{j=1}^{m}\frac{d_{j}}{2}\frac{A_{j}'}{A_{j}}\right)h-\frac{1}{\prod_{j=1}^{m}
A_{j}^{\frac{d_{j}}{2}}}\int\left(\sum_{j=1}^{m}\frac{d_{j}}{A_{j}}-S\right)\prod_{j=1}^{m}
A_{j}^{\frac{d_{j}}{2}}dx
\end{equation*}
and
\begin{equation*}
h''=h'\left(-\sum_{j=1}^{m}\frac{d_{j}}{2}\frac{A_{j}'}{A_{j}}\right)+h\left[\sum_{j=1}^{m}
\frac{d_{j}}{2}\left(\frac{A_{j}'}{A_{j}}\right)^{2}-\sum_{j=1}^{m}\frac{d_{j}}{2}\frac{A_{j}''}
{A_{j}}\right]
\end{equation*}
\begin{equation*}
+\left(\sum_{j=1}^{m}\frac{d_{j}}{2}\frac{A_{j}'}{A_{j}}\right)\left[\frac{1}{\prod_{j=1}^{m}
A_{j}^{\frac{d_{j}}{2}}}\int\left(\sum_{j=1}^{m}\frac{d_{j}}{A_{j}}-S\right)\prod_{j=1}^{m}
A_{j}^{\frac{d_{j}}{2}}dx\right]+S-\sum_{j=1}^{m}\frac{d_{j}}{A_{j}}.
\end{equation*}
Substituting these equations into equation (\ref{WangScalar})
gives equation (\ref{sumofmus}).  This completes the proof.
\end{proof}

While equation (\ref{GenConstScal}) holds in the Einstein case, it
also holds more generally.

    One may ask what the geometric significance of equation (\ref{sumofmus}) is.  We
demonstrated above that the manifold can be endowed naturally with
a complex structure $J$.  Equation (\ref{WangEinstOne}) of the
Einstein condition is equivalent to the equation
\begin{equation*}
r\left(\frac{\partial}{\partial x},\frac{\partial}{\partial
x}\right)=\lambda g\left(\frac{\partial}{\partial
x},\frac{\partial}{\partial x}\right)
\end{equation*}
and equation (\ref{WangEinstTwo}) is equivalent to
\begin{equation*}
r\left(\frac{\partial}{\partial \phi},\frac{\partial}{\partial
\phi}\right)=\lambda g\left(\frac{\partial}{\partial
\phi},\frac{\partial}{\partial \phi}\right).
\end{equation*}
Equation (\ref{sumofmus}) states that these two equations are
equivalent.  Since the metric is already $J$-invariant, equation
(\ref{sumofmus}) is equivalent to the condition
\begin{equation}\label{ricciinvariant}
r\left(J\frac{\partial}{\partial \phi},J\frac{\partial}{\partial
\phi}\right)=r\left(\frac{\partial}{\partial
\phi},\frac{\partial}{\partial \phi}\right).
\end{equation}

Note that if (\ref{ricciinvariant}) holds, then the Ricci tensor
is completely $J$-invariant.  The preceding theorem states that
when the Ricci tensor is $J$-invariant, the scalar curvature
equation can be integrated to (\ref{GenConstScal}).

~

\section{Constant Scalar Curvature Hermitian Metrics with $J$-Invariant Ricci Tensor}

    In the previous section, we were able to integrate the scalar curvature
equation when the Ricci tensor of the Hermitian manifold was
$J$-invariant.  Before using this integration to construct new
non-Kahler Hermitian metrics with $J$-invariant Ricci tensor and
constant scalar curvature, we wish to discuss briefly why we
believe this condition to be both natural and of interest,
particularly in dimension four.  To that end, we will try to
motivate the condition as a natural generalization of both the
Hermitian-Einstein and constant scalar curvature Kahler
conditions.

    Let $(M,g,J)$ be a compact Hermitian manifold which is not necessarily
Kahler.  By the definition of an Hermitian manifold, the metric,
$g$, is automatically $J$-invariant,
\begin{equation*}
g(J\cdot,J\cdot)=g(\cdot,\cdot).
\end{equation*}
The Ricci tensor, $r$, is \textit{not} automatically
$J$-invariant. It is, however, invariant if either of two
important conditions hold.  If the metric is Kahler (i.e. for
$\omega(\cdot,\cdot)=g(J\cdot,\cdot)$, $d\omega=0$), it is
well-known that the Ricci tensor is automatically $J$-invariant:
\begin{equation*}
\textrm{Kahler}\Rightarrow r(J\cdot,J\cdot)=r(\cdot,\cdot).
\end{equation*}
If the metric, $g$, is Hermitian-Einstein, then
$r(\cdot,\cdot)=\lambda g(\cdot,\cdot)$ where $\lambda$ is a
constant.  Since $g$ is $J$-invariant, the Einstein condition
implies that the Ricci tensor is also $J$-invariant,
\begin{equation*}
\textrm{Einstein}\Rightarrow r(J\cdot,J\cdot)=r(\cdot,\cdot).
\end{equation*}
Hermitian manifolds with $J$-invariant Ricci tensor provide a
natural generalization of both the Einstein and the Kahler
conditions.  However, Kahler metrics, and by extension Hermitian
metrics with $J$-invariant Ricci tensor, are very plentiful
(though not every complex manifold admits a Kahler metric).  We
would like to place an additional condition on Hermitian metrics
with $J$-invariant Ricci tensor.  We propose that an interesting
additional condition is to demand that the metrics be of constant
scalar curvature.

    Hermitian metrics of constant scalar curvature with $J$-invariant
Ricci tensor include both Hermitian-Einstein metrics (Einstein
metrics are automatically of constant scalar curvature) and
constant scalar curvature Kahler metrics, both of which are of
great interest to differential geometers.  Manifolds of this type
are also a subclass of the space of Hermitian manifolds with
constant scalar curvature.  While the classification of these
types of manifolds is an interesting problem in all dimensions, it
is particularly intriguing in dimension four.

    In dimension four, much is known about constant scalar curvature
Kahler metrics and Hermitian-Einstein metrics.  The compact
four-manifolds admitting Kahler-Einstein metrics have been
completely classified by the work of Yau, Tian, and others (see
LeBrun's article in \cite{Survey} for a discussion).  Many open
problems still exist in the classification of constant scalar
curvature Kahler metrics.  Much work has been done on this topic
by Tian, Donaldson, and others relating the existence of constant
scalar curvature metrics to various notions of stability, but no
classification as yet exists in dimension four or in higher
dimensions.  More is known about the existence of non-Kahler
Hermitian-Einstein metrics.  As mentioned above LeBrun
demonstrated that, in addition to the manifold
$\mathbb{C}P^{2}\sharp\overline{\mathbb{C}P^{2}}$ which is known
to admit a non-Kahler Hermitian-Einstein metric, the only other
compact four manifolds which \textit{could} admit a non-Kahler
Hermitian-Einstein metric are
$\mathbb{C}P^{2}\sharp2\overline{\mathbb{C}P^{2}}$ and
$\mathbb{C}P^{2}\sharp3\overline{\mathbb{C}P^{2}}$.

    We now turn our attention to the classification of constant scalar
curvature Hermitian metrics with $J$-invariant Ricci tensor in
dimension four.  As yet no classification exists; in fact, to the
author's knowledge this particular condition has not been studied.
However, some work has been done on Hermitian metrics with
$J$-invariant Ricci tensor on complex surfaces.  For example, see
\cite{Armstrong}.  We now review some of that work, beginning with
the following theorem due to Apostolov and Gauduchon:
\begin{theorem}\cite{Apostolov}
Suppose $(M,g,J)$ is a compact Hermitian four manifold with
$J$-invariant Ricci tensor.  The metric $g$ is locally conformal
to a Kahler metric.  It is globally conformal to a Kahler metric
$\bar{g}=f^{2}g$ if and only if the first Betti number $b_{1}(M)$
is even.  In that case, the vector field $J(\nabla_{g}f)$ is a
Killing vector field of both $g$ and $\bar{g}$.
\end{theorem}
In the case of $b_{1}$ even, finding an Hermitian metric with
$J$-invariant Ricci tensor is equivalent to finding a Kahler
metric with a Killing vector field.  In the case of $b_{1}$ odd,
there are no known examples in four dimensions of a non-Kahler
Hermitian metric with $J$-invariant Ricci tensor (Apostolov,
personal communication).

    To the author's knowledge, other than the Page metric, no non-Kahler Hermitian
metrics of constant scalar curvature with $J$-invariant Ricci
tensor on a compact four manifold have been constructed
previously.  In the following section, we construct a
one-parameter family of such metrics on
$\mathbb{C}P^{1}\times\mathbb{C}P^{1}$ and
$\mathbb{C}P^{2}\sharp\overline{\mathbb{C}P^{2}}$.  It would be
interesting to classify all of the metrics of this type in
dimension four.  In addition, one could ask whether \textit{all}
constant scalar curvature Hermitian metrics \textit{must} have
$J$-invariant Ricci tensor.  Unfortunately, we are unable to give
an answer to either of these questions.

\section{New Constant Scalar Curvature Hermitian Metrics in Dimension Four}
    Consider the case in which $m=1$ and $d_{1}=2$.  That is, each principal
orbit is an $S^{1}$-bundle over $\mathbb{C}P^{1}$ endowed with its
canonical Kahler-Einstein metric.  In this case,
\begin{equation*}
g=\begin{pmatrix}
\frac{1}{h} & & \\
 & h & \\
 & & A_{1}Id_{2} \\
\end{pmatrix}.
\end{equation*}
The Kahler case was discussed above in the section on
four-dimensional extremal Kahler metrics.  In that section, we
showed that, setting $b_{1}=b$,
\begin{equation*}
b=-\frac{q}{2}
\end{equation*}
where $q$ is an integer and each value of $q$ determines a
distinct principal orbit type.  Up to a change in orientation $q$
and $-q$ give the same principal orbit type.  Without loss of
generality, we can assume that $b\geq0$ and $q\leq0$.

    If the metric has two special orbits, which via translation and rescaling
we can fix at $x=-1$ and $x=1$, and that $A_{1}>0$ for
$x\in[-1,1]$ then the resulting manifold will be a Hirzebruch
surface
\begin{equation*}
\mathbb{P}(\mathcal{O}_{\mathbb{C}P^{1}}\oplus\mathcal{O}_{\mathbb{C}P^{1}}(q)).
\end{equation*}
    We want to find constant scalar curvature metrics when the Ricci tensor is
$J$-invariant, i.e. when equation (\ref{sumofmus}) holds.  Since
$m=1$, equation (\ref{sumofmus}) is equivalent to
\begin{equation*}
2A_{1}A_{1}''-(A_{1}')^{2}+b_{j}^{2}=0.
\end{equation*}
As seen above this equation admits two solutions; either
\begin{equation*}
A_{1}=bx+a
\end{equation*}
or
\begin{equation*}
A_{1}=e(x+c)^{2}-\frac{1}{4}\frac{b^{2}}{e}.
\end{equation*}
If the first solution holds then the metric is Kahler.  If the
second solution holds then the metric will be Hermitian non-Kahler
which is the case we will be concerned with in this section.

To simplify the calculations we make a slight change in variables
and set $l=2ec$ and
$k=ec^{2}-\frac{1}{4}\frac{b_{1}^{2}}{e}=\frac{1}{4e}(l^{2}-b^{2})$.
We now have
\begin{equation*}
A_{1}=ex^{2}+lx+k.
\end{equation*}

For the metric to have scalar curvature $S=\beta$ where $\beta$ is
a constant, the function $h$ must satisfy
\begin{equation*}
h=\frac{1}{ex^{2}+lx+k}\int\left[\int\left(\frac{2}{ex^{2}+lx+k}-\beta\right)(ex^{2}
+lx+k)dx\right]dx.
\end{equation*}

Furthermore, for the metric to be smooth at the special orbits
$x=-1$ and $x=1$, we require that $h(-1)=h(1)=0$, $h'(-1)=2$, and
$h'(1)=-2$.  Performing this double integral and imposing the
boundary conditions implies that
\begin{equation*}
S=6,
\end{equation*}
\begin{equation*}
k=\frac{1}{2}=\frac{1}{4e}(l^{2}-b^{2}),
\end{equation*}
and
\begin{equation*}
h=\frac{(1-x)(1+x)(ex^{2}+2lx+1+e)}{2ex^{2}+2lx+1}.
\end{equation*}
The variable $e$ is completely determined by $b_{1}$ and $l$, we
can therefore write $h$ as
\begin{equation*}
h=\frac{(1-x)(1+x)((l^{2}-b^{2})x^{2}+4lx+2+l^{2}-b^{2})}{2((l+b)x+1)((l-b)x+1)}.
\end{equation*}
Furthermore, we find that
\begin{equation*}
A_{1}=\frac{l^{2}-b^{2}}{2}x^{2}+lx+\frac{1}{2}=\frac{1}{2}((l+b)x+1)((l-b)x+1).
\end{equation*}
For the metric to be a positive definite metric on a Hirzebruch
surface, we must have $h(x)>0$ for $x\in(-1,1)$ and $A_{i}>0$ for
$x\in[-1,1]$. The condition $A_{i}>0$ implies that
\begin{equation*}
((l+b)x+1)((l-b)x+1)>0.
\end{equation*}
For this to hold on the entire interval $x\in[-1,1]$ it is clear
that both
\begin{equation}\label{a}
(l+b)x+1>0
\end{equation}
and
\begin{equation}\label{b}
(l-b)x+1>0.
\end{equation}
Adding these two conditions together implies that $2lx+2>0$ for
$x\in[-1,1]$.  We deduce that
\begin{equation*}
|l|<1.
\end{equation*}
As mentioned at the beginning of this section, we can assume that
$b\geq0$.   Evaluating inequality (\ref{a}) at $x=-1$ and
inequality (\ref{b}) at $x=1$ implies that
\begin{equation*}
1-l>b
\end{equation*}
and
\begin{equation*}
1+l>b
\end{equation*}
respectively.  Adding these two equations together implies that
\begin{equation*}
b<1
\end{equation*}
which implies that
\begin{equation*}
-q<2.
\end{equation*}
As $q$ is assumed to be a nonpositive integer, this leaves only
two possibilities: $q=0,$ or $-1$.

~

\textbf{Case $1$: $q=0$}

~

If $q=0$ then $b=0$ as well and the manifold is
$\mathbb{C}P^{1}\times\mathbb{C}P^{1}$.  The function $h$ reduces
to
\begin{equation*}
h=\frac{(1-x)(1+x)(l^{2}x^{2}+4lx+2+l^{2})}{2(lx+1)^{2}}
\end{equation*}
and $A_{1}$ becomes
\begin{equation*}
A_{1}=\frac{1}{2}(lx+1)^{2}.
\end{equation*}
This gives a constant scalar curvature Hermitian metric for all
values of $|l|<1$.  Note that when $l=0$, we recover a product
metric on $\mathbb{C}P^{1}\times\mathbb{C}P^{1}$.  However, it is
interesting to note that this is not the Kahler-Einstein metric as
the two spheres are of different volumes.

~

\textbf{Case $2$: $q=-1$}

~

If $q=-1$ then $b=\frac{1}{2}$ and the manifold is the blow-up of
$\mathbb{C}P^{2}$ at one point, i.e.
$\mathbb{C}P^{2}\sharp\overline{\mathbb{C}P^{2}}$.  The function
$h$ reduces to
\begin{equation*}
h=\frac{(1-x)(1+x)((4l^{2}-1)x^{2}+16lx+7+4l^{2})}{2((2l+1)x+2)((2l-1)x+2)}
\end{equation*}
and $A_{1}$ becomes
\begin{equation*}
A_{1}=\frac{1}{8}((2l+1)x+2)((2l-1)x+2).
\end{equation*}
By inspection, we see that we obtain a constant scalar curvature
Hermitian metric on
$\mathbb{C}P^{2}\sharp\overline{\mathbb{C}P^{2}}$ when
$|l|<\frac{1}{2}$.

When $l=0$, we recover the Hermitian-Einstein metric (the Page
metric) on $\mathbb{C}P^{2}\sharp\overline{\mathbb{C}P^{2}}$
\cite{WangWang}.

The above results prove the following proposition.

\begin{prop}
Let $(M,J)$ be a Hirzebruch surface,
$\mathbb{P}(\mathcal{O}_{\mathbb{C}P^{1}}\oplus\mathcal{O}_{\mathbb{C}P^{1}}(q))$
then $M$ admits a $U(2)$-invariant compatible Hermitian metric,
$g$, whose Ricci tensor is $J$-invariant and whose scalar
curvature is constant if and only if $q=0$ or $-1$.  In both
cases, the manifold admits a one-parameter family of such metrics.
Furthermore, only one of these metrics, the product metric on
$\mathbb{C}P^{1}\times\mathbb{C}P^{1}$, is Kahler.
\end{prop}

\section{A Curiosity}

    We end this section by noting an intriguing similarity between a
cohomogeneity-one Hermitian-Einstein metric in dimension four and
a cohomogeneity-one Kahler-Einstein metric in dimension six.
Above, in our discussion of Extremal Kahler Fiberwise Toric
manifolds in dimension six, we saw that the Kahler-Einstein metric
found by Sakane on the manifold
$\mathbb{P}(\mathcal{O}_{\mathbb{C}P^{1}\times\mathbb{C}P^{1}}\oplus\mathcal{O}_{\mathbb{C}P^{1}
\times\mathbb{C}P^{1}}(1,-1))$ sits inside a two-parameter family
of Extremal Kahler metrics on that manifold.  In the previous
section, we proved that the non-Kahler Hermitian-Einstein found by
Page on
$\mathbb{C}P^{2}\sharp\overline{\mathbb{C}P^{2}}=\mathbb{P}(\mathcal{O}_{\mathbb{C}P^{1}}
\oplus\mathcal{O}_{\mathbb{C}P^{1}}(1))$ sits inside a
one-parameter family of non-Kahler constant scalar curvature
Hermitian metrics with $J$-invariant Ricci tensor. Here, we want
to note the striking similarity between these two metrics.

    The Kahler-Einstein metric on $\mathbb{P}(\mathcal{O}_{\mathbb{C}P^{1}\times\mathbb{C}P^{1}}
    \oplus\mathcal{O}_{\mathbb{C}P^{1}\times\mathbb{C}P^{1}}(1,-1))$ can be written as
\begin{equation*}
g=\begin{pmatrix}
\frac{1}{h} & & & \\
 & h & & \\
 & & \frac{1}{2}(2+x)Id_{2} & \\
 & & & \frac{1}{2}(2-x)Id_{2} \\
\end{pmatrix}
\end{equation*}

where
\begin{equation*}
h=\frac{(1-x^{2})(7-x^{2})}{2(4-x^{2})}.
\end{equation*}

    The non-Kahler Hermitian-Einstein metric on $\mathbb{P}(\mathcal{O}_{\mathbb{C}P^{1}}
    \oplus\mathcal{O}_{\mathbb{C}P^{1}}(1))$ can be written as
\begin{equation*}
g=\begin{pmatrix}
\frac{1}{h} & & \\
 & h & \\
 & & \frac{1}{8}(2+x)(2-x)Id_{2} \\
\end{pmatrix}
\end{equation*}

where
\begin{equation*}
h=\frac{(1-x^{2})(7-x^{2})}{2(4-x^{2})}.
\end{equation*}

    It would be interesting to know if this is merely a coincidence or
whether it is indicative of a more widespread correspondence
between Kahler-Einstein and non-Kahler Hermitian-Einstein metrics.

\section{Non-Kahler Hermitian Metrics with Constant Scalar Curvature in Higher Dimension}

    Above, we demonstrated that if the Ricci tensor of the complex
manifolds is $J$-invariant on the manifolds being considered, i.e.
if equation (\ref{sumofmus}) holds then the scalar curvature
equation can be integrated to
\begin{equation*}
h=\frac{1}{\prod_{j=1}^{m}A_{j}^{\frac{d_{j}}{2}}}\int\left[\int\left(\sum_{j=1}^{m}
\frac{d_{j}}{A_{j}}-S\right){\prod_{j=1}^{m}A_{j}^{\frac{d_{j}}{2}}}dx\right]dx.
\end{equation*}
We want to find solutions to this equation in the case when
$S=\beta$ where $\beta$ is a constant.  Unfortunately, we are not
able to find all of the solutions to this problem because we do
not know \textit{all} of the solutions to equation
(\ref{sumofmus}).  We can, however, give explicit solution for a
subfamily of solutions to equation (\ref{sumofmus}) obtained by
adding the additional condition that the Riemannian curvature
tensor is $J$-invariant in the following way:
$R(J\cdot,J\cdot,J\cdot,J\cdot)=R(\cdot,\cdot,\cdot,\cdot)$.  When
the Riemannian curvature exhibits this invariance the functions
$A_{i}$ must satisfy $A_{i}=b_{i}x+a_{i}$ or
$A_{i}=e_{j}(x+c_{j})^{2}-\frac{1}{4}\frac{b_{j}^{2}}{e_{j}}$ as
discussed above.  In that case, we are able to perform the double
integral and obtain explicit solutions as we did in the four
dimensional case.  (Note that in the case $m=1$, equation
(\ref{sumofmus}) is equivalent to the Riemannian tensor being
$J$-invariant in this way, and that in dimension four $m=1$
automatically.)  In those solutions, the function $h$ will be
rational in the variable $x$.

    If all of the $A_{i}$ are of the form $A_{i}=b_{i}x+a_{i}$ then
the metric will be Kahler.  As this case was considered in the
section on fiberwise Kahler toric manifolds, we restrict our
attention to that case in which at least one of the functions
$A_{i}$ is quadratic in $x$.  That is, there is some $k>1$ such
that $A_{j}$ is quadratic for $1\leq j\leq k$ and $A_{j}$ is
linear if $k+1\leq j\leq m$.  When $A_{j}$ is quadratic, it will
sometimes be convenient to write $A_{j}$ as

\begin{equation*}
A_{j}=e_{j}(x+c_{j})^{2}-\frac{1}{4}\frac{b_{j}^{2}}{e_{j}}=e_{j}(x+c_{j}+\frac{1}{2}\frac{b_{j}}{e_{j}})
(x+c_{j}-\frac{1}{2}\frac{b_{j}}{e_{j}})=e_{j}x^{2}+l_{j}x+t_{j}
\end{equation*}
where $l_{j}=2e_{j}c_{j}$ and
$t_{j}=e_{j}c_{j}^{2}-\frac{1}{4}\frac{b_{j}^{2}}{e_{j}}$.

\subsection{Non-Compact Metrics of Constant Scalar Curvature}

    In this section, we restrict our attention to the search for non-Kahler
Hermitian metrics with constant scalar curvature and $J$-invariant
Riemannian tensor.  We construct metrics of this type on the
interval $[0,\infty)$.  We restrict our attention to the case in
which $A_{j}(0)>0$ for all $j$.  In that case, the manifold will
be a complex line bundle over the product manifold
$M_{1}\times...\times M_{m}$. For the metric to be smooth and
positive the following must hold
\begin{itemize}
\item $h(0)=0$ and $h'(0)=2$
\item $A_{j}=e_{j}x^{2}+l_{j}x+t_{j}>0$ for $x\in[0,\infty)$ which implies that $e_{j}>0$,
$l_{j}>0$, and $t_{j}>0$ for all $1\leq j\leq k$
\item $A_{j}=b_{j}x+a_{j}>0$ for $x\in[0,\infty)$ which implies that $b_{j}\geq0$ and
$a_{j}\geq0$ for all $k+1\leq j\leq m$ and
\item $h(x)>0$ for $x\in(0,\infty)$.
\end{itemize}
(To obtain the third condition, note that $e_{j}>0$ for the metric
to be positive for large values of $x$, $t_{j}>0$ for the
$A_{j}(0)>0$.  Also, because the $A_{j}$ can have only negative
roots,  $c_{j}>|\frac{1}{2}\frac{b_{j}^{2}}{e_{j}}|$.  From the
definition of $l_{j}$ we deduce that $l_{j}>0$.)

    In the constant scalar curvature case, $h$ must satisfy (\ref{GenConstScal}) with $S=\beta$
where $\beta$ is a constant.  We set $h(x)=\frac{P(x)}{Q(x)}$
where
\begin{equation*}
P(x)=\frac{1}{\prod_{j=1}^{k}(e_{j}x^{2}+l_{j}x+t_{j})^{\frac{d_{j}}{2}}\prod_{j=k+1}^{m}(b_{j}x
+a_{j})^{\frac{d_{j}}{2}}}\times
\end{equation*}

\begin{equation*}
\int\left[\int\left(\sum_{j=1}^{k}\frac{d_{j}}{e_{j}x^{2}+l_{j}x+t_{j}}+\sum_{j=k+1}^{m}
\frac{d_{j}}{b_{j}x+a_{j}}-\beta\right)\prod_{j=1}^{k}(e_{j}x^{2}+l_{j}x+t_{j})^{\frac{d_{j}}{2}}
\prod_{j=k+1}^{m}(b_{j}x+a_{j})^{\frac{d_{j}}{2}}dx\right]dx
\end{equation*}
and
\begin{equation*}
Q(x)=\prod_{j=1}^{k}(e_{j}x^{2}+l_{j}x+t_{j})^{\frac{d_{j}}{2}}\prod_{j=k+1}^{m}(b_{j}x
+a_{j})^{\frac{d_{j}}{2}}.
\end{equation*}
Note that because we have assumed that $A_{j}(x)>0$ for all $j$
when $x\geq0$ we see that $Q(x)>0$ for $x\in[0,\infty)$.  Both $P$
and $Q$ are polynomials in $x$ because each of the $d_{j}$ is
even.

    After performing the double integration, the polynomial $P$ can be written schematically as
\begin{equation*}
P=f+wx+...-\beta\frac{\prod_{j=1}^{k}e_{j}^{\frac{d_{j}}{2}}\prod_{j=k+1}^{m}b_{j}^{\frac{d_{j}}
{2}}}{(r+2)(r+1)}x^{r+2}
\end{equation*}
where $f$ and $w$ are the constants of integration and
$r=\sum_{j=1}^{k}d_{j}+\sum_{k+1}^{m}\frac{d_{j}}{2}$.  The
smoothness conditions at $x=0$ determine the two constants of
integration.  To this, we first note that $h(0)=0$ implies that
$P(0)=0$ as $Q(0)>0$.  Therefore,
\begin{equation*}
h(0)=0\Rightarrow P(0)=0\Rightarrow f=0.
\end{equation*}
Next, because $P(0)=0$, we have that $h'(0)=\frac{P'(0)}{Q(0)}$.
Therefore,
\begin{equation*}
h'(0)=2\Rightarrow\frac{P'(0)}{Q(0)}=2\Rightarrow
w=2\prod_{j=1}^{k}t_{j}^{\frac{d_{j}}{2}}\prod_{j=k+1}^{m}a_{j}^{\frac{d_{j}}{2}}.
\end{equation*}
This implies that
\begin{equation*}
P=2\Rightarrow
w=2\prod_{j=1}^{k}t_{j}^{\frac{d_{j}}{2}}\prod_{j=k+1}^{m}a_{j}^{\frac{d_{j}}{2}}x
+...-\beta\frac{\prod_{j=1}^{k}e_{j}^{\frac{d_{j}}{2}}\prod_{j=k+1}^{m}b_{j}^{\frac{d_{j}}{2}}}
{(r+2)(r+1)}x^{r+2}.
\end{equation*}
It remains only to check that $h(x)>0$ for $x\in(0,\infty)$.  This
is equivalent to demanding that $P(x)>0$ for $x\in(0,\infty)$.
Now, as $x\rightarrow\infty$
\begin{equation*}
P(x)\rightarrow-\beta\frac{\prod_{j=1}^{k}e_{j}^{\frac{d_{j}}{2}}\prod_{j=k+1}^{m}
b_{j}^{\frac{d_{j}}{2}}}{(r+2)(r+1)}x^{r+2}.
\end{equation*}
If $\beta>0$, then $P(x)\rightarrow-\infty$ as
$x\rightarrow\infty$ and the metric will fail to be positive for
large values of $x$.  Therefore, $\beta\leq0.$
    If $\beta\leq0,$ then we see that in the polynomial $P(x)$, all
of the coefficients of the $x^{k}$ are positive as $b_{j}\geq0$,
$a_{j}>0$. $e_{j}>0$, $l_{j}>0,$ $t_{j}>0,$ and $d_{j}>0$ for all
$j$.  Therefore, $P(x)>0$ for $x\in(0,\infty)$.  All of the
smoothness conditions and the positivity condition are satisfied
for all values of $a_{j}>0,$ $e_{j}>0$, and
$c_{j}>|\frac{1}{2}\frac{b_{j}^{2}}{e_{j}}|$.  As discussed above,
it is the $b_{i}$'s  which determine the principal orbit type.  So
given a collection of $b_{i}$, we that there exists a scalar-flat
metric for all values of $a_{j},$ $e_{j},$ and $c_{j}$, as well as
a metric with scalar curvature, say, $\beta=-1$.  This gives the
following theorem:
\begin{theorem}
Let $(M,g)$ be complex line-bundle over $M_{1}\times...\times
M_{m}$ endowed with a non-Kahler Hermitian metric as described
above with
$A_{j}=e_{j}(x+c_{j})^{2}-\frac{1}{4}\frac{b_{j}^{2}}{e_{j}}$ for
$1\leq j\leq k$ and $A_{j}=b_{j}x+a_{j}$ for $k+1\leq j\leq m$
with $k>1$.  If $b_{j}>0$ for $k+1\leq j \leq m$, every such
manifold admits a scalar-flat metric (S=0) as well as a negative
scalar curvature metric (S=-1) for every value of $a_{j}>0,$
$e_{j}>0,$ and $c_{j}>\frac{1}{2}\frac{b_{j}}{e_{j}}.$  That is,
each such manifold admits an $m+k$-parameter family of metrics of
scalar-flat metrics and an $m+k$-parameter family of metrics with
scalar curvature equal to $-1$.
\end{theorem}
    When $A_{j}(0)=0$ for some $j$, additional smoothness conditions
must be satisfied.  As in the Kahler case discussed above, we let
$A_{1}(0)=0$ and $A_{i}(0)>0$ for $2\leq i\leq m$.  As in the
Kahler case, we know that $\mathbb{C}P^{\frac{d_{1}}{2}}$.  If
$A_{1}$ is linear then $A_{1}=b_{1}x$.  If, on the other hand
$A_{1}$ is quadratic, then
$A_{1}=e_{1}(x+c_{1})^{2}-\frac{1}{4}\frac{b_{1}^{2}}{e}$.  The
condition $A_{1}=0$ determines the constant $c_{1}$ to give
\begin{equation*}
A_{1}=e_{1}x^{2}+b_{1}x.
\end{equation*}
For the function $A_{1}$ to be positive on the interval
$x\in(0,\infty)$, both $b_{1}$ and $e_{1}$ must be greater than
zero.

    In both the linear and the quadratic cases, the derivative
of $A_{1}$ at $0$ implies that $A'(0)=b_{1}$.  Since the first
derivatives of $A_{1}$ at $x=0$ agree in both the linear and
quadratic cases, the additional smoothness condition must be the
same.  As in the Kahler case, we have
\begin{equation*}
b_{1}=\frac{2}{d_{1}+2}.
\end{equation*}

    In either case, the total manifold will be a $\mathbb{C}^{\frac{d_{1}+2}{2}}$-bundle
over $M_{2}\times...\times M_{m}$.

\begin{theorem}
Let $(M,g)$ be a $\mathbb{C}^{\frac{d_{1}+2}{2}}$-bundle over
$M_{2}\times...\times M_{m}$ with
$A_{j}=e_{j}(x+c_{j})^{2}-\frac{1}{4}\frac{b_{i}^{2}}{e_{j}}$ for
$2\leq j\leq k,$ $A_{j}=b_{j}x+a_{j}$ with $b_{j}>0$ or all
$k+1\leq j\leq k$, and $A_{1}$ satisfies either
$A_{1}=\frac{2}{d_{1}+2}x$ or $A_{1}=ex^{2}+\frac{2}{d_{1}+2}$.
Every such manifold admits a scalar-flat metric ($S=0$) as well as
a constant negative scalar curvature metric (S=-1) for every value
of $a_{j}>0$, $e_{j}>0$ and $c_{j}>\frac{1}{2}\frac{b_{j}}{e_{j}}$
for all $j\neq 1$.
\end{theorem}
\begin{proof}
The proof of this theorem is almost identical to the proof of the
previous theorem.  We therefore omit the details.
\end{proof}

\subsection{Compact Metrics of Constant Scalar Curvature on $\mathbb{C}P^{1}$-Bundles}

    Here, we consider the case in which the metric is defined on the interval $x\in[0,1]$
(with $h$ vanishing at the endpoints) and $A_{i}(x)>0$ for
$x\in[0,1]$ for all $1\leq i\leq m$.  In this case, the manifold
$M$ is a $\mathbb{C}P^{1}$-bundle over $M_{1}\times...\times
M_{m}$.  Above, we proved the existence of constant scalar
curvature Kahler metrics on such manifolds when at least one of
the $b_{i}$ is positive and at least one is negative.  Therefore,
we restrict our attention to the case in which $b_{i}\geq0$ (we
could have equivalently considered the case in which all of the
$b_{i}$ are nonpositive) and at least one of the $b_{i}$, say
$b_{1},$ is positive. We also assume that $m\geq2$.  We make this
last assumption because, as we saw in our analysis of the
four-dimensional case, there exist manifolds of this type which do
not admit constant scalar curvature Hermitian metrics with
$J$-invariant Riemannian curvature when $m=1$.

\begin{theorem}
Let $(M,g,J)$ be a $\mathbb{C}P^{1}$-bundle over
$M_{1}\times...\times M_{m}$ with $m\geq2$ as described above.  If
$b_{i}\geq0$ for $1\leq i\leq m$ and $b_{1}>0$, then $M$ admits at
least a one-parameter family of non-Kahler Hermitian metrics of
constant scalar curvature with $J$-invariant Riemannian curvature
tensor.
\end{theorem}
\begin{proof}
To simplify the calculations, we set $A_{i}=b_{i}x+a_{i}$ for
$1\leq i\leq m-1$ and $A_{m}=e_{m}x^{2}+l_{m}x+t_{m}=ex^{2}+lx+t$.
Because one of the terms is quadratic, the metric will be
non-Kahler.  Furthermore, we can assume without loss of generality
that $b_{j}>0$ for $2\leq j\leq m-1$ since the metric would be a
product if one of these $b_{j}$ where to equal zero.  For the
metric to be smooth and positive definite the following must hold:
\begin{itemize}
\item $h(0)=0$ and $h'(0)=2$
\item $h(1)=0$ and $h'(1)=-2$
\item $A_{i}=b_{i}x+a_{i}>0$ for all $1\leq i\leq m$ for $x\in[0,1]$ which implies
that $a_{i}>0$
\item $A_{m}=ex^{2}+lx+t>0$ for $x\in[0,1]$ and
\item $h(x)>0$ for $x\in(0,1)$.
\end{itemize}
Let $c_{i}=\frac{a_{i}}{b_{i}}>0$ for $1\leq i\leq m-1$. Note that
$c_{i}>0$ because $b_{i}$ is assumed to be greater than zero.  We
have
\begin{equation*}
\prod_{j=1}^{m}A_{i}^{\frac{d_{j}}{2}}=\left(\prod_{j=1}^{m-1}b_{j}^{\frac{d_{j}}{2}}\right)
\left(\prod_{j=1}^{m-1}(x+c_{j})^{\frac{d_{j}}{2}}\right)(ex^{2}+lx+t)^{\frac{d_{m}}{2}}
\end{equation*}
and
\begin{equation*}
\sum_{j=1}^{m}\frac{d_{j}}{A_{j}}=\left(\sum_{j=1}^{m-1}\frac{\frac{d_{j}}{b_{j}}}{x+c_{j}}
\right)+\frac{d_{m}}{ex^{2}+lx+t}.
\end{equation*}

    By changing the value of the $a_{i}$ we can set
\begin{equation*}
c_{j}=c>0
\end{equation*}
for all $1\leq j\leq m-1$.  This implies that, in the constant
scalar curvature case $S=\beta$,
\begin{equation*}
h=\frac{1}{(x+c)^{p}(ex^{2}+lx+t)^{q}}\int\left[\int\left(\frac{r}{x+c}+\frac{2q}{ex^{2}
+lx+t}-\beta\right)(x+c)^{p}(ex^{2}+lx+t)^{q}dx\right]dx
\end{equation*}
where $p=\sum_{j=1}^{m-1}\frac{d_{j}}{2}$, $q=\frac{d_{m}}{2},$
and $r=\sum_{j=1}^{m-1}\frac{d_{j}}{b_{j}}$ which are all greater
than zero.

    Set $h(x)=\frac{T(x)}{L(x)}$ where
\begin{equation*}
T(x)=\int\left[\int\left(\frac{r}{x+c}+\frac{2q}{ex^{2}+lx+t}-\beta\right)(x+c)^{p}(ex^{2}
+lx+t)^{q}dx\right]dx
\end{equation*}
and
\begin{equation*}
L(x)=(x+c)^{p}(ex^{2}+lx+t)^{q}.
\end{equation*}

    Schematically, we can write $T$ as
\begin{equation*}
T=F_{p+q+4}x^{p+q+4}+...+F_{w}x^{w}+...+Ex+F
\end{equation*}
where $E$ and $F$ are the constants of integration.  It is
straightforward to check that the smoothness conditions $h(0)=0$
and $h'(0)=2$ are satisfied when $F=0$ and $E=2c^{p}t^{q}$.

    We will see that the two smoothness conditions at $x=1$ will
determine $\beta$ and will place a condition on the constants $c,$
$l,$ and $t$.  To obtain the constant scalar curvature metrics, it
will be helpful to group the terms by the order of $c$ and $t$.
Doing this, we see that
\begin{equation*}
T=[2c^{p}t^{q}x+c^{p-1}t^{q-1}(rt+2qc)\frac{x^{2}}{2}+...]-\beta[c^{p}t^{q}\frac{x^{2}}{2}
+c^{p-1}t^{q-1}(pt+lqc)\frac{x^{3}}{6}+...]
\end{equation*}
and
\begin{equation*}
L=c^{p}t^{q}+c^{p-1}t^{q-1}(pt+lqc)x+...
\end{equation*}
The condition $h(1)=0$ implies that
\begin{equation*}
\beta=\frac{2c^{p}t^{q}+\frac{1}{2}c^{p-1}t^{q-1}(rt+2qc)+...}{\frac{1}{2}c^{p}t^{q}
+\frac{1}{6}c^{p-1}t^{q-1}(pt+lqc)+...}.
\end{equation*}
As $c\rightarrow\infty$ and $t\rightarrow\infty$, we see that
$\beta\rightarrow4$.

    The final smoothness condition states that $h'(1)=-2$ which is equivalent
to $\frac{T'(1)}{L(1)}=-2$.  This implies that
\begin{equation*}
\beta=\frac{4c^{p}t^{q}+c^{p-1}t^{q-1}(rt+2qc)+2c^{p-1}t^{q-1}(pt+lqc)+...}{c^{p}t^{q}
+\frac{1}{2}c^{p-1}t^{q-1}(pt+lqc)+...}.
\end{equation*}
Equating these two expressions for $\beta$ implies that
\begin{equation*}
\frac{2}{3}c^{2p-1}t^{2q-1}(pt+lqc)+K(c,t)=0
\end{equation*}
where $K=K(c,t)$ is a polynomial in $c$ and $t$ in which each of
the terms in $K$ is of the form $K_{\eta\xi}c^{\eta}t^{\xi}$ such
that $\eta\leq2p,$ $\xi\leq2q,$ $\eta+\xi<2p+2q-1,$ and the
$K_{\eta\xi}$ are constants depending on the other terms.  This
implies that
\begin{equation}\label{B(ct)}
B(c,t)\equiv pt+lqc+\frac{3}{2}\frac{K(c,t)}{c^{2p-1}t^{2q-1}}=0.
\end{equation}

Note that as $c$ and $t$ approach positive infinity,
$\frac{3}{2}\frac{K(c,t)}{c^{2p-1}t^{2q-1}}$ approaches a
constant.  Recall that
$A_{m}=e_{m}x^{2}+2e_{m}c_{m}x+e_{m}c_{m}^{2}-\frac{1}{4}\frac{b_{j}^{2}}{e_{m}}=ex^{2}+lx+t$.
Letting $e_{m}=e$ be fixed positive number, equation (\ref{B(ct)})
implies that
\begin{equation*}
e[p(c_{m}^{2}-\frac{1}{4}\frac{b_{m}^{2}}{e^{2}})+2qcc_{m}]+\frac{3}{2}\frac{K(c,t)}
{c^{2p-1}t^{2q-1}}=0.
\end{equation*}
As $c_{m}\rightarrow-\infty$ and $c\rightarrow\infty$, this
equation, after dividing by $ec_{m}$, implies that
\begin{equation*}
-pc_{m}-2qc+C(c,c_{m})=0
\end{equation*}
where $C(c,c_{m})\rightarrow0$ as $c\rightarrow\infty$ and
$c_{m}\rightarrow-\infty$.  As $p$ and $q$ are positive, the above
equation must (by an intermediate value argument similar to those
used above) have a solution for $-c_{m}$ and $c$ arbitrarily
large.  As $c_{m}\rightarrow\infty$ and $c\rightarrow\infty$,
$A_{m}\rightarrow\infty$ and $A_{i}\rightarrow\infty$ for $1\leq
i\leq m-1$.  Furthermore, $h$ will be positive in $(0,1)$ for
large values of $c$ and $c_{m}^{2}$.

    What we have demonstrated is that we can always solve the smoothness
conditions from $-c_{m}$ and $c$ arbitrarily large.  Moreover, the
resulting function $h$ will be positive on $x\in(0,1)$.  This
proves the theorem.
\end{proof}

\subsection{Compact Metrics of Constant Scalar Curvature on $\mathbb{C}P^{\frac{d_{1}+2}{2}}$-
Bundles}

    We now consider the case in which $A_{1}(0)=0$.  For the metric
to be smooth, we know that
$M_{1}\cong\mathbb{C}P^{\frac{d_{1}}{2}},$
$b_{1}=\frac{2}{d_{1}+2}$ and $A_{1}=e_{1}x^{2}+b_{1}x$ or
$A_{1}=b_{1}x$.  In this case the manifold $M$ is a
$\mathbb{C}P^{\frac{d_{1}+2}{2}}$-bundle over
$M_{2}\times...\times M_{m}$.

    Above, we demonstrated that such manifolds admit constant scalar
curvature Kahler metrics when $m=1$ or when at least one of the
$b_{i}$ is positive and at least one is negative for $2\leq i\leq
m$.  In this section we ask which of the remaining manifolds
admits a constant scalar curvature Hermitian metric with
$J$-invariant Riemannian curvature.
\begin{theorem}
Let $(M,g,J)$ be a $\mathbb{C}P^{\frac{d_{1}+2}{2}}$-bundle over
$M_{2}\times...\times M_{m}$ with $m\geq3$ as described above.  If
$b_{i}\geq0$ for $3\leq i\leq m$ and $b_{2}>0$, then $M$ admits at
least a one-parameter family of non-Kahler Hermitian constant
scalar curvature metrics with $J$-invariant Ricci tensor.  The
following theorem provides a partial answer to this question.
\end{theorem}
\begin{proof}
As in the proof of the theorem of the previous section we set
$A_{i}=b_{i}x+a_{i}$ for $2\leq j\leq m-1$ and set
$A_{m}=e_{m}(x+c_{m})^{2}-\frac{1}{4}\frac{b_{m}^{2}}{e_{m}}$. The
presence of a quadratic term prevents the metric from being
Kahler.  We take $A_{1}$ to be either quadratic or linear.  The
details of this proof follow almost identically to those of the
previous theorem.  For brevity, we omit the details.
\end{proof}

\chapter{Integrating the Scalar Curvature Equation on a Cohomogeneity-One Manifold}

    In this thesis, we were able to integrate the scalar
curvature equation on a Kahler manifold of cohomogeneity-one under
a semisimple Lie group in the monotypic case.  We were also able
to integrate the scalar curvature equation in a non-Kahler
Hermitian case similar to the cohomogeneity-one Kahler case.  In
this section, we describe an integration of the scalar curvature
equation which applies to a large class of cohomogeneity-one
manifolds in the absence of a complex structure or special
holonomy.

    Let $(M^{n},g)$ be a manifold of cohomogeneity-one under
the action of a compact, simply-connected Lie group $G$ with
connected isotropy subgroup $K$.  This implies that $G/K$ must be
simply-connected.  The Lie algebra $\mathfrak{g}$ decomposes as
\begin{equation*}
\mathfrak{g}=\mathfrak{k}\oplus\mathfrak{p}_{0}\oplus...\oplus\mathfrak{p}_{m}
\end{equation*}
where the $\mathfrak{p}_{i}$ are irreducible
$\textrm{Ad}(K)$-invariant subspaces.  Let $d_{i}$ denote the
dimension of the summand $\mathfrak{p}_{i}$.  We begin with two
assumptions.

~

\textbf{Assumption One:}

~

~~~~~~~~~~~~~~~~~~~~~~~~~~~~~~~~~~~~~~~~~~~~~~~~~~~~~The
$\mathfrak{p}_{i}$ are distinct.

~

We restrict ourselves to the monotypic case because the metric, as
it is automatically diagonalizable, is more easily described.

~

~

~

\textbf{Assumption Two:}

~

~~~~~~~~~~~~~~~~~~~~~~~~~~~~~~~~~~~~~~~~~~~~~~~~~~~~~~~~~~dim$\mathfrak{p}_{0}=1.$

~

This is a more serious restriction; however, assumption two does
hold in a variety of interesting cases.  For example, if the
cohomogeneity-one metric admits a compatible $G$-invariant
almost-complex structure then there must be an invariant
one-dimensional subspace in the isotropy representation.  This can
be seen by setting $\frac{\partial}{\partial t}$ to be the unit
tangent vector to the horizontal distribution.  The image of the
vector field under the almost-complex structure,
$J\frac{\partial}{\partial t}$, belongs to the tangent space to
the vertical principal orbits.  The vector field
$J\frac{\partial}{\partial t}$ is invariant under the action of
$G$ because both $\frac{\partial}{\partial t}$ and $J$ are
invariant under the group action.  Therefore,
$J\frac{\partial}{\partial t}$ spans an invariant one-dimensional
subspace.  Of course, assumption two can hold even in the absence
of an almost-complex structure.

    Because $K$ is connected and $\mathfrak{p}_{0}$ is one-dimensional
we see that $\mathfrak{p}_{0}$ is in fact a trivial
representation.
\begin{lemma}\label{trivialsummand}Let $(M,g)$ be a cohomogeneity-one manifold under
the action of a compact connected Lie group $G$ with connected isotropy subgroup $K$
satisfying assumptions one and two.  On such a manifold,
\begin{equation}\label{p_{0}}
[\mathfrak{p}_{0},\mathfrak{p}_{i}]\subseteq\mathfrak{p}_{i}
\end{equation}
for all $i$.
\end{lemma}
\begin{proof}
Clearly, $[\mathfrak{p}_{0},\mathfrak{p}_{0}]=0$ as
$\mathfrak{p}_{0}$ is one-dimensional.  Furthermore, since $K$
acts trivially on $\mathfrak{p}_{0}$, we see that
$[\mathfrak{k},\mathfrak{p}_{0}]=0$.  This implies that
$\textrm{ad}(\mathfrak{p}_{0})$ acts on
$\mathfrak{p}_{1}\oplus...\oplus\mathfrak{p}_{m}$ by
$K$-morphisms.  Because we have assumed that all of the
$\mathfrak{p}_{i}$ are distinct, Schur's Lemma implies that
$\textrm{ad}(\mathfrak{p}_{0}):\mathfrak{p}_{i}\rightarrow\mathfrak{p}_{i}$
for all $1\leq i\leq m$.
\end{proof}

Under these assumptions we can write the metric, $g$, as
\begin{equation*}
g=dt^{2}+hQ|_{\mathfrak{p}_{0}}+\bigoplus_{i=1}^{m}A_{i}Q|_{\mathfrak{p}_{i}}
\end{equation*}
where $Q$ is some $\textrm{Ad}(K)$-invariant background metric on
$\mathfrak{p}=\mathfrak{p}_{0}\oplus...\oplus\mathfrak{p}_{m}$.
Let $X_{0}$ be an orthonormal basis for $Q|_{\mathfrak{p}_{0}}$
and let $\{Y^{i}_{j}\}_{j=1,...,d_{i}}$ be an orthonormal basis
for $Q|_{\mathfrak{p}_{i}}$.  Note also that $h$ and the $A_{i}$
are positive functions of $t$.

    Let $x=x(t)$ be defined by the equation
\begin{equation*}
\frac{1}{h}dx^{2}=dt^{2}.
\end{equation*}
This implies that $\frac{dx}{dt}=\sqrt{h}$.  We can now write the
metric as
\begin{equation*}
g=\begin{pmatrix}
\frac{1}{h} & & & & \\
& h & & & \\
& & A_{1}Id_{d_{1}} & & \\
& & & ... & \\
& & & & A_{m}Id_{d_{m}} \\
\end{pmatrix}
\end{equation*}
relative to basis $(\frac{\partial}{\partial x},
X_{0},Y^{1}_{1},...,Y^{m}_{d_{m}})$.

    The Einstein condition on this metric can be written schematically as
\begin{equation*}
r\left(\frac{\partial}{\partial x},\frac{\partial}{\partial
x}\right)=\lambda g\left(\frac{\partial}{\partial
x},\frac{\partial}{\partial x}\right),
\end{equation*}
\begin{equation*}
r(X_{0},X_{0})=\lambda g(X_{0},X_{0}),
\end{equation*}
and
\begin{equation*}
r(Y^{i}_{k},Y^{i}_{k})=\lambda g(Y^{i}_{k},Y^{i}_{k})
\end{equation*}
for all $i$ where $\lambda$ is the Einstein constant.  Note that
which $1\leq k\leq d_{i}$ we choose is immaterial.  In terms of
the functions $h$ and $A_{i}$ it is straightforward to calculate
that these equations become
\begin{equation}\label{Einscalarint1}
-\frac{1}{2}h''-h'\left(\sum_{j=1}^{m}\frac{d_{j}}{4}\frac{A_{j}'}{A_{j}}\right)+h
\left[\sum_{j=1}^{m}\frac{d_{j}}{4}\left(\frac{A_{j}'}{A_{j}}\right)^{2}-\sum_{j=1}^{m}
\frac{d_{j}}{2}\frac{A_{j}''}{A_{j}}\right]=\lambda,
\end{equation}
\begin{equation}\label{Einscalarint2}
-\frac{1}{2}h''-h'\left(\sum_{j=1}^{m}\frac{d_{j}}{4}\frac{A_{j}'}{A_{j}}\right)+\frac{1}{h}
\hat{r}(X_{0},X_{0})=\lambda,
\end{equation}
and
\begin{equation}\label{Einscalarint3}
-\frac{1}{2}h'\frac{A_{i}'}{A_{i}}+h\left[-\frac{1}{2}\frac{A_{i}''}{A_{i}}+\frac{1}{2}
\left(\frac{A_{i}'}{A_{i}}\right)^{2}-\frac{A_{i}'}{A_{i}}\left(\sum_{j=1}^{m}\frac{d_{j}}{4}
\frac{A_{j}'}{A_{j}}\right)\right]+\frac{1}{A_{i}}\hat{r}(Y_{k}^{i},Y^{i}_{k})=\lambda
\end{equation}
for all $1\leq i\leq m$.  Recall that $\hat{r}$ is the Ricci
curvature of the principal orbits as described in Chapter $2$.
More precisely, given a point $x$ in the quotient space $M/G$,
$\hat{r}|_{x}$ is the Ricci tensor of the principal orbit located
at $x$.

    We note that equations (\ref{Einscalarint1}) and
(\ref{Einscalarint2}) together imply that
\begin{equation}\label{rinv}
\frac{1}{h^{2}}\hat{r}(X_{0},X_{0})=\sum_{j=1}^{m}\frac{d_{j}}{4}\left(\frac{A_{j}'}{A_{j}}
\right)^{2}-\sum_{j=1}^{m}\frac{d_{j}}{2}\frac{A_{j}''}{A_{j}}.
\end{equation}

The scalar curvature equation can be written as
\begin{equation}\label{generalscalar}
S=-h''-h'\left(\sum_{j=1}^{m}d_{j}\frac{A_{j}'}{A_{j}}\right)+h\left[-\sum_{j=1}^{m}d_{j}
\frac{A_{j}''}{A_{j}}+\frac{3}{4}\sum_{j=1}^{m}d_{j}\left(\frac{A_{j}'}{A_{j}}\right)^{2}-
\left(\sum_{j=1}^{m}\frac{d_{j}}{2}\frac{A_{j}'}{A_{j}}\right)^{2}\right]+\hat{S}
\end{equation}
where
\begin{equation*}
\hat{S}=\frac{1}{h}\hat{r}(X_{0},X_{0})+\sum_{j=1}^{m}\frac{d_{j}}{A_{j}}\hat{r}(Y_{k}^{j},
Y_{k}^{j}).
\end{equation*}
is the scalar curvature of the homogeneous principal orbits and
$S$ is the scalar curvature of the total cohomogeneity-one
manifold $M$.  Specifically, for $x$ a point in the quotient
space, $\hat{S}|_{x}$ is the scalar curvature of the principal
orbit located at the point $x$.

    In \cite{WangZiller}, Wang and Ziller give a formula ($1.3$ on page
$181$) for the scalar curvature of a homogeneous manifold.  As
each of the principal orbits of $M$ is a homogeneous manifold,
$G/K$, we can use this formula to determine the form of $\hat{S}$.
From that formula we see that $\hat{S}$ can be written in terms of
the $h$ and the $A_{i}$.  More precisely, we can decompose
$\hat{S}$ by setting
\begin{equation*}
\hat{S}=S_{h}+S^{\ast}
\end{equation*}
where $S_{h}$ includes all of the terms involving $h$ and
$S^{\ast}$ does not include any terms with $h$ in them but only
terms involving the $A_{i}$. Because $\mathfrak{p}_{0}$ is
trivial, we have the following proposition:
\begin{prop}  For $(M,g)$ a cohomogeneity-one manifold satisfying assumptions one and two as
above, we have that
\begin{equation*}
S_{h}=-h\sum_{j=1}^{m}\frac{d_{j}}{4}\frac{b_{j}^{2}}{A_{j}^{2}}
\end{equation*}
for some constants $b_{j}$.
\end{prop}
\begin{proof}
Using lemma (\ref{trivialsummand}) and the fact that
$\mathfrak{p}_{0}$ is trivial, we see from the formula in
\cite{WangZiller} that the only factors of $\hat{S}$ which involve
the function $h$ are of the form
$-h\frac{d_{j}}{4}\frac{b_{j}^{2}}{A_{j}^{2}}$ for $1\leq j\leq
m$.
\end{proof}

The form of the function $S^{\ast}$ will vary depending on the
groups $G$ and $K$.  However, we stress that $S^{\ast}$ does
\textit{not} involve the function $h$.

    We next consider the term $\frac{1}{h}\hat{r}(X_{0},X_{0})$.

\begin{prop}\label{scalricci}
Let $(M,g)$ be a cohomogeneity-one manifold satisfying assumptions
one and two.  We have that
\begin{equation*}
\frac{1}{h}\hat{r}(X_{0},X_{0})=h\left(\sum_{j=1}^{m}\frac{d_{j}}{4}\frac{b_{j}^{2}}{A_{j}^{2}}
\right)
\end{equation*}
for some constants $b_{i}$.
\end{prop}
\begin{proof}
Again, this result follows from lemma (\ref{trivialsummand}) and
the fact the $\mathfrak{p}_{0}$ is trivial.  Given these two
facts, along with the discussion in \cite{WangZiller} and equation
$7.38$ of \cite{Besse}, it is straightforward to calculate that
this proposition holds.
\end{proof}
From the above we deduce immediately that
\begin{equation*}
S_{h}=-\frac{1}{h}\hat{r}(X_{0},X_{0}).
\end{equation*}

    Furthermore, proposition (\ref{scalricci}) and equation (\ref{rinv})
imply that
\begin{equation*}
\sum_{j=1}^{m}d_{j}\mu_{j}=0
\end{equation*}
where
\begin{equation*}
\mu_{j}=2\frac{A_{j}''}{A_{j}}-\left(\frac{A_{j}'}{A_{j}}\right)^{2}+\frac{b_{j}^{2}}{
A_{j}^{2}}.
\end{equation*}
    We pause now to discuss the significance of equation (\ref{rinv}).  We
define a map $\tilde{J}$ by the action
$\tilde{J}\frac{\partial}{\partial x}=\frac{1}{h}X_{0}$ and
$\tilde{J}X_{0}=-h\frac{\partial}{\partial x}$.  The map
$\tilde{J}$ can be written as
\begin{equation*}
\tilde{J}=\begin{pmatrix}
& -h & & & \\
\frac{1}{h} & & & & \\
& & 0 & & \\
& & & ... & \\
& & & & 0 \\
\end{pmatrix}
\end{equation*}
relative to basis $(\frac{\partial}{\partial x},
X_{0},Y_{1}^{1},...,Y_{d_{m}}^{m})$.  Clearly, the metric $g$ is
invariant under the action of $\tilde{J}$, i.e.
$g(\cdot,\cdot)=g\left(\tilde{J}\cdot,\tilde{J}\cdot\right)$.
However, the Ricci tensor, $r$, is \textit{not} automatically
invariant under the action of $\tilde{J}$.  In fact,
\begin{equation*}
r(\cdot,\cdot)|_{\langle\frac{\partial}{\partial
x}\rangle\oplus\mathfrak{p}_{0}}=r\left(\tilde{J}\cdot,\tilde{J}\cdot\right)|_{\langle\frac{
\partial}{\partial x}\rangle\oplus\mathfrak{p}_{0}}\Leftrightarrow
(\ref{rinv})~\textrm{holds}\Leftrightarrow\sum_{j=1}^{m}d_{j}\mu_{j}=0.
\end{equation*}

    We also note that the metric $g$ and the map $\tilde{J}$ define a $2$-form of
rank two such that $\tilde{\omega}=g\tilde{J}$, i.e.
$\tilde{\omega}(\cdot,\cdot)=g(\tilde{J}\cdot,\cdot)$.  This can
be written as
\begin{equation*}
\tilde{\omega}=\begin{pmatrix}
& -1 & & & \\
1 & & & & \\
& & 0 & & \\
& & & ... & \\
& & & & 0 \\
\end{pmatrix}
\end{equation*}
relative to basis $(\frac{\partial}{\partial x},
X_{0},Y_{1}^{1},...,Y_{d_{m}}^{m})$.  Alternatively, we could
write $\tilde{\omega}$ as
\begin{equation*}
\tilde{\omega}=dx\wedge\theta_{0}
\end{equation*}
where $\theta_{0}$ is the invariant dual $1$-form to the vector
field $X_{0}$ such that $\theta_{0}(X_{0})=1$.

We are now prepared to state and prove this section's main
theorem.
\begin{theorem}
Let $(M,g)$ be a cohomogeneity-one manifold satisfying assumptions
one and two.  The function $h$ can be written as
\begin{equation}\label{hdefinition}
h=\frac{1}{V^{\frac{1}{2}}}\int\left[\int(S^{\ast}-S)V^{\frac{1}{2}}dx\right]dx,
\end{equation}
where $V^{\frac{1}{2}}=\prod_{j=1}^{m}A_{j}^{\frac{d_{j}}{2}}$, if
and only if
\begin{equation*}
\sum_{j=1}^{m}d_{j}\mu_{j}=0.
\end{equation*}
\end{theorem}
\begin{proof}
Differentiating equation (\ref{hdefinition}) implies that
\begin{equation*}
S=-h''-h'\left(\sum_{j=1}^{m}d_{j}\frac{A_{j}'}{A_{j}}\right)+h\left[-\sum_{j=1}^{m}\frac{d_{j}}
{2}\frac{A_{j}''}{A_{j}}+\sum_{j=1}^{m}\frac{d_{j}}{2}\left(\frac{A_{j}'}{A_{j}}\right)^{2}-
\left(\sum_{j=1}^{m}\frac{d_{j}}{2}\frac{A_{j}'}{A_{j}}\right)^{2}\right]+S^{\ast}.
\end{equation*}
Comparing this to equation (\ref{generalscalar}) implies that
\begin{equation*}
S_{h}=h\left[-\sum_{j=1}^{m}\frac{d_{j}}{4}\left(\frac{A_{j}'}{A_{j}}\right)^{2}+\sum_{j=1}^{m}
\frac{d_{j}}{2}\frac{A_{j}''}{A_{j}}\right].
\end{equation*}
Since $S_{h}=-\frac{1}{h}\hat{r}(X_{0},X_{0})$, this is equivalent
to $\sum_{j=1}^{m}d_{j}\mu_{j}=0$.
\end{proof}

Put another way, the scalar curvature equation can be integrated
to give (\ref{hdefinition}) if and only if the Ricci tensor on
$\langle\frac{\partial}{\partial x}\rangle\oplus\mathfrak{p}_{0}$
is invariant under the map $\tilde{J}$. In the Einstein case, the
Ricci tensor must be so invariant and, therefore, equation
(\ref{hdefinition}) \textit{must} hold in the Einstein case.

\bibliographystyle{plain}

\end{document}